\magnification=1100
\font\bigbf=cmbx10 scaled \magstep2
\font\medbf=cmbx10 scaled \magstep1 
\hfuzz=50 pt
\input psfig.sty
\pageno=0
\def\la{\big\langle}
\def\ra{\big\rangle}
\font\medbf=cmbx10 at 13pt
\def\i{\item}
\def\n{\noindent}
\def\L{{\bf L}}
\def\dint{\int\!\!\int}
\def\forall{\hbox{~~for all}~~}
\def\wto{\rightharpoonup}
\def\ve{\varepsilon}

\def\dag{\diamondsuit}
\def\A{{\cal A}}
\def\D{{\cal D}}
\def\vp{\varphi}
\def\bfv{{\bf v}}
\def\b{\bullet}
\def\M{{\cal M}}
\def\N{{\cal N}}
\def\I{{\cal I}}
\def\Le{{\cal L}}
\def\O{{\cal O}}

\def\sgn{{\rm sign}}
\def\C{{\cal C}}
\def\R{{I\!\!R}} 
\def\c{\centerline}
\def\Hat{\widehat}

\def\dint{\int\!\!\int}
\def\T{{\cal T}}
\def\sgn{\hbox{sign}}

\def\tla{\tilde \lambda}
\def\hr{\hat r}
\def\vth{\eta}
\def\tr{\tilde r}

\def\tv{\hbox{Tot.Var.}}
\def\vs{\vskip 2em}
\def\vsk{\vskip 4em}
\def\v{\vskip 1em}
\def\sqr#1#2{\vbox{\hrule height .#2pt
\hbox{\vrule width .#2pt height #1pt \kern #1pt
\vrule width .#2pt}\hrule height .#2pt }}
\def\square{\sqr74}
\def\endproof{\hphantom{MM}\hfill\llap{$\square$}\goodbreak}
\def\Hat{\widehat}

\def\dint{\int\!\!\int}
\def\T{{\cal T}}
\def\sgn{\hbox{sign}}

\def\tv{\hbox{Tot.Var.}}
\def\osc{\hbox{Osc.}}
\def\vs{\vskip 2em}
\def\vsk{\vskip 4em}
\def\v{\vskip 1em}
\def\tla{\tilde\lambda}
\def\vars{\varsigma}
\def\hth{\hat\theta}
\def\vth{\eta}
\def\hr{\hat r}
\def\tr{\tilde r}
\def\hga{\hat\gamma}

\def\hvt{\hat\eta}
\def\io{\Upsilon}
\def\iot{g} 
\def\h{h}
\def\vvl{VVL}
\def\Ll{{\bf L}^1_{loc}}
\def\U{{\cal U}}
\def\hr{\hat r}
\def\hla{\hat\lambda}
\null

\c{\bigbf Vanishing Viscosity Solutions of}
\v
\c{\bigbf  Nonlinear Hyperbolic Systems}
\vsk
\centerline{\it Stefano Bianchini$\,^{(*)}$ and Alberto Bressan$^{(**)}$}
\vs
\centerline{(*)~~I.A.C. - Viale del Policlinico 137, Roma 00161, Italy}
\centerline{e-mail: bianchin@iac.rm.cnr.it}
\v
\centerline{(**)~~S.I.S.S.A. - Via Beirut 4, Trieste 34014, Italy}
\centerline{e-mail: bressan@sissa.it}
\vs
\c{(Dedicated to Prof.~Constantin Dafermos in the occasion of his 60-th
birthday)}
\vsk

\centerline{\bigbf Contents}
\vs
\halign{\qquad\quad#\hfil&\hfil\qquad\qquad\qquad\qquad#\cr
~~~~1 - Introduction & 1\cr
~~~~2 - Parabolic estimates& 7\cr
~~~~3 - Outline of the BV estimates& 12\cr
~~~~4 - A center manifold of viscous travelling waves& 15\cr
~~~~5 - Gradient decomposition& 19\cr
~~~~6 - Bounds on the source terms& 23\cr
~~~~7 - Transversal wave interactions& 32\cr
~~~~8 - Functionals related to shortening curves& 35\cr
~~~~9 - Energy estimates& 40\cr
~~~~10 - Proof of the BV estimates& 45\cr
~~~~11 - Stability estimates& 47\cr
~~~~12 - Propagation speed& 55\cr
~~~~13 - The vanishing viscosity limit& 58\cr
~~~~14 - The non-conservative Riemann problem& 63\cr
~~~~15 - Viscosity solutions and uniqueness of the semigroup& 71\cr
~~~~16 - Dependence on parameters and large time asymptotics& 74\cr
~~~~& \cr
~~~~Appendix A & 76\cr

~~~~Appendix B & 83\cr
~~~~Appendix C & 93\cr
~~~~Appendix D & 96\cr
~~~~References & 97\cr}
\vfill\eject

\n{\medbf 1 - Introduction}
\v
The Cauchy problem for
a system of conservation laws in one space dimension takes the form
$$u_t+f(u)_x=0,\eqno(1.1)$$
$$u(0,x)=\bar u(x).\eqno(1.2)$$
Here $u=(u_1,\ldots,u_n)$ is the vector of {\it conserved quantities},
while the components of $f=(f_1,\ldots,f_n)$ are the {\it fluxes}.
We assume that the flux function $f:\R^n\mapsto\R^n$ is smooth and that
the system is {\it
strictly hyperbolic}, i.~e., at each point $u$ 
the Jacobian matrix $A(u)=Df(u)$ has $n$ real, distinct eigenvalues
$$\lambda_1(u)<\cdots <\lambda_n(u).\eqno(1.3)$$
One can then select bases of right and left
eigenvectors $r_i(u)$, $l_i(u)$, normalized so that
$$ |r_i|\equiv 1\,\qquad\qquad l_i\cdot r_j=\cases{1\quad &if\quad $i=j$,\cr
0\quad &if\quad $i\not=j$.\cr}\eqno(1.4)$$

Several fundamental laws of physics take the form of
a conservation equation.
For the relevance of hyperbolic
conservation laws in continuum physics we refer to the 
recent book of Dafermos [D].

A distinguished feature of nonlinear hyperbolic systems is
the possible loss of regularity.  Even with smooth initial data,
it is well known that the solution can develop shocks within finite time.
Therefore, global solutions can only be
constructed within a space of discontinuous functions.
The equation (1.1) must then be interpreted in distributional sense.
A vector valued function $u=u(t,x)$ is a {\it weak solution} of
(1.1) if
$$\dint \big[u\,\phi_t+f(u) \,\phi_x\big]\,dxdt=0\eqno(1.5)$$
for every test function $\phi\in\C^1_c$, continuously differentiable with
compact support.  When discontinuities are present,
weak solutions may not be unique.
To single out a unique ``good'' solution of 
the Cauchy problem, additional {\it entropy
conditions} must be imposed along shocks [Lx], [L1]. 
These are often motivated
by physical considerations [D].

Toward a rigorous mathematical analysis of solutions,
the lack of regularity has always been a considerable source of
difficulties.
For discontinuous solutions, most of the standard tools
of differential calculus do not apply. Moreover, for general $n\times n$
systems, the powerful techniques of functional analysis
cannot be used. In particular, solutions cannot be obtained as fixed
points of a nonlinear transformation, or in variational form
as critical
points of a suitable functional.  Dealing with 
vector valued functions, comparison arguments
based on upper and lower solutions do not apply either.
Up to now, the theory of conservation laws has thus progressed
largely by developing {\it ad hoc} methods. In particular, 
a basic building block is the so-called {\it Riemann problem},
where the initial data is piecewise constant with a single 
jump at the origin:
$$u(0,x)=\cases{u^-\quad &if\quad $x<0\,$,\cr
u^+\quad &if\quad $x>0\,$.\cr}$$

Weak solutions to the 
Cauchy problem (1.1)-(1.2) were constructed in the 
celebrated paper of Glimm [G].
This global existence result is valid
for small $BV$ initial data and
under the additional assumption
\v
\i{(H)} For each $i\in\{1,\ldots,n\}$, the $i$-th characteristic field
is either {\it linearly degenerate}, so that
$$D\lambda_i(u)\cdot r_i(u)= 0\qquad\forall u\,,\eqno(1.6)$$
or else it is {\it genuinely nonlinear}, i.~e.
$$D\lambda_i(u)\cdot r_i(u)>0\qquad\forall u\,.\eqno(1.7)$$
\v
\n In [G], an approximate solution of the general Cauchy problem
is obtained by piecing together
solutions of several Riemann problems, with a restarting procedure
based on random sampling.  The key step in Glimm's proof
is an a priori estimate on the total variation of the approximate solutions,
obtained by introducing a {\it wave interaction potential}.
In turn, the control of the total variation yields the compactness of the
family of approximate solutions, and hence the existence of
a strongly convergent subsequence.  
Alternative constructions of 
approximate solutions, based on front tracking approximations, were
subsequently developed in [DP1], [B2], [Ri], [BaJ].

The above existence results are all based on a compactness 
argument which, by itself,
does not guarantee the uniqueness of solutions.  The well posedness of
the Cauchy problem has now been established in a series of papers
[B3], [BC1], [BCP], [BLY], [BLF], [BG], [BLe].  
The main results can be summarized as follows:
\v
\i{-} The solutions obtained as limits of Glimm or
front tracking approximations
are unique and depend Lipschitz continuously on the initial data,
in the $\L^1$ norm.
\v
\i{-} Every small $BV$ solution of the Cauchy problem
(1.1)-(1.2) which satisfies the Lax entropy conditions
coincides with the unique limit of front tracking approximations.
\v
\n For comprehensive account of the recent uniqueness 
and stability theory we refer to [B5].
\v
A long standing conjecture is that the entropic 
solutions of
the hyperbolic system (1.1) actually coincide with the
limits of solutions to the parabolic system
$$u_t+f(u)_x=\ve\,u_{xx}\,,\eqno(1.8)_\ve$$
letting the viscosity coefficient $\ve\to 0$.
In view of the recent uniqueness results, it 
looks indeed very plausible that the vanishing viscosity limit
should single out the unique ``good'' solution of the Cauchy problem,
satisfying the appropriate entropy conditions.  
In earlier literature, results 
in this direction were based on three main techniques:
\v
\n{\bf 1 - Comparison principles for parabolic equations.} 
For a {\it scalar} conservation law, the 
existence, uniqueness and global stability of vanishing
viscosity solutions has been established in a famous paper by
Kruzhkov [K]. The result is valid
within the more general class of $\L^\infty$ solutions,
also in several space dimensions.
For an alternative approach based on nonlinear semigroup theory, see
also [Cr].
\v
\n{\bf 2 - Singular perturbations.} Let $u$ be
a piecewise smooth solution of the $n\times n$ system (1.1),
with finitely many non-interacting, entropy admissible shocks.
In this special case, using a singular perturbation technique,
Goodman and Xin [GX] were able to construct 
a sequence of solutions $u^\ve$ to (1.8)$_\ve$, with $u^\ve\to u$
as $\ve\to 0$.  See also [Yu] for further results in this direction.
\v
\n{\bf 3 - Compensated compactness.} 
If, instead of a $BV$ bound, only
a uniform bound on the
$\L^\infty$ norm of solutions of (1.8)$_\ve$ is available, 
one can still construct
a weakly convergent subsequence $u^\ve\wto u$.
In general, we cannot expect
that this weak limit satisfies the nonlinear equations 
(1.5). However, for a class of $2\times 2$
systems, in [DP2] DiPerna showed
that this limit $u$ is indeed a weak solution of (1.1). The proof
relies on a compensated compactness argument, based on the
representation of the weak limit in terms of Young measures, which
must reduce to a Dirac mass due
to the presence of a large family of entropies.
We remark that the solution is here found
in the space $\L^\infty$.  Since the
known uniqueness results apply only to $BV$ solutions,
the uniqueness of solutions obtained by the compensated compactness method
remains a difficult open problem.   
\vs
In our point of view, to develop a satisfactory theory of
vanishing viscosity limits, the heart of the matter 
is to establish a priori $BV$ bounds on solutions $u(t,\cdot)$
of (1.8)$_\ve$, uniformly valid for all $t\in [0,\,\infty[\,$ and 
$\ve>0$.
This is indeed what we
will accomplish in the present paper.
Our results apply, more generally, to
strictly hyperbolic $n\times n$ systems with viscosity, 
not necessarily in conservation form:
$$u_t+A(u)u_x=\ve\, u_{xx}\,.\eqno(1.9)_\ve$$
As a preliminary, we observe that the rescaling
of coordinates
$s=t/\ve,~ y=x/\ve$ transforms the
Cauchy problem (1.9)$_\ve$, (1.2) into
$$u_s+A(u)u_y=u_{yy},\qquad\qquad u(0,y)=\bar u^\ve(y)\doteq \bar
u(\ve y)\,.$$
Clearly, the total variation of the initial data 
$\bar u^\ve$ does not change with $\ve$.
To obtain a priori $BV$ bounds and stability estimates 
for solutions of (1.9)$_\ve$, it thus suffices to
consider the system 
$$u_t+A(u)u_x=u_{xx}\,,\eqno(1.10)$$
and derive estimates uniformly valid for all 
times $t\geq 0\,$, depending only on the total variation
of the initial data $\bar u$.
\v
The first step in our proof
is a decomposition of the gradient $u_x=\sum v_i\tilde r_i$ into scalar
components.  In the purely hyperbolic case without viscosity, it is
natural to decompose $u_x$ along a basis $\{r_1,\ldots,r_n\}$
of eigenvectors of the
matrix $A(u)$. Remarkably, this choice does not work here.
Instead, we will decompose $u_x$ as a sum of gradients of 
viscous travelling waves,
selected by a center manifold technique.

As a second step, we study the evolution of each component $v_i$,
which is governed by a scalar conservation law with a source term,
accounting for nonlinear wave interactions.
Uniform  bounds on these source terms are achieved by means of
a {\it transversal interaction} functional, controlling the
interaction between waves of different families, and suitable
{\it swept area} and  {\it curve length} functionals, controlling the 
interaction of waves of the same family. 
All these can be regarded as
``viscous'' counterparts of the {\it wave interaction potential}, 
introduced by Glimm [G] in the purely hyperbolic case. 
Finally, on regions where the
diffusion is dominant, the strength of the source term is
bounded by an {\it energy} functional.
All together, these estimates yield the desired a priori bound on
$\big\|u_x(t,\cdot)\big\|_{\L^1}$,
independent of $t\in [0,\infty[\,$. 

Similar techniques can also be applied to a solution $z=z(t,x)$ of the
variational equation
$$z_t+\big[DA(u)\cdot z\big]u_x+A(u)z_x= z_{xx}\,,\eqno(1.11)$$
which describes the evolution of a first order perturbation to 
a solution $u$ of (1.10). Assuming that 
the total variation of $u$ remains small, we shall establish
an estimate of the form
$$\big\|z(t,\cdot)\big\|_{\L^1}\leq L\,
\big\|z(0,\cdot)\big\|_{\L^1}\qquad\forall t\geq 0\,,\eqno(1.12)$$
valid for all solutions of (1.11).
As soon as this estimate is proved, 
as in [B1] a standard homotopy argument yields the
Lipschitz continuity
of the flow of (1.10) w.r.t.~the initial data, uniformly
in time. 

By the simple rescaling of coordinates $t\mapsto \ve t$, $x\mapsto\ve x$,
all of the above estimates remain valid
for solutions $u^\ve$ of the system (1.9)$_\ve$.  
By a compactness argument, these $BV$ bounds imply the
existence of a strong limit $u^{\ve_m}\to u$ in $\L^1_{\rm loc}$,
at least for some subsequence $\ve_m\to 0$.
In the conservative case where $A=Df$, it is now easy to show that this
limit $u$ provides a weak solution to
the Cauchy problem (1.1)-(1.2).

At this intermediate stage of the analysis, since we are using a 
compactness argument, it is not yet clear whether the vanishing viscosity
limit is unique. In principle, different subsequences $\ve_m\to 0$
may yield different limits.  Toward a uniqueness result, 
in [B3] the second author introduced a definition of 
{\it viscosity solution} for the hyperbolic system of conservation laws
(1.1), based on local integral estimates.
Roughly speaking, a function $u$ is a {\it viscosity solution} if
\v
\i{$\bullet$} In a forward neighborhood of each point of jump,
the function $u$ is well approximated by the self-similar solution
of the corresponding Riemann problem.
\v
\i{$\bullet$} On a region where its total variation is small,
$u$ can be accurately approximated by the solution of a linear system
with constant coefficients.
\v
For a strictly hyperbolic system
of conservation laws satisfying the standard assumptions (H),
the analysis in [B3] proved that the viscosity solution of a 
Cauchy problem is unique, and coincides with the limit of
Glimm and front tracking approximations.
The definition given in [B3] was motivated
by a natural conjecture. Namely, the
viscosity solutions (characterized in terms of local integral estimates)
should coincide precisely 
with the limits of vanishing viscosity approximations.

In the present paper we adopt an entirely similar definition
of viscosity solutions and prove that the above conjecture
is indeed true.
Our results apply to the more general case of 
(possibly non-conservative) quasilinear strictly hyperbolic systems. 
In particular, we obtain the uniqueness of the vanishing viscosity limit.

As in [B3], [BLFP], the underlying idea is that a semigroup is entirely
determined by its local behavior on piecewise constant initial data.
Namely, if two semigroups yield the same solution to each
Riemann problem, then they coincide. 
In our proof of uniqueness, a basic step is thus the analysis
of the vanishing viscosity solution to a general Riemann problem.
The construction given here extends the previous results
by Lax and by Liu to general, non-conservative hyperbolic systems.
As in the cases considered in [Lx], [L1],
for a given left state $u^-$ there
exists a Lipschitz 
continuous curve of right states $u^+$ which can be connected
to $u^-$ by $i$-waves.  These right states are here obtained
by looking at the fixed point of a suitable contractive transformation.
Remarkably, our center manifold plays again a key role, in defining
this transformation.
\v
Our main results are as follows.
\vs
\n{\bf Theorem 1.} {\it 
Consider the Cauchy problem for the hyperbolic
system with viscosity
$$u_t+A(u)u_x=\ve\,u_{xx}\qquad\qquad u(0,x)=\bar u(x)\,.\eqno(1.13)_\ve$$
Assume that the matrices $A(u)$ are strictly hyperbolic, 
smoothly depending on $u$ in a
neighborhood of a compact set $K\subset\R^n$.
Then there exist constants $C,L,L'$ and $\delta>0$ such that the 
following holds.  If
$$\tv\{\bar u\}<\delta\,,\qquad\qquad 
\lim_{x\to -\infty}\bar u(x)\in K\,,
\eqno(1.14)$$
then for each $\ve>0$ the Cauchy problem (1.13)$_\ve$ 
has a unique solution $u^\ve$, defined for all $t\geq 0$.
Adopting a semigroup notation, this will be written as
$t\mapsto
u^\ve(t,\cdot)\doteq S^\ve_t\bar u$.
In addition, one has:
\v
$${\bf BV~ bounds :}\qquad\qquad\qquad\qquad
\tv\big\{S_t^\ve \bar u\big\}\leq C\,\tv\{\bar u\}\,.
\qquad\qquad\qquad\qquad\quad\eqno(1.15)$$
\v
$$\L^1  ~{\bf stability :}\qquad\qquad\qquad\quad
\qquad\big\|S^\ve_t\bar u-S^\ve_t\bar v\big\|_{\L^1}\leq L\,
\big\|\bar u-\bar v\big\|_{\L^1}\,,
\qquad\qquad\qquad\qquad\quad\eqno(1.16)$$
$$\big\|S^\ve_t\bar u-S^\ve_s\bar u\big\|_{\L^1}\leq L'\,
\Big(|t-s|+\big|\sqrt {\ve t}-\sqrt {\ve s}\,\big|\Big)\,.\eqno(1.17)
$$
\v
\n {\bf Convergence:}~
As $\ve\to 0+$, the solutions $u^\ve$ converge to the trajectories of 
a semigroup $S$ such that
$$\big\|S_t\bar u-S_s\bar v\big\|_{\L^1}\leq L\, \|\bar u-\bar v\|_{\L^1}
+L'\,|t-s|\,.\eqno(1.18)$$
These vanishing viscosity limits 
can be regarded as the unique {\bf vanishing viscosity solutions}
of the hyperbolic Cauchy problem
$$u_t+A(u)u_x=0,\qquad\qquad u(0,x)=\bar u(x)\,.\eqno(1.19)$$
\v
In the conservative case $A(u)=Df(u)$, every vanishing viscosity solution
is a weak solution of 
$$u_t+f(u)_x=0,\qquad\qquad u(0,x)=\bar u(x)\,,\eqno(1.20)$$
satisfying the Liu admissibility conditions.
\v
Assuming, in addition, that 
each field is genuinely nonlinear or 
linearly degenerate,
the vanishing viscosity solutions coincide with the unique limits
of Glimm and front tracking approximations.
\i{}
}
\v
Notice that in the above theorem the only key assumptions are
the strict hyperbolicity of the system and the small total variation
of the initial data. 
It is interesting to compare this result with
previous literature.
\v
\n{\bf 1.} Concerning the global existence of 
weak solutions, Glimm's proof requires the additional
assumption (H) of genuine nonlinearity or 
linear degeneracy of each characteristic field.
This assumption has been greatly relaxed in subsequent
works by Liu [L2] and Liu and Yang [LY], but never entirely removed.
The underlying technical reason is the following.
In all papers based on the Glimm scheme (or front tracking),
the construction of approximate solutions as well as the
$BV$ estimates rely on a careful analysis of the Riemann problem.
In this connection,
the hypothesis (H) is a simplifying assumption,
which guarantees that every Riemann problem
can be solved in terms of $n$ elementary waves
(shocks, centered rarefactions or contact discontinuities), 
one for for each characteristic
field $i=1,\ldots,n$.  
At the price of considerable technicalities, this assumption
can be replaced by
some other condition, implying that all
solutions of the Riemann problem can obtained by piecing together
a {\it finite} (but possibly large) number of elementary waves.  

On the other hand, our present approach is based on vanishing 
viscosity limits and does not make any reference to
Riemann problems. Global existence is obtained for the 
whole class of strictly hyperbolic systems.
\v
\n{\bf 2.} Concerning the uniform stability of entropy weak solutions,
the results previously available for $n\times n$ hyperbolic systems 
[BC1], [BCP], [BLY] 
always required the
assumption (H).  For $2\times 2$ systems, this
condition was somewhat relaxed in [AM].  Again, we remark that
the present result
makes no reference to the assumption (H).
\v
\n{\bf 3.} For the viscous system (1.10),
previous results in [L3], [SX], [SZ], [Yu] have
established the stability of special types of solutions,
such as travelling viscous shocks or viscous rarefactions,
w.r.t.~suitably small perturbations.
Taking $\ve=1$,
our present theorem yields at once the 
uniform Lipschitz stability
of {\it all} viscous solutions with sufficiently small total variation,
w.r.t.~the $\L^1$ distance.
\vs
\n{\bf Remark 1.1.} It remains an important open problem to establish
the convergence of vanishing viscosity approximations of the form
$$u_t+A(u)u_x=\ve\big(B(u)u_x\big)_x\eqno(1.21)_\ve$$
for more general viscosity matrices $B$.  In the present paper
we are exclusively concerned with the case where $B$
is the identity matrix. For systems which are not in conservative form,
we expect that the limit of solutions of (1.21)$_\ve$, as $\ve\to 0$,
will be heavily dependent on the choice of the matrix $B$.
\vs
The plan of the paper is as follows.
Section~2 collects those estimates
which can be obtained by standard parabolic techniques.
In particular, we show that the solution of (1.10)
with initial data $\bar u\in BV$ is well 
defined on an initial time interval $[0,\hat t]$ 
where the $\L^\infty$ norms of all derivatives decay rapidly.
Moreover, for large times, as soon as an estimate on the total variation
is available, one immediately obtains a bound on the $\L^1$ norms of
all higher order derivatives.
Our basic strategy for obtaining the BV estimate is outlined in
Section~3.
The decomposition of $u_x$ as a sum of 
gradients of viscous travelling profiles is performed in Section~5.  
This decomposition will depend pointwise on the second order jet
$(u_x,u_{xx})$, involving $2n$ scalar parameters.  
To fit these data, we must first select $n$
smooth families of viscous travelling waves, 
each depending on 2 parameters.  This preliminary
construction is achieved in Section 4, relying on the center
manifold theorem.
In Section 6 we derive the evolution equation satisfied by the
gradient components and analize the form of the various source terms.
As in [G], our point of view is that these 
source terms are the result of interactions between
viscous waves, and can thus be controlled by
suitable {\it interaction functionals}.
In Sections 7 to 9 we introduce various Lyapounov functionals,
which eventually allow us to
estimate the integral of all source terms.
The proof of the uniform BV bounds is then completed in Section 10.

In Section 11 we study the linearized evolution equation (1.11)
for an infinitesimal perturbation $z$, and derive the key estimate (1.12).
In turn, this yields the Lipschitz continuity of the flow,
stated in (1.16). Some of the estimates
here require lengthy calculations, which are postponed
to the Appendices.
Section 12 contains an additional estimate for solutions of 
(1.11), showing that, even in the parabolic case, the bulk of 
a perturbation propagates at a finite speed.
This estimate is crucial because, passing to the limit $\ve\to 0$,
it implies that the values of a vanishing viscosity solution
$u(t,\cdot)$ on an interval $[a,b]$ depend only
on the values of the initial data $u(0,\cdot)$ on 
a bounded interval $[a-\beta t,\,b+\beta t]$.
In Section 13 we study the existence and various properties of
a semigroup obtained as vanishing viscosity limit: $S=\lim S^{\ve_m}$.
At this stage, we only know that the limit exists for a suitable 
subsequence $\ve_m\to 0$.  In the
case of a system of conservation laws satisfying the standard assumptions
(H), we can show that every limit solution
satisfies the {\it Lax shock conditions}
and the {\it tame oscillation property}.  
Hence, by the uniqueness
theorem in [BG], the limit is unique and does not depend on the 
subsequence $\{\ve_m\}$.  This already achieves a proof
of Theorem 1 valid for this special case.

Toward a proof of uniqueness in the general case,
in Section 14 we construct a self-similar solution
$\omega(t,x)=\tilde\omega(x/t)$ to the non-conservative
Riemann problem, and show that it provides the unique 
vanishing viscosity limit.
A definition of {\it viscosity solution} in terms of local integral estimates
is introduced in Section 15.
By a minor modification of the arguments in [B3], [B5] we prove
that these viscosity solutions are unique and coincide with
the trajectories of any semigroup $S=\lim S^{\ve_m}$ obtained
as limit of vanishing viscosity approximations.
Since this result is independent of the subsequence $\{\ve_m\}$,
we obtain the convergence to a unique limit
of the whole family of viscous approximations
$S^\ve_t\bar u\to S_t\bar u$, over all real values of $\ve$.
This completes the proof of Theorem 1.

Finally, in Section 16
we derive two easy estimates. One is concerned with the dependence of the
the limit semigroup $S$ on the coefficients of the matrix $A$
in (1.19).  The other estimate describes
the asymptotic limit of solutions of the
parabolic system (1.10) as $t\to\infty$.
\vsk
\n{\medbf 2 - Parabolic estimates}
\v
In classical textbooks, the local existence and regularity
of solutions to the parabolic system (1.10)
are derived by regarding the hyperbolic term $A(u)u_x$ as a 
first order perturbation of the heat equation.
This leads to the definition of {\it mild solutions},
characterized by the representation
$$u(t)=G(t)*u(0)+\int_0^t G(t-s)*A(u(s)\big)u_x(s)\,ds$$
in terms of convolutions with the standard heat kernel $G$.

In this initial section we collect all the relevant estimates which
can be achieved by this approach.  In particular, we
prove various decay and regularity results for solutions of
(1.10) as well as (1.11).
Given a $BV$ solution $u=u(t,x)$
of (1.10), consider the state
$$u^*\doteq \lim_{x\to -\infty} u(t,x)\,,
\eqno(2.1)$$
which is clearly independent of time.  We then define
the matrix $A^*\doteq A(u^*)$ and let $\lambda_i^*$,
$r_i^*$, $l_i^*$ be the corresponding eigenvalues and
right and left eigenvectors,
normalized as in (1.4).
It will be convenient to  use ``$\bullet$'' 
to denote a directional derivative,
so that $z\bullet A(u)\doteq DA(u)\cdot z$ 
indicates the derivative of the matrix valued 
function $u\mapsto A(u)$ in the direction of the vector $z$.
We can now rewrite the systems (1.10) and (1.11) respectively as
$$\eqalignno{u_t+A^*u_x-u_{xx}&=\big(A^*-A(u)\big)u_x\,,&(2.2)\cr
z_t+A^*z_x-z_{xx}&=\big(A^*-A(u)\big)z_x -\big(z\b A(u)\big)u_x\,.
&(2.3)\cr}$$
In both cases, we
regard the right hand side as a perturbation of the
linear parabolic system with constant coefficients
$$w_t+A^*w_x-w_{xx}=0\,.\eqno(2.4)$$
We denote by $G^*$ the Green kernel for 
(2.4), so that
$$w(t,x)=\int G^*(t,\, x- y)\, w(0,y)\,dy\,.$$
The matrix valued function $G^*$ is easily computed. Indeed,
if $w$ solves (2.4), then its
$i$-th component $w_i\doteq
l^*_i\cdot w$ satisfies the scalar equation
$$w_{i,t}+\lambda_i^*w_{i,x}-w_{i,xx}=0\,.$$
Therefore
$w_i(t)=G_i^*(t)* w_i(0)$,
where
$$G_i^*(t,x)= {1\over 2 \sqrt{\pi t}} \exp \left\{ -
{( x - \lambda_i^*t)^2\over 4t} \right\}.
$$
Looking at the explicit form of its components, 
it is clear that the Green kernel
$G^*=G^*(t,x)$ satisfies the bounds
$$
\big\|G^*(t)\big\|_{\L^1}\leq \kappa\,,\qquad\qquad
\big\| G^*_x(t) \big\|_{\L^1} \leq {\kappa\over\sqrt{t}}\,, \qquad
\qquad
\big\| G^*_{xx}(t) \big\|_{\L^1} \leq {\kappa\over t}\,,
\eqno(2.5)$$
for some constant $\kappa$ and all $t>0$.
It is important to observe that, if $u$ is a solution of (2.2), then
$z=u_x$ is a particular solution 
of the variational equation (2.3). 
Hence all the estimates 
proved for $z_x,\,z_{xx}$ are certainly valid also for the corresponding 
derivatives $u_{xx},\,u_{xxx}$.
Assuming that the initial data $u(0,\cdot)$ has small total variation, we now
derive some estimates on higher derivatives. In particular, we will
show that
\v
\i{$\bullet$} The solution is well defined on some initial interval
$[0,\,\hat t]$, where the $\L^\infty$ norm of all derivatives
decays rapidly.
\v
\i{$\bullet$} As long as the total variation remains small,
the solution can be prolonged in time. In this case, all  
higher order derivatives remain small.  Indeed, waiting a long enough time,
one has
$\big\|u_{xxx}(t)\big\|_{\L^1}<<
\big\|u_{xx}(t)\big\|_{\L^1}<<\big\|u_x(t)\big\|_{\L^1}
=\tv\big\{u(t)\big\}$.
\vs
\n{\bf Proposition 2.1.} {\it
Let $u,z$ be solutions of the systems (2.2)-(2.3), 
satisfying the bounds
$$\big\| u_x(t) \big\|_{\L^1} \leq  \delta_0, \qquad
\qquad \big\|z(t) \big\|_{\L^1} \leq  \delta_0,\eqno(2.6)$$
for some constant $\delta_0<1$ and all $t\in[0,\hat t\,]$,
where 
$$\hat t\doteq \left({1\over 400 \kappa\,\kappa_A\,\delta_0}\right)^2
,\qquad\quad\kappa_A
\doteq \sup_u \big(\|DA\|+\|D^2A\|\big)\eqno(2.7)$$
and $\kappa$ is the constant in (2.5).
Then for $t\in [0,\hat t]$ 
the following estimates hold:}
$$\big\|u_{xx}(t) \big\|_{\L^1}\,,
~\big\|z_x(t) \big\|_{\L^1} \leq {2\kappa\delta_0\over \sqrt  t}\,,
\eqno(2.8)$$
$$\big\|u_{xxx}(t) \big\|_{\L^1}\,,~\big\|z_{xx}(t) \big\|_{\L^1}
\leq {5\kappa^2\delta_0\over t}\,,\eqno(2.9)$$
$$\big\|u_{xxx}(t) \big\|_{\L^\infty}\,,
~\big\|z_{xx}(t) \big\|_{\L^\infty}
\leq {16 \kappa^3 \delta_0\over t \sqrt t}
\,.\eqno(2.10)$$
\v
\n{\bf Proof.} 
The function $z_x$ can be represented as
$$
z_x(t) =G_x^*(t)* z(0) + \int_0^t
G_x^*(t-s) *\Big\{\big( A(u) - A^* \big) z_x(s)- \big( z \bullet A(u) \big) 
u_x(s)\Big\} \,ds\,.\eqno(2.11)$$
Using (2.5) and (2.6) we obtain
$$\eqalign{
&\left\| \int_0^t
G_x^*(t-s) *\Big\{ \big( A^*-A(u)\big) z_x(s)-
\big( z \bullet A(u) \big) u_x(s) \Big\}\, ds \right\|_{L^1} \cr
& \qquad \leq \int_0^t \big\|G^*_x(t-s)\big\|_{\L^1}
\cdot \Big\{ \big\| u_x(s)\big\|_{\L^1}
\big\| DA \big\|_{\L^\infty} \bigl\| z_x(s) \bigr\|_{L^1}
+ \big\|z(s)\big\|_{\L^\infty}\big\| DA \big\|_{\L^\infty} 
\bigl\| u_x (s)\bigr\|_{\L^1} \Big\}\,ds\cr
& \qquad \leq 
2\delta_0\kappa\,\big\| DA \big\|_{\L^\infty} \cdot\int_0^t
{1\over\sqrt{t-s}}\,\big\|z_x(s)\big\|_{\L^1}\,ds\,.\cr}$$
{}Consider first the case of smooth initial data.
We shall argue by contradiction.
Assume that there exists a first time $\tau<\hat t$ such that
the equality in (2.8) holds.
Then, observing that
$$\int_0^t{1\over\sqrt{s(t-s)}}\,ds=\int_0^1{1\over 
\sqrt{\sigma(1-\sigma)}}\,d\sigma
=\pi<4\,$$
we compute
$$\eqalign{\big\|z_x(\tau)\big\|_{\L^1}
&\leq{\kappa\over \sqrt\tau}\delta_0+
2\kappa\delta_0\,\big\| DA \big\|_{\L^\infty} \cdot\int_0^\tau
{1\over\sqrt{\tau-s}}\,{2\delta_0\kappa\over\sqrt s}\,ds\cr
&< {\kappa\delta_0\over \sqrt\tau}
+16\kappa^2\kappa_A\delta_0^2 \big\| DA \big\|_{\L^\infty}~
\leq~{2\kappa\delta_0\over\sqrt\tau}\,,\cr}$$
reaching a contradiction. 
Hence, (2.8) is satisfied as a strict inequality for all
$t\in [0,\hat t]$. 
Observing that this estimate depends only on the $\L^1$ norms of
$u_x$ and $z$, by an approximation argument we obtain the same bound for
general initial data, not necessarily smooth.
Since $z\doteq u_x$ is a particular solution of (2.3), the 
bounds (2.8) certainly apply also to $z_x=u_{xx}$.
\v
A similar technique is used to establish (2.9).
Indeed, we can write
$$\eqalign{z_{xx}(t)&=G_x^*(t/2)*z_x(t/2)-
\int_{t/2}^t G_x^*(t-s)\Big\{\big(z\b A(u)
\big) u_x(s) + \big( A(u) - A^* \big) z_x(s) \Big\}_x ds\,.\cr}
\eqno(2.12)$$
We will prove (2.9)
first in the case $z_{xx}=u_{xxx}$, then in the general case.
If (2.9) is satisfied as an equality at a first time
$\tau<\hat t$, using (2.12) and recalling the definitions (2.7)
we compute
$$\eqalign{
\big\| z_{xx}(\tau) \big\|_{\L^1} &\leq   {\kappa\over \sqrt{\tau/2}} \cdot
{2\kappa\delta_0\over \sqrt{\tau /2}}+ \int_{\tau /2}^\tau 
 {\kappa\over \sqrt{\tau -s}} \cdot \bigg\{ \big\| z_x \bullet A(u) u_x(s)
\big\|_{\L^1} + \big\| z \b \big( u_x \bullet A(u) \big) u_x(s)
\bigr\|_{\L^1} \cr
&~ \qquad + \big\| z \b A(u) u_{xx}(s) \big\|_{\L^1} +
\big\| u_x \bullet A(u) z_x(s) \bigr\|_{\L^1} + \Big\| \big( A(u) -
A^* \big) z_{xx}(s)
\Big\|_{\L^1} \bigg\}\, ds \cr
&\leq   {2\kappa^2\delta_0 \over \tau /2} + \int_{\tau /2}^\tau 
{\kappa\over \sqrt{\tau -s}} \cdot \bigg\{  \delta_0 \| DA
\|_{\L^\infty} \big\| z_{xx}(s) \big\|_{\L^1}
+ \delta_0 \|
D^2 A \|_{\L^\infty} \big\| u_{xx}(s) \big\|_{\L^1}^2 \cr
&~ \qquad +  \delta_0 \| DA
\|_{\L^\infty} \big\| u_{xxx}(s) \big\|_{\L^1}
+  \delta_0 \| DA \|_{\L^\infty} \big\| z_{xx}(s)
\big\|_{\L^1} +  \delta_0 \| DA \|_{\L^\infty} \bigl\|
z_{xx}(s) \big\|_{\L^1} \bigg\} ds \cr
&\leq {4\kappa^2\delta_0 \over \tau}
+\kappa\delta_0\big(4\kappa^2\delta_0^2\|D^2A\|_{\L^\infty}
+20\kappa^2\delta_0\|DA\|_{\L^\infty}\big)
\int_{\tau /2}^\tau  
{1\over s\sqrt {\tau -s}}\,ds\,,\cr
&< {4\kappa^2\delta_0 \over \tau} +
20\kappa^3\kappa_A\delta_0^2\cdot {4\over\sqrt{\tau/2}}~<~
{5\kappa^2\delta_0 \over \tau} \,,\cr}
$$
reaching a contradiction.   
\v
Finally, using (2.12) and (2.8)-(2.9), 
the bounds in (2.10) are proved by the estimate
$$\eqalign{
\big\| z_{xx}(\tau) \big\|_{L^\infty} &\leq  {\kappa\over \sqrt{\tau/2}} \cdot
{5\kappa^2 \delta_0\over  \tau/2} + \int_{\tau/2}^\tau
 {\kappa\over \sqrt{\tau-s}} \cdot \bigg\{ \big\| z_x \bullet A(u) u_x(s)
\big\|_{L^\infty} + \big\| z \bullet ( u_x \bullet A(u) ) u_x(s)
\big\|_{L^\infty} \cr
&~ \qquad \qquad + \big\| z \bullet A(u) u_{xx}(s) \big\|_{L^\infty} +
\big\| u_x \bullet A(u) z_x(s) \big\|_{L^\infty} + \big\| ( A(u) -
A_0 ) z_{xx}(s)
\big\|_{L^\infty} \bigg\} \, ds \cr
&\leq  {10 \sqrt 2 \kappa^3 \delta_0\over \tau \sqrt \tau} +
\big(8\kappa^4\delta_0^3 \|D^2A\|+46\kappa^4\delta_0^2\|DA\|\big)
\int_{\tau/2}^\tau {1\over s^{3/2}\sqrt{\tau-s}}\,ds\cr
&\leq  
{15 \kappa^3 \delta_0\over \tau \sqrt \tau} +
46\kappa^4\kappa_A\delta_0^2\cdot {4\over \tau/2}~<~
{16 \kappa^3 \delta_0\over \tau \sqrt \tau}\,.
\cr}$$
\endproof
\v
\n{\bf Corollary 2.2.} {\it In the same setting as Proposition 5.1,
assume that the bounds (2.6) hold on a larger interval
$[0,T]$.  Then for all $t\in [\hat t,\, T]$ there holds}
$$\eqalignno{\big\|u_{xx}(t) \big\|_{\L^1}\,,~\big\|u_x(t) 
\big\|_{\L^\infty}\,,
~\big\|z_x(t) \big\|_{\L^1} ~&=~
\O(1)\cdot \delta^2_0\,,
&(2.13)\cr
&&\cr
\big\|u_{xxx}(t) \big\|_{\L^1}\,,~\big\|u_{xx}(t) \big\|_{\L^\infty}\,,
~\big\|z_{xx}(t) \big\|_{\L^1}
~&=~
\O(1)\cdot \delta^3_0\,,&(2.14)\cr
&&\cr\big\|u_{xxx}(t) \big\|_{\L^\infty}\,,
~\big\|z_{xx}(t) \big\|_{\L^\infty}
~&=~
\O(1)\cdot \delta^4_0
\,.&(2.15)\cr}$$
\v
\n{\bf Proof.} It suffices to apply Proposition 2.1
on the interval $[t-\hat t,~t]$.
\endproof
\v 
\n{\bf Proposition 2.3.}
{\it 
Let $u=u(t,x)$, $z=z(t,x)$ be solutions of (2.2), (2.3) respectively,
such that
$$\tv\big\{u(0,\cdot)\big\}\leq {\delta_0\over 4\kappa}\,,\qquad\qquad
\big\|z(0)\big\|_{\L^1}\leq {\delta_0\over 4\kappa}\,.\eqno(2.16)$$
Then $u,z$ are well defined on the whole interval
$[0,\hat t]$ in (2.7), and satisfy}
$$\big\|u_x(t)\big\|_{\L^1}\leq {\delta_0\over 2}\,,
\qquad\qquad \big\|z(t)\big\|_{\L^1}\leq {\delta_0\over 2}\qquad\qquad 
t\in [0,\,\hat t]\,.
\eqno(2.17)$$
\v
\n{\bf Proof.} We have the identity
$$z(t)=G^*(t)z(0)-\int_0^t G^*(t-s)
\Big( z \bullet A(u) u_x(s) - 
\big( A(u) - A_0 \big) z_x(s) \Big) ds\,.\eqno(2.18)$$
As before, we first establish the result for $z=u_x$, then for a
general solution $z$ of (2.3).
Assume that there exists a first time
$\tau<\hat t$ where the
bound in (2.17) is satisfied as an equality. 
Estimating the right hand side of (2.18) by means of
(2.5) and (2.8), we obtain
$$\eqalign{\big\|z(\tau)\big\|_{\L^1}&\leq 
{\kappa\delta_0\over 4\kappa}+\int_0^\tau
{2\kappa\delta_0^2\over\sqrt s}\,\|DA\|_{\L^\infty}\,ds\cr
&\leq 
{\delta_0\over 4}+4\kappa\,\kappa_A\delta_0^2\sqrt\tau ~<~{\delta_0\over 2}
\,,\cr}$$
reaching a contradiction.
\endproof
\v
To simplify the proofs, in all previous
results we used the same hypotheses on the functions $u_x$ and $z$.
However, observing that $z$ solves a linear homogeneous
equation, similar estimates can be immediately derived
without any restriction on the initial size $\big\|z(0)\big\|_{\L^1}$.
In particular, from Proposition 2.3 it follows
\v
\n{\bf Corollary 2.4.} 
{\it 
Let $u=u(t,x)$, $z=z(t,x)$ be solutions of (2.2), (2.3) respectively,
such that
$\big\|u_x(0)\big\|_{\L^1}\leq \delta_0/ 4\kappa$.
Then $u,z$ are well defined on the whole interval
$[0,\hat t]$ in (2.7), and satisfy}
$$\big\|u_x(t)\big\|_{\L^1}\leq 2\kappa \big\|u_x(0)\big\|_{\L^1}\,,
\qquad\qquad \big\|z(t)\big\|_{\L^1}\leq 
2\kappa \big\|z(0)\big\|_{\L^1}\,,\qquad\qquad 
t\in [0,\,\hat t]\,.
\eqno(2.19)$$
\midinsert
\vskip 10pt
\centerline{\hbox{\psfig{figure=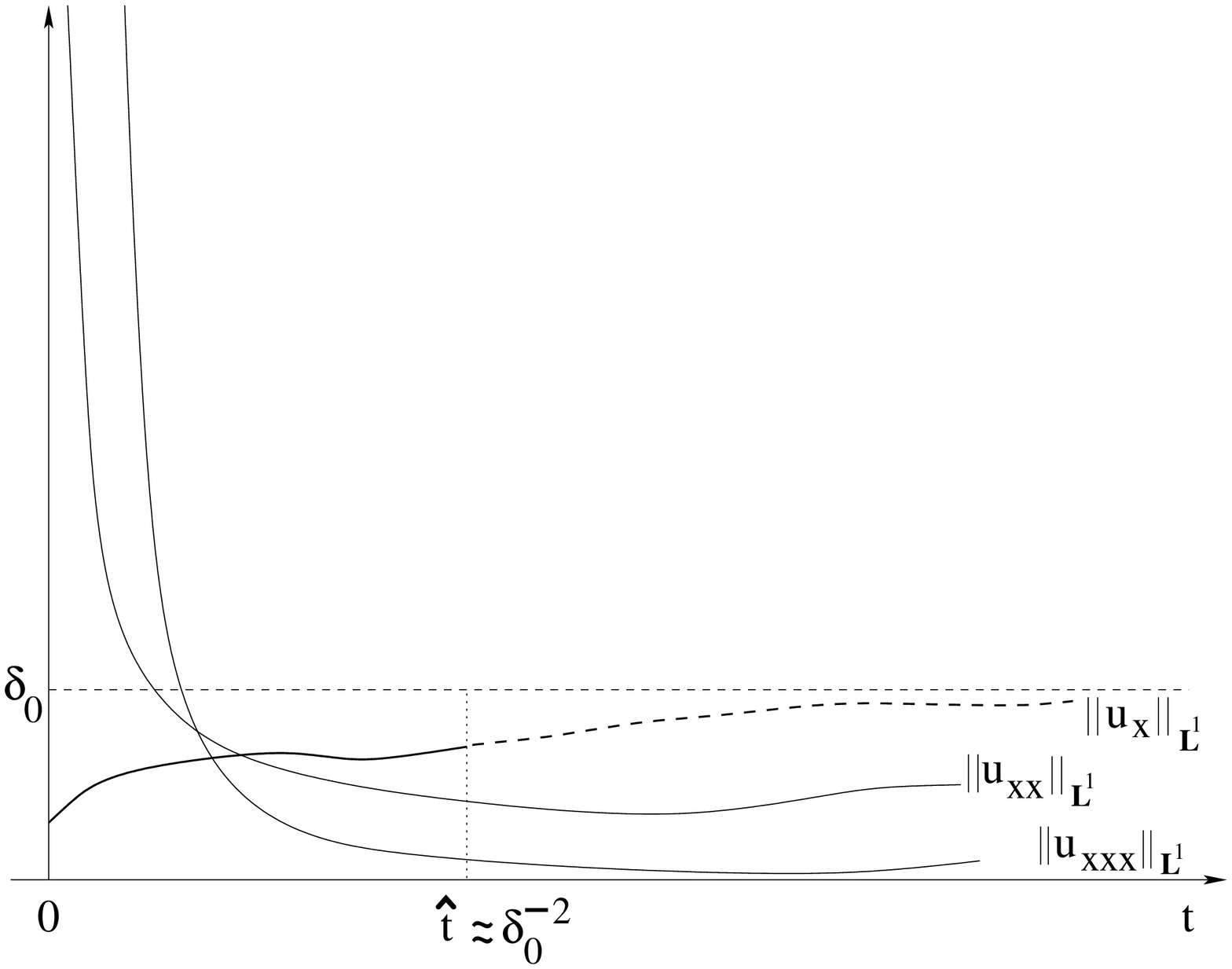,width=12cm}}}
\centerline{\hbox{figure 1}}
\vskip 10pt
\endinsert

A summary of the main estimates is illustrated in
fig.~1.   On the initial interval $t\in [0,\,\hat t]$, with $\hat t\approx
1/\delta_0^2$
we have 
$$ \big\|u_x(t)\big\|_{\L^1}\leq \delta_0,\eqno(2.20)
$$
while the norms of the higher derivatives decay:
$$\big\|u_{xx}\big\|_{\L^1}=\O(1)\cdot {\delta_0/ \sqrt t}\,,
\qquad\qquad
\big\|u_{xxx}\big\|_{\L^1}=\O(1)\cdot {\delta_0/ t}\,.$$
On the other hand, for $t\geq \hat t$, as long as (2.20) remains valid
we also have
$$\big\|u_{xx}\big\|_{\L^1}=\O(1)\cdot \delta^2_0\,,
\qquad\qquad
\big\|u_{xxx}\big\|_{\L^1}=\O(1)\cdot \delta_0^3\,.$$
These bounds (the solid lines in fig.~1) 
were obtained in the present section by standard parabolic-type estimates.
The most difficult part of the proof is to obtain the estimate
(2.20) for large times $t\in [\hat t,\,\infty[\,$ 
(the broken line in fig.~1).
This will require hyperbolic-type estimates, based on the local decomposition
of the gradient $u_x$ as a sum of travelling waves,
and on a careful analysis of all interaction terms.
\vsk
\n{\medbf 3 - Outline of the BV estimates}
\v
It is our aim to derive global a priori bounds on the total variation
of solutions of
$$u_t+A(u)u_x=u_{xx}\eqno(3.1)$$
for small initial data.
We always assume that the system is strictly hyperbolic, 
so that each matrix $A(u)$ has 
real distinct eigenvalues $\lambda_i(u)$ as in (1.3), and
and dual bases of right and left eigenvectors 
$r_i(u)$, ~$l_i(u)$~normalized as in (1.4).
The directional derivative of a function $\phi=\phi(u)$
in the direction of the vector $\bfv$ is written
$$\bfv\b \phi(u)\doteq D\phi\cdot\bfv=
\lim_{\epsilon\to 0} {\phi( u+\epsilon\bfv
)-
\phi(u)\over \epsilon}\,,\eqno(3.2)$$
while 
$$[r_j,r_k]\doteq r_j\b r_k-r_k\b r_j$$ denotes a Lie bracket. 
In order to obtain uniform bounds on $\tv\big\{ u(t,\cdot)\big\}$
for all $t>0$,
our basic strategy is as follows.
We choose $\delta_0>0$ sufficiently small and consider an
initial data $u(0,\cdot)=\bar u$ satisfying the first inequality in
(2.16).  By Proposition 2.3, 
the corresponding solution is well defined
on the initial time interval $[0,\hat t]$ and its total variation
remains bounded, according to (2.17). 
The main task is to establish $BV$ estimates
on the remaining interval $[\hat t,\,\infty[\,$.
For this purpose,
we decompose the gradient $u_x$ along a suitable
basis of unit vectors $\tilde r_1,\ldots,\tilde r_n$, say
$$u_x=\sum_{i=1}^n v_i \tilde r_i\,.\eqno(3.3)$$
Differentiating (3.1), we obtain a system of $n$ evolution equations
for these scalar components
$$v_{i,t}+(\tilde\lambda_i v_i)_x-v_{i,xx}=\phi_i
\qquad\qquad i=1,\ldots,n\,.
\eqno(3.4)$$
Since the left hand side is in conservation form, (3.4) implies
$$\big\|v_i(t,\cdot)\big\|_{\L^1}\leq 
\big\|v_i(\hat t,\cdot)\big\|_{\L^1}+\int_{\hat t}^\infty
\!\int \big|\phi_i(t,x)\big|
\,dxdt
\eqno(3.5)$$
for all $t\geq \hat t$.
By (3.3) it follows
$$\tv\big\{ u(t,\cdot)\big\}=\big\|u_x(t,\cdot)\big\|_{\L^1}\leq
\sum_i \big\|v_i(t,\cdot)\big\|_{\L^1} \,.\eqno(3.6)$$
In order to obtain a uniform bound on the total variation,
the key step is thus to construct the basis of unit vectors 
$\{\tilde r_1,\ldots,\tilde r_n\}$ in (3.3) in a clever way, 
so that
the functions $\phi_i$ on the right hand side of (3.4) become
integrable on the half plane $\{t>\hat t,~x\in\R\}$.

As a preliminary, we observe that the choice
$\tilde r_i\doteq r_i(u)$, the $i$-th eigenvector of the matrix
$A(u)$, seems quite natural.  
This choice was indeed adopted in [BiB1], where
the authors proved Theorem 1 restricted to
the special class of systems where
all Rankine-Hugoniot curves are straight lines.
Unfortunately, for general $n\times n$ hyperbolic
systems it does not work.
To understand why, let us write
$$u^i_x\doteq l_i(u)\cdot u_x\eqno(3.7)$$
for the $i$-th component of $u_x$ in this basis of 
eigenvectors. 
As shown in [BiB1], these components satisfy the system of
evolution equations
$$\eqalign{(u_x^i)_t+&(\lambda_i u_x^i)_x-(u^i_x)_{xx}\cr
&=
l_i\cdot\bigg\{\sum_{j\not= k} \lambda_j[r_j,r_k]u_x^j u_x^k+
2\sum_{j,k} (r_k\b r_j)(u_x^j)_x u_x^k +\sum_{j,k,\ell}[r_\ell,~r_k\b r_j]
u_x^j u_x^k u_x^\ell\bigg\}\cr 
&\doteq \phi_i\,.\cr}\eqno(3.8)$$
Assume that the $i$-th characteristic field is genuinely nonlinear,
with shock and rarefaction curves not coinciding, and
consider a travelling wave solution $u(t,x)=U(x-\lambda t)$,
representing a viscous $i$-shock.  It is then easy to see that the
right hand side of (3.8) is not identically zero. 
Since it corresponds to a travelling wave,
the integral
$$\int\big|\phi_i(t,x)\big|\,dx\,\not=\, 0$$
is constant in time.  Hence $\phi_i$ is certainly not integrable
over the half plane  $\{t>\hat t,~x\in\R\}$.
\v
\midinsert
\vskip 10pt
\centerline{\hbox{\psfig{figure=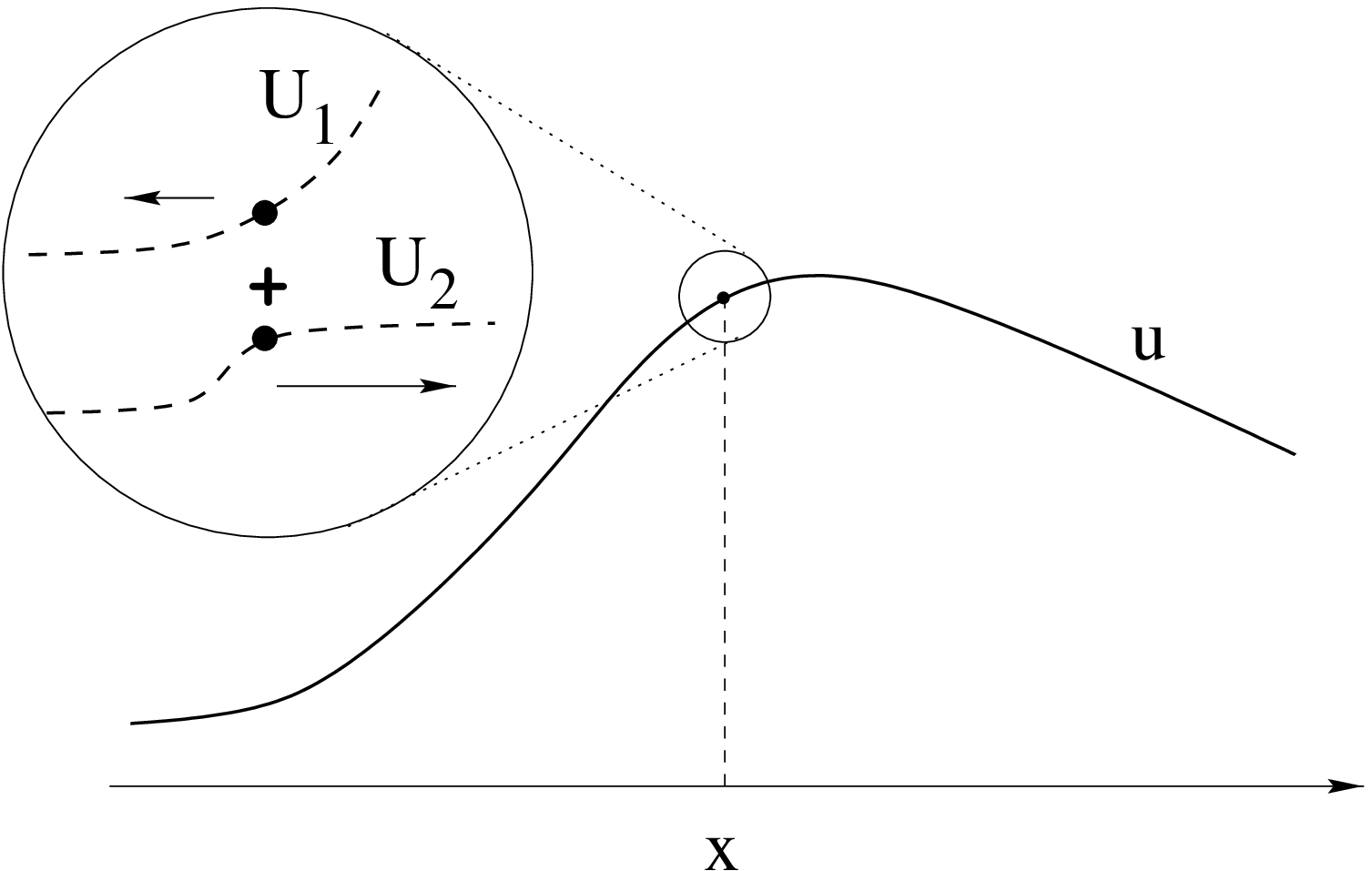,width=8cm}}}
\centerline{\hbox{figure 2}}
\vskip 10pt
\endinsert
The previous example clearly points out a basic requirement
for our decomposition (3.3). Namely, in connection with
a viscous travelling wave, the source terms $\phi_i$
in (3.4) should vanish identically. 
To achieve this goal, we shall seek a decomposition of
$u_x$ not along eigenvectors of the matrix $A(u)$, but as a {\bf sum
of gradients of viscous travelling waves.}
More precisely, consider a smooth function $u:\R\mapsto\R^n$.
At each point $x$, depending on the second order jet
$(u,u_x,u_{xx})$, 
we shall uniquely determine $n$ travelling
waves $U_1,\ldots,U_n$ passing through $u(x)$ (fig.~2). We then write
$u_x$ in the form (3.3), as the sum of the gradients of these waves.
As a guideline, we shall try to achieve the following
relations:
$$U_i(x)=u(x)\qquad\qquad i=1,\ldots,n,\eqno(3.9)$$
$$\sum_i U'_i(x)=u_x(x)\,,\qquad\qquad \sum_i U''_i(x)=u_{xx}(x)\,.
\eqno(3.10)$$
Details of this construction will be worked out 
in the next two sections.
\v
\n{\bf Remark 3.1.} For each $i\in\{1,\ldots,n\}$, one can find
an $(n+1)$-parameter family of viscous travelling waves $U_i$ passing through
a given state $u\in\R^n$.  Indeed, one can assign the speed 
$\sigma_i\approx \lambda_i^*$ and the first derivative
$U_i'\in\R^n$ arbitrarily, and then solve the second order O.D.E.
$$U_i''=\big(A(U_i)-\sigma_i\big)U_i'.$$
In all, this would give us $n(n+1)$ scalar parameters to determine.
Far too many, since
(3.10) is a system of only $2n$ equations.
In an attempt to fix this problem, we could
restrict ourselves only to {\it globally bounded} travelling wave profiles.
Assuming that the $i$-th field is genuinely nonlinear,
the analysis in [F] or [Se] shows that there exists
a 2-parameter family of viscous $i$-shock profiles trough any given point $u$.
Indeed, one can arbitrarily assign $\alpha_-,\alpha_+\leq 0$ and find 
unique asymptotic states
$u^+,u^-\in\R^n$ which are connected by a 
viscous shock profile $U_i(\cdot)$ passing through
$u$ (fig.~3), and such that the following inner products take the 
prescribed values:
$$\la r_i^*,~u-u^-\ra=\alpha_-\,,\qquad\qquad
\la r_i^*,~u^+-u\ra=\alpha_+\,.$$
In this case, summing over $n$ families,
we would end up with the right number of parameters
to fit the data, namely $2n$.  Unfortunately, viscous shock profiles
yield only negative values of gradient components:
$\la r_i^*,~U_i'\ra < 0$.  At a point $x$ where
$\la r_i^*,\,u_x\ra>0$, to achieve the decomposition 
(3.10) we have to consider also viscous
rarefaction profiles (which are not globally bounded).
In the next section, a suitable family of viscous travelling waves 
will be selected by the center manifold theorem.
\midinsert
\vskip 10pt
\centerline{\hbox{\psfig{figure=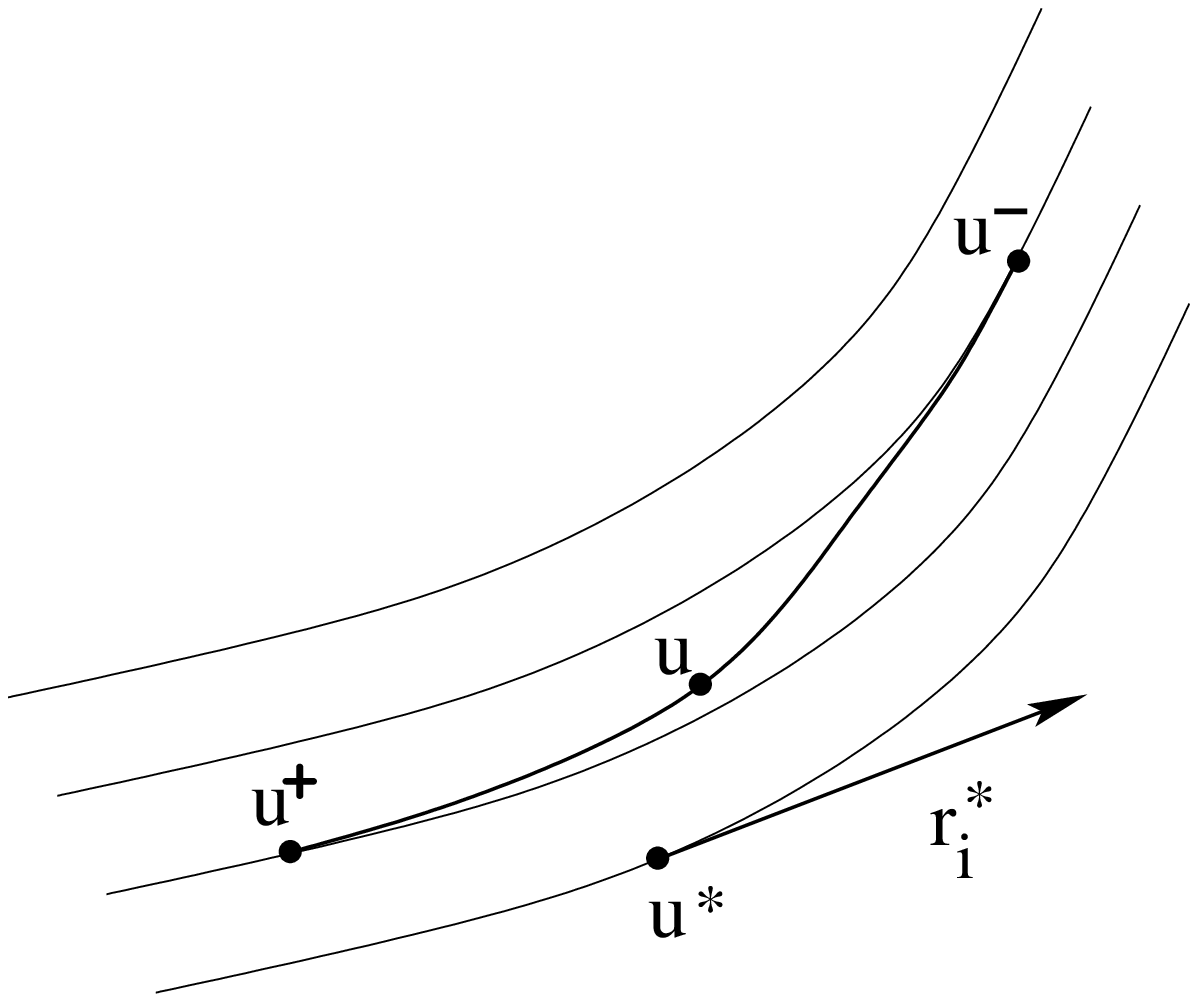,width=6cm}}}
\centerline{\hbox{figure 3}}
\endinsert

\vfill\eject
\n{\medbf 4 - A center manifold of viscous travelling waves}
\vs
To carry out our program, we must first select certain families of 
travelling waves,
depending on the correct number of parameters to fit the data. 
Given a state $u\in\R^n$, a second order jet $(u_x,\,u_{xx})$ determines 
$2n$ scalar parameters.  In order to
uniquely satisfy the equations (3.10), we thus need to 
construct $n$ families of travelling wave profiles through $u$, 
each depending on two scalar parameters. This will be 
achieved by an application of the center manifold theorem.

Travelling waves for the viscous hyperbolic system (3.1) correspond to
(possibly unbounded) solutions of 
$$\big(A(U)-\sigma\big) U'=U''.\eqno(4.1)$$ 
We write (4.1) as a first order system on the space $\R^n\times\R^n\times\R$:
$$\left\{\eqalign{\dot u&=v\,,\cr
\dot v&=\big(A(u)-\sigma\big)v\,,\cr
\dot \sigma&=0\,.\cr}\right.\eqno(4.2)$$
Let a state $u^*$ be given and fix an index $i\in \{1,\ldots,n\}$.
Linearizing (4.2) at the equilibrium point
$P^*\doteq \big(u^*,\,0,\,\lambda_i(u^*)\big)$
we obtain the linear system
$$\left\{\eqalign{\dot u&=v\,,\cr
\dot v&=\big(A(u^*)-\lambda_i(u^*)\big)v\,,\cr
\dot \sigma&=0\,.\cr}\right.\eqno(4.3)$$

\midinsert
\vskip 10pt
\centerline{\hbox{\psfig{figure=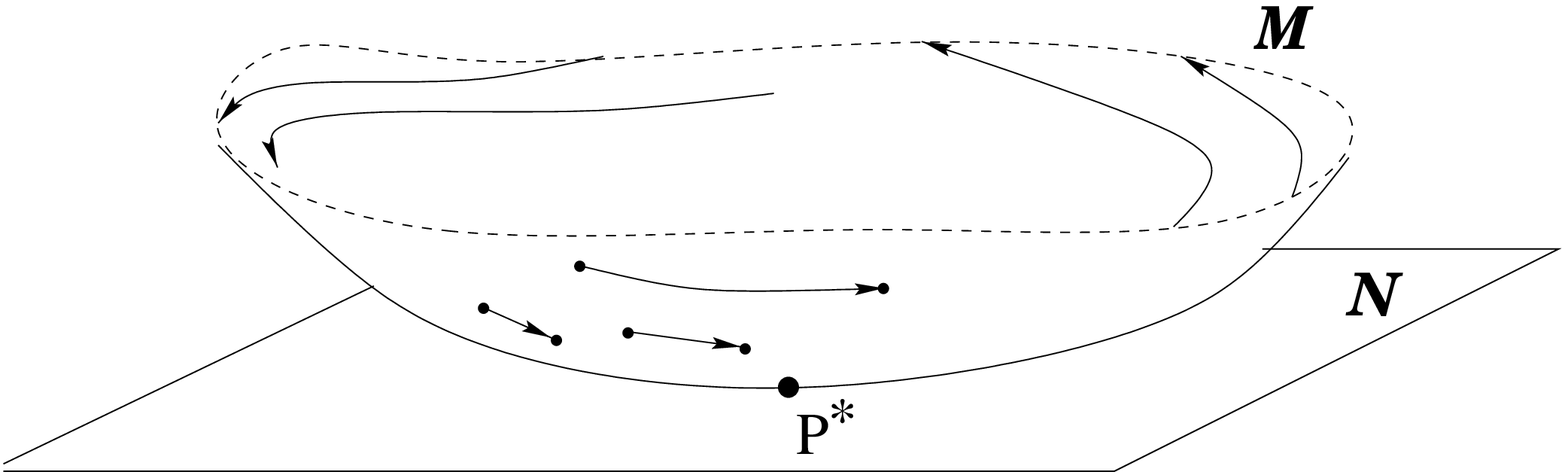,width=9cm}}}
\centerline{\hbox{figure 4}}
\vskip 10pt
\endinsert

\n Let $\{r_1^*,\ldots,r_n^*\}$ and $\{l_1^*,\ldots,l_n^*\}$ be
dual bases of right and left eigenvectors of $A(u^*)$ normalized as in (1.4).
We call $(V_1,\ldots,V_n)$ the coordinates of a vector $v\in\R^n$ 
w.r.t.~this basis, so that
$$|r_i^*|=1,\qquad\qquad
v=\sum_j V_jr^*_j,\qquad\qquad V_j\doteq l^*_j\cdot v\,.$$
The null space ${\cal N}$ for (4.3) consists of all vectors $(u,v,\sigma)\in 
\R^n\times\R^n\times\R$ such that
$$V_j=0\qquad\qquad \hbox{for all}~~j\not= i,\eqno(4.4)$$
and therefore has dimension $n+2$.
By the center manifold theorem [V], there exists a smooth manifold
$\M\subset\R^{n+n+1}$, tangent
to ${\cal N}$ at the stationary
point $P^*$ (fig.~4), which is locally invariant under the flow of (4.2).
This manifold has dimension $n+2$ and can be locally
defined by the $n-1$ equations
$$V_j=\varphi_j(u,V_i,\sigma)\qquad\qquad\qquad j\not= i.\eqno(4.5)$$
We can assume that the $n-1$ smooth scalar functions $\varphi_j$
are defined on the domain
$$\D\doteq\Big\{ |u-u^*|<\epsilon,\qquad |V_i|<\epsilon,\qquad \big|\sigma
-\lambda_i(u^*)|<\epsilon\Big\}.$$
Moreover, the tangency condition implies
$$\varphi_j(u,V_i,\sigma)=\O(1)\cdot \Big(|u-u^*|^2+|V_i|^2+\big|\sigma
-\lambda_i(u^*)\big|^2\Big).\eqno(4.6)$$
\v
We now take a closer look at the flow on this center manifold.
By construction, every trajectory 
$$t\mapsto P(t)\doteq\big(u(t),\,v(t),\,\sigma(t)\big)$$
of (4.2),
which remains within a small neighborhood of the point
$P^*\doteq \big(u^*,0,\lambda_i(u^*)\big)$ for all $t\in\R$,
must lie entirely on the manifold $\M$.
In particular, $\M$ contains all viscous $i$-shock profiles joining
a pair of states $u^-,u^+$ sufficiently close to $u^*$.
Moreover, all equilibrium points
$(u,0,\sigma)$ with $|u-u^*|<\epsilon$ and 
$\big|\sigma-\lambda_i(u^*)\big|<\epsilon$ must lie on $\M$.
Hence
$$\varphi_j(u,0,\sigma)=0\qquad\qquad \forall j\not=i.\eqno(4.7)$$
By (4.7) and the smoothness of the functions $\varphi_j$, we can 
``factor out'' the component $V_i$ and write
$$\varphi_j(u,V_i,\sigma)=\psi_j(u,V_i,\sigma)\cdot V_i,$$
for suitable smooth funtions $\psi_j$. From (4.6) it follows
$$\psi_j\to 0\qquad\hbox{as}\qquad (u,V_i,\sigma)\to
\big(u^*,0,\lambda_i(u^*)\big).\eqno(4.8)$$
On the manifold $\M$ we thus have
$$v=\sum_k V_k r_k^*=V_i\cdot \left(r_i^*+\sum_{j\not=i}
\psi_j(u,V_i,\sigma)\,r_j^*\right)
\doteq V_i \,r_i^\sharp(u,V_i,\sigma).\eqno(4.9)$$
By (4.8), the function $r^\sharp$ defined by the last equality in (4.9)
satisfies
$$r_i^\sharp(u,V_i,\sigma)\to r_i^*\qquad \hbox{as}\qquad (u,V_i,\sigma)
\to (u^*,0,\lambda_i(u^*)). \eqno(4.10)$$
\v
\n{\bf Remark 4.1.} Trajectories on the center manifold correspond to 
the profiles
of viscous travelling
$i$-waves.  
We thus expect that the derivative $\dot u=v$ should be a vector
``almost parallel'' to the eigenvector $r_i^*\doteq 
r_i(u^*)$.  This is indeed confirmed by (4.10).
\v
We can now define the new variable
$$v_i=v_i(u,V_i,\sigma)\doteq V_i\cdot \big|r_i^\sharp(u,V_i,\sigma)\big|.
\eqno(4.11)$$
As $(u,V_i,\sigma)$ range in a small neighborhood of $(u^*,0,\lambda_i(u^*))$,
by (4.10) the vector $r_i^\sharp$ remains close to the eigenvector
$r_i^*$. In particular, its norm remains uniformly positive.
Therefore, the transformation $V_i\longleftrightarrow v_i$ is
invertible and smooth.
We can thus reparametrize the center manifold $\M$ in terms of the variables
$(u,v_i,\sigma)\in\R^n\times\R\times\R$.
Moreover, we define the unit vector
$$\tilde r_i(u,v_i,\sigma)\doteq {r_i^\sharp\over |r_i^\sharp|}\,.
\eqno(4.12)$$
Observe that  $\tilde r_i$ is also a smooth function of its arguments. 
With the above definitions, instead of (4.5)
we can write the manifold $\M$ in terms of the equation
$$v=v_i \tilde r_i.\eqno(4.13)$$
\v
The above construction of a center manifold can be repeated
for every $i=1,\ldots,n$.   We thus obtain $n$ center manifolds
$\M_i\subset \R^{2n+1}$ 
and vector functions $\tilde r_i=\tilde r_i(u,v_i,\sigma_i)$ 
such that
$$|\tilde r_i|\equiv 1,\eqno(4.14)$$
$$\M_i=\big\{(u,v,\sigma_i)~;~~v=v_i\,\tilde r_i(u,v_i,\sigma_i)\big\}\,,
\eqno(4.15)$$
as $(u,v_i,\sigma_i)\in\R^n\times\R\times\R$ ranges in a neighborhood
of $\big(u^*,\,0,\,\lambda_i(u^*)\big)$.
\vs
We derive here some useful identities, for later use. The partial
derivatives of $\tilde r_i=\tr_i(u,v_i,\sigma_i)$ w.r.t.~its arguments 
will be written as
$$\tr_{i,u}\doteq{\partial\over\partial u}\tr_i\,,
\qquad\quad
\tr_{i,v}\doteq{\partial\over\partial v_i}\tr_i\,,
\qquad\quad
\tr_{i,\sigma}\doteq{\partial\over\partial \sigma_i}\tr_i\,.
$$
Clearly, $\tr_{i,u}$ is an $n\times n$ matrix, while $\tr_{i,v}$,
$\tr_{i,\sigma}$ are $n$-vectors.  Higher order derivatives
are denoted as $\tr_{i,u\sigma}$, $\tr_{i,\sigma\sigma}\ldots$~ 
We claim that
$$\tr_i(u,0,\sigma_i)=r_i(u)\qquad\qquad\forall u,\sigma_i\,.\eqno(4.16)$$
Indeed, consider again the equations for a viscous travelling 
$i$-wave:
$$u_{xx}=\big(A(u)-\sigma_i\big)u_x.\eqno(4.17)$$
For a solution contained in the center manifold, taking the derivative 
w.r.t.~$x$ of
$$u_x=v=v_i\tr_i(u,v_i,\sigma_i)\eqno(4.18)$$
and using (4.17) we obtain
$$v_{i,x}\tr_i+v_i\tr_{i,x}=\big(A(u)-\sigma_i\big)v_i\tr_i.\eqno(4.19)$$
Since $|\tr_i|\equiv 1$, the vector $\tr_i$ is perpendicular to
all of its  derivatives. Taking the inner product of (4.19) with $\tr_i$
we thus obtain
$$v_{i,x}=(\tilde\lambda_i-\sigma_i)v_i\,,\eqno(4.20)$$
where we defined the speed
$\tilde\lambda_i=\tilde\lambda_i(u,v_i,\sigma_i)$ as the inner product
$$\tla_i\doteq \la \tr_i\,,~A(u)\tr_i\ra\,.\eqno(4.21)$$
Using (4.20) in (4.19) and dividing by $v_i$ we finally obtain
$$(\tilde\lambda_i-\sigma_i)v_i\tr_i+v_i\big(\tr_{i,u}\tr_iv_i+\tr_{i,v}
(\tilde\lambda_i-\sigma_i)v_i\big)=\big(A(u)-\sigma_i\big)v_i\tr_i\,,
\eqno(4.22)$$
$$
v_i\big(\tr_{i,u}\tr_i+\tr_{i,v}
(\tilde\lambda_i-\sigma_i)\big)=\big(A(u)-\tla_i\big)\tr_i\,.\eqno(4.23)$$
By (4.23), as $v_i\to 0$, the unit vector $\tr_i(u,v_i,\sigma_i)$
approaches an eigenvector of the matrix $A(u)$, while $\tla_i$
approaches the corresponding eigenvalue.  
By continuity, this establishes (4.16).

In turn, by the smoothness of the vector field $\tr_i$ we also have
$$
\eqalign{
\tr_i(u,v_i,\sigma_i)-r_i(u)&=\O(1)\cdot v_i,\cr
\tr_{i,u\sigma}&=\O(1)\cdot v_i,\cr}\qquad\qquad
\eqalign{\tr_{i,\sigma}&=\O(1)\cdot v_i,\cr
\tr_{i,\sigma\sigma}&=\O(1)\cdot v_i.\cr}\eqno(4.24)$$
Observing that the vectors $\tr_{i,v}$
and $\tr_{i,\sigma}$ are both perpendicular to $\tr_i$, from
(4.24) we deduce
$$\big|\tilde\lambda_i(u,v_i,\sigma_i)-\lambda_i(u)\big|=\O(1)\cdot
v_i,\qquad \tilde\lambda_{i,v}=\O(1)\cdot v_i,\qquad
\tilde\lambda_{i,\sigma}=\O(1)\cdot v_i^2.\eqno(4.25)$$

A further identity will be of use.
Differentiating (4.19) one finds
$$v_{i,xx}\tr_i+2v_{i,x}\tr_{i,x}+v_i\tr_{i,xx}=
\big(A(u)v_i\tr_i\big)_x-\sigma_i v_{i,x}\tr_i-\sigma_i v_i\tr_{i,x}\,.
\eqno(4.26)$$
{}From the identities
$$\la \tr_i,\,\tr_{i,x}\ra=0\,,\qquad\qquad\la \tr_i,\,\tr_{i,xx}\ra=
-\la \tr_{i,x},\,\tr_{i,x}\ra,$$
taking the inner product of (4.19) with $\tr_{i,x}$
we obtain
$$\la\tr_i,\,\tr_{i,xx}\ra 
v_i=-\la \tr_{i,x},\,A(u)\tr_i\ra v_i\,.\eqno(4.27)$$
Taking now the inner product of (4.26) with $\tr_i$ we find
$$v_{i,xx}+\la \tr_i,\tr_{i,xx}\ra v_i=
\la \tr_i,~(A(u)\tr_i v_i)_x\ra -\sigma_i v_{i,x}\,.$$
Since $v_{i,t}+\sigma_i v_{i,x}=0$, using the identity (4.27) we conclude
$$v_{i,t}+(\tilde\lambda_iv_i)_x-v_{i,xx}=0,\eqno(4.28)$$
where $\tilde\lambda_i$ is the speed at (4.21). 
\v
\n{\bf Remark 4.2.} It is important to appreciate the difference between
the identities
$$\big(A(u)-\lambda_i\big) r_i=0\,,\qquad\qquad \big(A(u)-\tla_i\big)\tr_i
=v_i\big(\tr_{i,u}\tr_i+\tr_{i,v}
(\tilde\lambda_i-\sigma_i)\big),\eqno(4.29)$$
satisfied respectively by an eigenvector $r_i$ and by 
a unit vector $\tr_i$ parallel to the gradient of a travelling wave. 
Decomposing $u_x$ along the eigenvectors $r_i$ one obtains 
the evolution equations (3.8), 
with non-integrable source terms on the right hand side.
When a similar computation is performed in connection with the
vectors $\tr_i$, thanks to the presence of  
the additional terms on the right hand side in (4.29)
a crucial cancellation is achieved. 
In this case, we will show that
the source terms $\phi_i$ in (3.4) are integrable
over the half plane $x\in\R$, $t>\hat t$.
\vfill\eject
\n{\medbf 5 - Gradient decomposition}
\v
Let $u:\R\mapsto\R^n$ be a smooth function with small total variation.
At each point $x$, we seek a decomposition of the gradient $u_x$ in the form
(3.3), where $\tilde r_i=\tilde r_i(u,v_i,\sigma_i)$ are the
vectors defining the center manifold in (4.15). 
To uniquely determine the $\tilde r_i$, 
we should first define the
wave strengths $v_i$ and speeds $\sigma_i$ in terms
of $u$, $u_x$, $u_{xx}$.  

Consider first the special case where $u$ is precisely
the profile of a viscous travelling wave of the $j$-th family
(contained in the center manifold $\M_j$). 
In this case, our decomposition should clearly contain one single
component:
$$u_x=v_j\tilde r_j(u,v_j,\sigma_j)\,.\eqno(5.1)$$
It is easy to guess what $v_j,\sigma_j$ in (5.1)
should be. Indeed, since by construction $|\tilde r_j|=1 $,
the quantity
$$v_j=\pm |u_x|$$
is the signed strength of the wave. 
Notice also that for a travelling wave the vectors $u_x$ and $u_t$ are
always parallel, since $u_t=-\sigma_j u_x$ where $\sigma_j$ is 
the speed of the wave.
We can thus write
$$u_t=u_{xx}-A(u)u_x=\omega_j\tilde r_j(u,v_j,\sigma_j)\eqno(5.2)$$
for some scalar $\omega_j$.   The speed of the wave is now obtained as
$\sigma_j=-\omega_j/ v_j$.
\v
Motivated by the previous analysis, as a first attempt we define
$$u_t=u_{xx}-A(u)u_x\eqno(5.3)$$ 
and try to find scalar quantities $v_i,\omega_i$ such that
$$\left\{\eqalign{ u_x&=\sum_i v_i\,\tilde r_i(u,v_i,\sigma_i),\cr
u_t&=\sum_i \omega_i\,\tilde r_i(u,v_i,\sigma_i),\cr}\right.
\qquad\qquad \sigma_i=-{\omega_i\over v_i}\,.\eqno(5.4)$$
The trouble with (5.4) is that the vectors $\tilde r_i$ are
defined only for speeds $\sigma_i$ close to the $i$-th characteristic
speed $\lambda_i^*\doteq \lambda_i(u^*)$.  However, when $u_x\approx 0$
one has $v_i\approx 0$ and the ratio $\omega_i/v_i$ may become
arbitrarily large.

\midinsert
\vskip 10pt
\centerline{\hbox{\psfig{figure=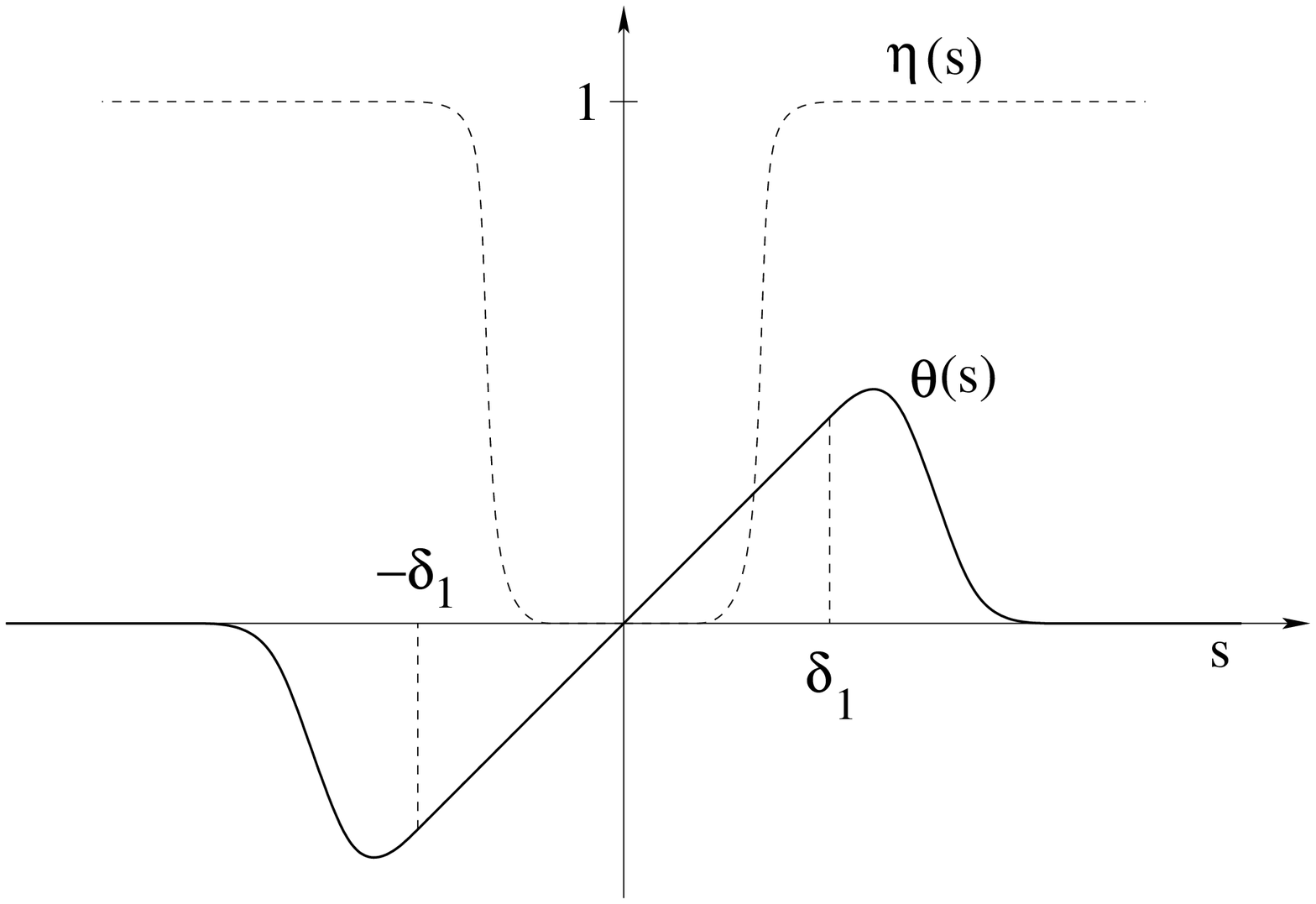,width=9cm}}}
\centerline{\hbox{figure 5}}
\vskip 10pt
\endinsert

To overcome this problem, we introduce a cutoff function (fig.~5).
Fix $\delta_1>0$ sufficiently small. Define a smooth odd function 
$\theta:\R\mapsto [-2\delta_1,\,2\delta_1]$ such that
$$\theta(s)=\cases{ s\qquad &if\qquad $|s|\leq \delta_1$\cr
0\qquad &if\qquad $|s|\geq 3\delta_1$\cr}\qquad\qquad |\theta'|\leq 1,
\qquad |\theta''|\leq 4/\delta_1\,.\eqno(5.5)$$
We now rewrite (5.4) in terms of the new variable $w_i$, related to $\omega_i$
by $\omega_i\doteq
w_i-\lambda_i^*v_i$. We
require that $\sigma_i$ coincides with $-\omega_i/v_i$ only when this
ratio is sufficiently close to $\lambda_i^*\doteq\lambda_i(u^*)$.
Our basic equations thus take the form
$$\left\{\eqalign{ u_x&=\sum_i v_i\,\tilde r_i(u,v_i,\sigma_i),\cr
u_t&=\sum_i (w_i-\lambda_i^*v_i)\,\tilde r_i(u,v_i,\sigma_i),\cr}\right.
\eqno(5.6)$$
where
$$u_t=u_{xx}-A(u)u_x\,,\qquad\qquad
\sigma_i=\lambda_i^*-\theta\left(w_i\over v_i\right)\,.\eqno(5.7)$$
Notice that $\sigma_i$ is not well defined when $v_i=w_i=0$. However,
recalling (4.16), in this case we have $\tilde r_i=r_i(u)$, regardless
of $\sigma_i$. Hence the two equations
in (5.6) are still meaningful.
\v
\n{\bf Remark 5.1.}  The decomposition (5.6) corresponds to 
viscous travelling waves $U_i$ such that
$$U_i(x)=u(x),\qquad\quad
U_i'(x)=v_i\tilde r_i\,,\qquad\quad
U''_i=\big(A(u)-\sigma_i\big)U'_i\,.$$
From the first equation in (5.6) it follows
$$u_x(x)=\sum_i U_i'(x)\,.$$
If $\sigma_i=\lambda_i^*-w_i/v_i$ for all $i=1,\ldots,n$, i.e.~if
none of the cutoff functions is active, then
$$\eqalign{u_{xx}(x)&=u_t+A(u)u_x\cr
&=\sum_i(w_i-\lambda_i^*v_i)\tilde r_i+A(u) \sum_i v_i\tr_i\cr
&=\sum_i \big(A(u)-\sigma_i\big)v_i\tr_i\cr
&=\sum_i U''_i(x)\,.\cr}$$
In this case, both of the equalities in (3.10) hold. 
Notice however that the second equality in (3.10)
may fail if $|w_i/v_i|>\delta_1$
for some $i$.
\v
\n{\bf Lemma 5.2.} 
{\it For $|u-u^*|$, $|u_x|$ and $|u_{xx}|$ sufficiently small,
the system of $2n$ equations (5.6) has a unique solution
$(v,w)=(v_1,\ldots,v_n,\,w_1,\ldots,w_n)$.
The map $(u,u_x,u_{xx})\mapsto (v,w)$ is smooth outside the
$n$ manifolds $\N_i\doteq\{v_i=w_i=0\}$;
moreover it is $\C^{1,1}$, i.e.~continuously differentiable 
with Lipschitz continuous derivatives
on a whole neighborhood of the point $(u^*,0,0)$.
}
\v
\n{\bf Proof.}~ Given $(v,w)$ in a neighborhood of $(0,0)\in\R^{2n}$,
the vectors $u_x,u_t$ are uniquely determined.
Hence the solution of (5.6)-(5.7) is certainly unique.
To prove its existence, 
consider the mapping $\Lambda:\R^n\times
\R^{n}\times\R^n\mapsto\R^{2n}$ defined by
$$\Lambda(u,v,w)\doteq\sum_{i=1}^n\Lambda_i(u,v_i,w_i),\eqno(5.8)$$
$$\Lambda_i(u,v_i,w_i)\doteq \pmatrix{ v_i\,\tilde r_i
\big(u,~v_i,~\lambda_i^*-\theta(w_i/v_i)\big)\cr 
(w_i-\lambda_i^*v_i)\,\tilde r_i
\big(u,~v_i,~\lambda_i^*-\theta(w_i/v_i)\big)\cr}.\eqno(5.9)
$$
This map is well defined and continuous
also when $v_i=0$, because in this case
(5.16) implies $\tilde r_i=r_i(u)$.
Computing the Jacobian matrix of partial derivatives w.r.t.~$(v_i,w_i)$ we find
$$\eqalign{&{\partial\Lambda_i\over\partial(v_i,w_i)}=
\pmatrix{\tilde r_i& 0\cr
-\lambda_i^* \tilde r_i & \tilde r_i\cr}
\cr
&\quad+
\pmatrix{ v_i \tilde r_{i,v}+(w_i/v_i)\theta'_i\tr_{i,\sigma} & 
-\theta'_i\tr_{i,\sigma}\cr
w_i\tr_{i,v}-\lambda_i^*v_i\tr_{i,v}-\lambda_i^*(w_i/v_i)\theta'_i
\tr_{i,\sigma}+(w_i/v_i)^2\theta'_i\tr_{i,\sigma}
&\lambda_i^*\theta'_i\tr_{i,\sigma}
-(w_i/v_i)\theta'_i\tr_{i,\sigma}\cr}\,.\cr}\eqno(5.10)$$
Here and throughout the following, by $\theta_i$, $\theta'_i$ we denote
the function $\theta$ and its derivative, evaluated at 
the point $s=w_i/v_i$.  By (5.10) we can write
$${\partial\Lambda\over\partial(v,w)}=B_0(u,v,w)+B_1(u,v,w),\eqno(5.11)$$
Because of (4.24), the matrix functions $B_0,B_1$
are well defined and continuous also when $v_i=0$.
Moreover, for $(v,w)$ small, $B_0$ has a uniformly bounded inverse 
and $B_1\to 0$ as $(v,w)\to 0$.  Since $\Lambda(u,0,0)=0\in\R^{2n}$,
we conclude that the
map $(v,w)\mapsto\Lambda(u;\,v,w)$ is $\C^1$ and invertible 
in a neighborhood of the origin. Therefore, 
given $(u,u_x,u_{xx})$, there exist
unique values of $(v,w)$ such that 
$$\Lambda(u,v,w)=\big(u_x,~u_{xx}-A(u)u_x\big).\eqno(5.12)$$
\v
The inverse of the map $\Lambda$ w.r.t.~the variables
$v,w$ will be denoted by $\Lambda^{-1}(u;\,p,q)$. In other words,
$$\Lambda^{-1}(u;\,p,q)=(v,w)\qquad\hbox{iff}\qquad \Lambda(u;\,v,w)=(p,q)\,.
$$
Since $\tr_i(u,0,\sigma_i)=r_i(u)$, we have
$$\Lambda(u,0,w)=\Big(~0~,~\sum_i w_i r_i(u)\Big).$$
Therefore,
$$\Lambda^{-1}(u,0,q)=(0,w)\qquad\hbox{where}
\qquad w_i=l_i(u)\cdot q\,.$$
In particular, $\Psi(u,0,0)=(0,0)\in \R^{2n}$.
Concerning first derivatives (which we regard here as linear operators), 
we have
$${\partial\Lambda(u;0,w)\over\partial (v,w)}\cdot (\hat v,\hat w)=
B_0(u;0,w)\cdot (\hat v,\hat w)=\left(\sum_i \hat v_i r_i(u)~,~\sum_i
\big(\hat w_i-\lambda_i(u)\hat v_i\big)r_i(u)\right),\eqno(5.13)$$
$${\partial\Lambda^{-1}(u;0,q)\over\partial (p,q)}\cdot (\hat p,\hat q)=
(\hat v,\hat w)\qquad\hbox{where}
\qquad \hat v_i=l_i(u)\cdot \hat p_i\,,\qquad
\hat w_i=l_i(u)\cdot \hat q_i\,.\eqno(5.14)$$

We shall not compute the second derivatives explicitly. 
However, one easily checks that
$${\partial^2\Lambda\over\partial v_i\partial v_j}=
{\partial^2\Lambda\over\partial v_i\partial w_j}=
{\partial^2\Lambda\over\partial w_i\partial w_j}=0\qquad\hbox{if}\quad
i\not= j\,.\eqno(5.15)$$
Moreover, recalling (4.24) and (5.5), we have the estimate
$${\partial^2\Lambda\over\partial v_i^2}\,,~
{\partial^2\Lambda\over\partial v_i\partial w_i}\,,~
{\partial^2\Lambda\over\partial w_i^2}~=~
\O(1)\cdot {1\over\delta_1}\,.\eqno(5.16)$$
Since the cutoff function $\theta$ vanishes for $|s|\geq 3\delta_1$,
it is clear that each $\Lambda_i$ is smooth outside the
manifold $\N_i\doteq\big\{(v,w)\,;~v_i=w_i=0\big\}$,
having codimension 2.
Since all second derivatives are uniformly bounded outside 
the $n$ manifolds $\N_i$, we conclude that $\Lambda$ is 
continuously differentiable
with Lipschitz continuous first derivatives
on a whole neighborhood of the point $(u^*,0,0)$.
Hence the same holds for $\Lambda^{-1}$.
\endproof 
\v
\n{\bf Remark 5.3.} By possibly performing a linear transformation of
variables, we can assume that the
matrix $A(u^*)$ is diagonal, hence its eigenvectors $r^*_1,\ldots r^*_n$
form an orthonormal basis:
$$\la r^*_i,\,r^*_j\ra =\delta_{ij}\,.\eqno(5.17)$$
Observing that
$$\big|\tilde r_i(u,v_i,\sigma_i)-r^*_i\big|=\O(1)\cdot 
\big(|u-u^*|+|v_i|\big),\eqno(5.18)$$
from (4.16) and the above assumption we deduce
$$\eqalign{\la \tilde r_i(u,v_i,\sigma_i),\,\tilde r_j(u,v_j,\sigma_j)
\ra&=
\delta_{ij}+\O(1)\cdot 
\big(|u-u^*|+|v_i|+|v_j|\big)\cr
&=\delta_{ij}+\O(1)\cdot\delta_0\,,\cr}\eqno(5.19)$$
$$\la\tr_i,\,A(u)\tr_j\ra=\O(1)\cdot\delta_0\qquad\qquad j\not= i\,.
\eqno(5.20)$$
Another useful consequence of (5.17)-(5.18) is the following.
Choosing $\delta_0>0$ small enough, 
the decomposition (5.6) will satisfy
$$|u_x|\leq\sum_i |v_i|\leq
2\sqrt n |u_x|\,.\eqno(5.21)$$
\v
We conclude this section by deriving estimates corresponding to
(2.13)-(2.15), valid for the components $v_i,w_i$.
In the following, given a solution $u=u(t,x)$ of (3.1) with
small total variation, we consider the
decomposition (5.6) of $u_x$ in terms of gradients of travelling waves.
It is understood that the vectors $\tr_i$ are constructed as in Section 4,
taking $P_i^*\doteq \big(u^*,0,\lambda_i(u^*)\big)$
as basic points in the construction of the center manifolds $\M_i$.
Here $u^*\doteq u(t,-\infty)$ is the constant state in (2.1).
\v
\n{\bf Lemma 5.4.} {\it In the same setting as Proposition 2.1,
assume that the bounds (2.6) hold on a larger interval
$[0,T]$.  
Then for all $t\in [\hat t,\, T]$, the decomposition
(5.6) is well defined.  The components
$v_i,w_i$ satisfy the estimates}
$$\eqalignno{\big\|v_i(t) \big\|_{\L^1}\,,
~\big\|w_i(t) \big\|_{\L^1} \,&=~
\O(1)\cdot \delta_0\,,
&(5.22)\cr
\big\|v_{i}(t) \big\|_{\L^\infty}\,,~\big\|w_{i}(t) \big\|_{\L^\infty}\,,
\big\|v_{i,x}(t) \big\|_{\L^1}\,,~~\big\|w_{i,x}(t) \big\|_{\L^1}
\,&=~
\O(1)\cdot \delta^2_0\,,&(5.23)\cr\big\|v_{i,x}(t) \big\|_{\L^\infty}\,,~
\big\|w_{i,x}(t) \big\|_{\L^\infty}
\,&=~
\O(1)\cdot \delta^3_0
\,.&(5.24)\cr}$$
\v
\n{\bf Proof.}
By Lemma 5.2, in a neighborhood
of the origin the map $(v,w)\mapsto\Lambda(u,v,w)$ in (5.8) is 
well defined, locally invertible, and continuously differentiable
with Lipschitz continuous derivatives.  Hence, for $\delta_0>0$ suitably
small, the $\L^\infty$ bounds in (2.13) and (2.14) 
guarantee that the decomposition (5.6) is well defined.
From the identity
(5.12) it now follows
$$v_i, w_i=\O(1)\cdot \big(|u_x|+|u_{xx}|\big).$$
By (2.6) and (2.13)-(2.14)
this yields the $\L^1$ bounds in (5.22) and the $\L^\infty$
bounds in (5.23).
Differentiating (5.12) w.r.t.~$x$ we obtain
$${\partial\Lambda\over\partial u}u_x+
{\partial\Lambda\over\partial (v,w)}(v_x,w_x)=
\Big(u_{xx},~u_{xxx}-A(u)u_{xx}-
\big(u_x\bullet A(u)\big)u_x\Big).\eqno(5.25)$$
Using the estimate
$${\partial\Lambda\over\partial u}=\O(1)\cdot \big(|v|+|w|\big),$$
since the derivative $\partial\Lambda/\partial (v,w)$
has bounded inverse, from (5.25) we deduce
$$(v_x,w_x)=\O(1)\cdot 
\Big(|u_{xx}|+|u_{xxx}|+|u_x|^2+|u_x|\big(|v|+|w|\big)
\Big).$$
This yields the remaining estimates in (5.23) and (5.24).
\endproof
\vsk
\n{\medbf 6 - Bounds on the source terms}
\v
We now consider a smooth solution $u=u(t,x)$ of 
(3.1) and let $v_i$, $w_i$, be the corresponding
components in the decomposition (5.6), which are well defined 
in view of Lemma 5.2.  
The equations governing the evolution of these $2n$ 
components can be written in the form
$$\left\{\eqalign{
v_{i,t}+(\tilde\lambda_i v_i)_x-v_{i,xx}&=
\phi_i\,,\cr
w_{i,t}+(\tilde\lambda_i w_i)_x-w_{i,xx}&=
\psi_i\,.\cr
}\right.
\eqno(6.1)$$
As in (4.21), we define here the speed 
$\tla_i\doteq \la \tr_i\,,~A(u)\tr_i\ra\,$.
The source terms $\phi_i,\psi_i$
can be computed by differentiating (3.1) and using
the implicit relations (5.6).  
However, it is not necessary to carry out in detail
all these computations.  Indeed, we are interested not in the 
exact form of these terms, but only in an upper bound for 
the norms $\|\phi_i\|_{\L^1}\,$ and
$\|\psi_i\|_{\L^1}\,$.  

Before giving these estimates, we provide an intuitive
explanation of how the source terms arise.
Consider first
the special case where $u$ is precisely one of the 
travelling
wave profiles on the center manifold (fig.~6a), 
say $u(t,x)=U_j(x-\sigma_j t)$.
We then have
$$u_x=v_j\tr_j\,,\qquad u_t=(w_j-\lambda_j^*v_j)\tr_j\,,
\qquad\quad v_i=w_i=0\quad
\hbox{for}~i\not= j\,,$$
and therefore
$$\left\{\eqalign{
v_{i,t}+(\tilde\lambda_i v_i)_x-v_{i,xx}&=
0\,,\cr
w_{i,t}+(\tilde\lambda_i w_i)_x-w_{i,xx}&=
0\,.\cr
}\right.
\eqno(6.2)
$$
Indeed, this is obvious when $i\not= j$. The identity $\phi_j=0$
follows from (4.28), while the relation $w_j=(\lambda_j^*-\sigma_j) v_j$
implies $\psi_j=0$.

\midinsert
\vskip 10pt
\centerline{\hbox{\psfig{figure=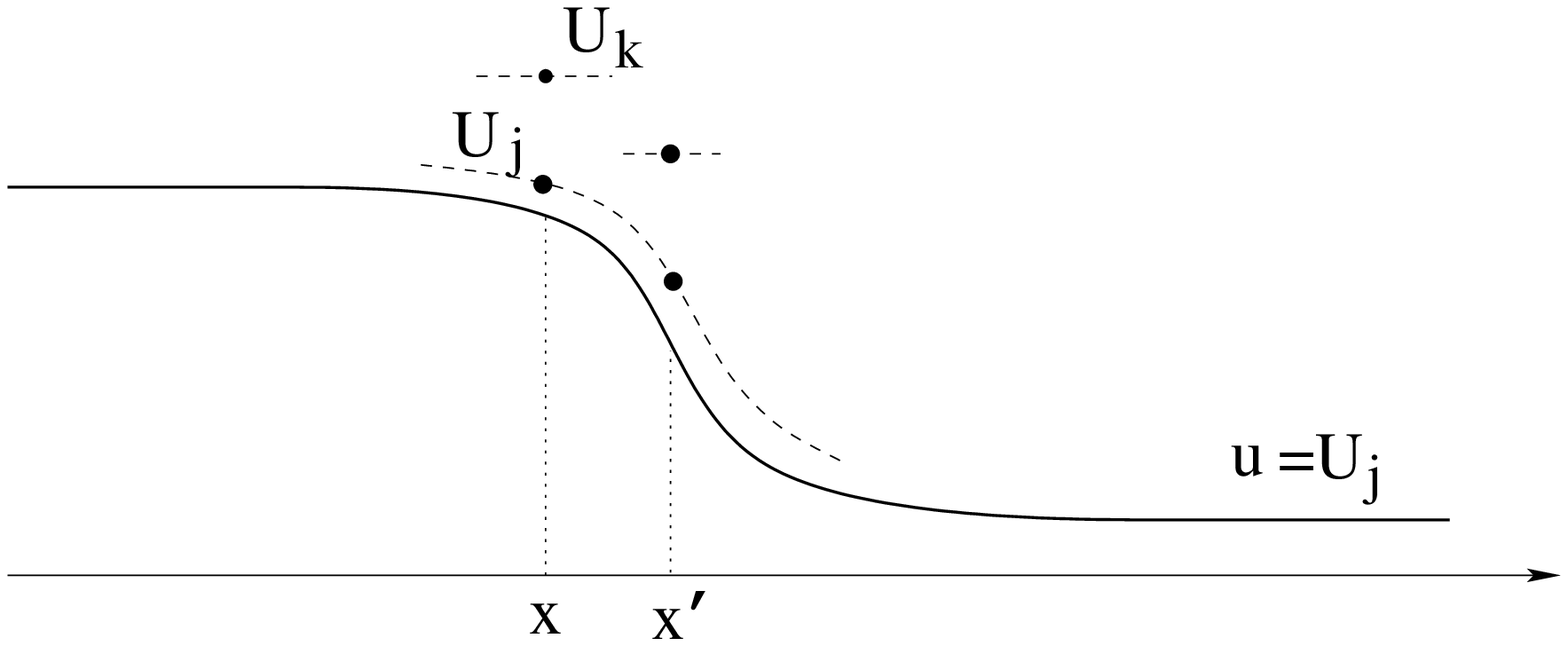,width=9cm}}}
\centerline{\hbox{figure 6a}}
\endinsert
\midinsert
\centerline{\hbox{\psfig{figure=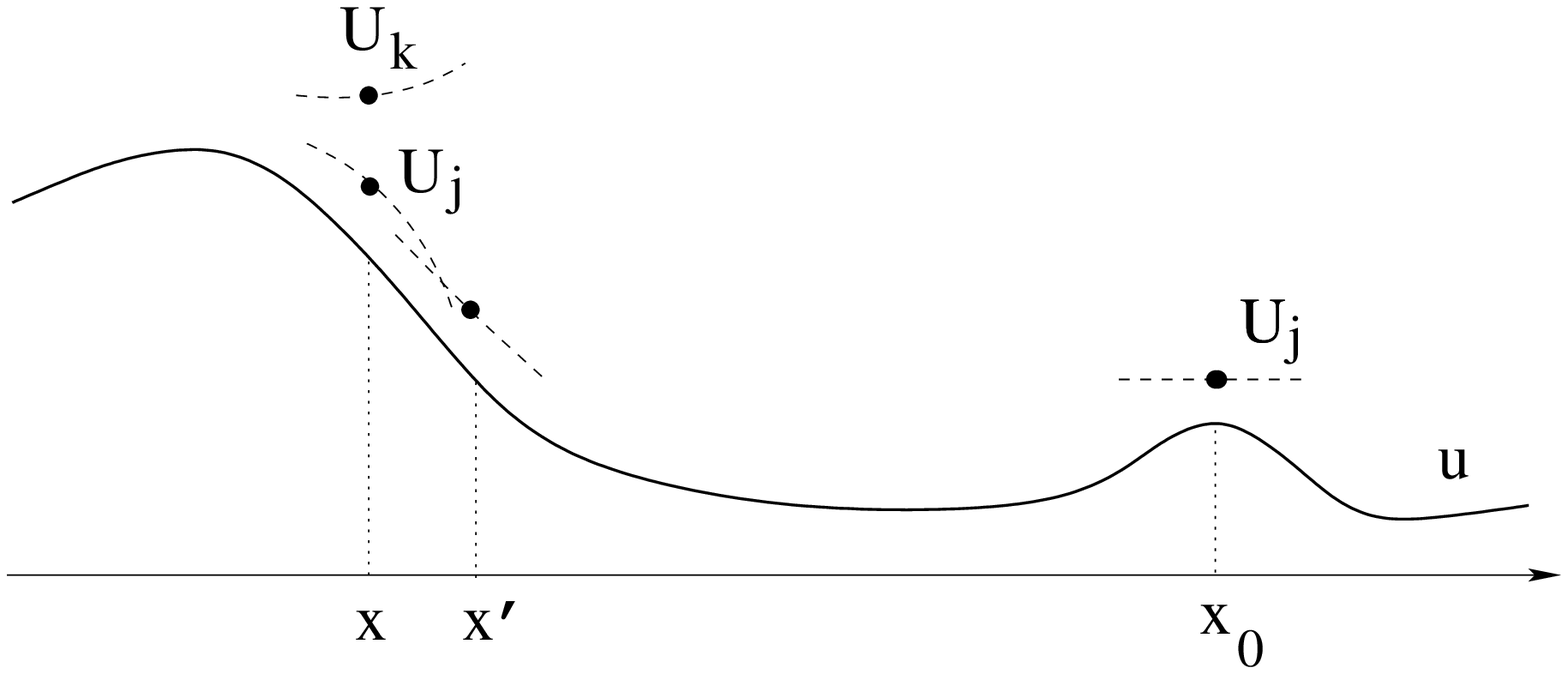,width=9cm}}}
\centerline{\hbox{figure 6b}}
\vskip 10pt
\endinsert

Next, consider the case of a general solution $u=u(t,x)$.
The sources on the right hand sides of (6.1) arise for 
three different  reasons (fig.~6b).
\v
\n{\bf 1.} The ratio $|w_j/v_j|$ is large and hence the
cutoff function $\theta$ in (5.7) is active. Typically,
this will happen near a 
point $x_0$ where $u_x=0$ but $u_t=u_{xx}\not= 0$.
In this case the identity (4.28) fails because of a ``wrong'' 
choice of the speed:
$\sigma_j\not= \lambda_j^*-(w_j/v_j) $.
\v
\n{\bf 2.}  Waves of two different families $j\not= k$
are present at a given point $x$.  These will produce
quadratic source terms, 
due to transversal interactions.
\v
\n{\bf 3.}  Since the decomposition (3.10) is defined pointwise,
it may well happen that
the travelling $j$-wave profile $U_j$ at
a point $x$ is not the same as 
the profile $U_j$ at a nearby point $x'$. 
Indeed, these two travelling waves may have slightly
different speeds.  It is the rate of change in this speed,
i.e.~$\sigma_{j,x}$, that determines
the infinitesimal interaction between 
nearby waves of the same family.
A detailed analysis will show that
the corresponding source terms
can only be linear
or quadratic w.r.t.~$\sigma_{j,x}$, with the square of the strength
of the wave always appearing as a factor.
These terms can thus be estimated as 
$\O(1)\cdot v_j^2\sigma_{j,x}+\O(1)\cdot v_j^2\sigma^2_{j,x}$.
\v
\n{\bf Lemma 6.1.} {\it The source terms in (6.1) satisfy the bounds
$$\eqalign{\phi_i,~\psi_i=\,&\O(1)\cdot \sum_j
\big( |v_{j,x}|+|w_{j,x}|\big)\cdot
|w_j-\theta_j v_j|~~\qquad\qquad
~\hbox{\bf (wrong speed)}\cr
&+\O(1)\cdot \sum_j |v_{j,x}w_j-v_jw_{j,x}|\qquad\qquad\qquad
~~~~~~~\hbox{\bf (change in speed, linear)}
\cr
&+\O(1)\cdot \sum_j \left|v_j\, \left(w_j\over v_j
\right)_x\right|^2\cdot
\chi_{\strut \big\{|w_j/v_j|<3\delta_1\big\}}\qquad~~~
\hbox{\bf (change in speed, quadratic)}
\cr
&+\O(1)\cdot\sum_{j\not= k}\big(|v_jv_k|+|v_{j,x}v_k|+
|v_jw_k|+|v_{j,x}w_k|+|v_jw_{k,x}|+|w_jw_k|\big)\cr
&\qquad\qquad\qquad\qquad\quad\qquad~~\hbox{\bf (interaction of waves of
different families)}\cr}$$
}
\v
From a direct inspection of the equations (6.1), it will be clear that
the source terms depend only on the third order jet $(u,\, u_x,\,
u_{xx},\, u_{xxx})$.  Since all functions
$\phi_i$, $\psi_i$  vanish in the case of a travelling wave,
for a general solution $u$  their size
can be estimated in terms of the distance between the third order jet
of $u$ and the (nearest) jet of some travelling wave.
This is indeed the strategy adopted in the following proof.
An alternative proof, based on more direct calculations,
will be given in Appendix~A.
\v
\n{\bf Proof of Lemma 6.1.} The conclusion will be reached in several steps.
\v
\n{\bf 1.} The vector
$(u_x,u_t)=\Lambda(u,v,w)$ satisfies the evolution equation
$$\pmatrix{u_x\cr u_t\cr}_t+\left(\left[\matrix{A(u) &0\cr
0& A(u)\cr}\right]\pmatrix{ u_x\cr u_t\cr}\right)_x-
\pmatrix{ u_x\cr u_t\cr}_{xx}=
\pmatrix{0\cr \big(u_x\bullet A(u)\big) u_t-
\big(u_t\bullet A(u)\big)u_x\cr}\,.
\eqno(6.3)
$$
Observe that, in the conservative case $A(u)=Df(u)$, the right hand
side vanishes because
$$\big(u_x\bullet A(u)\big) u_t=
\big(u_t\bullet A(u)\big)u_x=D^2f(u)\,(u_x\otimes u_t)\,.$$
In the general case, recalling (5.6)
we deduce
$$\big(u_x\bullet A(u)\big) u_t-
\big(u_t\bullet A(u)\big)u_x=\O(1)\cdot\sum_{j\not= k}
\big(|v_jv_k|+|v_jw_k|\big).\eqno(6.4)$$
\v
\n{\bf 2.}
For notational convenience, we introduce the variable
$z\doteq (v,w)$ and write 
$\tilde\lambda$ for the $2n\times 2n$ diagonal
matrix with entries $\tilde \lambda_i$ defined at (4.21):
$$
\tilde \lambda \doteq \left( \matrix{
\hbox{diag} ( \tilde \lambda_i ) & 0 \cr
0 & \hbox{diag} ( \tilde \lambda_i ) } \right).
$$
{}From (6.3) it now follows
$$
\eqalign{{\partial\Lambda\over\partial u}u_t +
{\partial\Lambda\over\partial z}&\pmatrix{v\cr w\cr}_t
+
\left(\left[\matrix{A(u) &0\cr
0& A(u)\cr}\right]\,\Lambda\right)_x
-{\partial\Lambda\over\partial z}
\pmatrix{v\cr w\cr}_{xx}-{\partial\Lambda\over\partial u}u_{xx}
\cr
&-{\partial^2\Lambda\over\partial u^{[2]}}(u_x\otimes u_x)-
{\partial^2\Lambda\over\partial z^{[2]}}\cdot
\pmatrix{v_x\cr w_x\cr}\otimes\pmatrix{v_x\cr w_x\cr}
-2{\partial^2\Lambda\over\partial u\,\partial z}\,u_x\otimes 
\pmatrix{v_x\cr w_x\cr}\cr
&\qquad=
\pmatrix{0\cr \big(u_x\bullet A(u)\big) u_t-
\big(u_t\bullet A(u)\big)u_x\cr}\,.\cr}
$$
Therefore,
$$\eqalign{{\partial\Lambda\over\partial z}&\left[
\pmatrix{v\cr w\cr}_t+\left(\tilde\lambda \pmatrix{v\cr w\cr}
\right)_x-\pmatrix{v\cr w\cr}_{xx}\right]\cr
&=
{\partial\Lambda\over\partial z}\left(\tilde\lambda 
\pmatrix{v\cr w\cr}\right)_x -\left(\left[\matrix{A(u) &0\cr
0& A(u)\cr}\right]\,\Lambda\right)_x+
{\partial\Lambda\over\partial u}A(u)u_x
+\pmatrix{0\cr \big(u_x\bullet A(u)\big) u_t-
\big(u_t\bullet A(u)\big)u_x\cr}\cr
&\qquad 
+{\partial^2\Lambda\over\partial u^{[2]}}(u_x\otimes u_x)+
{\partial^2\Lambda\over\partial z^{[2]}}\cdot
\pmatrix{v_x\cr w_x\cr}\otimes\pmatrix{v_x\cr w_x\cr}
+2{\partial^2\Lambda\over\partial u\,\partial z}\,u_x\otimes 
\pmatrix{v_x\cr w_x\cr}\cr
&\doteq E}\eqno(6.5)$$
Since the differential $\partial\Lambda/\partial z$ has 
uniformly bounded inverse, the right hand sides
in (6.1) clearly satisfy the bounds
$$\phi_i=\O(1)\cdot E\,,\qquad\qquad
\psi_i=\O(1)\cdot E,\qquad\qquad i=1,\ldots,n\,.\eqno(6.6)$$
\v
\n{\bf 3.} To estimate the quantity $E$ in (6.5), it is convenient to 
introduce the function
$$\Lambda_i(u,v_i,w_i,\sigma_i)\doteq\pmatrix{v_i\tr_i(u,v_i,\sigma_i)
\cr
(w_i-\lambda_i^*v_i)\tr_i(u,v_i,\sigma_i)\cr}\,,\eqno(6.7)$$
so that
$\Lambda = \sum\Lambda_i$ and $E = \sum E_i$, where
$$\eqalign{
E_i =&~ {\partial\Lambda_i\over\partial z_i}
\pmatrix{\tla_i v_i\cr \tla_i w_i\cr}_x -\left(\left[\matrix{A(u) &0\cr
0& A(u)\cr}\right]\,\Lambda_i\right)_x +
{\partial\Lambda_i\over\partial \sigma_i}
\,\tilde\lambda_i\,\sigma_{i,x} + {\partial
\Lambda_i \over \partial u} \sum_j A(u)v_j
\tilde r_j \cr
& + \pmatrix{0\cr\pi_i\big((u_x\bullet
A(u)) u_t - (u_t\bullet A(u))u_x\big) \cr} + {\partial^2 \Lambda_i \over
\partial u^{[2]}} u_x \otimes u_x
+ {\partial^2\Lambda_i\over\partial z_i^{[2]}} \cdot
\pmatrix{v_{i,x}\cr w_{i,x}\cr}\otimes\pmatrix{v_{i,x}\cr w_{i,x}\cr} \cr
& + {\partial^2\Lambda_i\over\partial \sigma_i^2}\,
\sigma_{i,x}^2 + 2 {\partial^2\Lambda_i\over\partial \sigma_i\,\partial z_i}\,
\pmatrix{\sigma_{i,x}v_{i,x}\cr \sigma_{i,x}w_{i,x}\cr} +
2 {\partial^2\Lambda_i\over\partial u\,\partial z_i}\,
\pmatrix{v_{i,x}\cr w_{i,x}\cr}\otimes u_x
+ 2 {\partial^2\Lambda_i\over\partial u\,\partial
\sigma_i}\,\sigma_{i,x} u_x \cr
& + {\partial \Lambda_i \over \partial \sigma_i} \left( {\partial^2
\sigma_i \over \partial v_i^2} v_{i,x}^2 + 2 {\partial^2
\sigma_i \over \partial v_i \partial w_i} v_{i,x} w_{i,x} + {\partial^2
\sigma_i \over \partial w_i^2}  w_{i,x}^2 \right).
} \eqno(6.8)$$
Notice that in (6.5) we regarded $\Lambda$ as a function of the
three independent variables $(u,v,w)$, while in (6.8) we think
$\Lambda$ as a function of the four independent variables
$(u,v,w,\sigma)$.
Regarding the $\sigma_i$ as independent variables, 
one has the advantage that
the maps $\Lambda_i=\Lambda_i(u,v_i,w_i,\sigma_i)$ are now smooth, while 
$\Lambda_i=\Lambda_i(u,v,w)$ in (5.8) was only $\C^{1,1}$, because of 
the singularities of the map $(u,v_i,w_i)\mapsto\sigma_i$ in (5.7).
The last term in (6.8) is due to the nonlinear dependence of $\sigma_i$
w.r.t.~$v_i,w_i$. 
By $\pi_i({\bf v})$ we denoted the $i$-th component of a vector ${\bf v}$~
w.r.t.~the
basis $\{r_1^*,\ldots,r_n^*\}$.
Also notice that
in the previous computation
we used the identity
$$\eqalign{
{\partial \Lambda_i \over \partial \sigma_i}& {\partial \sigma_i\over
\partial v_i}
\tla_{i,x}v_i 
+{\partial \Lambda_i \over \partial \sigma_i} {\partial \sigma_i\over
\partial w_i}
\tla_{i,x}w_i\cr
& =\tla_{i,x} \left( \matrix{w_i /v_i \cdot  
\theta_i'\tilde r_{i,\sigma} & -
 \theta_i'\tilde r_{i,\sigma} \cr
w_i/v_i ( w_i/v_i - \lambda_i^* ) \theta_i' \tilde r_{i,\sigma} & -
( w_i/v_i - \lambda_i^* )  \theta_i'\tilde r_{i,\sigma} } \right)
\pmatrix{ v_i \cr w_i \cr} = 
\pmatrix{0\cr 0\cr}\,.\cr}
$$
\v
\n{\bf 4.} By Lemma 5.2, 
the inverse map $\Lambda^{-1}$ sets a one to one correspondence
$$(u,u_x,u_{xx})\mapsto (u, v,w)$$
between two neighborhoods of the point $(u^*,0,0)\in\R^{3n}$.
This map is $\C^1$ with Lipschitz continuous derivative.
It can be prolonged to a map
$$(u,u_x,u_{xx},u_{xxx})\mapsto (u, v,w,\sigma,v_x,w_x,\sigma_x)
\eqno(6.9)$$
which is one-to-one, but of course
not onto.
Indeed, (5.6) and the identity $u_t+A(u)u_x=u_{xx}$ together imply
$$\eqalign{&\sum_i w_i\tr_i+\sum_i\big(A(u)-\lambda_i^*)v_i\tr_i
\cr
&\qquad =\sum_iv_{i,x}\tr_i+\sum_{ij}v_i\tr_{i,u}\,v_j\tr_j
+\sum_i v_i\tr_{i,v}v_{i,x}+\sum_i v_i
\tr_{i,\sigma}\sigma_{i,x}\,.\cr}\eqno(6.10)$$
A vector $(u, v,w,\sigma, v_x,w_x,\sigma_x)\in\R^{7n}$ 
corresponds to some third order jet $(u,u_x,u_{xx},u_{xxx})$
provided that it satisfies the vector equation (6.10), 
together with
$$\sigma_i=\lambda_i^*-\theta\left(w_i\over v_i\right)\,,
\qquad\qquad\sigma_{i,x}={w_i v_{i,x}-w_{i,x}v_i\over v_i^2}\,
\theta'\left(w_i\over v_i\right)\qquad\qquad i=1,\ldots,n\,.\eqno(6.11)$$
\v
\n{\bf 5.} By the analysis at (6.2), $E_i(u,v^\dag, w^\dag,\sigma^\dag,
v_x^\dag,w_x^\dag,
\sigma_x^\dag)=0$
whenever the argument 
corresponds to the third order jet of a viscous
travelling $i$-wave. This is the case
if
$$v_j^\dag=w_j^\dag=v^\dag_{j,x}=w^\dag_{j,x}=0, \qquad\qquad
\sigma^\dag_{j,x}=0,\qquad\qquad \hbox{for all}~~j\not= i,\eqno(6.12)$$
$$\left\{ \eqalign{v^\dag_{i,x}&=(\tla_i-\sigma^\dag_i)v^\dag_i\,,\cr
w^\dag_{i,x}&=(\tla_i-\sigma^\dag_i)w^\dag_i\,,\cr}\right.\qquad\qquad
\left|w^\dag_i\over v^\dag_i\right|<3\delta_1\,,\qquad
\sigma^\dag_i=\lambda_i^*-{w^\dag_i\over v^\dag_i}\,,\qquad
\sigma^\dag_{i,x}=0\,.\eqno(6.13)$$
In order to estimate $E_i(u,v,w,\sigma ,v_x,w_x,\sigma_x)$ 
we proceed as follows.
We introduce a new vector 
$(u,v^\dag, w^\dag,\sigma^\dag,
v_x^\dag,w_x^\dag,\sigma_x^\dag)$ 
corresponding to the jet of a travelling
$i$-wave, by setting
$$v_i^\dag=v_i\,,\qquad w_i^\dag=\theta\left(w_i\over v_i\right)v_i,
\qquad \sigma^\dag_i=\sigma_i=\lambda_i^*-{w^\dag_i\over v^\dag_i}\,.
\eqno(6.14)$$
The quantities $v^\dag_{i,x}\,,w^\dag_{i,x},\,\sigma^\dag_{i,x}$
are then defined according to (6.13), while
the components $j\not= i$ are as in
(6.12).  
The above construction implies $E^\dag_i\doteq
E_i(u,v^\dag, w^\dag,\sigma^\dag,
v_x^\dag,w_x^\dag,\sigma_x^\dag)=0$.
Hence $E_i=E_i-E^\dag_i$.
\v
\n{\bf 6.} 
Taking the inner product of 
(6.10) with $\tr_i$, recalling that $\tr_i$ 
has unit norm and is thus orthogonal to its 
derivatives, we obtain
$$w_i+(\tilde\lambda_i-\lambda^*_i)v_i=v_{i,x}+\Theta,$$
where
$$\eqalign{\Theta&\doteq
\sum_{j\not= i} \Big\langle \tr_i,~\big(\lambda_j^*-A(u)\big)
\tr_j\Big\rangle v_j+\sum_{j\not= i}
\sum_k\la \tr_i,~\tr_{j,u}\tr_k\ra v_jv_k\cr
&\qquad+\sum_{j\not= i} \la r_i,\,r_{j,v}\ra v_jv_{j,x}
+\sum_{j\not= i} \la r_i,\,r_{j,\sigma}\ra v_j\sigma_{j,x}\cr
&=
\O(1)\cdot\delta_0
\sum_{j\not= i}|v_j|\,.\cr}\eqno(6.15)$$
The above estimate on $\Theta$ is obtained using (5.20)
together with the $\L^\infty$ bounds in (5.23)-(5.24) and the
bound on $\tr_{j,\sigma}$ in (4.24).
We can now write
$$\eqalign{v_{i,x}&=w_i+(\tilde\lambda_i-\lambda^*_i)v_i+\O(1)\cdot
\delta_0\sum_{j\not= i}|v_j|\cr
&=(\tla_i-\sigma_i)v_i+(w_i-\theta_iv_i)+
\O(1)\cdot\delta_0\sum_{j\not= i}|v_j|\,.\cr}\eqno(6.16)$$
We recall that $\theta_i\doteq\theta(w_i/v_i)$.
The first equality in (6.16) yields the implications
$$|w_i|<3\delta_1|v_i|\qquad\Longrightarrow \qquad v_{i,x}=
\O(1)\cdot v_i+\O(1)\cdot\delta_0\sum_{j\not= i}|v_j|
\,.\eqno(6.17)$$
$$|w_i|>\delta_1|v_i|\qquad\Longrightarrow \qquad v_i=
\O(1)\cdot v_{i,x}+\O(1)\cdot\delta_0\sum_{j\not= i}|v_j|
\,.\eqno(6.18)$$
Moreover, from the second equality in (6.16) we deduce
$$\eqalign{(\tla_i-\sigma_i)w_i&=(\tla_i-\sigma_i)
\big[v_{i,x}-(\tla_i-\lambda_i^*)v_i\big]+\O(1)\cdot
\delta_0\sum_{j\not= i}|v_j|\cr
&=(\tla_i-\sigma_i)v_{i,x}-(\tla_i-\lambda_i^*)\big(v_{i,x}
-(w_i-\theta_iv_i)\big)+\O(1)\cdot\delta_0\sum_{j\not= i}|v_j|
\cr
&={w_i\over v_i}v_{i,x}-\left({w_i\over v_i}-\theta_i\right)v_{i,x}
+(\tla_i-\lambda_i^*)(w_i-\theta_iv_i)+\O(1)\cdot\delta_0\sum_{j\not= i}|v_j|
\,,\cr}$$
and hence, by (6.18),
$$w_{i,x}-(\tla_i-\sigma_i)w_i={w_{i,x}v_i-w_iv_{i,x}\over v_i}
+\O(1)\cdot |w_i-\theta_iv_i|+\O(1)\cdot\delta_0\sum_{j\not= i}|v_j|
\,.$$
{}From the definitions (6.13)-(6.14), using the above estimates
we obtain
$$\eqalign{|w_i-w_i^\dag|&=|w_i-\theta_iv_i|+\O(1)\cdot
\delta_0\sum_{j\not= i}|v_j|\,,\cr
|v_{i,x}-v^\dag_{i,x}|&=\O(1)\cdot |w_i-\theta_iv_i|+\O(1)\cdot
\delta_0\sum_{j\not= i}|v_j|\,,\cr
|w_{i,x}-w^\dag_{i,x}|&=\left|{w_{i,x}v_i-w_iv_{i,x}\over v_i}
\right|+
\O(1)\cdot |w_i-\theta_iv_i|+\O(1)\cdot\delta_0\sum_{j\not= i}|v_j|
\,.\cr}\eqno(6.19)$$
\v
\n{\bf 7.}  We now compute
$$\eqalign{
E_i =&~ E_i -  E^\dag_i \cr
=&~ {\partial\Lambda_i\over\partial z_i} 
\pmatrix{\tla_i v_i\cr \tla_i w_i\cr}_x -\left(\left[\matrix{A(u) &0\cr
0& A(u)\cr}\right]\,\Lambda_i\right)_x - {\partial
\Lambda^\dag_i\over\partial z_i} 
\pmatrix{\tla_i^\dag v_i^\dag\cr \tla_i^\dag w_i^\dag\cr}_x
-\left(\left[\matrix{A(u) &0\cr
0& A(u)\cr}\right]\,\Lambda^\dag_i\right)_x 
\cr
&+ 
{\partial\Lambda_i\over\partial \sigma_i}
\,\tilde\lambda_i\,\sigma_{i,x} - {\partial \Lambda_i \over \partial
u} \sum_{j\not= i} A(u)v_j\tr_j
\cr
&+ 
\pmatrix{0\cr\pi_i\big((u_x\bullet
A(u)) u_t - (u_t\bullet A(u))u_x\big) \cr} + \left(
{\partial^2 \Lambda_i \over \partial u^{[2]}} u_x \otimes u_x -
{\partial^2 \Lambda^\dag_i \over \partial u^{[2]}}\, 
v_i \tilde r_i^\dag \otimes
v_i\tilde r_i^\dag \right) \cr
& + {\partial^2\Lambda_i\over\partial z_i^{[2]}}\cdot
\pmatrix{v_{i,x}\cr w_{i,x}\cr}\otimes\pmatrix{v_{i,x}\cr w_{i,x}\cr}
- {\partial^2 \Lambda^\dag_i\over\partial z_i^{[2]}}\cdot
\pmatrix{v^\dag_{i,x}\cr w^\dag_{i,x}\cr}\otimes\pmatrix{v^\dag_{i,x}\cr
w^\dag_{i,x}\cr} \cr
& + {\partial^2\Lambda_i\over\partial \sigma_i^2}\,
\sigma_{i,x}^2 + 2 {\partial^2\Lambda_i\over\partial u\partial z_i}\, 
\pmatrix{v_{i,x}\cr w_{i,x}\cr}\otimes u_x
- 2 {\partial^2 \Lambda^\dag_i\over\partial u\partial z_i}\,
\pmatrix{v^\dag_{i,x}\cr w^\dag_{i,x}\cr}\otimes v_i \tilde r_i^\dag
+ 2 {\partial^2\Lambda_i\over\partial u\,\partial
\sigma_i}\,\sigma_{i,x}u_x \cr
& + 2 {\partial^2\Lambda_i\over\partial z_i\partial\sigma_i}\,
\pmatrix{\sigma_{i,x}v_{i,x}\cr \sigma_{i,x}w_{i,x}\cr}
+ {\partial \Lambda_i \over \partial \sigma_i} \left( {\partial^2
\sigma_i \over \partial^2 v_i} v_{i,x}^2 + 2{\partial^2
\sigma_i \over \partial v_i \partial w_i} v_{i,x} w_{i,x} + {\partial^2
\sigma_i \over \partial^2 w_i}w_{i,x}^2 \right).
} \eqno(6.20) $$
Observe that the
quantities $u,v_i,\sigma_i,\tilde r_i,\tla_i$ remain the same in
the computations of $E_i$ and $E_i^\dag$.
Moreover, all the terms involving derivatives w.r.t.~$\sigma_i$
vanish when we compute $E_i^\dag$.

In the remaining steps, we will examine the various terms on the right
hand side of (6.20) and show that they can all be bounded according
to the lemma.  As a preliminary, we observe that by (6.7) and (4.24) 
the derivatives of the smooth function 
$\Lambda_i=\Lambda_i(u,z_i,\sigma_i)$ satisfy
$${\partial \Lambda_i \over \partial u}\,,~~ {\partial^2 \Lambda_i \over
\partial u^{[2]}}\,,~~ {\partial^2 \Lambda_i \over \partial z_i
\partial \sigma_i} = \O(1)\cdot\big(|v_i|+|w_i|\big), \eqno(6.21)$$
$${\partial \Lambda_i \over
\partial \sigma_i},~~{\partial^2 \Lambda_i \over
\partial u \partial \sigma_i}\,,~~ {\partial^2 \Lambda_i \over
\partial \sigma_i^2} = \O(1) \cdot\big(|v_i^2|+|v_iw_i|\big)\,.\eqno(6.22)
$$
\v
\n{\bf 8.}
We start by collecting some transversal terms.
Using (6.4) and (6.21)-(6.22) we obtain
$$
\eqalign{
- {\partial \Lambda_i \over \partial u}& \sum_{j \not= i} 
A(u)v_j \tilde r_j + \pmatrix{0\cr\pi_i\big((u_x\bullet
A(u)) u_t - (u_t\bullet A(u))u_x\big) \cr} \cr
&\qquad\qquad + {\partial^2
\Lambda_i \over \partial u^{[2]}} \left( \sum_{j,k} v_j v_k \tilde r_j
\otimes \tilde r_k - v_i^2 \tilde r_i \otimes \tilde r_i \right) +
2 {\partial^2 \Lambda_i \over \partial u \partial z_i} \sum_{j \not= i}
\left( \matrix{
v_{i,x} \cr w_{i,x} } \right) \otimes v_j \tilde r_j\cr
&= \O(1)\cdot \sum_{j \not= k} \big(|v_j v_k| + 
|w_j v_k|\big) +\O(1)\cdot \sum_{j \not= i} \big(|w_j w_i| + |v_j v_{i,x}|
+|v_j w_{i,x}|\big).
}\eqno(6.23)
$$
Here and in the following, by ``transversal terms'' we mean
terms whose size is bounded by products of distinct components $j\not=k$,
as in (6.23).
\v
\n{\bf 9.} We now look at terms involving
derivatives w.r.t.~$\sigma_i$. One should here keep in mind that,
if $\sigma_{i,x}\not= 0$, then both sides of the implication (6.16)
hold true.
Using (6.21) we obtain
$$
\eqalign{
{\partial^2\Lambda_i\over\partial z_i\partial \sigma_i}\,
&\pmatrix{\sigma_{i,x}v_{i,x}\cr \sigma_{i,x}w_{i,x}\cr} =\O(1)\cdot
v_i\big(|v_{i,x}|+|w_{i,x}|\big)\sigma_{i,x}\cr
&\qquad=\O(1)\cdot v_i^2\sigma_{i,x}+\O(1)\cdot  (w_i v_{i,x} - w_{i,x}
v_i) \sigma_{i,x} +\hbox{transversal terms}\cr
&\qquad=\O(1)\cdot |w_i v_{i,x} - w_{i,x}
v_i| + \left| v_i \cdot \left( w_i\over v_i\right)_x\right|^2
\cdot\chi_{\strut
\{|w_i/v_i|<3\delta_1\}}+\hbox{transversal terms}\,.\cr}\eqno(6.24)$$
An application of (6.22) yields
$$\eqalign{{\partial\Lambda_i\over\partial \sigma_i}\tla_i&\sigma_{i,x}+
{\partial^2\Lambda_i\over\partial \sigma_i^2}\sigma_{i,x}^2
+2{\partial^2\Lambda_i\over\partial u\partial \sigma_i}\sigma_{i,x}u_x=
\O(1)\cdot v_i^2\big(|\sigma_{i,x}|+|\sigma_{i,x}^2|\big)
+\hbox{transversal terms}\cr
&=\O(1)\cdot |v_iw_{i,x}-w_iv_{i,x}|+\O(1)\cdot
\left|v_i\cdot \left( w_i\over v_i\right)_x\right|^2\cdot\chi_{\strut
\{|w_i/v_i|<3\delta_1\}}
+\hbox{transversal terms}\,.\cr}\eqno(6.25)$$
Next, we observe that the quantity
$${\partial^2
\sigma_i \over \partial^2 v_i}v_{i,x}^2 + 2{\partial^2
\sigma_i \over \partial v_i \partial w_i} v_{i,x} w_{i,x} + {\partial^2
\sigma_i \over \partial^2 w_i} w_{i,x}^2$$
vanishes in the special case where $w_{i,x}=(w_i/v_i)v_{i,x}$.
In general, recalling (6.16) a direct computation yields
$$\eqalign{
{\partial \Lambda_i \over \partial \sigma_i} &\left( {\partial^2
\sigma_i \over \partial^2 v_i}v_{i,x}^2 + 2{\partial^2
\sigma_i \over \partial v_i \partial w_i} v_{i,x} w_{i,x} + {\partial^2
\sigma_i \over \partial^2 w_i} w_{i,x}^2 \right)\cr
&={\partial \Lambda_i \over \partial \sigma_i}\left\{
\left(-\theta_i''{w^2_i\over v_i^4}-2\theta_i'{w_i\over v_i^3}\right)
v_{i,x}^2+2\left(\theta_i''{w^2_i\over v_i^4}+\theta_i'{w_i\over v_i^3}
\right)v_{i,x}w_{i,x}-\theta_i''{w^2_i\over v_i^4}w_{i,x}^2\right\}\cr
&={\partial \Lambda_i \over \partial \sigma_i}\left\{
2\left(\theta_i''{w^2_i\over v_i^4}+\theta_i'{w_i\over v_i^3}
\right)v_{i,x}\left(w_{i,x}-{w_iv_{i,x}\over v_i}\right)
-\theta_i''{w^2_i\over v_i^4}\left(w_{i,x}-{w_iv_{i,x}\over v_i}\right)^2
\right\}\cr
&=\O(1)\cdot |v_iw_{i,x}-w_iv_{i,x}|+\O(1)\cdot
\left|v_i\left(w_i\over v_i\right)_x\right|^2\cdot
\chi_{\strut \{|w_j/v_j|<3\delta_1\}}+\hbox{transversal terms}
\,.\cr}
\eqno(6.26)
$$
\v
\n{\bf 10.} We now complete the analysis of the remaining terms.
As a preliminary, we observe that 
the only difference between $\Lambda_i^\dag$ and $\Lambda_i$
is due to the fact that one may have $w_i^\dag\not= w_i$.
The first estimate in (6.19)
thus implies
$$\eqalign{\big|\Lambda_i^\dag-\Lambda_i\big|\,,~~
\big|D\Lambda_i^\dag-D\Lambda_i\big|\,,~~
\big|D^2\Lambda_i^\dag-D^2\Lambda_i\big|~=~\O(1)\cdot
|w_i-\theta_iv_i|+\O(1)\cdot
\delta_0\sum_{j\not=i}|v_j|\,.\cr}
\eqno(6.27)$$
By (6.17) and (6.27),
if  we compute  $\Lambda_i$ or its partial derivatives 
at the point $(u,v_i,w_i,\sigma_i)$ instead of 
$(u,v^\dag_i,w^\dag_i,\sigma^\dag_i)=(u,v_i,w_i^\dag,\sigma_i)$, 
the difference in each of the corresponding
terms in (6.20) will have magnitude
$$\O(1)\cdot |w_i-\theta_iv_i|\cdot\big(|v_{i,x}|+|w_{i,x}|\big)
+\hbox{transversal terms}\,.$$
For example,
$$\eqalign{\left({\partial^2\Lambda_i\over\partial u^{[2]}}-
{\partial^2\Lambda_i^\dag\over\partial u^{[2]}}\right)(v_i\tr_i\otimes
v_i\tr_i)
&=\O(1)\cdot |w_i-\theta_iv_i|\,v^2_i
+\hbox{transversal terms}\cr
&=\O(1)\cdot |w_i-\theta_iv_i|\,v^2_{i,x}
+\hbox{transversal terms}\,.\cr}\eqno(6.28)$$
Indeed, if $w_i\not=\theta_iv_i$, then both sides of the implication
(6.17) hold true.
\v
Observing that $\partial^2\Lambda_i/\partial w_i^2=0$ and recalling
(6.17), we have
$$
\eqalign{
{\partial^2 \Lambda_i \over \partial z_i^{[2]}}& \bigg[\left( \matrix{
v_{i,x} \cr w_{i,x} } \right) \otimes \left( \matrix{
v_{i,x} \cr w_{i,x} } \right) - \left( \matrix{
v_{i,x}^\dag \cr w_{i,x}^\dag } \right) \otimes \left( \matrix{
v_{i,x}^\dag \cr w_{i,x}^\dag } \right) \bigg]\cr
&= \O(1)\cdot \big(
v_{i,x}^2 - ( v_{i,x}^\dag )^2 \big)  + {\cal O}(1)\cdot \big(
v_{i,x} w_{i,x} - v_{i,x}^\dag w_{i,x}^\dag ) \cr
&=\O(1)\cdot v_{i,x} | w_i - \theta_i v_i |
+ \O(1)\cdot w_{i,x} | w_i - \theta_i v_i |
+ \O(1) | w_i v_{i,x} - v_i w_{i,x} |
+ \hbox{transversal terms}\,.\cr }
\eqno(6.29)
$$
In a similar way, using (6.17) and (6.19) one derives the estimate
$$
\eqalign{
{\partial^2 \Lambda_i \over \partial u \partial z_i}
\left( \matrix{ v_{i,x} - v_{i,x}^\dag \cr
w_{i,x} - w_{i,x}^\dag } \right)\otimes v_i\tr_i &=\O(1)\cdot 
v_{i,x} ( w_i -
\theta_i v_i)
+ \O(1)\cdot w_{i,x} | w_i - \theta_i v_i | \cr
&\qquad+ \O(1)\cdot | w_i v_{i,x} - v_i w_{i,x} |
+ \hbox{transversal terms}\,. }\eqno(6.30)
$$
Using the identity (4.23),
we now compute
$$\eqalign{
{\partial\Lambda_i\over\partial z_i} &
\pmatrix{\tla_i v_i\cr \tla_i w_i\cr}_x -\left(\left[\matrix{A(u) &0\cr
0& A(u)\cr}\right]\,\Lambda_i\right)_x = \left[ \matrix{ I - A(u) & 0
\cr
0 & I - A(u) } \right] \left[\matrix{v_i
\tilde r_{i,v} &0\cr
( w_i - \lambda^*_i v_i ) \tilde r_{i,v} & 0\cr}\right]
\pmatrix{ \tla_i v_i\cr \tla_i w_i\cr}_x \cr
&+ \left[\matrix{\tilde r_i &0\cr
- \lambda_i^* \tilde r_i &  \tilde r_i \cr} \right]
\pmatrix{ \tla_{i,x} v_i\cr \tla_{i,x} w_i\cr} - v_i \left[\matrix{ 
\tilde r_{i,u} \tilde r_i +  ( \tla_i - \sigma_i )
\tilde r_{i,v} & 0 \cr
- \lambda_i^* ( \tilde r_{i,u} \tilde r_i +  ( \tla_i -
  \sigma_i )
\tilde r_{i,v} ) & \tilde r_{i,u} \tilde r_i +  ( \tla_i -
  \sigma_i ) \tilde r_{i,v} \cr} \right]
\pmatrix{ v_{i,x}\cr w_{i,x}\cr} \cr
&- \left[\matrix{A(u) &0\cr 0& A(u)\cr}\right] {\partial
\Lambda_i \over \partial \sigma_i} \sigma_{i,x} - \sum_j \left[ \matrix{
DA(u) \tilde r_j & 0 \cr
0 & DA(u) \tilde r_j } \right] v_j \Lambda_i - \sum_j \left[ \matrix{
A(u) & 0 \cr
0 & A(u) } \right] {\partial \Lambda_i \over \partial u} v_j \tilde r_j\,.
} $$
With similar arguments as above, we obtain
$$ \eqalign{&
\left( \matrix{ ( \tilde \lambda_i v_i)_x v_i
\tilde r_{i,v} \cr
( \tilde \lambda_i v_i
)_x v_i ( w_i - \lambda^*_i v_i ) \tilde r_{i,v} }
\right) - \left( \matrix{ ( \tilde \lambda_i^\dag v_i^\dag
)_x v_i^\dag
\tilde r_{i,v}^\dag \cr
( \tilde \lambda_i^\dag v_i^\dag
)_x v_i^\dag ( w_i^\dag - \lambda^*_i v_i^\dag ) \tilde r_{i,v}^\dag }
\right) \cr
&\qquad\qquad = \O(1)\cdot | v_i w_{i,x} - v_{i,x} w_i |\cdot
\chi_{\strut\{ | w_i/v_i | \leq 3 \delta_1 \}}
+\O(1)\cdot v_{i,x} | w_i - \theta_i v_i |+ \hbox{transversal terms}\,,
} $$
$$
\pmatrix{ \tla_{i,x} v_i\cr \tla_{i,x} w_i\cr} - \pmatrix{
\tla_{i,x}^\dag v_i^\dag \cr \tla_{i,x}^\dag w_i^\dag \cr} = \O(1)
\cdot |v_i w_{i,x} - v_{i,x} w_i | + \O(1)\cdot v_{i,x} | w_i -
\theta_i v_i | + \hbox{transversal terms}\,,
$$
$$
\pmatrix{ v_i v_{i,x} \cr v_i w_{i,x} \cr} - \pmatrix{ v_i^\dag v_{i,x}^\dag
\cr v_i^\dag w_{i,x}^\dag \cr} = \O(1)
\cdot |v_i w_{i,x} - v_{i,x} w_i | + \O(1)\cdot v_{i,x} | w_i -
\theta_i v_i | + \hbox{transversal terms}\,.
$$
The above estimates together imply
$$\eqalign{
{\partial\Lambda_i\over\partial z_i} &
\pmatrix{\tla_i v_i\cr \tla_i w_i\cr}_x -\left(\left[\matrix{A(u) &0\cr
0& A(u)\cr}\right]\,\Lambda_i\right)_x -
{\partial\Lambda_i^\dag\over\partial z_i}
\pmatrix{\tla_i^\dag v_i^\dag \cr \tla_i^\dag w_i^\dag\cr}_x
-\left(\left[\matrix{A(u) &0\cr
0& A(u)\cr}\right]\,\Lambda_i^\dag\right)_x \cr
&\quad =\O(1) \cdot| w_i -
\theta_i v_i |\big(|v_{i,x}|+|w_{i,x}|\big)
+ \O(1) \cdot | w_i v_{i,x} - v_i w_{i,x} |
+ \hbox{transversal terms}\,.\cr }
\eqno(6.31)$$
This completes the proof of Lemma 6.1.
\endproof
\vsk
\n{\medbf 7 - Transversal wave interactions}
\v
The goal of this section is to establish an a priori bound on the
total amount of interactions between waves of different families.
More precisely, let $u=u(t,x)$ be a solution of the
parabolic system (3.1) and assume that
$$\big\|u_x(t)\big\|_{\L^1}\leq\delta_0\qquad\qquad t\in [0,T]\,.\eqno(7.1)$$
In this case, for $t\geq \hat t$, by Corollary 2.2
all higher derivatives will
be suitably small and  
we can thus
define the components 
$v_i,w_i$ according to (5.6)-(5.7).
These will satisfy the linear evolution equation (6.1), with source terms
$\phi_i,\psi_i$ described in Lemma 6.1.
Assuming that
$$\int_{\hat t}^T\!\int  \big|\phi_i(t,x)\big|+
\big|\psi_i(t,x)\big|\,dxdt\leq\delta_0\,,\qquad\qquad i=1,\ldots,n,
\eqno(7.2)$$
and relying on the bounds (5.22)--(5.24),
we shall prove the estimate
$$\int_{\hat t}^T\!
\int \sum_{j\not= k} \big( |v_jv_k|+|v_{j,x}v_k|+
|v_jw_k|+|v_{j,x}w_k|+|v_j w_{k,x}|+|w_jw_k|\big)\,dxdt
=\O(1)\cdot\delta_0^2\,.\eqno(7.3)$$
\v
As a preliminary, we establish a more general estimate
on solutions of two independent linear parabolic
equations, with strictly different drifts (fig.~7).
\v
\n{\bf Lemma 7.1.} {\it Let $z,z^\sharp $ be solutions of the
two independent scalar equations
$$\eqalign{z_t+\big(\lambda(t,x)z\big)_x-z_{xx}&=\vp(t,x)\,,\cr
z^\sharp _t+\big(\lambda^\sharp (t,x)z^\sharp \big)_x
-z^\sharp _{xx}&=\vp^\sharp(t,x)\,,\cr}\eqno(7.4)$$
defined for $t\in [0,T]$. Assume that
$$\inf_{t,x}\lambda^\sharp(t,x)~-~\sup_{t,x}\lambda(t,x)~\geq~ c~>~0\,.
\eqno(7.5)$$
Then}
$$\eqalign{\int_0^T &\!\!\int\big|z(t,x)\big|\,\big|z^\sharp(t,x)\big|
\,dxdt\cr
&\leq {1\over c}
\left(\int  \big|z(0,x)\big|\,dx+\int_0^T\!\int
\big|\vp(t,x)\big|\,dxdt\right)
\left(\int  \big|z^\sharp(0,x)\big|\,dx+\int_0^T\!\int
\big|\vp^\sharp(t,x)\big|\,dxdt\right).\cr}\eqno(7.6)$$

\midinsert
\vskip 10pt
\centerline{\hbox{\psfig{figure=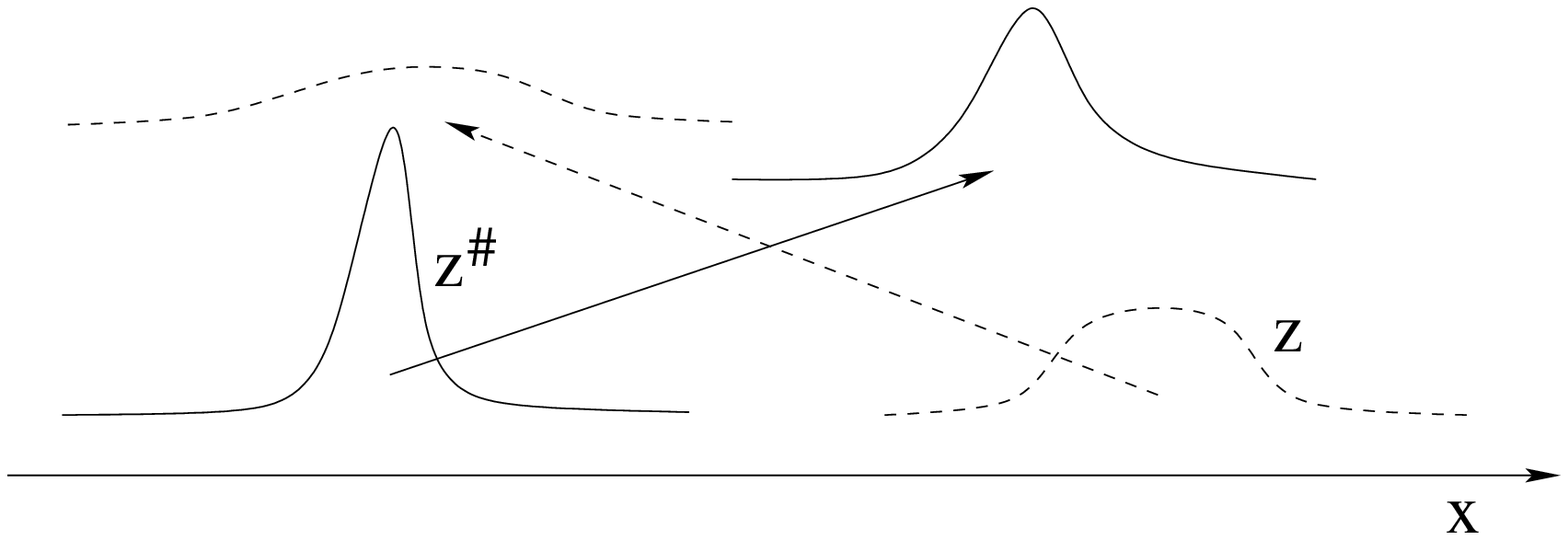,width=10cm}}}
\centerline{\hbox{figure 7}}
\vskip 10pt
\endinsert

\n{\bf Proof.} We consider first the homogeneous case, where
$\vp=\vp^\sharp=0$. Define
the interaction potential
$$Q(z,z^\sharp)\doteq\dint K(x-y)\,\big| z(x) \big|
\, \big| z^\sharp(y)
\big| \,dx dy\,,\eqno(7.7)$$
with
$$K(s)\doteq\cases{
1/c\qquad &if\qquad $s \geq 0$,\cr
1/c\cdot e^{cs/2}\qquad &if\qquad $s < 0$.\cr}
\eqno(7.8)$$
Computing the
distributional derivatives of the kernel $K$ we find that
$cK'-2K''$ is precisely the Dirac distribution, i.e.~a
unit mass at the origin.
A direct computations now yields
$${d\over dt}Q\big(z(t),z^\sharp(t)\big)
\leq -\int\big| z(t,x) \big|
\, \big| z^\sharp(t,x)
\big| \,dx\,.$$
Therefore
$$
\int_0^T\!\!\int \big| z(t,x)\big|\,\big| z^\sharp(t,x) \big|\, dx dt
\leq Q\big(z(0),\,z^\sharp(0)\big)
\leq {1\over c}\,\big\| z(0) \big\|_{\L^1} \big\| z^\sharp(0) \big\|_{\L^1}\,.
\eqno(7.9)$$
proving the lemma in the homogeneous case.

To handle the general case, call $\Gamma,\Gamma^\sharp$ the Green functions
for the corresponding linear homogenous systems.
The general solution of (7.4) can thus be written in the form
$$\eqalign{z(t,x)&=\int \Gamma(t,x,0,y) z(0,y) dy + \int_0^t\!\!\int
\Gamma(t,x,s,y) \vp(s,y) dyds,\cr
z^\sharp(t,x)&=\int \Gamma^\sharp(t,x,0,y) z^\sharp(0,y) dy + \int_0^t\!\!\int
\Gamma^\sharp(t,x,s,y) \vp^\sharp(s,y) dyds.\cr}\eqno(7.10)$$
From (7.9) it follows
$$\int_{\max\{s,s'\}}^T\int  \Gamma(t,x,s,y)\cdot 
\Gamma^\sharp(t,x,s',y')\,dxdt
\leq {1\over c}\eqno(7.11)$$
for every couple of initial points $(s,y)$ and $(s',y')$.  
The estimate (7.6) now follows from (7.11) and the representation 
formula (7.10).
\endproof
\v
\n{\bf Remark 7.2.}  Exactly the same estimate (7.6) 
would be true also for a system without viscosity.
In particular, if
$$z_t+\big(\lambda(t,x)z\big)_x=0,\qquad\qquad
z^\sharp _t+\big(\lambda^\sharp (t,x)z^\sharp \big)_x=0,$$
and if the speeds 
satisfy the gap condition (7.5), then
$$
{d\over dt}\left[ {1\over c}\dint_{x<y} 
\big|z^\sharp(t,x) z(t,y)\big|\,dxdy\right]
\leq -\int\big| z(t,x) \big|
\, \big| z^\sharp(t,x)
\big| \,dx\,.$$
In the case where viscosity is present, our
definition (7.7)-(7.8) thus provides a natural counterpart
to the Glimm interaction potential between waves of different families,
introduced in [G] for strictly hyperbolic systems. 
\v
Lemma 7.1 allows us to
estimate the integral of the terms $|v_iv_k|$,
$|v_jw_k|$ and $|w_jw_k|$ in (7.3).
We now work toward an estimate of the remaining terms
$|v_{j,x}v_k|$, $|v_{j,x}w_k|$ and $|v_j w_{k,x}|$,
containing one derivative w.r.t.~$x$. 
\v
\n{\bf Lemma 7.3.}  {\it Let $z,z^\sharp$ be solutions of (7.4)
and assume that (7.5) holds, together with the estimates
$$\int_0^T\!\int \big|\vp(t,x)\big|\,dxdt\leq\delta_0,\qquad\qquad
\int_0^T\!\int \big|\vp^\sharp(t,x)\big|\,dxdt\leq\delta_0\,,
\eqno(7.12)$$
$$\big\|z(t)\big\|_{\L^1}\,,~\big\|z^\sharp(t)\big\|_{\L^1}\leq \delta_0\,,
\qquad\qquad
\big\|z_x(t)\big\|_{\L^1}\,,~\big\|z^\sharp(t)\big\|_{\L^\infty}
\leq C^*\delta_0^2\,,
\eqno(7.13)$$
$$\big\|\lambda_x(t)\big\|_{\L^1}\leq 
C^*\delta_0\,,
\qquad\qquad
\lim_{x\to -\infty} \lambda(t,x)=0\eqno(7.14)$$
for all $t\in [0,T]$. 
Then one has the bound}
$$\int_0^T\!\int \big|z_x(t,x)\big|\,
\big|z^\sharp(t,x)\big|\,dxdt=\O(1)\cdot\delta_0^2.
\eqno(7.15)$$
\v
\n{\bf Proof.} 
The left hand side of (7.15) is clearly bounded by the quantity
$$\I(T)\doteq\sup_{(\tau,\xi)\in [0,T]\times\R}\int_0^{T-\tau}
\!\int  \big|z_x(t,x)z^\sharp(t+\tau,x+\xi)\big|\,dxdt\leq
(C^*\delta_0^2)^2\cdot T\,,$$
the last inequality being a consequence of (7.13). 
For $t>1$ we can write $z_x$ in the form
$$z_x(t,x)=\int G_x(1,y)z(t-1,\,x-y)\,dy
+\int_0^1\!\int G_x(s,y) 
\big[\vp-(\lambda z)_x\big](t-s,\,x-y)\,dyds\,,
$$
where $G(t,x)\doteq \exp\{-x^2/4t\}/2\sqrt{\pi t}$
is the standard heat kernel.
Using (7.6)
we obtain
$$\eqalign{&\int_1^{T-\tau}\!\int 
\big|z_x(t,x)\,z^\sharp(t+\tau,x+\xi)\big|\,dxdt
\cr
&~\leq
\int_1^{T-\tau}\!\dint \Big|G_x(1,y)\,z(t-1,\,x-y)
\,z^\sharp(t+\tau,x+\xi)\Big|\,dydxdt\cr
&\qquad +\int_1^{T-\tau}\!\int \int_0^1\int
\big\|\lambda_x\|_{\L^\infty}\,\Big|G_x(s, y)\,z(t-s,\,x-y)
\,z^\sharp(t+\tau,x+\xi)\Big|\,dydsdxdt\cr
&\qquad +\int_0^{T-\tau}\!\int \int_0^1\int
\big\|\lambda\|_{\L^\infty}\,\Big|G_x(s, y)\,z_x(t-s,\,x-y)
\,z^\sharp(t+\tau,x+\xi)\Big|\,dydsdxdt\cr
&\qquad +\int_1^{T-\tau}\!\int
\int_{t-1}^t\int \Big|G_x(t-s,\,x-y) \,\vp(s,y) \,
z^\sharp(t+\tau,x+\xi)\Big|\,dydsdxdt
\cr
&~\leq
\left(\int\big|G_x(1,y)\big|\,dy+\|\lambda_x\|_{\L^\infty}\int_0^1\!\int
\big|G_x(s,y)\big|\,dyds\right)\cr
&\qquad\qquad\cdot
\sup_{s,y,\tau,\xi}\,\left(\int_1^{T-\tau}\!\int \big|z(t-s,\,x-y)\big|\,
\big|z^\sharp(t+\tau,x+\xi)\big|\,dxdt\right)\cr
&\qquad+
\left(\|\lambda\|_{\L^\infty}\cdot\int_0^1\!\int
\big|G_x(s,y)\big|\,dyds\right)
\cdot \left(\sup_{s,y,\tau,\xi}\int_1^{T-\tau}\!\int \big|z_x(t-s,\,x-y)\big|\,
\big|z^\sharp(t+\tau,x+\xi)\big|\,dxdt\right)\cr
&\qquad +\|z^\sharp\|_{\L^\infty}\cdot\int_0^1\int\big|G_x(s,y)\big|\,dsdy\cdot
\int_0^T\int\big|\vp(t,x)\big|\,dxdt\cr
&~\leq
2\cdot {1\over c} 2\delta_0\cdot 2\cdot\delta_0
+2\|\lambda\|_{\L^\infty} \I(T)+C^*\delta_0^2\cdot 2\delta_0\,.\cr
}\eqno(7.16)
$$
On the initial time interval $[0,1]$, by (7.13) one has 
$$\int_0^1\int  \big|z_x(t,x)\,z^\sharp(t+\tau,\,x+\xi)\big|\,dxdt\leq
\int_0^1\big\|z_x(t)\big\|_{\L^1}\,\big\|z^\sharp(t+\tau)\big\|_{\L^\infty}\,dt
\leq (C^*\delta_0^2)^2\,.\eqno(7.17)$$
Because of (7.14) we have
$$2\|\lambda\|_{\L^\infty}\leq 2\|\lambda\|_{\L^1}\leq 2C^*\delta_0<\!<1\,.$$
{}From (7.16) and (7.17) it thus follows
$$\I(T)\leq (C^*\delta_0^2)^2+{8\delta_0^2\over c}+2C^*\delta_0\cdot \I(T)+
2C^*\delta_0^3\,.$$
For $\delta_0$ sufficiently small, this implies $\I(T)\leq 9\delta_0^2/c$,
proving the lemma.
\endproof
\v
Using the two previous lemmas we now prove the estimate
(7.3). 
Setting $z\doteq v_j$, $z^\sharp\doteq
v_k$,
$\lambda\doteq\tilde\lambda_j$, $\lambda^\sharp\doteq\tilde\lambda_k$,
an application of Lemma 7.1 yields the desired bound on the  
integral of $|v_jv_k|$.  Moreover, Lemma 7.3 allows us to estimate
the integral of $|v_{j,x}v_k|$. Notice that the assumptions
(7.13)-(7.14) are a consequence of (5.22)-(5.23). 
The simplifying condition $\lambda(t,-\infty)=0$ in (7.14)
can be easily achieved, using a new space coordinate
$x'\doteq x-\lambda_j^*t$.

The other terms $|v_jw_k|$, $|w_jw_k|$, $|v_{j,x}w_k|$ and $|v_jw_{k,x}|$
are handled similarly.
\vsk
\n{\medbf 8 - Functionals related to shortening curves}
\v
We now study the interaction of viscous waves of the same family.
As in the previous section, 
let $u=u(t,x)$ be a solution of the
parabolic system (3.1)
whose total variation remains bounded according to (7.1).
Assume that the components
$v_i,w_i$ satisfy the evolution equation (6.1), with source terms
$\phi_i,\psi_i$ bounded as in (7.2).
Relying on the bounds (5.22)--(5.24),
for each $i=1,\ldots,n$ we shall prove the estimates
$$\int_{\hat t}^T\!\int |w_{i,x}v_i-w_iv_{i,x}|\,dxdt=\O(1)\cdot\delta_0^2,
\eqno(8.1)$$
$$
\int_{\hat t}^T\!\int_{|w_i/v_i|<3\delta_1}
|v_i|^2\left| \left(w_i\over v_i
\right)_x\right|^2\,dxdt=\O(1)\cdot \delta_0^3.\eqno(8.2)$$
\v
The above integrals will be controlled in terms of two functionals, 
related to shortening curves.
Consider a parametrized curve in the plane
$\gamma:\R\mapsto\R^2$.   Assuming that $\gamma$ is
sufficiently smooth, its {\bf length} is computed by
$$\Le(\gamma)\doteq\int \big|\gamma_x(x)\big|\,dx\,.\eqno(8.3)$$
Following [BiB2], we also define the {\bf area} functional
as the integral of a wedge product:
$$\A(\gamma)\doteq {1\over 2}
\dint_{x<y}  \big| \gamma_x(x)\wedge\gamma_x(y)\big|
\,dxdy\,.\eqno(8.4)$$
To understand its geometrical meaning, observe that if $\gamma$ 
is a closed curve, the integral 
$${1\over 2}\int\gamma(y)\wedge\gamma_x(y)\,dy
={1\over 2}\dint_{x<y} \gamma_x(x)\wedge
\gamma_x(y)\,dx\,dy$$
yields the sum of the areas of the regions enclosed by the curve $\gamma$,
multiplied by the corresponding winding number (fig.~8).
In general, the quantity $\A(\gamma)$ provides un upper bound
for the area of the convex hull of $\gamma$.

\midinsert
\vskip 10pt
\centerline{\hbox{\psfig{figure=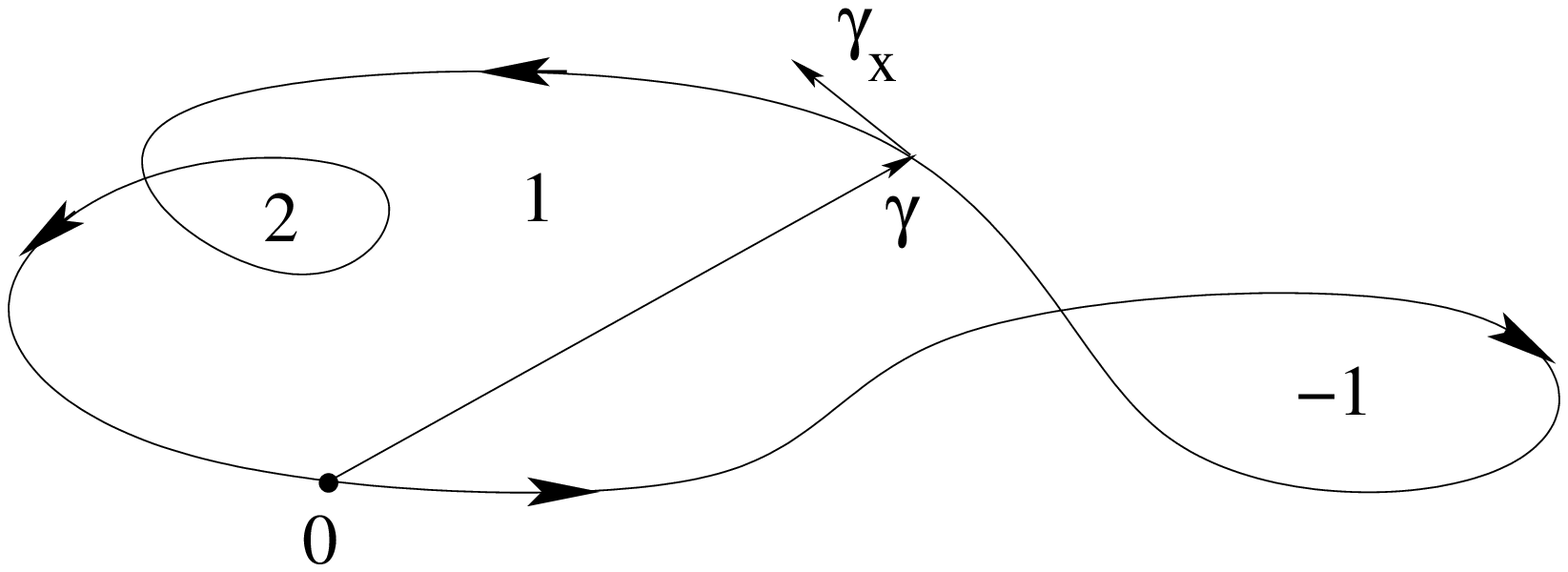,width=8cm}}}
\centerline{\hbox{figure 8}}
\vskip 10pt
\endinsert

Let now $\gamma=\gamma(t,x)$ be a planar curve which evolves in time,
according to the vector equation
$$\gamma_t+\lambda\gamma_x=\gamma_{xx}\,.\eqno(8.5)$$
Here $\lambda=\lambda(t,x)$ is a sufficiently smooth scalar function.
It is then clear that the length $\Le\big(\gamma(t)\big)$ of the curve
is a decreasing function of time. It was shown in [BiB2]
that also the area functional $\A\big(\gamma(t)\big)$
is monotonically decreasing.  Moreover,
the amount of decrease
dominates the area swept by the curve during its motion.
An intutive way to see this is the following.
In the special case where $\gamma$ is a polygonal line,
with vertices at the points $P_0,\ldots,P_m$,
the integral in (8.4) reduces to a sum:
$$\A(\gamma)=
{1\over 2}\sum_{i<j}\big|\bfv_i\wedge\bfv_j\big|\,,
\qquad\qquad \bfv_i\doteq P_i-P_{i-1}\,.$$
If we now replace $\gamma$ by a new curve $\gamma'$ obtained by
replacing two consecutive edges $\bfv_h$, $\bfv_k$
by one single edge (fig.~9b), the area between $\gamma$ and $\gamma'$
is precisely $|\bfv_h\wedge \bfv_k|/2$, while an easy computation
yields
$$\A(\gamma')\leq \A(\gamma)-{1\over 2}|\bfv_h\wedge \bfv_k|\,.$$
The estimate on the area swept
by a smooth curve (fig.~9a) is now obtained by approximating
a shortening curve $\gamma$ by a sequence of polygonals,
each obtained from the previous one by replacing two consecutive edges
by a single segment.

\midinsert
\vskip 10pt
\centerline{\hbox{\psfig{figure=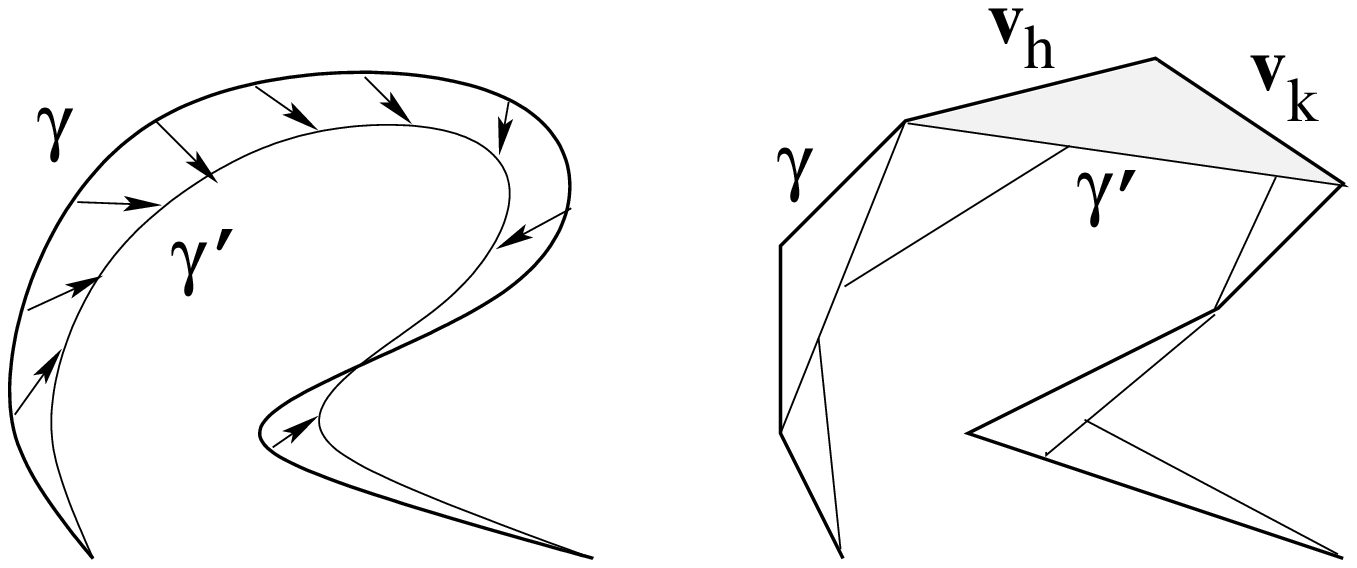,width=9cm}}}
\centerline{\hbox{figure 9a~~~~~~~~~~~~~~~~~~~~~~~~~~figure 9b}}
\vskip 10pt
\endinsert

We shall apply the previous geometric considerations 
toward a proof of the estimates of (8.1)-(8.2).
Let $v,w$ be two scalar functions, satisfying
$$\eqalign{v_t+(\lambda v)_x-v_{xx}&=\phi\,,\cr
w_t+(\lambda w)_x-w_{xx}&=\psi\,.\cr}\eqno(8.6)$$
Define the planar curve $\gamma$ by setting
$$
\gamma(t,x) = \left( \int_{-\infty}^x v(t,y) dy, ~~\int_{-\infty}^x
w(t,y) dy \right).\eqno(8.7)$$
Integrating (8.6) w.r.t.~$x$, one finds 
the corresponding evolution equation for $\gamma$ :
$$\gamma_t + \lambda \gamma_x - \gamma_{xx} = \left(
\int_{-\infty}^x \phi(t,y) dy, ~~\int_{-\infty}^x \psi(t,y) dy
\right).\eqno(8.8)$$
In particular, if no sources were present, the 
motion of the curve would reduce to (8.5).
At each fixed time $t$, we now define 
the {\bf Length Functional} as
$$\Le(t)=\Le\big(\gamma(t)\big)=\int\sqrt{v^2(t,x)+w^2(t,x)
}~dx\eqno(8.9)$$
and
the {\bf Area Functional} as
$$\A(t)=\A\big(\gamma(t)\big)= {1\over 2}
\dint_{x<y}  \big| v(t,x) w(t,y) -
v(t,y)w(t,x) \big| \,dxdy\,.\eqno(8.10)$$
\v
We now estimate the time derivative of the above functionals, in the
general case when sources are present.
\v
\n{\bf Lemma 8.1} {\it Let $v,w$ be solutions of (8.3), defined for
$t\in [0,T]$.  For each $t$, assume that
the maps $x\mapsto v(t,x)$
$x\mapsto w(t,x)$ and $x\mapsto \lambda(t,x)$ are $\C^{1,1}$, 
i.e.~continuously differentiable with Lipschitz derivative.
Then the corresponding area functional (8.10) satisfies}
$${d\over dt}\A(t)\leq - \int \Big| v_x(t,x) w(t,x) - v(t,x)
w_x(t,x) \Big| dx + \big\| v(t) \big\|_{L^1} \big\| \psi(t)
\big\|_{L^1} + \big\| w(t) \big\|_{L^1} \big\| \phi(t)
\big\|_{L^1}\,.\eqno(8.11)$$
\v
\n{\bf Proof.} 
Differentiating (8.10) w.r.t. time we obtain
$$\eqalign{
&{d\A\over dt}={1\over 2} \int\!\!\int_{x<y} \sgn \bigl( v(x)
w(y) - w(x) v(y) \bigr) \cdot \Bigl\{ v_t(x) w(y) + v(x) w_t(y) -
v_t(y) w(x) - v(y) w_t(x) \Bigr\} dxdy \cr
&= \int\!\!\int_{x<y} \sgn \bigl( v(x) w(y) - w(x) v(y) \bigr)
\cdot \Bigl\{ \bigl( v_x(x) - \lambda(x) v(x) \bigr)_x w(y) - v(y)
\bigl( w_x(x) - \lambda(x) w(x) \bigr)_x \Bigr\} dxdy \cr
& \qquad + \int\!\!\int_{x<y} \sgn \bigl( v(x)
w(y) - w(x) v(y) \bigr) \cdot \Bigl\{ \phi(x) w(y) + v(x) \psi(y)
\Bigr\} dxdy \cr
&= \int\!\!\int_{x<y} \sgn \bigl( v(x) w(y) - w(x) v(y) \bigr)
\cdot {\partial\over\partial x}\Bigl\{ \bigl( v_x(x) - \lambda(x)
v(x) \bigr) w(y) - v(y)
\bigl( w_x(x) - \lambda(x) w(x) \bigr) \Bigr\} dxdy \cr
& \qquad + \int\!\!\int_{x<y} \sgn \big( v(x)
w(y) - w(x) v(y) \big) \cdot \Big\{ \phi(x) w(y) + v(x) \psi(y)
\Big\} dxdy\,. \cr}\eqno(8.12)
$$
To simplify the first term on the right hand side of (8.12),
for a fixed $y$ we let $\big\{x_\alpha(y)\big\}$ be the set of
points $\leq y$ such that
$v(x_\alpha)w(y)-w(x_\alpha)v(y)=0$.
Relying on an approximation argument, we can assume that these points
are in finite number,  say $x_N(y) <\cdots< x_2(y) < x_1(y) < x_0(y)=y$,
and that the function $x\mapsto v(x) w(y) - w(x) v(y)$
changes sign across each $x_\alpha(y)$.  
For convenience, we define the additional point $x_{N+1}\doteq -\infty$
and let $\iota(y)=\pm 1$ be the sign of 
$v(x)w(y)-v(y)w(x)$ on the last interval, i.e.~when $x\in [x_1,\,y]$.
We now can write
$$\eqalign{
\int_{-\infty}^y &\sgn \bigl( v(x) w(y) - w(x) v(y) \bigr)
\cdot {\partial\over\partial x}\Bigl\{ \bigl( v_x(x) - \lambda(x)
v(x) \bigr) w(y) - v(y)
\bigl( w_x(x) - \lambda(x) w(x) \bigr) \Bigr\} dx \cr
&=\iota(y)\cdot\sum_{\alpha=0}^N 
( - 1 )^\alpha \int_{x_{\alpha+1}(y)}^{x_\alpha(y)} 
{\partial\over\partial x}\Bigl\{ \bigl( v_x(x) - \lambda(x)
v(x) \bigr) w(y) - v(y) \bigl( w_x(x) - \lambda(x) w(x) \bigr) \Bigr\}
dx \cr
&=2\, \iota(y)\cdot\left(\sum_{\alpha=1}^N( - 1)^\alpha 
\Bigl( v_x \bigl( x_\alpha(y)
\bigr) w(y) - w_x \bigl(x_\alpha(y) \bigr) v(y) \Bigr) +
\big( v_x(y) w(y) - w_x(y) v(y) \big)\right).
\cr}
$$
By direct inspection, one checks that the factor
$\iota(y)\cdot (-1)^\alpha$
has always the opposite
sign of $v_x(x_\alpha) w(y) - w_x(x_\alpha) v(y)$.  Therefore
$$\eqalign{
\int\!\!\int_{x<y} &\sgn \bigl( v(x) w(y) - w(x) v(y) \bigr)
\cdot {\partial\over\partial x}\Bigl\{ \bigl( v_x(x) - \lambda(x)
v(x) \bigr) w(y) - v(y)
\bigl( w_x(x) - \lambda(x) w(x) \bigr) \Bigr\} dxdy\cr
&=~-  \int \Bigl| v_x(y) w(y) - v(y) w_x(y) \Bigr| dy - 2 \sum_{\alpha=1}
^N
\int \Bigl| v_x\big(x_\alpha(y)\big) w(y) - v(y) w_x\big(x_\alpha(y)
\big) \Bigr| dy\,.
\cr}\eqno(8.13)
$$
The bound (8.11) is now an immediate consequence of (8.12) and (8.13).
\endproof
\v
\n{\bf Lemma 8.2.} {\it Together with the hypotheses of Lemma 8.1, 
at a fixed time $t$ assume that $\gamma_x(t,x)\not=0$ for every $x$.
Then}
$${d\over dt}\Le(t)\leq 
-{1\over( 1 + 9 \delta_1^2)^{3/2}} \int_{|w/v| \leq 3
\delta_1} \big|v(t)\big|
\,\left| \left(w(t)\over v(t)
\right)_x\right|^2 dx+\big\| \phi(t) \big\|_{\L^1}
+ \big\| \psi(t) \big\|_{\L^1}\,.\eqno(8.14)$$
\v
\n{\bf Proof.} 
As a preliminary, we derive the identity
$$\eqalign{| \gamma_{xx}|^2 \,| \gamma_x|^2 
- \la \gamma_{x}, \,\gamma_{xx} \ra ^2 &= \big( v_x^2 + w_{xx}^2 \big)
( v^2 + w_x^2) -\big( v v_x + w_x w_{xx}
\big)^2 \cr
&= \big( v w_{xx}+ v_x w_x \big)^2 = v^4
\big|( w/v)_x\big|^2.
\cr}$$
Thanks to the assumption that $\gamma_x$ never vanishes,
we can now integrate by parts and obtain
$$\eqalign{
{d\over dt}\Le(t)&=
\int {\la \gamma_x, \,\gamma_{xt}
\ra\over\sqrt{\la \gamma_x,\, \gamma_x \ra }}~ dx
= \int\left\{ {\la \gamma_x, \,\gamma_{xxx}
\ra\over |\gamma_x|}
- {\la \gamma_x, \,( \lambda \gamma_x)_x
\ra\over |\gamma_x|}
+ {\la \gamma_x, \,( \phi, \psi) \ra\over |\gamma_x|}\right\}\,dx
\cr
&= \int \left\{ | \gamma_x |_{xx} - \big(
\lambda \,| \gamma_x | \big)_x
- { | \gamma_{xx} |^2 - \la \gamma_x/| \gamma_x|, \,
\gamma_{xx} \ra^2\over |\gamma_x|}\right\} dx + \int
{\la \gamma_x, \,( \phi, \psi) \ra\over |\gamma_x|}
\,dx\cr
&\leq - \int { |v| \,\big|( w/v)_x\big|^2\over\big( 1 +
( w/v )^2
\big)^{3/2}} \,dx + 
\big\| \phi(t) \big\|_{\L^1}
+ \big\| \psi(t) \big\|_{\L^1}\,.\cr}
$$
Since the integrand is non-negative, the last inequality clearly implies
(8.14).
\endproof
\v
\n{\bf Remark 8.3.}  Let $u=u(t,x)$ be a solution to a scalar, viscous
conservation law
$$u_t+f(u)_x-u_{xx}=0\,,$$
and consider the planar curve $\gamma\doteq \big( u,~f(u)-u_x\big)$ 
whose components are
respectively the conserved quantity and the flux (fig.~10).
If $\lambda\doteq f'$, the components $v\doteq u_x$ and $w\doteq -u_t$
evolve according to (8.6), with $\phi=\psi=0$, hence
$\gamma_t+\lambda \gamma_x-\gamma_{xx}=0$.
Defining the speed (fig.~11)
$s(x)\doteq - {u_t(x)/ u_x(x)}$,
the area functional $\A(\gamma)$ in (8.4) can now be written as
$$\eqalign{\A(\gamma)
&={1\over 2} \dint_{x<y} \big|u_x(x)u_t(y)-u_t(x)u_x(y)\big|
\,dxdy\cr
&={1\over 2} \dint_{x<y} \big|u_x(x)\,dx\big|\cdot  \big|u_x(y)\,dy\big|
\cdot \big|s(x)-s(y)\big|\cr
&={1\over 2} \dint_{x<y} [\hbox{wave at}~x]\times
[\hbox{wave at}~y]\times [\hbox{difference in speeds}]\,.\cr}$$
It now becomes clear that the area functional can be
regarded as an interaction
potential between waves of the same family. In the case where
viscosity is present, this provides a
counterpart to the interaction functional introduced in [L2]
in connection with strictly hyperbolic systems.

\midinsert
\vskip 10pt
\centerline{\hbox{\psfig{figure=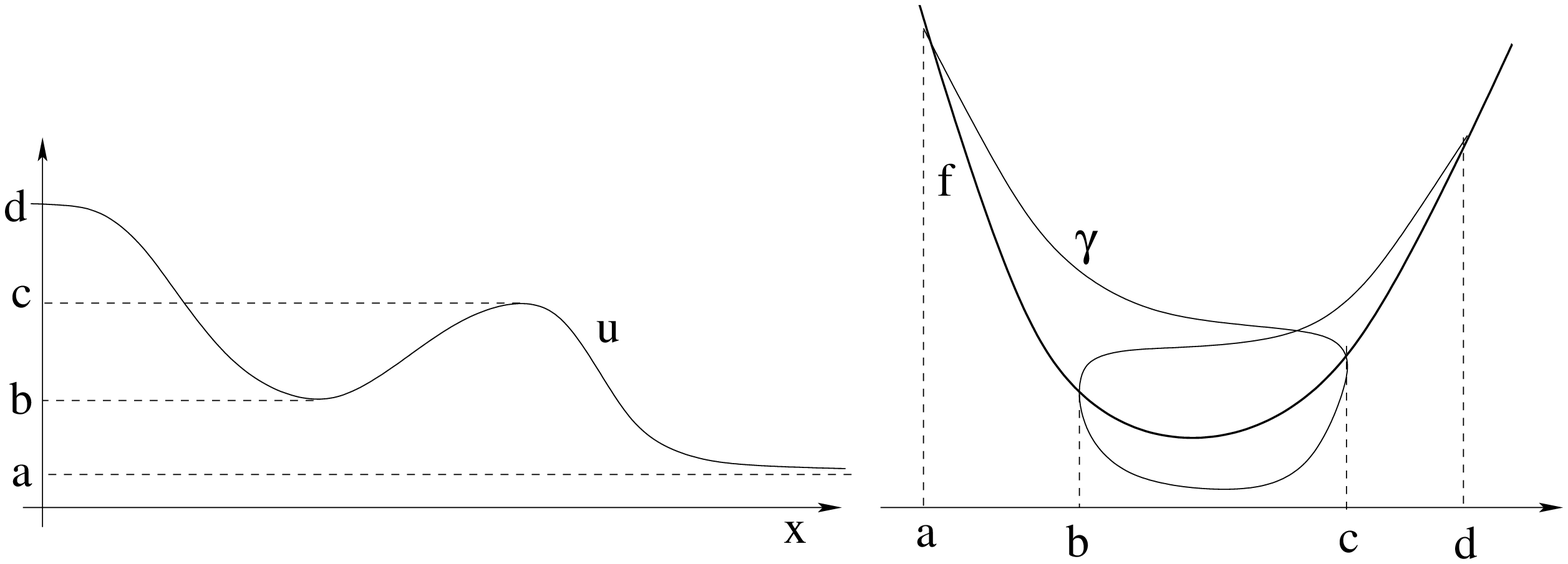,width=12cm}}}
\c{\hbox{figure 10}}
\vskip 10pt
\endinsert

\midinsert
\vskip 10pt
\centerline{\hbox{\psfig{figure=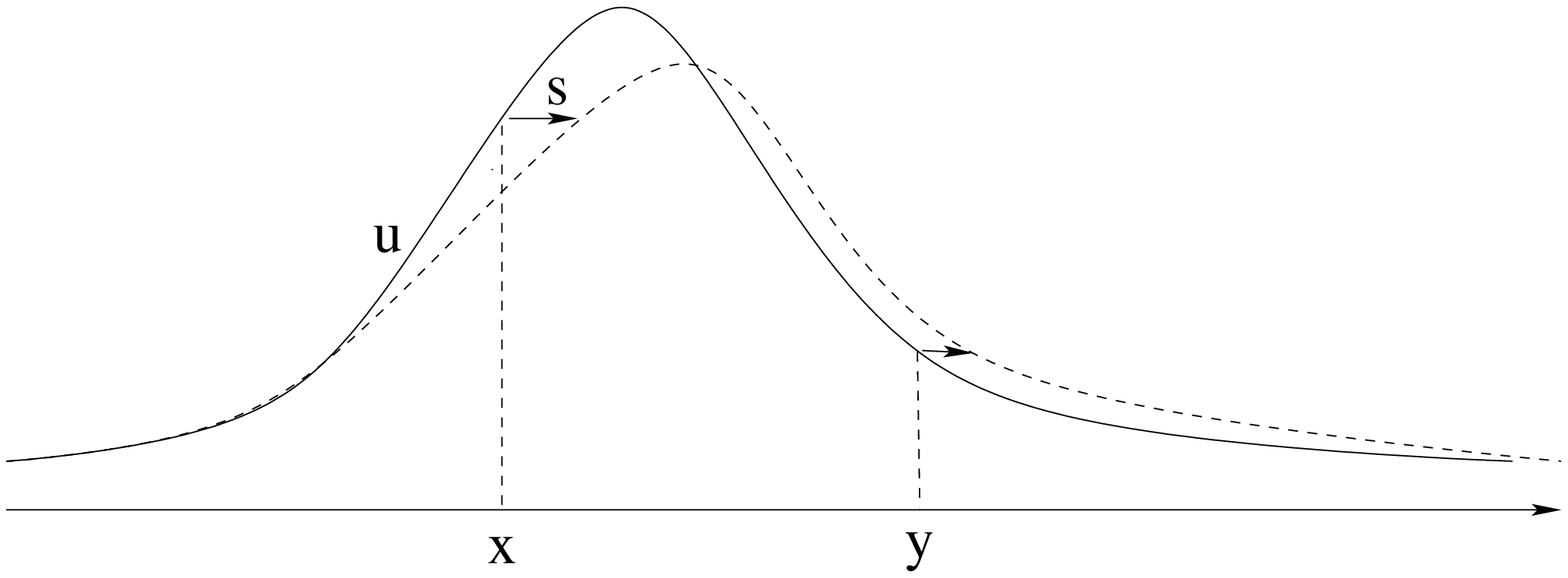,width=10cm}}}
\centerline{\hbox{figure 11}}
\vskip 10pt
\endinsert

\v
Recalling that the components $v_i,w_i$ satisfy the equations (6.1), 
we can apply the previous lemmas with 
$v\doteq v_i$, $w\doteq w_i$, $\lambda\doteq \tilde\lambda_i$,
$\phi\doteq\phi_i$, $\psi\doteq\psi_i$, 
calling $\Le_i$ and $\A_i$ the corresponding length and area functionals.
For $t\in [\hat t,\,T]$,
the bounds (5.22)--(5.23) yield
$$\A_i(t)\leq \big\|v_i(t)\big\|_{\L^\infty}\cdot\big\|w_i(t)\big\|_{\L^1}
=\O(1)\cdot\delta_0^3,\eqno(8.15)$$
$$\Le_i(t)\leq \big\|v_i(t)\big\|_{\L^1}+\big\|w_i(t)\big\|_{\L^1}
=
\O(1)\cdot\delta_0.\eqno(8.16)$$
Using (8.11) we now obtain
$$\eqalign{\int_{\hat t}^T 
\!\int \big|w_{i,x}v_i&-w_iv_{i,x}\big|\,dxdt
\leq
\int_{\hat t}^T \left|{d\over dt}\A_i(t)\right|\,dt+\int_
{\hat t}^T\Big( \big\| v_i(t) \big\|_{L^1} \big\| \psi_i(t)
\big\|_{L^1} + \big\| w_i(t) \big\|_{L^1} \big\| \phi_i(t)
\big\|_{L^1}\Big)\,dt\cr
&\leq \A_i(\hat t)+\sup_{t\in [\hat t,T]}
\Big(\big\|v_i(t)\big\|_{\L^1}+\big\|w_i(t)\big\|_{\L^1}\Big)
\cdot 
\int_{\hat t}^T\!\int \Big( \big|\phi_i(t,x)\big|+
\big|\psi_i(t,x)\big|\Big)\,dxdt
\cr
&=\O(1)\cdot \delta_0^2\,,\cr}\eqno(8.17)$$
proving (8.1).
To establish (8.2), 
we first observe that, by an approximation argument, it is
not restrictive to assume that 
the set of points in the $t$-$x$ plane where
$v_{i,x}(t,x)=w_{i,x}(t,x)=0$ is at most 
countable. In this case, for almost every $t\in [\hat t,T]$
the inequality (8.14) is valid, and hence
$$\eqalign{
\int_{\hat t}^T\!\int_{|w_i/v_i|<3\delta_1}
|v_i|\,\left| \left(w_i\over v_i
\right)_x\right|^2\,dxdt
&\leq \int_{\hat t}^T\left|{d\over dt}\Le_i(t)\right|dt
+\int_{\hat t}^T\Big(\big\| \phi_i(t)
\big\|_{L^1} + \big\| \psi_i(t)
\big\|_{L^1}\Big)\,dt\cr
&\leq \Le_i(\hat t)+
\int_{\hat t}^T\!\int \Big( \big|\phi_i(t,x)\big|+\big|\psi_i(t,x)\big|
\Big)dxdt\cr
&=\O(1)\cdot \delta_0.\cr}\eqno(8.18)$$
Using the bound (5.23) on $\|v_i\|_{\L^\infty}$, from (8.18)
we deduce (8.2).
\vsk
\n{\medbf 9 - Energy estimates}
\v
In the same setting as the two previous sections, we shall
now prove the estimate
$$ \int_{\hat t}^T\!\int 
\big( |v_{i,x}|+|w_{i,x}|\big)\,|w_i-\theta_i v_i|\,dxdt
=\O(1)\cdot \delta_0^2\,.\eqno(9.1)$$
We recall that
$\theta_i\doteq \theta(w_i/v_i)$, where $\theta$ is the cut-off
function introduced in (5.5).  Notice that the integrand 
can be $\not= 0$ only when $|w_i/v_i|>\delta_1$.
\v
Consider another cut-off function $\vth:\R\mapsto [0,1]$
such that (fig.~5)
$$
\vth(s)=\cases{0\qquad &if\qquad $|s| \leq 3 \delta_1/5\,$,\cr
1\qquad &if\qquad $|s| \geq 4 \delta_1/5\,$.\cr}\eqno(9.2)$$ 
We can assume that $\vth$ is a smooth even function,
such that 
$$| \vth'|\leq 21/\delta_1\,,\qquad\qquad
|\vth''|\leq 101/\delta_1^2\,.$$
For convenience, we shall write $\vth_i\doteq\vth(w_i/v_i)$.
As a preliminary, we prove some simple estimates relating the
sizes of $v_i$, $w_i$ and $v_{i,x}$. 
It is here useful to
keep in mind the bounds
$$\tilde\lambda_i-\lambda_i^*=\O(1)\cdot |\tr_i-r_i^*|=
\O(1)\cdot\delta_0\,,\qquad\qquad
|v_i|,\,|w_i|=\O(1)\cdot \delta_0^2\,,\eqno(9.3)$$
valid for $t\geq \hat t$ and $i=1,\ldots,n$.
Recall also our choice of the constants
$$0<\delta_0<\!<\delta_1\leq {1\over 3}\,.\eqno(9.4)$$
\v
\n{\bf Lemma 9.1.} {\it If $|w_i/v_i|\geq 3\delta_1/5$, then
$$|w_i|\leq 2|v_{i,x}|+\O(1)\cdot\delta_0\sum_{j\not= i} |v_j|,\qquad\qquad
|v_i|\leq {5\over 2\delta_1} |v_{i,x}|+\O(1)\cdot
\delta_0\sum_{j\not= i}|v_j|\,.
\eqno(9.5)$$
On the other hand, if  $|w_i/v_i|\leq 4\delta_1/5$,
then}
$$|v_{i,x}|\leq \delta_1 |v_i|+\O(1)\cdot
\delta_0\sum_{j\not= i}|v_j|\,.\eqno(9.6)$$
\v
\n{\bf Proof.} We recall the first estimate in (6.16): 
$$\eqalign{v_{i,x}&=w_i+(\tilde\lambda_i-\lambda^*_i)v_i+\Theta\cr
&=w_i+(\tilde\lambda_i-\lambda^*_i)v_i+\O(1)\cdot\delta_0
\sum_{j\not= i}|v_j|\,,\cr}\eqno(9.7)$$
with $\Theta$ defined as in (6.15).
By (9.3)-(9.4), 
from the condition $|w_i/v_i|\geq 3\delta_1/5$ two cases can arise.
On one hand, if
$$|\Theta|\leq {\delta_1\over 10}|v_i|\,,\eqno(9.8)$$
then
$$|v_{i,x}|\geq {3\delta_1\over 5}|v_i|-\O(1)\cdot \delta_0 |v_i|
-{\delta_1\over 10}|v_i|\geq {2\delta_1\over 5}|v_i|\,,$$
and hence
$$|v_i|\leq {5\over 2\delta_1} |v_{i,x}|,\qquad\qquad
|w_i|\leq  |v_{i,x}|+{\delta_1\over 5}|v_i|\leq 2|v_{i,x}|\,.
\eqno(9.9)$$
On the other hand, if (9.8) fails, then by (9.3) we have
$$\big|\tla_i-\lambda^*_i\big|\,|v_i|=
\O(1)\cdot\delta_0^2\sum_{j\not= i} |v_j|\,.\eqno(9.10)$$
In both cases, the estimates in (9.5) hold.
\v
Next, if $|w_i/v_i|\leq 4\delta_1/5$, from (9.7) we deduce
$$|v_{i,x}|\leq { 4\delta_1\over 5} |v_i|+{\delta_1\over 10}|v_i|
+\O(1)\cdot\delta_0\sum_{j\not= i} |v_j|\,.\eqno(9.11)$$
If (9.8) holds, then $|v_{i,x}|\leq\delta_1|v_i|$.
If (9.8) fails, then (9.10) is valid.  In both cases 
we have (9.6).
\endproof
\v
Toward a proof of the estimate (9.1),
we first reduce the integrand to
a more tractable expression.
Since the term $|w_i-\theta_iv_i|$ vanishes when $|w_i/v_i|\leq\delta_1$,
and is $\leq |w_i|$ otherwise, by (9.5) we always have the bound
$$ |w_i-\theta_iv_i|\leq |\vth_i w_i|\leq
\vth_i\left(2|v_{i,x}|+\O(1)\cdot\delta_0\sum_{j\not= i}|v_j|\right).$$
Therefore
$$\eqalign{
\big( |v_{i,x}|+|w_{i,x}|\big)\cdot|w_i-\theta_iv_i|
&
\leq 
\big( |v_{i,x}|+|w_{i,x}|\big)
\vth_i\left(2|v_{i,x}|+\O(1)\cdot\delta_0\sum_{j\not= i}|v_j|\right)\cr
&\leq
2\vth_iv_{i,x}^2+2\vth_i|v_{i,x}w_{i,x}|+\sum_{j\not= i}
\big( |v_jv_{i,x}|+|v_jw_{i,x}|\big)\cr
&\leq
3\vth_iv_{i,x}^2+\vth_i w_{i,x}^2+\sum_{j\not= i}
\big( |v_jv_{i,x}|+|v_jw_{i,x}|\big)\cr
\,.\cr}\eqno(9.12)$$
Since we already proved the bounds (7.3)
on the integrals of transversal terms, to prove
(9.1) we only need to consider the integrals of $v_{i,x}^2$ and
$w_{i,x}^2$, in the region where $\vth_i\not=0$. 
In both cases, energy type estimates will be used. 
\v
We start with $v_{i,x}^2$. 
Multiplying the first equation in (6.1) by $\vth_iv_i$ and
integrating by parts, we obtain
$$\eqalign{\int\vth_iv_i\phi_i\,dx&=
\int\Big\{\vth_iv_i v_{i,t}+\vth_i
v_i(\tla_iv_i)_x-\vth_iv_iv_{i,xx}
\Big\}\,dx\cr
&=\int\Big\{\vth_i
(v_i^2/ 2)_t
-\vth_i\tla_i v_iv_{i,x}-\vth_{i,x}\tla_iv_i^2+
\vth_i v_{i,x}^2+
\vth_{i,x}v_{i,x}v_i\Big\}\,dx\cr
&=\int\Big\{\big(\vth_iv_i^2/2 \big)_t+(\tla_i\vth_i)_x(v_i^2/2)
-\big(\vth_{i,t}
+2\tla_i\vth_{i,x}-\vth_{i,xx}\big)(v_i^2/2) +\vth_i v_{i,x}^2
\Big\}\,dx\,.\cr}$$
Therefore
$$\eqalign{ \int \vth_i v_{i,x}^2\,dx&=
-{d\over dt}\left[\int\vth_i v_i^2/2\,dx\right]
+\int \big(\vth_{i,t}
+\tla_i\vth_{i,x}-\vth_{i,xx}\big)(v_i^2/2)\,dx\cr
&\qquad -\int \tla_{i,x}\vth_i(v_i^2/2)\,dx+
\int\vth_iv_i\phi_i\,dx\,.\cr}
\eqno(9.13)$$
A direct computation yields
$$\eqalign{\vth_{i,t}
+\tla_i\vth_{i,x}-\vth_{i,xx}
&=\vth_i'\left({w_{i,t}\over v_i}-{v_{i,t}w_i\over v_i^2}
\right)+\tilde\lambda_i\vth_i'\left({w_{i,x}\over v_i}-
{v_{i,x}w_i\over v_i^2}\right)\cr
&\qquad -\vth''_i\left({w_i\over v_i}\right)_x^2
-\vth'\left({w_{i,xx}\over v_i}-{v_{i,xx}w_i\over v_i^2}
-2{v_{i,x}w_{i,x}\over v_i^2}+2{v_{i,x}^2w_i\over v_i^3}\right)\cr
&= \Big[ \vth_i' \big( w_{i,t} + ( \tilde
\lambda_i w_i)_x - w_{i,xx} \big) / v_i - \vth'_i w_i \big(
v_{i,t} + ( \tilde \lambda_i v_i)_x - v_{i,xx} \big)/v_i^2 \Big] \cr
&\qquad + 2 v_{i,x} \vth_i'/v_i \cdot ( w_i/v_i)_x 
- \vth''_i ( w_i / v_i )_x^2 \cr
&=\vth'_i \left( {\psi_i\over v_i} - {w_i\over v_i}
{\phi_i\over v_i} \right) + 2 \vth'_i {v_{i,x}\over v_i} \left(
{w_i\over v_i} \right)_x
- \vth''_i \left( {w_i\over v_i} \right)_x^2\,.\cr}
\eqno(9.14)$$
Since $\tla_{i,x}=(\tla_i-\lambda_i^*)_x$,
integrating by parts and using the second estimate in (9.5) one obtains
$$\eqalign{&\left|\int \tla_{i,x}\vth_i(v_i^2/2)\,dx\right|=
\left|\int (\tla_i-\lambda_i^*) \big(\vth_{i,x}v_i^2/2+
\vth_iv_iv_{i,x}
\big)\,dx
\right|\cr
&\qquad\leq \|\tla_i-\lambda_i^*\|_{\L^\infty}\cdot\left\{{1\over 2}\int
\big|\vth_i'\big|\,|w_{i,x}v_i-v_{i,x}w_i|\,dx+{5\over 2\delta_1} 
\int \vth_i v_{i,x}^2\,dx+\O(1)\cdot
\delta_0\int\sum_{j\not= i}|v_iv_j|\,dx\right\}\cr
&\qquad\leq \int
|w_{i,x}v_i-v_{i,x}w_i|\,dx+{1\over 2}\int \vth_i v_{i,x}^2\,dx+
\delta_0\int\sum_{j\not= i}|v_iv_j|\,dx\,.
\cr}
\eqno(9.15)$$
Indeed, by (9.3)-(9.4), $|\tla_i-\lambda_i^*|=\O(1)\cdot\delta_0<\!<\delta_1$.
Using (9.14)-(9.15) in (9.13) we now obtain
$$\eqalign{
{1\over 2}&\int \vth_i\,v_{i,x}^2 \, dx ~\leq~ -
{d\over dt}\left[ \int {\vth_i 
v_i^2\over 2} \,dx\right] 
+{1\over 2}\int \big|\vth_i'\big|\big(|v_i\psi_i|+|w_i\phi_i|\big)\,dx
+\int \left| \vth_i' v_iv_{i,x}\left(w_i\over v_i\right)_x
\right|\,dx
\cr
&+{1\over 2}\int\left|\vth_i''v_i^2
\left(w_i\over v_i\right)_x^2\right|\,dx
+\int
\big| w_{i,x}
v_i -w_i v_{i,x} \big|\, dx+\delta_0\int\sum_{j\not=i}|v_iv_j|\,dx
+\int |v_i\phi_i|\,dx\,.\cr}\eqno(9.16)$$
Recalling the definition of $\vth_i$,
on regions where $\vth_i'\not=0$ one has $|w_i/v_i|\leq 4\delta_1/5$,
hence the bounds (9.6) hold. In turn, they imply
$$\eqalign{
\left|\vth_i' v_iv_{i,x}\left(w_i\over v_i\right)_x
\right|&\leq
\left|\delta_1 \vth_i' v_i^2\left(w_i\over v_i\right)_x
\right|+\O(1)\cdot
\delta_0\sum_{j\not= i}\left| \vth_i' v_iv_j\left(w_i\over v_i\right)_x
\right|\cr
&\leq \big|\delta_1\vth_i'\big|\, |w_{i,x}v_i-w_iv_{i,x}|
+\O(1)\cdot
\delta_0|\vth_i'|\sum_{j\not= i}\left(|v_jw_{i,x}|+|v_jv_{i,x}|\left|
w_i\over v_i\right|\right)\,.\cr}\eqno(9.17)$$
Using the bounds (5.22)-(5.23). (7.2)-(7.3) and (8.1)-(8.2), from (9.16) we
conclude
$$\eqalign{
\int_{\hat t}^T \!&\int  \vth_i\,v_{i,x}^2 \, dxdt
\leq \int\vth_i v_i^2(\hat t,x)\,dx
+\O(1)\cdot\int_{\hat t}^T\!\int \big(|v_i\psi_i|+|w_i\phi_i|\big)\,dxdt
\cr
&\qquad
+\O(1)\cdot\int_{\hat t}^T\!\int|w_{i,x}v_i-w_iv_{i,x}|\,dxdt
+\O(1)\cdot\delta_0\int_{\hat t}^T\!\int 
\sum_{j\not= i}\big(|v_jw_{i,x}|+|v_jv_{i,x}|\big)\,dxdt
\cr
&\qquad+\O(1)\cdot\int_{\hat t}^T\!\int_{|w_i/v_i|<\delta_1}
\big|v_i(w_i/ v_i)_x\big|^2\,dxdt
+2\delta_0\int_{\hat t}^T\!\int\sum_{j\not=i}|v_iv_j|\,dxdt
+2\int_{\hat t}^T\!\int |v_i\phi_i|\,dxdt\cr
&=\O(1)\cdot \delta_0^2\,.\cr}\eqno(9.18)$$
\v
We now perform a similar computation for $w_{i,x}^2$.
Multiplying the second equation in (6.1) by $\vth_iw_i$ and
integrating by parts, we obtain
$$\int\vth_iw_i\psi_i\,dx
=\int\Big\{\big(\vth_iw_i^2/2 \big)_t+(\tla_i\vth_i)_x(w_i^2/2)
-\big(\vth_{i,t}
+2\tla_i\vth_{i,x}-\vth_{i,xx}\big)(w_i^2/2) +\vth_i w_{i,x}^2
\Big\}\,dx\,.$$
Therefore, the identity (9.13) still holds, with
$v_i,\phi_i$ replaced by $w_i,\psi_i$, respectively:
$$\eqalign{ \int \vth_i w_{i,x}^2\,dx&=
-{d\over dt}\left[\int\vth_i w_i^2/2\,dx\right]
+\int \big(\vth_{i,t}
+\tla_i\vth_{i,x}-\vth_{i,xx}\big)(w_i^2/2)\,dx\cr
&\qquad -\int \tla_{i,x}\vth_i(w_i^2/2)\,dx+
\int\vth_i w_i\psi_i\,dx\,.\cr}
\eqno(9.19)$$
The equality (9.14) can again be used.  To obtain a suitable
replacement for (9.15), we observe that, if $\vth_i\not= 0$
then (9.5) implies
$$|w_iw_{i,x}|\leq 2|v_{i,x}w_{i,x}|+\delta_0 \sum_{j\not= i}
|v_jw_{i,x}\big|\leq v^2_{i,x} + w_{i,x}^2
+\O(1)\cdot\delta_0 \sum_{j\not= i}
|v_jw_{i,x}|\,.$$
Integrating by parts we thus obtain
$$\eqalign{\left|\int \tla_{i,x}\vth_i(w_i^2/2)\,dx\right|&=
\left|\int (\tla_i-\lambda_i^*) \big(\vth_{i,x}w_i^2/2+
\vth_iw_iw_{i,x}
\big)\,dx
\right|\cr
&\leq \|\tla_i-\lambda_i^*\|_{\L^\infty}\cdot\bigg\{\int
\big|\vth_i'\big|\,|w_{i,x}v_i-v_{i,x}w_i|\,\left|w_i^2\over v_i^2
\right|\,dx
+\int \vth_i v_{i,x}^2\,dx\cr
&\qquad+
\int \vth_i w_{i,x}^2\,dx
+\O(1)\cdot\delta_0\int\sum_{j\not= i}|v_j w_{i,x}|\,dx\bigg\}\cr
&\leq 
\int |w_{i,x}v_i-v_{i,x}w_i|\,dx
+{1\over 2}\int \vth_i v_{i,x}^2\,dx
+{1\over 2}\int \vth_i w_{i,x}^2\,dx
+\delta_0\int\sum_{j\not= i}|v_j w_{i,x}|\,dx
\,.\cr}
\eqno(9.20)$$
Using (9.14) and (9.19) in (9.18) and observing that
$|w_i^2/ v_i^2|\leq \delta_1^2$ on the region where $\vth_i'\not=0$, 
we now obtain an estimate similar to (9.16):
$$\eqalign{
{1\over 2}\int \vth_i\,w_{i,x}^2 \, dx &\leq -
{d\over dt}\left[ \int {\vth_i 
w_i^2\over 2} \,dx\right] 
+{\delta_1^2\over 2}\int 
\big|\vth_i'\big|\big(|v_i\psi_i|+|w_i\phi_i|\big)\,dx
+\delta_1^2\int \left| \vth_i' v_iv_{i,x}\left(w_i\over v_i\right)_x
\right|\,dx
\cr
&\qquad +{\delta_1^2\over 2}\int\left|\vth_i''v_i^2
\left(w_i\over v_i\right)_x^2\right|\,dx
+\int\big| w_{i,x}v_i -w_i v_{i,x} \big|\, dx\cr
&\qquad +
{1\over 2}\int \vth_i v_{i,x}^2\,dx
+\delta_0\int\sum_{j\not= i}|v_jw_{i,x}|\,dx
+\int |w_i\psi_i|\,dx\,.\cr}\eqno(9.21)$$
Using the bounds (5.22)-(5.23), (7.2)-(7.3), (8.1)-(8.2) and (9.17)-(9.18), 
from (9.21) we
conclude
$$\eqalign{
\int_{\hat t}^T \!&\int \vth_i\,w_{i,x}^2 \, dxdt
\leq \int\vth_i w_i^2(\hat t,x)\,dx
+\O(1)\cdot\int_{\hat t}^T\!\int 
\big(|v_i\psi_i|+|w_i\phi_i|\big)\,dxdt
\cr
&\qquad
+\O(1)\cdot\int_{\hat t}^T\!\int|w_{i,x}v_i-w_iv_{i,x}|\,dxdt
+\O(1)\cdot\delta_0\int_{\hat t}^T\!\int 
\sum_{j\not= i}\big(|v_jw_{i,x}|+|v_jv_{i,x}|\big)\,dxdt
\cr
&\qquad+\O(1)\cdot\int_{\hat t}^T\!\int_{|w_i/v_i|<\delta_1} 
\big|v_i(w_i/v_i)_x\big|^2\,dxdt
+\int_{\hat t}^T\!\int \vth_i\,v_{i,x}^2 \, dxdt
\cr
&\qquad
+\delta_0\int_{\hat t}^T\!\int\sum_{j\not=i}|w_{i,x}v_j|\,dxdt
+2\int_{\hat t}^T\!\int |w_i\psi_i|\,dxdt\cr
&=\O(1)\cdot \delta_0^2\,.\cr}\eqno(9.22)$$
Using (9.18) and (9.22) in (9.12), we obtain the desired estimate
(9.1).
\vfill\eject

\n{\medbf 10 - Proof of the BV estimates}
\v
In this section we conclude the proof of the uniform BV bounds. 
Consider any initial data $\bar u:\R\mapsto\R^n$, 
with 
$$\tv\{\bar u\}\leq {\delta_0\over 8\sqrt n \,\kappa}\,,\qquad
\qquad\lim_{x\to
-\infty}\bar u(x)=u^*\in K\,.\eqno(10.1)$$ 
We recall that
$\kappa$ is the constant defined at (2.5), 
related to the Green kernel $G^*$
of the linearized equation (2.4).
This constant actually depends on the matrix $A(u^*)$,
but it is clear that it
remains uniformly bounded when $u^*$ varies in a compact set $K\subset\R^n$.

An application of Corollary 2.4 yields the existence of 
the solution to the Cauchy problem (1.10), (1.2) on an initial interval
$[0,\hat t]$, satisfying the bound
$$\big\|u_x(\hat t)\big\|_{\L^1}\leq {\delta_0
\over 4\sqrt n}\,.\eqno(10.2)$$
This solution can be prolonged in time as long as its total variation
remains small. 
Define the time
$$T\doteq\sup\left\{ \tau~;~~~\sum_i\int_{\hat t}^\tau\!\!\int \big|\phi_i(t,x)
\big|+\big|\psi_i(t,x)\big|\,dxdt\leq{\delta_0\over 2}
\right\}\,.\eqno(10.3)$$
If $T<\infty$, a contradiction is obtained as follows.
By (5.21) and (10.2), for all $t\in [\hat t,\,T]$ one has
$$\eqalign{\big\|u_x(t)\big\|_{\L^1}&\leq\sum_i\big\|v_i(t)\big\|_{\L^1}
\cr
&\leq \sum_i\left(\big\|v_i(\hat t)\big\|_{\L^1}+
\int_{\hat t}^T\!\int\big|\phi_i(t,x)\big|\,dxdt\right)\cr
&\leq 2\sqrt n\, \big\| u_x(\hat t)\big\|_{\L^1}
+{\delta_0\over 2}~\leq~\delta_0\,.\cr}\eqno(10.4)$$
Using Lemma 6.1 and the bounds
(7.3), (8.1), (8.2) and (9.1) we now obtain
$$\sum_i\int_{\hat t}^T\!\int
\big|\phi_i(t,x)\big|+\big|\psi_i(t,x)\big|\,dxdt
=\O(1)\cdot \delta_0^2~<~{\delta_0\over 2}\,,\eqno(10.5)$$
provided that $\delta_0$ was chosen suitably small.
Therefore $T$ cannot be a supremum.
This contradiction with (10.3) shows that the total variation remains
$<\delta_0$ for all $t\in[\hat t,\,\infty[\,$.  In particular, 
the solution $u$ is globally defined.
\v
\n{\bf Remark 10.1.} The estimates (8.1) and (9.1) were obtained under the
assumption (7.2) on the source terms.  A posteriori, 
by (10.5) the integral of the source terms is 
quadratic w.r.t.~$\delta_0$. Using (10.5) instead of (7.2) in 
the inequalities (8.17) and (9.18), (9.22), we now see that
the quantities in (8.1) and (9.1) are both $=\O(1)\cdot \delta_0^3$.
Recalling that $\delta_0$ is the order of magnitude of the total variation,
we see here another analogy with the  
the purely hyperbolic case [G].
Namely, the total amount of interactions between waves of different
families is of quadratic order w.r.t.~the total variation,
while the interaction between waves of the same
family is cubic.
\v
\n{\bf Remark 10.2.}
Within the previous proof, we constructed wave speeds
$\sigma_i\doteq\lambda_i^*-\theta(w_i/v_i)$
for which the following holds.
Decomposing the gradients $u_x, u_t$ according to
$$\left\{\eqalign{ u_x&=\sum_i v_i\,\tilde r_i(u,v_i,\sigma_i),\cr
u_t&=\sum_i (w_i-\lambda_i^*v_i)\,\tilde r_i(u,v_i,\sigma_i),\cr}\right.
\eqno(10.6)$$
the components $v_i, w_i$ then satisfy
$$\left\{\eqalign{v_{i,t}+(\tilde\lambda_iv_i)_x-v_{i,xx}&=\phi_i(t,x)\,,\cr
w_{i,t}+(\tilde\lambda_iw_i)_x-w_{i,xx}&=\psi_i(t,x)\,,\cr}\right.\eqno(10.7)$$
where all source terms $\phi_i$,$\psi_i$ are integrable:
$$\int_0^\infty\int\big|\phi_i(t,x)\big|\,dxdt<\delta_0\,,
\qquad\qquad
\int_0^\infty\int\big|\psi_i(t,x)\big|\,dxdt<\delta_0\,.\eqno(10.8)$$
In general, the speeds $\sigma_i$ defined at (5.7) are not
even continuous, as functions of $t,x$.  However,
by a suitable mollification we can find slightly different
speed functions $\sigma_i(t,x)$ which are smooth and such that
the corresponding decomposition (10.6) is achieved in terms
of (smooth) functions $v_i,w_i$ satisfying a system of the form
(10.7), with source terms again bounded as in (10.8).
\v
We conclude this section by studying the continuous dependence
w.r.t.~time of the solution $t\mapsto u(t,\cdot)$.
By (10.4) we have 
$$\tv\big\{u(t)\big\}=\big\|u_x(t)\big\|_{\L^1}\leq\delta_0\qquad
\forall t>0.\eqno(10.9)$$
By the estimate (2.8) in Proposition 2.1, the second derivative
satisfies
$$\big\|u_{xx}(t)\big\|_{\L^1}\leq \cases{
2\kappa\delta_0/\sqrt t\qquad
&if\quad $t<\hat t$,\cr
2\kappa\delta_0/\sqrt{\hat t}\qquad
&if\quad $t\geq\hat t$.\cr}\eqno(10.10)$$
Therefore, from (1.10) it easily follows
$$\big\|u_t(t)\big\|_{\L^1}\leq L'\Big(1+{1\over 2\sqrt t}\Big),$$
for some constant $L'$.
For any $t>s\geq 0$ we now have
$$\eqalign{\big\|u(t)-u(s)\big\|_{\L^1}&
\leq\int_s^t\big\|u_t(\tau)\big\|_{\L^1}\,d\tau\cr
&\leq L'\Big(|t-s|+\big|\sqrt t - \sqrt s\,\big|\Big).\cr}\eqno(10.11)$$
\v
\n {\bf Remark 10.3.} A more careful analysis shows that in (10.11) one can 
actually take 
$L'=\O(1)\cdot \tv\{\bar u\}$. However, this sharper estimate
will not be needed in the sequel.
\vfill\eject
\n{\medbf 11 - Stability estimates}
\v
Let $u=u(t,x)$ be any solution of (3.1) with small total variation.
The evolution of a first order perturbation $z=z(t,x)$ is then
governed by the linear equation
$$z_t+\big(A(u)z\big)_x-z_{xx}=
\big(u_x\bullet A(u)\big) z-\big(z\bullet A(u)\big) u_x\,.\eqno(11.1)$$
As usual, by ``$\bullet$'' we denote a directional derivative.
The primary goal of our
analysis is to establish the bound
$$\big\|z(t,\cdot)\big\|_{\L^1}\leq L\,
\big\|z(0,\cdot)\big\|_{\L^1}\qquad\forall t\geq 0\,,\eqno(11.2)$$
for some constant $L$.
By a standard homotopy argument [B1], [BiB1], this implies the uniform 
stability of solutions, w.r.t.~the $\L^1$ distance.  Indeed, consider two
initial data $\bar u, \bar v$ with suitably small
total variation.  We can assume that
$u^*\doteq \bar u(-\infty)=\bar v(-\infty)$, otherwise 
$\|\bar u-\bar v\|_{\L^1}=\infty$ and there is nothing to prove.
We construct the smooth path
$$\theta\mapsto \bar u^\theta\doteq \theta \bar u +(1-\theta) \bar v,
\qquad\qquad \theta\in [0,1].$$
Calling $t\mapsto u^\theta(t,\cdot)$ the
solution of (3.1) with initial data $\bar u^\theta$, 
for every $t\geq 0$ we have
$$\eqalign{\big\|u(t)-v(t)\big\|_{\L^1}&\leq
\int_0^1 \left\| du^\theta(t)\over d\theta\right\|_{\L^1}\,d\theta\cr
&\leq L\cdot
\int_0^1 \left\| du^\theta(0)\over d\theta\right\|_{\L^1}\,d\theta\cr
&=L\cdot \|\bar u-\bar v\|_{\L^1}\,.\cr}\eqno(11.3)$$
Indeed, the tangent vector
$$z^\theta(t,x)\doteq {du^\theta\over d\theta}(t,x)$$
is a solution of the linearized Cauchy problem
$$z^\theta_t+\big[DA(u^\theta)\cdot z^\theta\big]u^\theta_x+
A(u^\theta)z^\theta_x=z^\theta_{xx},$$
$$z^\theta(0,x)=\bar z^\theta(x) = \bar u(x)-\bar v(x),
$$
hence it satisfies (11.2) for every $\theta$.
The bound (11.3) 
provides the Lipschitz continuous dependence of solutions of (3.1)
w.r.t.~the initial data, with a Lipschitz constant independent of time.
In particular, it shows that all solutions with small total
variation are uniformly stable.
\v
\n{\bf Remark 11.1.}  In the hyperbolic case, a priori estimates
on first order tangent vectors for solutions with shocks
were first derived in [B2].  
However, even with the aid of these estimates,
controlling the $\L^1$ distance between any two solutions
remains a difficult task.
Indeed, a straightforward use of the homotopy argument
fails, due to lack of regularity.
These difficulties were eventually overcome in [BC] and [BCP],
at the price of heavy technicalities.
On the other hand, in the present case with viscosity,
all solutions are smooth and the homotopy argument goes through
without any effort.
\v
Throughout the following, we consider a reference solution 
$u=u(t,x)$ of (3.1) with small total variation. 
According to Remark 10.2, we can assume that there exist
smooth functions $v_i,w_i\sigma_i$ for which the decomposition (10.6)
holds, together with (10.7) and (10.8).

The techniques that we shall use to prove (11.2) are
similar to those used to control the total variation.
By (2.19) we already know that the desired estimate holds
on the initial time interval $[0,\hat t]$.
To obtain a uniform estimate valid for all $t>0$, we decompose
the vector $z$ along a basis of unit vectors and derive an evolution
equation for these scalar components. 
At first sight, it looks promising to write
$$z=\sum_i z_i \tilde r_i(u,v_i,\sigma_i),$$
where $\tr_1,\ldots,\tr_n$ are the same vectors
used in the decomposition of $u_x$ at (5.6).
Unfortunately, this choice would lead to non-integrable source terms.
Instead, we shall use a different basis of unit vectors
$\hr_1,\ldots,\hr_n$, depending not only
on the reference solution $u$ but also
on the perturbation $z$.

Toward this decomposition, we 
introduce the variable
$$
\io \doteq z_x - A(u) z\,,
$$
related to the flux of $z$.
By (11.1), this quantity evolves according to the equation
$$\eqalign{
\io_t + \bigl( A(u) \io \bigr)_x - \io_{xx} &= \Big[ \big(
u_x \bullet A(u)
\big) z - \big( z \bullet A(u) \big) u_x \Big]_x 
- A(u) \Big[ \big( u_x \bullet A(u)
\big) z - \big( z \bullet A(u) \big) u_x \Big]
\cr &\qquad
+ \big( u_x \bullet A(u) \big) \io - \big( u_t \bullet A(u)
\big) z \,.\cr}$$
We now decompose $z,\io$ according to
$$\left\{\eqalign{
z &= \sum_i \h_i \tilde r_i\big( u, v_i, \lambda_i^*-\theta(\iot_i/\h_i)
\big)\,,\cr
\io &= \sum_i (\iot_i- \lambda_i^* h_i) \tr_i
\big( u, v_i, \lambda_i^*-\theta(\iot_i/\h_i)\big)\,,\cr}\right.
\eqno(11.4)$$
where $\theta$ is the cutoff function introduced at (5.5).
In the following we shall write 
$$\hr_i\doteq\tilde r_i\big(u,v_i,\,
\lambda_i^*-\theta(\iot_i/\h_i)\big),$$ 
to distinguish these
unit vectors from the vectors $\tilde
r_i\big(u,v_i,\,\lambda_i^* - \theta(w_i/v_i)\big)$
previously 
used in the decomposition (5.6) of $u_x$.
Moreover we introduce the speed
$$
\hat \lambda_i \doteq \la \hat r_i,\, A(u) \hat r_i \ra\,,\eqno(11.5)
$$
and denote by
$$
\hth_i \doteq \theta ( \iot_i/\h_i)\eqno(11.6)
$$
the correction in the speed for the perturbation.
The next result, similar to Lemma 5.2, provides
the existence and regularity of the decomposition (11.4).
\v
\n{\bf Lemma 11.2} {\it 
Let $|u-u^*|$ and $|v|$ be sufficiently small.
Then for all $z,\io\in\R^n$ the system of $2n$ equations
(11.4) has a unique solution $(\h_1,\ldots,\h_n,\iot_1,\ldots,\iot_n)$.
The map $(z,\io)\mapsto(\h,\iot)$ is Lipschitz continuous.
Moreover, it is smooth outside the $n$ manifolds 
$\Hat\N_i\doteq\{\h_i=\iot_i=0\}$.}
\v
\n{\bf Proof.}~ 
The uniqueness of the decomposition is clear.
To prove the existence, 
consider the mapping $\Hat\Lambda:\R^{2n}\mapsto\R^{2n}$ defined by
$$\Hat\Lambda(\h,\iot)\doteq\sum_{i=1}^n\Hat
\Lambda_i(\h_i,\iot_i),\eqno(11.7)$$
$$\Hat\Lambda_i(\h_i,\iot_i)\doteq \pmatrix{ \h_i\,\tilde r_i
\big(u,\,v_i,\,\lambda_i^*-\theta(\iot_i/\h_i)\big)\cr 
(\iot_i-\lambda_i^*\h_i)\,\tilde r_i
\big(u,\,v_i,\,\lambda_i^*-\theta(\iot_i/\h_i)\big)\cr}.\eqno(11.8)
$$
Computing the Jacobian matrix of partial derivatives we find
$${\partial\Hat\Lambda_i\over\partial(\h_i,\iot_i)}=
\pmatrix{\hr_i+(\iot_i/\h_i)\hth'_i\hr_{i,\sigma} & 
-\hth'_i\hr_{i,\sigma}\cr
-\lambda_i^*\hr_i-\lambda_i^*(\iot_i/\h_i)\hth'_i
\hr_{i,\sigma}+(\iot_i/\h_i)^2\hth'_i\hr_{i,\sigma}
&\hr_i+\lambda_i^*\hth'_i\hr_{i,\sigma}
-(\iot_i/\h_i)\hth'_i\hr_{i,\sigma}\cr}\,.\eqno(11.9)$$
By (4.24), $\hr_{i,\sigma}=\O(1)\cdot v_i$.  
Hence, for $v_i$ small enough,
the differential $D\Hat\Lambda$ is invertible.
By the implicit function theorem, $\Hat\Lambda$ is a one-to-one
map whose range covers a whole neighborhood of the origin.
Observing that $\Hat\Lambda$ is positively homogeneous of degree 1,
we conclude that the decomposition is well defined and
Lipschitz continuous on the whole space $\R^{2n}$.
Outside the manifolds $\Hat\N_i$, $i=1,\ldots,n$,
the smoothness of the decomposition is clear.
\endproof
\v
Writing the identity $\io=z_x-A(u)z$ in terms of the decomposition
(11.4) we obtain
$$\eqalign{\sum_i (\iot_i-\lambda_i^*\h_i)\hr_i&=\sum_i h_{i,x}\hr_i
-\sum_i A(u)\h_i\hr_i+\sum_{ij}h_i\hr_{i,u}\,v_j\tr_j+\sum_i
h_i\hr_{i,v}v_{i,x}-\sum_i h_i\hr_{i,\sigma}\hth_{i,x}\,.\cr}\eqno(11.10)$$
Taking the inner product with $\hr_i$ and observing that $\hr_i$ 
is a unit vector and hence is perpendicular to its derivatives, we obtain
$$\iot_i=h_{i,x}+(\hla_i-\lambda_i^*)h_i+\Hat\Theta$$
with
$$\Hat\Theta=\sum_{j\not= i}
\la \hr_i\,,\,A(u)\hr_j\ra h_j+\sum_{j\not= i}\sum_k
\la \hr_i\,,~\hr_{j,u}\tr_k\ra h_j v_k+
\sum_{j\not= i}\la \hr_i\,,\,\hr_{j,v}\ra h_jv_{j,x} 
+\sum_{j\not= i}\la \hr_i\,,\,\hr_{j,\sigma}\ra h_j\hth_{j,x}\,.
\eqno(11.11)$$
Hence, by (5.20) and (4.24),
$$\iot_i=h_{i,x}+(\hla_i-\lambda_i^*)h_i
+\O(1)\cdot \delta_0\sum_{j\not= i}\big(|h_j|+|v_j|\big)\,.
\eqno(11.12)$$
A straightforward consequence of (11.12) is the following analogue of
Lemma 9.1.
\v
\n{\bf Corollary 11.3.} {\it 
If $|\iot_i/\h_i|\geq 3\delta_1/5$, then
$$|\iot_i|\leq 2|\h_{i,x}|+\O(1)\cdot\delta_0\sum_{j\not= i} 
\big(|v_j|+|\h_j|\big),\qquad\qquad
|\h_i|\leq {5\over 2\delta_1} |\h_{i,x}|+\O(1)\cdot
\delta_0\sum_{j\not= i}\big(|v_j|+|\h_j|\big)\,.
\eqno(11.13)$$
On the other hand, if  $|\iot_i/\h_i|\leq 4\delta_1/5$,
then}
$$|\h_{i,x}|\leq \delta_1 |\h_i|+\O(1)\cdot
\delta_0\sum_{j\not= i}\big(|v_j|+|\h_j|\big)\,.\eqno(11.14)$$
\v
Our eventual goal is to show that
the components $\h_i,\iot_i$ satisfy a system of 
evolution equations of the form
$$\left\{\eqalign{\h_{i,t}+(\tla_i\h_i)_x-\h_{i,xx}&=\hat\phi_i\,,\cr
\iot_{i,t}+(\tla_i\iot_i)_x-\iot_{i,xx}&=\hat\psi_i\,,\cr}
\right.\eqno(11.15)$$
where the source terms on the right hand sides are integrable
on $[\hat t,\,\infty[\,\times\R$.
Before embarking in calculations, 
we must first dispose of a technical difficulty due to the lack of regularity
of the equations (11.4).

Since our equations (3.1) and (11.1) are uniformly parabolic,
it is clear that for $t>0$ all solutions 
are smooth.   Moreover, by Remark 10.2, we can slightly
modify the speeds $\sigma_i$ occurring in the decomposition
of $u_x$, so that (10.6)--(10.8)
hold and the corresponding functions $v_i$ are now smooth. 
On the other hand, the map $\Hat\Lambda$ in (11.7) is only
Lipschitz continuous, hence the same is true in general
for the functions $\h_i=\h_i(t,x)$ and $\iot_i=\iot_i(t,x)$.
Indeed, at points where $\h_i=\iot_i=0$ for some index $i$,
the derivatives $\h_{i,x}$ or $\iot_{i,x}$ may well be discontinuous.
In this case, the equations (11.15) would make no sense.
To avoid this unpleasant situation, we observe that 
each manifold $\Hat\N_i$ has codimension 2.
Given the smooth functions $z,\io$ and $\epsilon>0$, by an arbitrarily
small perturbation we can 
construct new functions $\h^\sharp,\io^\sharp$ satisfying
$$\|z^\sharp-z\|_{\C^2}+\|\io^\sharp-\io\|_{\C^2}<\epsilon$$
and such that the corresponding decomposition (11.4)
is $\C^\infty$ outside a countable set of isolated points 
$(t_m,x_m)_{m\geq 1}$.  A further implementation of this technique
yields
\v
\n{\bf Lemma 11.4.}  {\it
Let $z$, $\io$ be solutions of (3.1) and (11.1) respectively. 
Then for 
any $\epsilon>0$
there exists smooth functions $z^\sharp$, $\io^\sharp$ such that the
corresponding coefficients in the decomposition (11.4) are
smooth except at countably many isolated points $(t_m,x_m)$, $m\geq 1$.
Moreover, these perturbed functions solve the system of equations
$$\eqalign{
z_t^\sharp + \big( A(u) z^\sharp \big)_x - z_{xx}^\sharp 
&= \big(u_x \bullet A(u)\big) z^\sharp - \big(z^\sharp \bullet
A(u) \big)u_x + e_1(t,x)\,,\cr
\io^\sharp_t + \bigl( A(u) \io^\sharp \bigr)_x - \io^\sharp_{xx} 
&= \Big[ \big(
u_x \bullet A(u)
\big) z^\sharp - \big( z^\sharp \bullet A(u) \big) u_x \Big]_x 
- A(u) \Big[ \big( u_x \bullet A(u)
\big) z^\sharp - \big( z^\sharp \bullet A(u) \big) u_x \Big]
\cr &\qquad
+ \big( u_x \bullet A(u) \big) \io^\sharp - \big( u_t \bullet A(u)
\big) z^\sharp+e_2(t,x) \,,\cr}$$
for some perturbations $e_1,e_2$ 
such that}
$$\int_{\hat t}^\infty\!\!\int \big|e_1(t,x)\big|+
\big|e_2(t,x)\big|\,dxdt <\epsilon\,.$$
\v
Thanks to this lemma, we can study the time evolution of the
components $\h_i,\iot_i$ 
by means of a second order parabolic system,
at the price of an arbitrarily small perturbation on the right hand side.
In the remainder of the paper, for simplicity we derive all the 
estimates in the case $e_1=e_2=0$.  
The general case easily follows by an approximation argument.
\v
In Section 6 we showed that the source terms in the equations (6.1)
could be reduced to four basic types.  
The following result is an analogue of Lemma 6.1,
providing an estimate for the source terms in the equations (11.15).
The proof, involving lengthy calculations,
will be given in Appendix~B.
\v
\n{\bf Lemma 11.5.} {\it The source terms in the equations (11.15)
satisfy the estimates}
$$\eqalign{
\hat \phi_i&(t,x),~\hat \psi_i(t,x)=\O(1)\cdot
\sum_j \Big( | \h_{j,x}| + | \h_j v_j| + | \iot_j v_j | + |
\iot_{j,x} | \Big) \,| w_j -\theta_j v_j | 
\cr
&+ \O(1)\cdot\sum_j \Big(| v_j \h_{j,x} - \h_j v_{j,x} | + 
| v_{j,x} \iot_j - \iot_{j,x} v_j | + | \h_j
w_{j,x} - w_j h_{j,x} | 
+  | \iot_j w_{j,x} - \iot_{j,x} w_j|\Big) \cr
&
+ \O(1)\cdot\sum_j \big(|v_j|+|\h_j|\big)
\left| \h_j \left(
{\iot_j\over\h_j} \right)_x^2 \right|\cdot\chi_{\strut 
\big\{|\iot_j/\h_j|<3\delta_1\big\}}
\cr
&
+ \O(1)\cdot\sum_{j \not= k}\Big( | \h_j v_k| + 
| \h_{j,x} v_k | +| \h_j v_{k,x} |
+ | \h_j w_k |
+ | \iot_j v_k | +  | \iot_{x,j} v_k | +|
\iot_j v_{k,x} | + | \h_j \h_k |+|\h_j\iot_k|\Big)
\cr
&+ \O(1)\cdot\sum_j \Big(| \h_j
\phi_j| + | \h_j \psi_j|+ | \iot_j \phi_j|
+ | \iot_j
\psi_j|\Big) . \cr}\eqno(11.16)$$
\v
The key step in establishing the bound (11.2) is to prove
\v
\n{\bf Lemma 11.6.} {\it 
Consider a solution $z$ of (11.1),
satisfying
$$\big\|z(t)\big\|_{\L^1}\leq\delta_0\qquad\quad\hbox{for all}~
t\in [0,T],\eqno(11.17)$$
and assume that the source terms in (11.4) satisfy
$$\int_{\hat t}^T\!\int\big|\hat\phi_i(t,x)\big|
+\big|\hat\psi_i(t,x)\big|\,dxdt\leq \delta_0\qquad\quad\qquad i=1,\ldots,n.
\eqno(11.18)$$
Then for each $i=1,\ldots,n$ one has the estimates}
$$\int_{\hat t}^T\!\int\big|\hat\phi_i(t,x)\big|\,dxdt=\O(1)\cdot\delta_0^2
\,,
\qquad\qquad
\int_{\hat t}^T\!\int\big|\hat\psi_i(t,x)\big|\,dxdt=\O(1)\cdot\delta_0^2\,.
\eqno(11.19)$$
\v
Assuming the validity of this lemma, we can easily recover the estimate
(11.2).   Indeed, since the equations (11.1)
are linear, it suffices to prove the estimate 
in the case where
$$\big\|z(0)\big\|_{\L^1}= {\delta_0\over 8\sqrt n\,\kappa}\,.
\eqno(11.20)$$
We recall that $\kappa$ is the constant 
defined at (3.5).
By Corollary 2.4,
on the initial interval $[0,\hat t]$ we have
$$
\big\|z(t)\big\|_{\L^1}\leq 
2\kappa \big\|z(0)\big\|_{\L^1}={\delta_0\over 4\sqrt n}\qquad\qquad 
t\in [0,\,\hat t]\,.
\eqno(11.21)$$
Define the time
$$T\doteq\sup\left\{ \tau~;~~~\sum_i\int_{\hat t}^\tau\!\!\int 
\big|\hat\phi_i(t,x)
\big|+\big|\hat\psi_i(t,x)\big|\,dxdt\leq{\delta_0\over 2}
\right\}\,.\eqno(11.22)$$
If $T<\infty$, a contradiction is obtained as follows.
First, we observe that the inequalities in (5.17) remain valid for
the decomposition of $z$, namely
$$|z|\leq\sum_i |\h_i|\leq
2\sqrt n\,|z|\,.\eqno(11.23)$$
For every $\tau\in [\hat t,\,T]$, by
(11.22) and (11.23) one has
$$\eqalign{\big\|z(\tau)\big\|_{\L^1}
&\leq\sum_i\big\|\h_i(\tau)\big\|_{\L^1}
\cr
&\leq \sum_i\left(\big\|\h_i(\hat t)\big\|_{\L^1}+
\int_{\hat t}^\tau\int\big|\hat \phi_i(t,x)\big|\,dxdt\right)\cr
&\leq 2\sqrt n\,\big\| z(\hat t)\big\|_{\L^1}
+{\delta_0\over 2}~\leq~\delta_0\,.\cr}\eqno(11.24)$$
We can thus use Lemma 11.5 and conclude
$$\sum_i\int_{\hat t}^T\!\int\big|\hat\phi_i(t,x)\big|
+\big|\hat\psi_i(t,x)\big|\,dxdt=\O(1)\cdot\delta_0^2~<~{\delta_0\over 2}
\,,\eqno(11.25)$$
provided that $\delta_0$ was chosen suitably small.
Therefore $T$ cannot be a supremum.
This contradiction shows that the bound (11.2)
holds for all $t\geq 0$ and $z\in\L^1$, with
$L\doteq 8\kappa\sqrt n$. 
The remainder of this section is aimed at establishing the estimates
(11.19).
\v
\n{\bf Proof of Lemma 11.6.}
By Corollary 2.2, for $t\in[\hat t,T]$, 
as long as $\big\|z(t)\big\|_{\L^1}\leq \delta_0$
we also have the bounds
$$\big\|z_x(t)\big\|_{\L^1}=\O(1)\cdot \delta_0^2\,,\qquad
\big\|z_{xx}(t)\big\|_{\L^1}=\O(1)\cdot \delta_0^3\,,
\qquad\big\|z_{xx}(t)\big\|_{\L^\infty}=\O(1)\cdot \delta_0^4\,.
$$
By Lemma 11.1, the map $(z,\io)\mapsto (\h,\iot)$ is uniformly
Lipschitz continuous. From the previous bounds, for every 
$t\in[\hat t,T]$ and all $j=1,\ldots,n$
it thus follows
$$\eqalignno{
\big\|\h_{j,x}(t)\big\|_{\L^1}\,,~~\big\|\iot_{j,x}(t)\big\|_{\L^1}\,,
~~\big\|\h_j(t)\big\|_{\L^\infty}\,,~~
\big\|\iot_j(t)\big\|_{\L^\infty}
~&=~\O(1)\cdot\delta_0^2\,,&(11.26)\cr&\cr
\big\|\h_{j,x}(t)\big\|_{\L^\infty}\,,~~
\big\|\iot_{j,x}(t)\big\|_{\L^\infty}~&=~\O(1)\cdot\delta_0^3\,.&(11.27)
\cr}$$
Recalling that $v_i,w_i,\h_i,\iot_i$
satisfy the systems of equations (10.7) and (11.15)
with source terms bounded by (10.8) and (11.18), 
we  now provide an estimate on the integrals of
all terms on the right hand side of (11.16).
\v
The same techniques used in Section 7 yield an estimate
on all transversal terms, with $j\not= k$:
$$\eqalign{\int_{\hat t}^T\!\!\int 
\Big( &| \h_j v_k| + 
| \h_{j,x} v_k | +| \h_j v_{k,x} |
+ | \h_j w_k |
+ | \iot_j v_k | +  | \iot_{j,x} v_k | +|
\iot_j v_{k,x} | + | \h_j \h_k |+|\h_j\iot_k|\Big)\,dxdt
\cr
&=\O(1)\cdot\delta_0^2\,.\cr}
\eqno(11.28)$$
\v
{}From (10.8) and (11.26) one easily obtains
$$\eqalign{\int_{\hat t}^T\!\!\int 
 \Big(| \h_j
\phi_j| + | \h_j \psi_j|+ | \iot_j \phi_j|
+ | \iot_j
\psi_j|\Big)\,dxdt&~\leq~\int_{\hat t}^T\!\!\int 
 \Big(\| \h_j\|_{\L^\infty}+\|\iot_j\|_{\L^\infty}\Big)\cdot\big(
|\phi_j| + |\psi_j|\big)\,dxdt\cr
&~=~\O(1)\cdot\delta_0^3\,.\cr}\eqno(11.29)$$
\v
A further set of terms will now be bounded using functionals
related to shortening curves, as in Section 8.
At each fixed time $t\in [\hat t,\,T]$, for $i=1,\ldots,n$
consider the curves
$$\gamma^{(v,\h)}_i(x)\doteq
\left( \int_{-\infty}^x v_i(t,y) dy, ~~\int_{-\infty}^x
h_i(t,y) dy \right).$$
By obvious meaning of notations, we also consider the curves 
$\gamma^{(v,\iot)}_i$, $\gamma_i^{(w,\h)}$, 
$\gamma_i^{(w,\iot)}$,
$\gamma_i^{(\h,\iot)}$.
By (6.1) and (11.15), the evolution of these curves is 
governed by vector equations similar to (8.8).
For example,
$$\gamma^{(v,\h)}_{i,t} + \tla \gamma^{(v,\h)}_{i,x} 
- \gamma^{(v,\h)}_{i,xx} = \left(
\int_{-\infty}^x \phi_i(t,y) dy, ~~\int_{-\infty}^x \hat\phi_i(t,y) dy
\right).$$
As in (8.9)-(8.10), we introduce the corresponding
{\it Length} and {\it Area Functionals}, by setting
$$\Le^{(v,\h)}_i(t)=\Le\big(\gamma_i^{(v,\h)}(t)\big)=
\int\sqrt{v_i^2(t,x)+\h_i^2(t,x)
}~dx\,,$$
$$\A_i^{(v,\h)}(t)=\A\big(\gamma_i^{(v,\h)}(t)\big)= {1\over 2}
\dint_{x<y}  \big| v_i(t,x) \h_i(t,y) -
v_i(t,y)\h_i(t,x) \big| \,dxdy\,.$$
Similarly we define $\Le_i^{(v,\iot)}(t)$,
$\A^{(v,\iot)}_i(t)$, etc$\ldots~~$
A computation entirely analogous to (8.17) now yields the
bounds
$$\eqalign{\int_{\hat t}^T\!\int
 \Big(| v_i \h_{i,x} &- \h_i v_{i,x} | \,+\, 
| v_{i,x} \iot_i - \iot_{i,x} v_i \big| \,+\, | \h_i
w_{i,x} - w_i h_{i,x} | \cr
&+\,  | \iot_i w_{i,x} - \iot_{i,x} w_i|
\,+\,| \iot_i \h_{i,x} - \h_i \iot_{i,x} |\Big)dxdt =\O(1)\cdot\delta_0^2\,.
\cr}\eqno(11.30)$$
Moreover, repeating the argument in (8.18)
we obtain
$$\int_{\hat t}^T\!\int_{|w_i/v_i|<3\delta_1}
\left|\h_i \left(\iot_i\over \h_i
\right)^2_x\right|\,dxdt=\O(1)\cdot\delta_0\,.\eqno(11.31)$$
Using the bounds (5.19) on $\|v_i\|_{\L^\infty}$
and (11.26) on $\|\h_i\|_{\L^\infty}$, from (11.31)
we deduce
$$\int_{\hat t}^T\!\int_{|w_i/v_i|<3\delta_1}
\big(|v_i|+|h_i|\big)\left| \h_i\left(\iot_i\over \h_i
\right)_x^2\right|\,dxdt=\O(1)\cdot\delta_0^3\,.\eqno(11.32)$$
\v
The integrals of the
remaining terms in (11.16) will be bounded by means of energy estimates.
For convenience, we write $\hvt_i\doteq\vth(\iot_i/h_i)$,
where $\vth$ is the cutoff function introduced at (9.2).
In Appendix~C we will prove the estimates
$$
\int_0^T\!\!\int \h_{i,x}^2 \,\hvt_i dx =
\O(1)\cdot \delta_0^2\,,\eqno(11.33)$$
$$
\int_0^T\!\!\int \iot^2_{i,x} \,\hvt_i dx =
\O(1)\cdot \delta_0^2\,.
\eqno(11.34)$$
Using (11.33)-(11.34) we now bound the terms containing 
the ``wrong speed'' $|w_i -
\theta_i v_i|$.  All these terms can be $\not= 0$ only when
$|w_i/v_i|>\delta_1$. Hence by 
(6.18) we can write
$$\eqalign{\big( | \h_i v_i| + | \iot_i v_i |\big) \,| w_i -\theta_i v_i |
&= \O(1)\cdot\big( |h_i|+|\iot_i|\big)\bigg(
\Big| v_{i,x} \,( w_i - \theta_i v_i) \Big| +
\sum_{j\not= i} \Big| v_j \,( w_i - \theta_i v_i ) \Big|\bigg).\cr}
$$
By (7.3), (9.1) and (11.26) it thus follows
$$\int_{\hat t}^T\!\int \big( | \h_i v_i| + | \iot_i v_i |\big) 
\,| w_i -\theta_i v_i |dxdt =\O(1)\cdot \delta_0^4\,.\eqno(11.35)$$
To estimate the remaining terms, we split the domain according to
the size of $|\iot_i/\h_i|$. 
\v
\n CASE 1: $|\iot_i/\h_i|> 4\delta_1/5$, $|w_i/v_i|>\delta_1$.  
Recalling (9.5) we then have
$$\eqalign{\big( | \h_{i,x}| + 
|\iot_{i,x} | \big) \,| w_i -\theta_i v_i |&\leq \big( | \h_{i,x}| + 
|\iot_{i,x} | \big)\, |w_i|
\cr
&=\big( | \h_{i,x}| + |\iot_{i,x} | \big) \Big( 2|v_{i,x}|
+\O(1)\cdot\sum_{j\not= i} |v_j|\Big)\cr
&\leq  \big(\h_{i,x}^2+\iot_{i,x}^2 +2v_{i,x}^2\big)
+
\O(1)\cdot\sum_{j\not= i} |\h_{i,x} v_j|+\O(1)\cdot\sum_{j\not= i} 
|\iot_{i,x} v_j|\cr
}$$
Using (11.33)-(11.34), (9.17) and (11.28), we conclude
$$\eqalign{\int_{\hat t}^T &\!\int_{|\iot_i/\h_i|>4\delta_1/5}
\big( | \h_{i,x}| + 
|\iot_{i,x} | \big) \,| w_i -\theta_i v_i|\,dxdt\cr
&\quad\leq
\int_{\hat t}^T\!\int  \big(\h_{i,x}^2 \,\hvt_i +\iot_{i,x}^2\,\hvt_i
+2v_{i,x}^2\vth_i\big)\,dxdt+
\O(1)\cdot\int_{\hat t}^T\!\int \sum_{j\not= i}\big(
|\h_{i,x} v_j|+|\iot_{i,x} v_j|\big)\,dxdt\cr
&\quad=\O(1)\cdot \delta_0^2\,.\cr}\eqno(11.36)$$
\v
\n CASE 2: $|\iot_i/\h_i|\leq 4\delta_1/5$, $|w_i/v_i|>\delta_1$.  
In this case we have
$$|\iot_i v_i|\leq {4\delta_1|\h_i|\over 5}\,{|w_i|\over\delta_1}
={4\over 5}|\h_iw_i|\,.\eqno(11.37)$$
Using (11.14) we can write
$$\eqalign{| \h_{i,x}| \,| w_i -\theta_i v_i |&\leq  | \h_{i,x}w_i|\cr
&=\left| \iot_i-(\hat\lambda_i-\lambda_i^*)\h_i+
\O(1)\cdot\delta_0\sum_{j\not= i}\big(|\h_j|+|v_j|\big)\right|\,|w_i|\cr
&\leq \delta_1|\h_iw_i|+\O(1)\cdot
\delta_0\sum_{j\not= i}\big(|\h_j|+|v_j|\big)|w_i|\,.\cr
}
\eqno(11.38)$$
By (11.28), the integral of the last terms on the right hand side of
(11.38) is $\O(1)\cdot\delta_0^3$.
Concerning the first term, using
(11.36) and then (6.16) and (11.14), 
we can write
$$\eqalign{{1\over 5}|\h_iw_i|&\leq   |\h_iw_i-\iot_iv_i|\cr
&\leq |\h_i v_{i,x}-\h_{i,x} v_i|+|h_i|\,\big|w_i-v_{i,x}\big|+
|v_i|\,\big|\iot_i-h_{i,x}\big|\cr
&=|\h_i v_{i,x}-\h_{i,x} v_i|+\O(1)\cdot \delta_0 |h_i|\Big( |w_i|+
\sum_{j\not= i}|v_j|\Big) +\O(1)\cdot\delta_0 |w_i|\Big(|h_i|+
\sum_{j\not= i}\big(|h_j|+|v_j|\big)\Big).\cr}$$
Hence, for $\delta_0$ small, one has
$$|\h_iw_i|\leq  6 |\h_iv_{i,x}-\h_{i,x}v_i|+\O(1)\cdot
\sum_{j\not= i}\big(|h_iv_j|+|w_ih_j|+|w_iv_j|\big)\,.$$
By (11.28) and (11.30) we conclude
$$\eqalign{\int_{\hat t}^T\!\!\int_{|\iot_i/\h_i|\leq 4\delta_1/5}
| \h_{i,x}|\,\big| w_i -\theta_i v_i\big|\,dxdt
&\leq\int_{\hat t}^T\!\!\int_{|\iot_i/\h_i|\leq 4\delta_1/5,~
|w_i/v_i|>\delta_1}
| \h_{i,x}w_i|\,dxdt\cr
&=\O(1)\cdot 
\delta_0^2\,.\cr}\eqno(11.39)$$
Recalling that $\delta_1\leq 1$,
the last remaining term can now be bounded as
$$\eqalign{| \iot_{i,x}|\,\big| w_i -\theta_i v_i\big|&\leq
|\iot_{i,x}w_i|\cr
&\leq |\iot_{i,x}w_i-\iot_i w_{i,x}|+|\iot_i w_{i,x}|\cr
&\leq |\iot_{i,x}w_i-\iot_i w_{i,x}|+{4\delta_1\over 5}|\h_i w_{i,x}|\cr
&\leq |\iot_{i,x}w_i-\iot_i w_{i,x}|+|\h_i w_{i,x}-\h_{i,x}w_i|
+|h_{i,x}w_i|\,.\cr}$$
By (11.30) and (11.39) we conclude
$$\int_{\hat t}^T\!\int_{|\iot_i/\h_i|\leq 4\delta_1/5}
| \iot_{i,x}|\,\big| w_i -\theta_i v_i\big|\,dxdt=\O(1)\cdot 
\delta_0^2\,.\eqno(11.40)$$
This completes the proof of Lemma 11.6.
\endproof
\vsk
\n{\medbf 12 - Propagation speed}
\v
Consider two solutions $u,v$ of the same
viscous system (1.10), whose initial data
coincide outside a bounded interval $[a,b]$.
Since the system is parabolic, at a given time $t>0$ one may
well have $u(t,x)\not= v(t,x)$ for all $x\in\R$.
Yet, we want to show that the bulk of the difference $|u-v|$
remains confined within a bounded interval $[a-\beta t,\, b+\beta t]$.
This result will be useful in the final section of the paper,
because it implies the finite propagation speed
of vanishing viscosity limits.
\v
\n{\bf Lemma 12.1}  {\it For some constants $\alpha,\beta>0$
the following holds. 
Let $u,v$ be solutions of (1.10) with
small total variation, whose initial data
satisfy 
$$u(0,x)=v(0,x)\qquad\qquad x\notin [a,b]\,.\eqno(12.1)$$
Then for all $x\in\R$, $t>0$ one has
$$\big|u(t,x)-v(t,x)\big|\,dx\leq \big\|u(0)-w(0)\big\|_{\L^\infty}\cdot
\min\Big\{ \alpha e^{\beta t-(x-b)},~\alpha e^{\beta t+(x-a)}\Big\}
\,.\eqno(12.2)$$
On the other hand, assuming that
$$u(0,x)=v(0,x)\qquad\qquad x\in [a,b]\,,\eqno(12.3)$$
one has}
$$\big|u(t,x)-v(t,x)\big|\,dx\leq \big\|u(0)-w(0)\big\|_{\L^\infty}\cdot
\Big(\alpha e^{\beta t-(x-a)}+\alpha e^{\beta t+(x-b)}\Big)
.\eqno(12.4)$$
\v
\n{\bf Proof.~~1.}  As a first step, we consider a solution $z$
of the linearized system 
$$z_t+\big[A(u)z\big]_x+\big[DA(u)\cdot z\big]u_x-\big[DA(u)\cdot u_x\big]z
=z_{xx}\,\eqno(12.5)$$
with initial data satisfying
$$\left\{\eqalign{\big|z(0,x)\big|&\leq 1\qquad\hbox{if}\quad x\leq 0\,,\cr
z(0,x)&=0\qquad\hbox{if}\quad x> 0\,.\cr}\right.$$
We will show that $z(t,x)$ becomes exponentially small on a domain
of the form $\{x>\beta t\}$.
More precisely,
let $B(t)$ be a continuous increasing function such that
$$B(t)\geq 1+
2 \|A\|_{\infty} \int_0^t\left(  {1\over \sqrt{t-s}}
+ \sqrt{\pi} \right)B(t)\,dt\,, \qquad\qquad B(0)=1.
$$
One can show that such a function exists, satisfying the additional
inequality $B(t)\leq 2e^{Ct}$, for some constant $C$ large enough and 
for all $t\geq 0$. 
We claim that
$$\big|z(t,x)\big|\leq E(t,x)\doteq 
B(t)\,\exp \left\{ 4 \|
DA \|_{L^\infty}
\int_0^t \bigl\| u_x(s) \bigr\|_{L^\infty} ds + t - x
\right\}\eqno(12.6)$$
for all $x\in\R$ and $t\geq 0$.
Indeed, any solution of (12.5) admits  
the integral representation
$$\eqalign{
z(t) &= G(t)*z(0) -
\int_0^t G_x(t-s)* \big[A(u)z\big](s) ds \cr
&\qquad+ \int_0^t G(t-s)* \Bigl[ \bigl( u_x \bullet A(u) \bigr)
z(s) - \bigl( h \bullet A(u) \bigr) u_x(s) \Bigr]ds
\,,\cr}$$
in terms of convolutions with the standard heat kernel
$G(t,x)\doteq e^{-x^2/4t}/2\sqrt{\pi t}$.
Therefore
$$\eqalign{
\big|z(t,x)\big| &\leq \int G(t,x-y) \bigl|z(0,y) \bigr| dy 
+ \| A \|_{L^\infty}
\int_0^t\!\!\int\Bigl| G_x(t-s,x-y) \Bigr| \bigl| z(s,y) \bigr|
dyds \cr
&\qquad+ 2 \| DA \|_{L^\infty} \int_0^t\!\!\int \bigl\| u_x(s)
\bigr\|_{L^\infty}
G(t-s,x-y) \bigl| z(s,y) \bigr| dyds\,.\cr}\eqno(12.7)$$
For every $t>0$ the following estimates hold (see Appendix D for details):
$$\int G(t,\,x-y)\,\big|z(0,y)\big|\,dy<
\int{e^{-(x-y)^2/4t}\over 2\sqrt{\pi t}}\,e^{-y} dy
=e^{t-x}\,,\eqno(12.8)
$$
$$
\|A\|_{L^\infty} \int_0^t\!\!\int \bigl| G_x(t-s,x-y) \bigr|
E(s,y) dyds \leq  {1\over 2} E(t,x)-  {1\over 2} e^{t-x},
\eqno(12.9)
$$
$$2 \|DA\|_{L^\infty}
\int_0^t \bigl\| u_x(s) \bigr\|_{L^\infty} \left(\int G(t-s,x-y)
E(s,y) dy\right)ds 
\leq  {1\over 2} E(t,x)-  {1\over 2} e^{t-x}.
\eqno(12.10)
$$
The bounds (12.7)--(12.10) show that, if (12.6) is satisfied
for all $t\in [0,\tau[$, then at time $t=\tau$ 
one always has a strict inequality:
$\big|z(\tau,x)\big|<E(\tau,x)$.  A simple argument now yields the validity
of (12.6) for all $t>0$ and $x\in\R$.
\v
\n{\bf 2.} Recalling (10.10) we have
$$\big\|u_x(s)\big\|_{\L^\infty}\leq \max\,\left\{
{2\kappa\delta_0\over\sqrt s}\,~{2\kappa\delta_0\over\sqrt {\hat t}}\right\}.$$
{}From the definition of $E$ at (12.6), for some constants $\alpha,\beta>0$
we now obtain
$$\big|z(t,x)\big|\leq E(t,x)\leq 2e^{Ct}\exp\Big\{ 4\|DA\|_{\L^\infty}
\cdot 2\kappa\delta_0\big(2\sqrt t +t/\sqrt{\hat t}~\big)\Big\}\,
e^{t-x}\leq\alpha e^{\beta t-x}.
\eqno(12.11)$$
\v
\n{\bf 3.} More generally, let now $z$ be a solution of (12.5) 
whose initial data satisfies
$$\left\{\eqalign{
\big|z(0,x)\big|&\leq \rho\qquad\hbox{if}\quad x\leq b\,,\cr
z(0,x)&=0\qquad\hbox{if}\quad x>b\,.\cr}\right.$$
By the linearity of the equations (11.1) and translation invariance,
a straightforward extension of the above arguments yields
$$\big|z(t,x)\big|\leq \rho\cdot \alpha e^{\beta t-(x-b)}.$$
On the other hand, if 
$$\left\{\eqalign{
\big|z(0,x)\big|&\leq \rho\qquad\hbox{if}\quad x\geq a\,,\cr
z(0,x)&=0\qquad\hbox{if}\quad x<a\,,\cr}\right.$$
then
$$\big|z(t,x)\big|\leq \rho\cdot \alpha e^{\beta t+(x-a)}.$$
\v
\n{\bf 4.} Having established the corresponding 
bounds on first order tangent vectors,
the estimates (12.2) and (12.4) can now be recovered by a simple homotopy
argument.
For each $\theta\in [0,1]$, let $u^\theta$ be the solution of (1.10)
with initial data 
$$u^\theta(0)=\theta u(0)+(1-\theta) v(0).$$
Moreover, call $z^\theta$ the solution of the linearized Cauchy
problem
$$z^\theta_t+\big[DA(u^\theta)\cdot z^\theta\big]u^\theta_x+
A(u^\theta)z^\theta_x=z^\theta_{xx},$$
$$z^\theta(0,x)= u(0,x)-v(0,x).
$$
If (12.1) holds, then
by the previous analysis all functions $z^\theta$ satisfy
the two inequalities
$$\eqalign{\big|z^\theta(t,x)\big|&\leq \big\|u(0)-v(0)\big\|_{\L^\infty}
\cdot\alpha e^{\beta t-(x-b)},
\cr
\big|z^\theta(t,x)\big|&\leq \big\|u(0)-v(0)\big\|_{\L^\infty}
\cdot\alpha e^{\beta t+(x-a)}.
\cr}
$$
Therefore
$$\eqalign{\big|u(t,x)-v(t,x)\big|&\leq
\int_0^1 \left|{du^\theta(t,x)\over d\theta}\right|\,d\theta
=\int_0^1 \big|z^\theta (t,x)\big|\,d\theta\cr
&\leq \big\|u(0)-v(0)\big\|_{\L^\infty}
\cdot\min\Big\{\alpha e^{\beta t-(x-b)},~\alpha
e^{\beta t-(a-x)}\Big\}.
\cr}$$
This proves (12.2).
On the other hand, if (12.3) holds, 
we consider a third solution 
$w$ of (1.10), with initial data
$$w(0,x)=\cases{u(0,x)\qquad &if\quad $x\leq b$,\cr
v(0,x)\qquad &if\quad $x\geq a$.\cr}
$$
For every $x\in\R$ and $t>0$, 
the previous arguments now yield
$$\eqalign{\big|u(t,x)-w(t,x)|&\leq \big\|u(0)-w(0)\big\|_{\L^\infty}\cdot
\alpha e^{\beta t-(x-a)},\cr
\big|v(t,x)-w(t,x)|&\leq \big\|w(0)-v(0)\big\|_{\L^\infty}\cdot
\alpha e^{\beta t+(x-b)}.\cr}$$
Combining these two inequalities we obtain (12.4).
\endproof
\vsk
\n{\medbf 13 - The vanishing viscosity limit}
\v
Up to now, all the analysis has been concerned with 
solutions of the parabolic system (1.10), with unit viscosity.
Our results, however, can be immediatly applied to the Cauchy problem
$$u^\ve_t+A(u^\ve)u^\ve_x=\ve\,u^\ve_{xx}\,,
\qquad\qquad u^\ve(0,x)=\bar u(x)\eqno(13.1)$$
for any $\ve>0$.  Indeed, as remarked in the Introduction,
a function $u^\ve$ is a solution of (13.1)
if and only if 
$$u^\ve(t,x)=u(t/\ve,\,x/\ve),\eqno(13.2)$$ where $u$ is
the solution of the Cauchy problem
$$u_t+A(u)u_x=u_{xx}\,,\qquad\qquad u(0,x)=
\bar u(\ve x)\,.\eqno(13.3)$$
Since the rescaling (13.2) does not change the total variation, 
from our earlier analysis we easily 
obtain the first part of Theorem 1.  Namely, for every
initial data $\bar u$ with sufficiently small total variation,
the corresponding solution $u^\ve(t)\doteq S^\ve_t\bar u$ 
is well defined for all times $t\geq 0$.  The bounds 
(1.15)--(1.17) follow from
$$\tv\big\{u^\ve(t)\big\}=\tv\big\{u(t/\ve)\big\}
\leq C\,\tv\{\bar u\},\eqno(13.4)$$
$$\big\|u^\ve(t)-v^\ve(t)\big\|_{\L^1}=
\ve \big\|u(t)-v(t)\big\|_{\L^1}\leq \ve L\,
\big\|u(0)-v(0)\big\|_{\L^1}=\ve L\,{1\over\ve}\|\bar u-\bar v\|_{\L^1}\,,
\eqno(13.5)$$
$$
\big\|u^\ve(t)-u^\ve(s)\big\|_{\L^1}\leq 
\ve \big\|u(t/\ve)-u(s/\ve)\big\|_{\L^1}\leq 
\ve L'\left(\Big|{t\over\ve}-{s\over\ve}\Big|+\left|\sqrt
{t\over\ve}-\sqrt{s\over\ve}\right|\right).\eqno(13.6)$$
Moreover, if 
$\bar u(x)=\bar v(x)$ for $x\in [a,b]$, then
(12.4) implies
$$\big|u^\ve(t,x)-v^\ve(t,x)\big|
\leq
\|\bar u-\bar v\|_{\L^\infty}\cdot
\left\{\alpha \exp\Big({\beta t-(x-a)\over\ve}\Big)+
\alpha \exp\Big({\beta t+(x-b)\over\ve}\Big)\right\}.\eqno(13.7)$$
\v
We now consider the vanishing viscosity limit.
Call $\U\subset\Ll$ the set of all functions $\bar u:\R\mapsto
\R^n$ with small total variation, satisfying (1.14).
For each $t\geq 0$ and every initial condition $\bar u\in\U$,
call $S^\ve_t\bar u\doteq u^\ve(t,\cdot)$ the corresponding
solution of (13.1).
Thanks to the uniform BV bounds (13.4), we can apply
Helly's compactness theorem and obtain a sequence $\ve_\nu\to 0$
such that $$\lim_{\nu\to\infty}
u^{\ve_\nu}(t,\cdot)= u(t,\cdot)\qquad\qquad\hbox{in}\quad\Ll\,.
\eqno(13.8)$$
holds for some BV function $u(t,\cdot)$.
By extracting further subsequences and then using a standard diagonalization 
procedure, we can assume that
the limit in (13.8) exists for all rational times $t$ and all
solutions $u^\ve$ with initial data in a countable dense set 
$\U^*\subset\U$.   
Adopting a semigroup notation, we thus define
$$S_t\bar u\doteq \lim_{m\to\infty} S^{\ve_m}_t \bar u\qquad
\qquad\hbox{in}~~\Ll\,,
\eqno(13.9)$$
for some particular subsequence $\ve_m\to 0$.
By the uniform continuity of the maps
$(t,\bar u)\mapsto u^\ve(t,\cdot)\doteq S^\ve_t\bar u\,,$
stated in (13.5)-(13.6), the set of couples $(t,\bar u)$ for which
the limit (13.9) exists must be closed in $\R_+\times\U$.
Therefore, this limit is well defined
for all $\bar u\in\U$ and $t\geq 0$.
\v
\n{\bf Remark 13.1.}  The function $u(t,\cdot) =S_t\bar u$
is here defined as a limit in $\Ll$. Since it has bounded variation,
we can remove any ambiguity concerning its pointwise values
by choosing, say,  a right continuous representative:
$$u(t,x)=\lim_{y\to x+}u(t,y).$$
With this choice, the function $u$ is certainly jointly measurable 
w.r.t. $t,x$ (see [B5], p.16).
\v
To complete the proof of Theorem 1, we need to show that
the map $S$ defined at (13.9) is a semigroup, satisfies the
continuity properties (1.18) and does not depend on the choice of
the subsequence $\{\ve_m\}$. 
These results will be achieved 
in several steps.
\v
\n{\bf 1. (Continuous dependence)} 
Let $S$ be the map defined by (13.9). Then
$$\big\|S_t\bar u-S_t\bar v\big\|_{\L^1}=\sup_{r>0}\int_{-r}^r
\Big|\big(S_t\bar u\big)(x)-\big(S_t\bar v\big)(x)\Big|\,
dx\,.$$
For every $r>0$, the convergence in $\Ll$ implies
$$\int_{-r}^r
\Big|\big(S_t\bar u\big)(x)-\big(S_t\bar v\big)(x)\Big|\,
dx=\lim_{m\to\infty}
\int_{-r}^r
\Big|\big(S^{\ve_m}_t\bar u\big)(x)-\big(S_t^{\ve_m}\bar v\big)(x)\Big|\,
dx\leq L\,\big\|\bar u-\bar v\big\|_{\L^1}\,.$$
because of (13.5). This yields the Lipschitz 
continuous dependence w.r.t.~the initial data:
$$\big\|S_t\bar u-S_t\bar v\big\|_{\L^1}
\leq L\,\big\|\bar u-\bar v\big\|_{\L^1}\,.\eqno(13.10)$$

The continuous dependence w.r.t.~time is proved in a similar way.
By (13.6), for every $r>0$ we have
$$\eqalign{\int_{-r}^r
\Big|\big(S_t\bar v\big)(x)-\big(S_s\bar v\big)(x)\Big|\,
dx&=\lim_{m\to\infty}
\int_{-r}^r
\Big|\big(S^{\ve_m}_t\bar v\big)(x)-\big(S_s^{\ve_m}\bar v\big)(x)\Big|\,
dx\cr
&=\lim_{m\to\infty} \ve_m L'  
\left(\Big|{t\over\ve_m}-{s\over\ve_m}\Big|+\left|\sqrt
{t\over\ve_m}-\sqrt{s\over\ve_m}\,\right|\,\right)\cr
&=L'|t-s|\,.\cr}$$
Hence
$$\big\|S_t\bar v-S_s\bar v\big\|_{\L^1}\leq L'|t-s|\,.
\eqno(13.11)$$

Together, (13.10) and (13.11) yield (1.18).
\v
\n{\bf 2. (Finite propagation speed)} 
Consider any interval $[a,b]$ and two initial data 
$\bar u,\bar v$, with 
$\bar u(x)=\bar v(x)$ for $x\in [a,b]$.
By (13.7), for every $t\geq 0$ and $x\in \,]a+\beta t,\, b-\beta t[\,$ 
one has
$$\eqalign{\Big|\big(S_t\bar u\big)(x)-\big(S_t\bar v\big)(x)
\Big|&\leq\limsup_{m\to\infty}
\Big|\big(S^{\ve_m}_t\bar u\big)(x)-\big(S^{\ve_m}_t\bar v\big)(x)
\Big|\cr
&\leq \lim_{m\to\infty} \|\bar u-\bar v\|_{\L^\infty}\cdot
\left\{\alpha \exp\Big({\beta t-(x-a)\over\ve_m}\Big)+
\alpha \exp\Big({\beta t+(x-b)\over\ve_m}\Big)\right\}\cr
&=0\,.\cr}\eqno(13.12)
$$
In other words, the restriction of the function $S_t\bar u\in\Ll$
to a given interval $[a',b']$ depends only on the 
values of the initial data $\bar u$ on the interval
$[a'-\beta t,\, b'+\beta t]$.
Using (13.12), we now prove a sharper version
of the continuous dependence estimate (13.10):
$$\int_a^b
\Big|\big(S_t\bar u\big)(x)-\big(S_t\bar v\big)(x)
\Big|\,dx\leq L\cdot\int_{a-\beta t}^{b+\beta t}\big|
\bar u(x)-\bar v(x)\big|\,dx\,.\eqno(13.13)$$
valid for every $\bar u,\bar v$ and $t\geq 0$.
Indeed, define the auxiliary
function 
$$\bar w(x)=\cases{\bar u(x)\qquad &if\quad $x\in [
a-\beta t,~b+\beta t]\,,$\cr
\bar v(x)\qquad &if\quad $x\notin [
a-\beta t,~b+\beta t]\,.$\cr}$$
{}Using the finite propagation speed, we now have
$$\eqalign{\int_a^b
\Big|\big(S_t\bar u\big)(x&)-\big( S_t\bar v\big)(x)
\Big|\,dx=
\int_a^b
\Big|\big(S_t\bar w\big)(x)-\big(S_t\bar v\big)(x)
\Big|\,dx\cr
&\leq L\,\|\bar w-\bar v\|_{\L^1}=
L\cdot\int_{a-\beta t}^{b+\beta t}\big|
\bar u(x)-\bar v(x)\big|\,dx\,.\cr}$$
\v
\n{\bf 3. (Semigroup property)}
We now show that the map $(t,\bar u)\mapsto S_t\bar u$ is a semigroup,
i.e.
$$S_0\bar u=\bar u\,,\qquad\qquad S_sS_t\bar u=S_{s+t}\bar u\,.\eqno(13.14)$$
Since every $S^\ve$ is a semigroup, 
the first equality in (13.14) is a trivial consequence of the
definition (13.9).
To prove the second equality, we observe that
$$S_{s+t}\bar u=\lim_{m\to\infty}S_s^{\ve_m}S_t^{\ve_m}\bar u\,,
\qquad\qquad S_sS_t\bar u=\lim_{m\to\infty}S_s^{\ve_m}S_t\bar u\,.
\eqno(13.15)$$
We can assume $s>0$. Fix any $r>0$ and
consider the function
$$\tilde u_m(x)\doteq\cases{\big(S_t\bar u\big)(x)\qquad 
&if\quad $|x|\leq r+2\beta s\,,$\cr
\big(S_t^{\ve_m}\bar u\big)(x)\qquad &if\quad $|x|> r+2\beta s\,.$\cr
}$$
Since $S_t^{\ve_m}\bar u\to S_t\bar u$ in $\Ll$,
using (13.7) and (13.5) one obtains
$$\eqalign{\limsup_{m\to\infty}\int_{-r}^r
&\Big|\big(S_s^{\ve_m}S_t^{\ve_m}\bar u\big)(x)
-\big(S_s^{\ve_m}S_t\bar u\big)(x)\Big|\,dx\cr
&\leq \lim_{m\to\infty} 2r\cdot \sup_{|x|<r}
\Big|\big(S_s^{\ve_m}S_t^{\ve_m}\bar u\big)(x)
-\big(S_s^{\ve_m}\tilde u_m\big)(x)\Big|+
\lim_{m\to\infty}\big\|S_s^{\ve_m}\tilde u_m-
S_s^{\ve_m}S_t^{\ve_m}\bar u\big\|_{\L^1}\cr
&\leq \lim_{m\to\infty} 2r\,
\big\|S_t^{\ve_m}\bar u-\tilde u_m\big\|_{\L^\infty}\cdot
2\alpha
e^{-\beta s/\ve_m}+
\lim_{m\to\infty} L\cdot\big\|\tilde u_m-S_t^{\ve_m}\bar u\big\|_{\L^1}\cr
&=0\,.\cr}$$
By (13.15), this proves the second identity in (13.14).
\v
\n{\bf 4. (Tame Oscillation)} We now exhibit
a regularity property which is shared by all semigroup trajectories.
This property, introduced in [BG], plays a key role in the  
proof of uniqueness.  We begin by recalling the main definitions.
Given $a<b$ and $\tau\geq 0$, we denote by
$\tv\big\{u(\tau)\,;~]a,b[\,\big\}$ the total variation
of $u(\tau,\cdot)$ over the open interval $]a,b[\,$.
Moreover, consider the triangle
$$\Delta_{a,b}^\tau\doteq\big\{ (t,x)\,;~~t>\tau,~
a+\beta t<x<b-\beta t\big\}\,.$$
The oscillation of $u$ over $\Delta_{a,b}^\tau$ will be denoted by
$$\osc\big\{u\,;~\Delta_{a,b}^\tau\big\}\doteq\sup
\Big\{\big|u(t,x)-u(t',x')\big|\,;~~
(t,x),(t',x')\in
\Delta_{a,b}^\tau\Big\}.$$
We claim that each function
$u(t,x)=\big(S_t\bar u\big)(x)$ satisfies the 
{\it tame oscillation property}: there exists a constant $C'$ such that,
for every $a<b$ and $\tau\geq 0$, one has
$$\osc\big\{u\,;~\Delta_{a,b}^\tau\big\}
\leq C'\cdot\tv\big\{u(\tau)\,;~]a,b[\,\big\}\,.\eqno(13.16)$$
Indeed, let $a,b,\tau$ be given, together with an initial data
$\bar u$. By the semigroup
property, it is not restrictive to assume $\tau =0$.
Consider the auxiliary initial condition 
$$\bar v(x)\doteq\cases{ \bar u(x)\qquad &if\quad $a<x<b$,\cr
\bar u(a+)\qquad &if\quad $x\leq a$,\cr
\bar u(b-)\qquad &if\quad $x\geq b$,\cr}\eqno(13.17)$$
and call $v(t,x)\doteq \big(S_t\bar v\big)(x)$ the corresponding
trajectory.
Observe that 
$$\lim_{x\to -\infty} v(t,x)=\bar u (a+)$$
for every $t\geq 0$.
Using (1.15) and the finite propagation speed, we can thus write 
$$\eqalign{\osc\big\{&u\,;~\Delta_{a,b}^\tau\big\}=
\osc\big\{v\,;~\Delta_{a,b}^\tau\big\}
\leq 2\sup_t \Big(\tv\big\{S_t\bar v\big\}\Big)\cr
&\leq 2 C \cdot\tv\{\bar v\}
=2C \cdot\tv\big\{u(\tau)\,;~]a,b[\,\big\}\,,\cr}$$
proving (13.16) with $C'=2C$.
\v
\n{\bf 5. (Conservation equations)}
Assume that 
the system (13.1) is in conservation form, i.e.~$A(u)=Df(u)$
for some flux function $f$. In this special case, we claim that
every vanishing viscosity limit is a weak solution of the 
system of conservation laws (1.1).  Indeed, 
with the usual notations, if $\phi$ is a
$\C^2$ function with compact support contained in the half
plane $\{x\in\R,\,t>0\}$, one can repeatedly integrate 
by parts and obtain
$$\eqalign{\int\!\!\int &\big[u\,\phi_t+f(u)\phi_x\big]\,dxdt~=~
\lim_{m\to\infty} 
\int\!\!\int\big[u^{\ve_m}\,\phi_t+f(u^{\ve_m})\phi_x\big]\,dxdt
\cr
&=-\lim_{m\to\infty} 
\int\!\!\int\big[u_t^{\ve_m}\,\phi+f(u^{\ve_m})_x\phi\big]\,dxdt
~=~-\lim_{m\to\infty} \int\!\!\int
\ve_m\,u^{\ve_m}_{xx}\phi\,dxdt
\cr
&=-\lim_{m\to\infty} \int\!\!\int
\ve_m\,u^{\ve_m}\phi_{xx}\,dxdt~=~0\,.\cr}$$
An easy approximation argument shows that 
the identity (1.5) holds more generally,
assuming only $\phi\in \C^1_c$.
\v
\n{\bf 6. (Approximate jumps)} From the uniform 
bound on the total variation and
the Lipschitz continuity w.r.t.~time, it
follows that each function $u(t,x)=\big(S_t\bar u\big)(x)$
is a BV function, jointly w.r.t.~the two variables 
$t,x$.  In particular, an application of Theorem 2.6 in [B5]
yields the existence of a set of times $\N\subset\R_+$ of measure zero
such that, for every $(\tau,\xi)\in\R_+\times\R$ with $\tau\notin\N$,
the following holds.
Calling
$$u^-\doteq\lim_{x\to\xi-} u(\tau,x),\qquad\qquad
u^+\doteq\lim_{x\to\xi+} u(\tau,x),\eqno(13.18)$$
there exists a finite speed $\lambda$ such that the function
$$U(t,x)\doteq\cases{u^-\quad &if\quad $x<\lambda t$,\cr
u^+\quad &if\quad $x>\lambda t$,\cr}\eqno(13.19)$$
for every constant $\kappa>0$ satisfies
$$\lim_{r\to 0+}
{1\over r^2}\int_{-r}^r\!\int_{-\kappa r}^{\kappa r}
\big| u(\tau+t,\,\xi+x)-U(t,x)\big|\,dxdt=0,\eqno(13.20)$$
$$\lim_{r\to 0+}
{1\over r}\int_{-\kappa r}^{\kappa r}
\big| u(\tau+r,\,\xi+x)-U(r,x)\big|\,dx=0.\eqno(13.21)$$

In the case where $u^-\not= u^+$, we say that $(\tau,\xi)$ is a
point of {\it approximate jump} for the function $u$.
On the other hand, if $u^-=u^+$ (and hence $\lambda$ can 
be chosen arbitrarily), we say that $u$ is {\it approximately
continuous} at $(\tau,\xi)$.
The above result can thus be restated as follows:
with the exception of a null set $\N$ of ``interaction times'',
the solution $u$ is either approximately continuous or has an
approximate jump discontinuity at each point $(\tau,\xi)$.
\v
\n{\bf 7. (Shock conditions)} Assume again that the system is in
conservation form.  Consider a semigroup
trajectory $u(t,\cdot)=S_t\bar u$ and a point $(\tau,\xi)$
where $u$ has an approximate jump.  Since $u$ is a weak solution,
the states $u^-,u^+$ and the speed $\lambda$ 
in (13.19) must satisfy the
Rankine-Hugoniot equations
$$\lambda\,(u^+-u^-)=f(u^+)-f(u^-).\eqno(13.22)$$
For a proof, see Theorem 4.1 in [B5].

If $u$ is a limit of vanishing viscosity approximations, the same is
true of the solution $U$ in (13.19).  In particular (see [MP] or [D]), the 
{\it Liu shock conditions} must hold.
More precisely, call $s\mapsto S_i(s)$ the 
parametrized shock curve through $u^-$ and let
$\lambda_i(s)$ be the speed of the corresponding shock.
If $u^+=S_i(s)$ for some $s$, then
$$\lambda_i(s')\leq \lambda_i(s')\qquad\hbox{for all}~s'\in [0,s]\,.
\eqno(13.23)$$

Under the additional assumption that each characterictic 
field is either linearly degenerate or genuinely nonlinear,
it is well known that the Liu conditions imply
the Lax shock conditions:
$$\lambda_i(u^+)\leq \lambda\leq\lambda_i(u^-).\eqno(13.24)$$
\v
\n{\bf 8. (Uniqueness in a special case)} 
Assume that the system is in conservation form and that
each  characterictic 
field is either linearly degenerate or genuinely nonlinear.
By the previous steps, the semigroup trajectory
$u(t,\cdot)=S_t\bar u$ provides a weak solution to the Cauchy problem
(1.1)-(1.2) which satisfies the Tame Oscillation and the
Lax shock conditions.  By a well known uniqueness theorem
in [BG], [B5], such a weak solution is unique and coincides with
the limit of front tracking approximations.
In particular, it does not depend
on the choice of the subsequence
$\{\ve_m\}$ :
$$S_t\bar u=\lim_{\ve\to 0+} S_t^\ve\bar u\,,$$
i.e.~the same limit actually holds over all real values of $\ve$.
\v
The above results already yield a proof of Theorem 1 in the special case
where the system is in conservation form and
satisfies the standard assumptions (H).   To handle the general 
(non-conservative) case, we shall need to understand 
first the solution of the Riemann problem.
\vsk
\n{\medbf 14 - The non-conservative Riemann problem}
\v
Aim of this section is to characterize the vanishing viscosity
limit for solutions $u^\ve$ of (13.1), in the case of
Riemann data
$$\bar u(x)=\cases{u^-\quad &if\quad $x<0$,\cr
u^+\quad &if\quad $x>0$.\cr}\eqno(14.1)$$
More precisely, we will show that, as $\ve\to 0+$, the solutions 
$u^\ve$ converge to a self-similar limit $\omega(t,x)=\tilde\omega(x/t)$.
We first describe a method for constructing this solution 
$\omega$.
\v
As a first step, given a left state $u^-$ and $i\in\{1,\ldots,n\}$, 
we seek a
one-parameter curve of right states
$u^+=\Psi_i(s)$ such that the non-conservative Riemann problem
$$\omega_t+A(u)\omega_x=0,\qquad\qquad \omega(0,x)=
\cases{u^-\quad &if\quad $x<0$\cr
u^+\quad &if\quad $x>0$\cr}\eqno(14.2)$$
admits a vanishing viscosity solution consisting only of $i$-waves.
In the case where the system is in conservation form and 
the $i$-th field is genuinely nonlinear, it is well known [Lx] that one
should take
$$\Psi_i(s)=\cases{R_i(s)\quad &if\quad $s\geq 0$,\cr
S_i(s)\quad &if\quad $s< 0$.\cr}$$
Here $R_i$ and $S_i$ are the $i$-th rarefaction and shock curves
through $u^-$, respectively.
We now describe a method for constructing such curve $\Psi_i$
in the general case.
\v
Fix $\epsilon,s>0$.  Consider the family $\Gamma\subset\C^0\big([0,s]\,;~
\R^n\times\R\times\R\big)$
of all continuous curves
$$\tau\mapsto\gamma(\tau)=\big(u(\tau),\,v_i(\tau),\,\sigma_i(\tau)\big),
\qquad\qquad \tau\in [0,s]\,,$$
with 
$$u(0)=u^-,\qquad \big|u(\tau)-u^-\big|\leq\epsilon\,,
\qquad \big|v_i(\tau)\big|\leq
\epsilon\,,\qquad \big|\sigma_i(\tau)-\lambda_i(u^-)\big|\leq
\epsilon\,.$$
In connection with a given curve $\gamma\in\Gamma$, 
define the scalar flux function
$$f_i(\gamma,\tau)\doteq\int_0^\tau\tla_i\big(u(\vars),\,v_i(\vars),\,
\sigma_i(\vars)\big)\,d\vars\qquad\qquad \tau\in [0,\,s]\,,\eqno(14.3)$$
where $\tla_i$ is the speed in (4.21).
Moreover, consider the lower convex envelope
$$\hbox{conv}\,f_i(\gamma,\tau)\doteq\inf\Big\{
\theta f_i(\gamma,\tau')+(1-\theta)f_i(\gamma,\tau'')\,;~~~
\theta\in [0,1]\,,~~\tau',\tau''\in [0,s]\,,~~
\tau=\theta \tau'+(1-\theta)\tau''\Big\}\,.$$
We now define a continuous mapping $\T_{i,s}:\Gamma\mapsto\Gamma$ by
setting
$\T_{i,s}\gamma=\hga=(\hat u,\hat v_i,\hat \sigma_i)$, where
$$\left\{
\eqalign{\hat u(\tau)&\doteq u^-+\int_0^\tau \tilde r_i
\big(u(\vars),\,v_i(\vars),\,\sigma_i(\vars)\big)\,d\vars\,,\cr
\hat v_i(\tau)&\doteq f_i(\gamma,\tau)-\hbox{conv}\,f_i(\gamma,\tau)\,,\cr
\hat\sigma_i(\tau)&\doteq {d\over d\tau}\, 
\hbox{conv}\,f_i(\gamma,\tau)\,.\cr}
\right.\eqno(14.4)$$
We recall that $\tr_i$ are the unit vectors that define the center manifold
in (4.13).
Because of the bounds
$$\big|\hat u(\tau)-u^-\big|\leq \tau\leq s\,,$$
$$\big|\hat \sigma_i(\tau)-\lambda_i(u^-)\big|=\O(1)\cdot \sup_{\vars\in
[0,s]}
\Big|\tla_i\big(u(\vars),\,v_i(\vars),\,
\sigma_i(\vars)\big)-\lambda_i(u^-)\Big|=\O(1)\cdot s\,,$$
$$\big|\hat v_i(\tau)\big|=\O(1)\cdot s\,
\sup_{\vars\in [0,s]}
\Big|\tla_i\big(u(\vars),\,v_i(\vars),\,
\sigma_i(\vars)\big)-\lambda_i(u^-)\Big|=\O(1)\cdot s^2\,,$$
it is clear that for $0\leq s<\!<\epsilon <\!<1$ 
the transformation
$\T_{i,s}$ maps $\Gamma$ into itself.
We now show that, in this same range of parameters,
$\T_{i,s}$ is a contraction
with respect to the weighted norm
$$\big\|(u,v_i,\sigma_i)\big\|_\dagger\doteq \|u\|_{\L^\infty}
+\|v_i\|_{\L^\infty}
+\epsilon\,\|\sigma_i\|_{\L^\infty}
\,.$$
Indeed, consider two curves $\gamma,\gamma'\in\Gamma$.
For each $\tau\in [0,s]$ one has
$$\eqalign{\|\hat u-\hat u'\|_{\L^\infty}&\leq
\int_0^s \Big|\tr_i(u,v_i,\sigma_i)-
\tr_i(u',v'_i,\sigma'_i)\Big|\,d\vars\cr
&=\O(1)\cdot s\,\Big(
\big\|u-u'\big\|_{\L^\infty}+\|v_i-v'_i\|_{\L^\infty}+
\|v_i\|_{\L^\infty}\|\sigma_i
-\sigma'_i\|_{\L^\infty}\Big),\cr}$$
$$\eqalign{\|\hat v_i-\hat v_i'\|&\leq
2\big\|f_i(\gamma)-f_i(\gamma')\big\|_{\L^\infty}\cr
&\leq
2\int_0^s \Big|\tla_i(u,v_i,\sigma_i)-
\tla_i(u',v'_i,\sigma'_i)\Big|\,d\vars\cr
&=\O(1)\cdot\,s\,\Big(
\big\|u-u'\big\|_{\L^\infty}+\|v_i-v'_i\|_{\L^\infty}+
\|v_i\|_{\L^\infty}\|\sigma_i
-\sigma'_i\|_{\L^\infty}\Big)\,,\cr}$$
$$\eqalign{\|\hat \sigma_i(\tau)-\hat \sigma_i'\|&
\leq\sup_{\tau\in [0,s]}\left|{d\over d\tau}\hbox{conv}\,f_i(\gamma,\tau)
-{d\over d\tau}\hbox{conv}\,f'_i(\gamma,\tau)\right|
\leq\sup_{\tau\in [0,s]}\left|{d\over d\tau}\,f_i(\gamma,\tau)
-{d\over d\tau}\,f'_i(\gamma,\tau)\right|\cr
&
\leq
\big\|\tla_i-\tla_i'\big\|_{\L^\infty}=\O(1)\cdot\Big(
\big\|u-u'\big\|_{\L^\infty}+\|v_i-v'_i\|_{\L^\infty}+
\|v_i\|_{\L^\infty}\|\sigma_i
-\sigma'_i\|_{\L^\infty}\Big)\,.\cr}$$
For some constant $C_0$, the previous estimates imply
$$\|\hat\gamma-\hat\gamma'\|_\dagger\leq
C_0\,\epsilon\,
\|\gamma-\gamma'\|_\dagger\leq {1\over 2} \|\gamma-\gamma'\|_\dagger\,,
\eqno(14.5)$$
provided that $\epsilon$ is sufficiently small.  Therefore,
by the contraction mapping principle, the map $\T_{i,s}$
admits a unique fixed point, i.e. a continuous curve $\gamma=(u,v_i,\sigma_i)$
such that
$$\left\{
\eqalign{u(\tau)&= u^-+\int_0^\tau \tilde r_i
\big(u(\vars),\,v_i(\vars),\,\sigma_i(\vars)\big)\,d\vars\,,\cr
v_i(\tau)&= f_i(\gamma,\tau)-\hbox{conv}\,f_i(\gamma,\tau)\,,\cr
\sigma_i(\tau)&= {d\over d\tau}\, 
\hbox{conv}\,f_i(\gamma,\tau)\,.\cr}
\right.\eqno(14.6)$$
Recalling the definition (14.3), from the continuity of $u,v_i,\sigma_i$
it follows that the maps $\tau\mapsto u(\tau)$, $\tau\mapsto v_i(\tau)$
and $\tau\mapsto f_i(\gamma,\tau)$
are continuously diferentiable. 
We now show that, taking $u^+=u(s)$ corresponding to the endpoint
of this curve $\gamma$, the Riemann problem (14.2) admits a
self-similar solution containing only $i$-waves.
\v
\n{\bf Lemma 14.1} {\it In the previous setting, let
$\gamma: \tau\mapsto \big(u(\tau),\,v_i(\tau),\,\sigma_i(\tau)\big)$
be the fixed point of the transformation $\T_{i,s}$. Define the right state
$u^+\doteq u(s)$.  Then the unique vanishing
viscosity solution
of the Riemann problem (14.2) is the function}
$$\omega(t,x)\doteq\left\{\eqalign{&u^-\cr
&u(\tau)\cr &u^+\cr}\qquad
\eqalign{&\hbox{if}\quad x/t\leq \sigma_i(0)\,,\cr
&\hbox{if}\quad x/t=\sigma(\tau)\,,\cr
&\hbox{if}\quad x/t\geq\sigma(s)\,.\cr}\right.\eqno(14.7)$$
\v
\n{\bf Proof.}
With the semigroup notations introduced in Theorem 1,
we will show that, for every $t\geq 0$,
$$\lim_{\ve\to 0+}\,\big\|\omega(t)-S^\ve_t \omega(0)\big\|_{\L^1}=0\,.
\eqno(14.8)$$
The proof will be given in several steps.
\v
\n{\bf 1.} Assume that we can construct a family $v^\ve$ of solutions
to
$$v_t+A(v)v_x=\ve\,v_{xx}\,,\eqno(14.9)_\ve$$
with
$$\lim_{\ve\to 0+}\big\|v^\ve(t)-\omega(t)\big\|_{\L^1}=0\eqno(14.10)$$
for all $t\in [0,1]$. Then (14.8) follows.
Indeed, by a simple rescaling we immediately have a family
of solutions $v^\ve$ such that (14.9)$_\ve$-(4.10) hold 
on any fixed interval $[0,T]$.
For every $t\in [0,T]$, since by assumption
$v^\ve(t)=S^\ve_t v^\ve(0)$, using (1.16) 
we obtain
$$\eqalign{\lim_{\ve\to 0+}\,\big\|\omega(t)-S^\ve_t \omega(0)\big\|_{\L^1}
&\leq \lim_{\ve\to 0+}\,\big\|\omega(t)-v^\ve(t)\big\|_{\L^1}
+\lim_{\ve\to 0+}\,\big\|v^\ve(t)-S^\ve_t \omega(0)\big\|_{\L^1}
\cr
&\leq 0 +L\cdot \lim_{\ve\to 0+} 
\,\big\|v^\ve(0)-\omega(0)\big\|_{\L^1}
~=~0\,.\cr}$$
\v
\n{\bf 2.} For notational convenience,
call $\vvl$ the set of all vanishing viscosity limits, i.e.~all
functions $v:[0,1]\times \R\mapsto\R^n$ 
such that 
$$\lim_{\ve\to 0+}\big\|v^\ve(t)-v(t)\big\|_{\L^1}=0\qquad\qquad
t\in [0,1]\eqno(14.11)$$
for some family of solutions $v^\ve$ of (14.9)$_\ve$ .
By Step 1, it suffices to show that the function
$\omega$ at (14.7) lies in $\vvl$.

Let us make some preliminary considerations.
Consider a  piecewise smooth function $v=v(t,x)$ 
which provides a classical solution to the quasilinear system
$$v_t+A(v)v_x=0\qquad\qquad t\in [0,1]\,,$$
outside a finite number of straight lines, say
$x=x_j(t)$, $j=1,\ldots,N$.
Assume that there exists $\delta>0$ and 
constant states
$u_j^-,u_j^+$ such that
$$u(t,x)=\cases{u_j^-\qquad &if\quad $x_j(t)-\delta\leq x<x_j(t)\,$,\cr
u_j^+\qquad &if\quad $x_j(t)<x\leq x_j(t)+\delta\,$,\cr}$$
Moreover, assume that each pair of
states $u_j^-,u_j^+$ can be
connected by a viscous travelling wave having speed $\dot x_j$.
Finally, let $v$ be constant on each of the two regions where
$x>r_0$ or $x<-r_0$, for some $r_0$ sufficiently large.
Under all of the above hypotheses, it is then clear that $v\in\vvl$.
Indeed, a family of viscous approximations $v^\ve$
can be constructed by a simplified version of the singular perturbation
technique used in [GX].

As a second observation, notice that if we have a sequence 
of functions $v_m\in \vvl$ with 
$$\lim_{m\to\infty}\big\|v_m(t)-v(t)\big\|_{\L^1}=0\qquad\qquad
t\in [0,1]\,,$$
then also $v\in\vvl$.
\v
\n{\bf 3.}
Consider first the (generic) case where the set of points
in which $f_i$ is disjoint from its convex envelope is a finite union of
open intervals (fig.~12), say
$$\Big\{ \tau\in [0,s]\,;~~f_i(\gamma,\tau)>\hbox{conv}\,f_i(\gamma,\tau)
\Big\}=\bigcup_{j=1}^N \,]a_j,\,b_j[\,.\eqno(14.12)$$
Our strategy is to prove that $\omega\in\vvl$
first in this special case.  Later we shall deal with the general case,
by an approximation argument.

\midinsert
\vskip 10pt
\centerline{\hbox{\psfig{figure=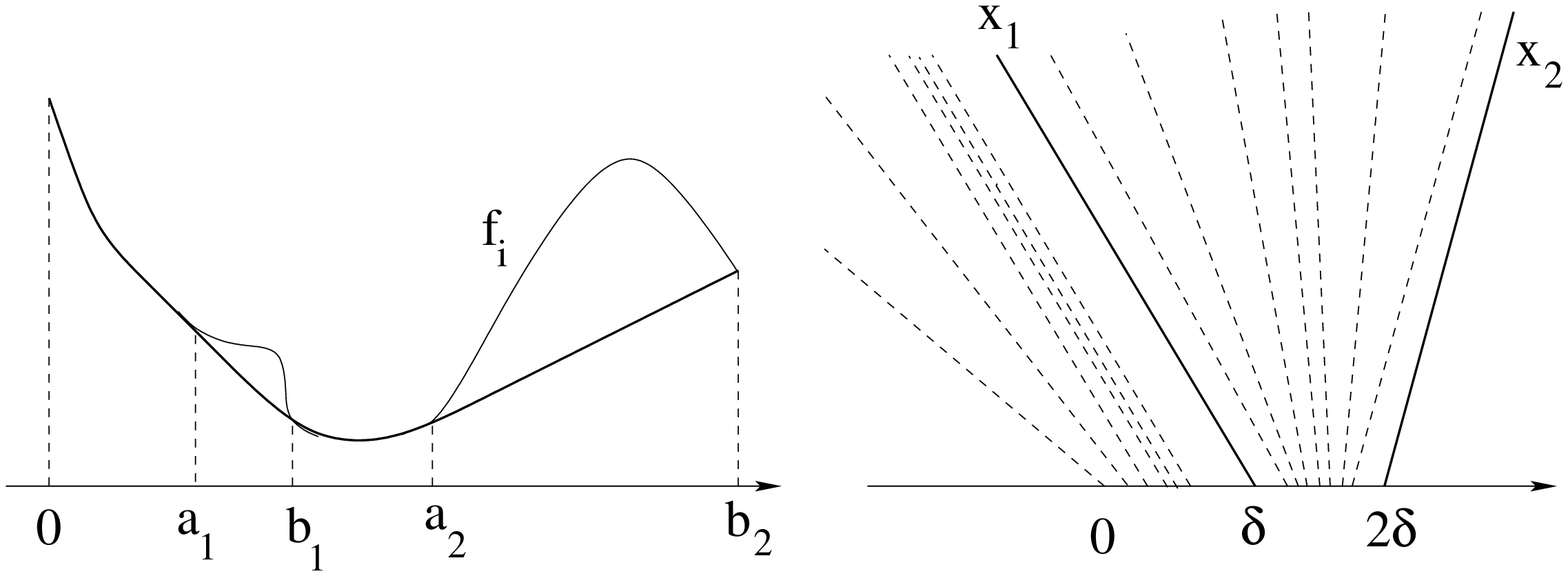,width=13cm}}}
\centerline{\hbox{figure 12~~~~~~~~~~~~~~~~~~~~~~~~~~~~~~~~~~~~~~~~~figure 13}}
\vskip 10pt
\endinsert

If (14.12) holds, we can make the two 
following observations.
\v
\n{\bf (i)}
For each $j=1,\ldots,N$, we claim that
the left and right states
$u(a_j)$, $u(b_j)$ are connected by a viscous travelling 
profile $U$ such that
$$U''=\big(A(U)-\sigma_{ij}\big)U',\qquad\qquad U(-\infty)=u(a_j),\quad 
U(\infty)=u(b_j).\eqno(14.13)$$
Here $\sigma_{ij}$ is the constant speed
$$\sigma_{ij}\doteq\sigma_i(\tau)=
{d\over d\tau}\, 
\hbox{conv}\,f_i(\gamma,\tau)\qquad\qquad \tau\in [a_j,b_j]\,.$$
To construct the function $U$,
consider the variable transformation $]a_j,\,b_j[\,\mapsto\R$,
say
$\tau\mapsto x(\tau)$, defined by
$$x\Big({a_j+b_j\over 2}\Big)=0\,,\qquad\qquad {dx(\tau)\over d\tau}=
{1\over v_i(\tau)}\,.$$
Let $\tau=\tau(x)$ be its inverse.
Then the function $U(x)\doteq u\big(\tau(x)\big)$ is the required
travelling wave profile.  Indeed, $U$ obviously
takes the correct limits at
$\pm\infty$.  Moreover,
$$U'={du\over d\tau}\,{d\tau\over dx}=v_i\tr_i,$$
$$\eqalign{U''&=v_{i,x}\tr_i+v_i\tr_{i,x}\cr
&=v_i\big(v_{i,\tau}\tr_i+v_i\tr_{i,\tau}\big)\cr
&=v_i(\tla_i-\sigma_{ij})\tr_i+v_i^2\big(\tr_{i,u}\tr_i
+\tr_{i,v}(\tla_i-\sigma_{ij})v_i\big).\cr}$$
Recalling the identity (4.22), we see that $U$ also satisfies the differental
equation in (14.13), thus proving our claim.
\v
\n{\bf (ii)} On the intervals where
$f_i(\gamma,\tau)=\hbox{conv}\,f_i(\gamma,\tau)$ we have 
$v_i(\tau)=0$. Hence, by the first equation in (14.6) and by (4.16), 
$u_\tau=\tr_i=r_i(u)$ is an $i$-eigenvector
of the matrix $A(u)$.
\v
\n{\bf 4.} In general, even if the condition 
(14.12) is satisfied, we do not expect that the function
$\omega$ has the regularity specified in Step 2.
However, we now show that it can be approximated in $\L^1$ by functions
$\omega_\delta$ satisfying all the required assumptions.
To fix the ideas, let
$$0\doteq b_0\leq a_1<b_1\leq a_2<b_2\leq \cdots 
\leq a_N<b_N\leq a_{N+1}\doteq s\,.$$
A piecewise smooth viscosity solution can be defined as follows
(fig.~13).
Fix a small $\delta>0$. For each $k=0,\ldots,N$, consider a smooth
non-decreasing map
$$\tau_k:\big[k\delta,\,(k+1)\delta\big]\mapsto [b_k,\,a_{k+1}]$$
such that
$$\tau_k(x)=\cases{b_k\qquad &if\quad $x\leq k\delta+\delta/3$,\cr
a_{k+1}\qquad &if\quad $x\geq k\delta+2\delta/3$.\cr
}$$
We then define the initial condition
$$\omega_\delta(0,x)\doteq 
\cases{u(0)\qquad &if \quad $x<0\,$,\cr
u\big(\tau_k(x)\big)
\qquad &if\quad $k\delta<x<(k+1)\delta\,,$\cr
u(s)\qquad &if \quad $x>(N+1)\delta\,$.\cr}$$
A corresponding solution of the Cauchy problem can then
be constructed by the method of characteristics:
$$\omega_\delta(t,x)\doteq 
\cases{u(0)=u^-\qquad &if\qquad $x<\sigma_i(0)\,t\,$,\cr
u(0)\qquad &if\qquad $x=y+\sigma_i\big(\tau_k(y)\big)t$~~for some~~
$y\in \,\big]k\delta,\,(k+1)\delta\big[\,$,\cr
u(s)=u^+\qquad &if\qquad $x>(N+1)\delta+\sigma_i(s)t\,$.\cr}$$
It is clear that the above function $\omega_\delta$ 
satisfies all of the 
assumptions considered in Step 2.  Hence $\omega_\delta\in VVL$.
Letting $\delta\to 0$ we have 
$\big\|\omega_\delta(t)-\omega(t)\big\|_{\L^1}\to 0$
for every $t$.  Therefore, by the last observation in 
Step 2 we conclude that also $\omega\in \vvl$.
\v
\n{\bf 5.}
To prove the Lemma in the general case, where the set in (14.12)
may by the union of infinitely many open intervals,
we use an approximation argument.
For each $\delta>0$, by slightly perturbing the 
values of $A$, we can construct a second
matrix valued function $A'$ with
$$\sup_u\big|A'(u)-A(u)\big|\leq \delta\,,
\eqno(14.14)
$$ such that the following properties hold.
For some right state $\tilde u^+$ with
$|\tilde u^+-u^+|\leq\delta$, the non-conservative Riemann problem
$$u_t+A'(u)u_x=0,\qquad\qquad u(0,x)=
\cases{u^-\quad &if\quad $x<0$,\cr
\tilde u^+\quad &if\quad $x>0$.\cr}\eqno(14.15)$$
admits a self-similar solution $\omega'$ which is limit of vanishing
viscosity approximations and satisfies
$$\int_{-3\beta}^{3\beta}\big|\omega'(t,x)-\omega(t,x)\big|\,dx\leq\delta
\qquad\qquad t\in [0,1]\,.\eqno(14.16)$$
Clearly, the fact that $\omega'$ is a limit of vanishing viscosity
approximations can be achieved by choosing $A'$ so that
a corresponding transformation $\T'_{s'}$ will admit as fixed point
some curve $\gamma':\tau\mapsto\big(u'(\tau)\,,
v_i'(\tau)\,,\sigma_i'(\tau)\big)$
for which
$u'(0)=u^-$, $u'(s')\approx u^+$ and
with $f_i(\gamma',\tau)$ differing from its convex envelope
on a finite number of open intervals.

Call $\omega^\ve$ the solution of the viscous Riemann problem
(13.1) with initial data (14.1).
Using (13.7) with $v^\ve\equiv u^+$, $a=0$, $b=\infty$,
for all $t\in [0,1]$ we obtain
$$\lim_{\ve\to 0}\int_{\beta}^\infty\big|\omega^\ve(t,x)-\omega(t,x)\big|
\,dx\leq \lim_{\ve\to 0}
\int_{\beta}^\infty |u^--u^+|\cdot \alpha e^{(\beta t-x)/\ve}\,dx=0\,.
\eqno(14.17)$$
Similarly,
$$\lim_{\ve\to 0}\int_{-\infty}^{-\beta}\big|\omega^\ve(t,x)-\omega(t,x)\big|
\,dx=0\,.\eqno(14.18)$$
To establish the convergence also on the interval $[-\beta,\beta]$,
call $v^\ve$ the solution of the Cauchy problem
$$v^\ve_t+A'(v^\ve)v^\ve_x=\ve v^\ve_{xx}\,,\qquad\qquad
v^\ve(0,x)=\cases{u^-\quad &if\quad $x<0\,$,\cr
\tilde u^+\quad &if\quad $0<x<3\beta\,$,\cr
u^+\quad &if\quad $x>3\beta\,.$\cr}$$
Clearly, 
$$\lim_{\ve\to 0}\int_{-\beta}^\beta\big|v^\ve(t,x)-\omega'(t,x)\big|
\,dx=0\,.\eqno(14.19)$$
because $\omega'$ is a vanishing viscosity limit and because of the
finite propagation speed.  Using the triangle inequality we can write
$$\eqalign{\limsup_{\ve\to 0}
\int_{-\beta}^\beta\big|\omega^\ve(t,x)-&\omega(t,x)\big|
\,dx\leq\limsup_{\ve\to 0}
\int_{-\beta}^\beta\big|\omega^\ve(t,x)-v^\ve(t,x)\big|\,dx
\cr
&+
\lim_{\ve \to 0}
\int_{-\beta}^\beta
\big|v^\ve(t,x)-\omega'(t,x)\big|\,dx
+\int_{-\beta}^\beta\big|\omega'(t,x)-\omega(t,x)\big|\,dx
\,.\cr}\eqno(14.20)$$
Since $t\mapsto \omega^\ve(t)=S^\ve_t\omega(0)$ is a trajectory 
of the Lipschitz semigroup $S^\ve$,
recalling (14.14) we have the estimate
$$\eqalign{\big\|v^\ve(t)-\omega^\ve(t)\big\|_{\L^1}
\,dx&\leq L\,\big\|v^\ve(0)-\omega^\ve(0)\big\|_{\L^1}
+L\cdot \int_0^t\left\{\lim_{h\to 0+}{\big\|v^\ve(s+h)-S^\ve_h v^\ve(s)
\big\|_{\L^1}\over h}
\right\}ds\cr
&\leq  3\beta L\,|\tilde u^+-u^+|
+L\cdot \int_0^t\left\{\int \Big|A(v^\ve(s,x))-A'(v^\ve(s,x))\Big|
\,\big|v^\ve_x(s,x)\big|\,dx
\right\}ds\cr
&\leq 3\beta L\delta+L\delta \,C\cdot\tv\big\{v^\ve(0)\big\}\cr
&\leq C''\delta\,,\cr}\eqno(14.21)$$
for some constant $C''$.
Estimating the right hand side of (14.20) by means of
(14.21), (14.19) and (14.16), we obtain
$$
\limsup_{\ve\to 0}
\int_{-\beta}^\beta\big|\omega^\ve(t,x)-\omega(t,x)\big|
\,dx\leq C''\delta+0+\delta.$$
Since $\delta>0$
can be arbitrarily small, together with (14.17)-(14.18) this yields
$$\lim_{\ve\to 0}\big\|\omega^\ve(t)-\omega(t)\big\|_{\L^1}=0
\qquad\qquad \hbox{for all}~t\in [0,1]\,,$$
completing the proof.
\endproof
\v
\n{\bf Remark 14.2.} The transformation
$\T_{i,s}$ defined at (14.4) depends on the vectors $\tr_i$, and
hence on
center manifold (which is not unique).
However, the curve $\gamma$ that we obtain as fixed point
of $\T_{i,s}$ involves only a concatenation of bounded travelling
profiles or stationary solutions.   These are bounded solutions
of (4.2), and will certainly be included in every center manifold.
For this reason, the curve $\gamma$
(and hence the solution of the Riemann problem) is independent of our
choice of the center manifold.
\v
For negative values of the parameter $s$, 
a right state $u^+=\Psi_i(s)$ can be constructed
exactly in the same way as before, except that one now takes the upper
concave envelope of $f_i$ :
$$\hbox{conc}\,f_i(\gamma,\tau)\doteq\sup\Big\{
\theta f_i(\gamma,\tau')+(1-\theta)f_i(\gamma,\tau'')\,;~~~
\theta\in [0,1]\,,~~\tau',\tau''\in [0,s]\,,~~
\tau=\theta \tau'+(1-\theta)\tau''\Big\}\,,$$
instead of the lower concave envelope.
\v
Our next step is to study the regularity of the curve 
of right states $u^+=\Psi_i(s)$.
\v
\n{\bf Lemma 14.3.} {\it 
Given a left state $u^-$ and $i\in \{1,\ldots,n\}$, the curve
of right states $s\mapsto\Psi_i(s)$ is Lipschitz continuous and satisfies}
$$\lim_{s\to 0}\, {d\Psi_i(s)\over ds}=r_i(u^-).\eqno(14.22)$$
\v
\n{\bf Proof.}  We assume $s>0$, the other case being entirely similar.
For sake of clarity, let us introduce some notations.
For fixed $i$ and $s>0$, let $\gamma^{i,s}=(u^{i,s}, v_i^s, \sigma_i^s)$ 
be  the fixed point of the transformation $\T_{i,s}$ in
(14.4).
Then by definition
$$\Psi_i(s)\doteq u^{i,s}(s)\,.$$
For $0<s'<s$, let $\gamma'\doteq (u', v_i',\sigma_i')$ be the restriction 
of $\gamma^{i,s}$ to the subinterval $[0,s']$.
Since $\T_{i,s'}$ is a strict contraction,
the distance of $\gamma'$ from the fixed point of $\T_{i,s'}$ 
is estimated as
$$\big\|\gamma'-\gamma^{i,s'}\big\|_\dagger
=\O(1)\cdot\big\|\gamma'-\T_{i,s'}\gamma'\big\|_\dagger
=\O(1)\cdot (s-s')\,s\,.$$
In particular,
$$\big|u^{i,s'}(s')-u^{i,s}(s')\big|=\O(1)\cdot (s-s')\,s\,.$$
Observing that
$$u^{i,s}(s)-u^{i,s}(s')=\int_{s'}^s \tr_i\big(
u^{i,s}(\vars),\,v_i^s(\vars),\,\sigma^s_i(\vars)\big)\,d\vars
=(s-s')\cdot r_i(u^-)+\O(1)\cdot (s-s')s\,,$$
we conclude
$$\big|\Psi_i(s')-\Psi_i(s')-(s-s')r_i(u^-)\big|=\O(1)\cdot (s-s')\,s\,.
\eqno(14.23)$$
By
(14.23), the map $s\mapsto \Psi_i(s)$ is Lipschitz continuous,
hence differentiable almost everywhere, by Rademacher's theorem.
The limit in (14.22) is again a consequence of (14.23).
\endproof
\v
Thanks to the previous analysis, the
solution of the general Riemann problem (14.2) can now 
be constructed following a standard procedure.
Given a left state $u^-$, call $s\mapsto \Psi_i(s)(u^-)$
the curve of right states that can be connected to $u^-$ by
$i$-waves.
Consider the composite mapping
$$\Psi:(s_1,\ldots,s_n)\mapsto \Psi_n(s_n)\circ\cdots\circ
\Psi_1(s_1)(u^-)\,.$$
By Lemma 14.3 and a version of the implicit function theorem
valid for Lipschitz continuous maps (see [Cl], p.253),
$\Psi$ is a one-to-one mapping
from a neighborhood of the origin in $\R^n$ onto a neighborhood of $u^-$.
Hence, for all $u^+$ sufficiently close to $u^-$, one can find
unique values $s_1,\ldots,s_n$ such that $\Psi(s_1,\ldots,s_n)=u^+$.
In turn, this yields 
intermediate states
$u_0=u^-, u_1,\ldots,u_n=u^+$ such that
each Riemann problem with data $u_{i-1}, u_i$ admits
a vanishing viscosity solution $\omega_i=\omega_i(t,x)$ 
consisting only of $i$-waves.
By strict hyperbolicity, we can now choose
intermediate speeds 
$$-\infty\doteq\lambda_0'<\lambda_1'<\lambda_2'
<\cdots<\lambda_{n-1}'<\lambda'_n=\infty$$
such that
all $i$-waves in the solution $\omega_i$ have speeds contained inside the 
interval $[\lambda_{i-1}',\,\lambda_i']$.
The general solution of the general Riemann problem (14.2)
is then given by
$$\omega(t,x)=\omega_i(t,x)\qquad\qquad \hbox{for}~\lambda_{i-1}'<
{x\over t}<\lambda_i'\,.\eqno(14.24)$$
Because of Lemma 14.2, it is clear that the function $\omega$
is the unique limit of viscous approximations:
$$\lim_{\ve\to 0+}\big\|\omega(t)-S^\ve_t\omega(0)\big\|_{\L^1}=0\qquad\quad
\hbox{for every}~t\geq 0\,.\eqno(14.25)$$
\vsk
\n{\medbf 15 - Viscosity solutions and uniqueness of the semigroup}
\v
In [B3], one of the authors introduced a definition of
{\it viscosity solution} for a system of conservation laws,
based on local integral estimates. 
Assuming the existence of a Lipschitz semigroup of
entropy weak solutions,
it was proved that such a semigroup is necessarily unique and every
viscosity solution coincides with a semigroup trajectory.
We shall follow here exactly the same approach, 
in order to prove the uniqueness
of the Lipschitz semigroup constructed in (13.9) as limit of 
vanishing viscosity approximations.

Toward the definition of a {\it viscosity solution} for the
general hyperbolic system
$$u_t+A(u)u_x=0,
\eqno(15.1)$$
we first introduce some notations.
Given a function $u=u(t,x)$ and a point $(\tau,\xi)$,
we denote by  $U^\sharp_{(u;\tau,\xi)}$ the solution
of the Riemann problem (14.1)
with initial data
$$ u^-=\lim_{x\to \xi-} u(\tau,x),\qquad\qquad u^+=\lim_{x\to \xi+}
u(\tau,x).\eqno(15.2)$$
Of course, we refer here to the vanishing viscosity solution
constructed in Section 14.
In addition, we define $U^\flat_{(u;\tau,\xi)}$ as the solution
of a linear hyperbolic Cauchy problem with constant coefficients:
$$w_t+\Hat A w_x=0,\qquad\qquad w(0,x)=u(\tau,x).\eqno(15.3)$$
Here $\Hat A\doteq A\big(u(\tau,\xi)\big)$.
Observe that (15.3) is obtained from the quasilinear system
(15.1)
by ``freezing'' the coefficients of the matrix $A(u)$ 
at the point $(\tau,\xi)$ and choosing $u(\tau)$ as initial data.

As in [B3], the notion of {\it viscosity solution}
is now defined
by locally comparing a function $u$ with 
the self-similar solution of a Riemann problem
and with the solution of a linear hyperbolic system with 
constant coefficients.
\v
\n{\bf Definition 15.1.} 
A function $u=u(t,x)$ is a {\bf viscosity solution} of the 
system (15.1) if $t\mapsto u(t,\cdot)$ is continuous
as a map with values into $\Ll$, and moreover the 
following integral estimates hold.
\v
\item{(i)}  At every point $(\tau,\xi)$, for every $\beta'>0$ one has
$$\lim_{h\to 0+} {1\over h}\int_{\xi-\beta' h}
^{\xi+\beta' h}\Big| u(\tau+h,~x)- U^\sharp
_{(u;\tau,\xi)}(h,~x-\xi)\Big|~dx~=~0.\eqno(15.4)$$
\v
\item{(ii)}  There
exist constants $C,\beta>0$ such that,
or every $\tau\geq 0$ and $a<\xi<b$, one has
$$\limsup_{h\to 0+}{1\over h}\int_{a+\beta h}
^{b-\beta h}\Big| u(\tau+h,~x)- U^\flat
_{(u;\tau,\xi)}(h,x)\Big|~dx\leq C\cdot \Big(\tv\big\{
u(\tau);~]a,\,b[ \ \big\}\Big)^2.\eqno(15.5)$$ 
\v
The main result of this section shows that the above viscosity
solutions coincide 
precisely with the limits
of vanishing viscosity approximations.
\v
\n{\bf Lemma 15.2.} {\it Let $S:\D\times[0,\infty[\,\mapsto \D$
be a semigroup of vanishing viscosity solutions, constructed
as limit of a sequence $S^{\ve_m}$ as in (13.9) and defined on
a domain $\D\subset\Ll$ of functions with small total variation. 
A map 
$u:[0,T]\mapsto \D$ satisfies
$$u(t)=S_t u(0)\qquad\qquad \hbox{for all}~t\in [0,T]\eqno(15.6)$$
if and only if $u$ is a viscosity solution of(15.1).}
\v
\n{\bf Proof.~~Necessity:}
Assume that (15.6) holds.  By (13.11), the map $t\mapsto u(t)$ is
continuous. 
Let any $\beta'>0$ be given and
let $L$, $\beta$ be the constants
in (13.13). Then, for any
$(\tau,\xi)$, an application of (13.13) yields
$$\eqalign{&\int_{\xi-\beta' h}
^{\xi+\beta' h}\Big| u(\tau+h,~x)- U^\sharp
_{(u;\tau,\xi)}(h,~x-\xi)\Big|\,dx\cr
&\quad\leq L\cdot\left\{\int_{\xi-(\beta+\beta') h}^\xi
\Big| u(\tau,x)- u(\tau,\xi-)\Big|\,dx+
\int_\xi^{\xi+(\beta+\beta') h}
\Big| u(\tau,x)- u(\tau,\xi+)\Big|\,dx\right\}\cr
&\quad\leq L(\beta+\beta') h\bigg\{
\sup_{\xi-(\beta+\beta') h\,<\,x\,<\,\xi}\,\big|
u(\tau,x)- u(\tau,\xi-)\big|~+~
\sup_{\xi\,<\,x\,<\,\xi+(\beta+\beta') h}\,\big|
u(\tau,x)- u(\tau,\xi+)\big|\bigg\}.\cr}$$
Hence (15.4) is clear.
\v
To prove the second estimate, fix $\tau$ and $a<\xi<b$.
Define the function
$$\bar v(x)\doteq\cases{u(\tau, a+)\qquad &if\quad $x\leq a\,,$\cr
u(\tau, x)\qquad &if\quad $a<x<b \,,$\cr
u(\tau, b-)\qquad &if\quad $x\geq b \,.$\cr
}$$
Call $v^\ve,w^\ve$ respectively the solutions of the
viscous systems
$$v^\ve_t+A(v^\ve) v^\ve_x=\ve\,v^\ve_{xx}\,,
\qquad\qquad w^\ve_t+\Hat A w^\ve_x=\ve\,w^\ve_{xx}\,,\eqno(15.7)$$
with the same initial data $v^\ve(0,x)=w^\ve(0,x)=\bar v(x)$.

Recalling that $S^\ve$
is a semigroup with Lipschitz constant $L$, 
as in [B3], [B5]
we can use the error formula
$$\eqalign{\big\|w^\ve(h)-v^\ve(h)\big\|_{\L^1}&=
\big\|w^\ve(h)-S^\ve_h \bar v\big\|_{\L^1}\cr
&\leq
L\cdot\int_0^h \left\{\liminf_{r\to 0+}\,{
\big\|w^\ve(t+r)-S^\ve_r w^\ve(s)\big\|_{\L^1}\over r}\right\}\,dt\cr
&\leq L\cdot\int_0^h \int 
\Big|\Hat A-A\big( w^\ve(t,x)\big)
\Big|\, \big|w^\ve_x(t,x)\big|\,dxdt\cr
&\leq L \, h\,\left(\sup_{t,x}
\Big|A\big( w^\ve(0,\xi)\big)-A\big( w^\ve(t,x)\big)\Big|\right)
\cdot
\sup_t \big\|w_x^\ve(t)\big\|_{\L^1}\cr
&\leq C\, h\,\Big(\tv\{\bar v\}\Big)^2,
\cr}$$
for some constant $C$. 
Letting $\ve\to 0$ and using the estimate (13.13) on the
finite speed of propagation,
we obtain
$$\eqalign{
{1\over h}\int_{a+\beta h}
^{b-\beta h}&\Big| u(\tau+h,~x)- U^\flat
_{(u;\tau,\xi)}(h,~x)\Big|~dx\leq {1\over h}\lim_{\ve\to 0} \int_{a+\beta h}
^{b-\beta h}\big| v^\ve(h,x)- w^\ve(h,x)\big|~dx\cr
&\leq C\,\Big(\tv\{\bar v\}\Big)^2~=~
C\,\Big(\tv\big\{\bar u\,;~]a,b[\,\big\}\Big)^2.\cr}$$
This proves (15.5), with $\beta$ the constant in (13.13).
\v
\n{\bf Sufficiency:} Let $u=u(t,x)$ be a viscosity
solution of (15.1).  By assumption, the map $t\mapsto
u(t)$ is continuous with values in a domain $\D\subset\Ll$
of functions with small total variation. From (15.5) and this uniform bound
on the total variation it follows that this map 
is actually Lipschitz continuous:
$$\big\|u(t)-u(s)\big\|_{\L^1}\leq L''\,|t-s|\,,\eqno(15.8)$$
for some constant $L''$ and all $s,t\in [0,T]$.
Let $L$ be the Lipschitz constant of the semigroup $S$, as in (13.13).
Given any interval $[a,b]$, thanks to (15.8) one has the error estimate
$$\eqalign{\int_{a+t\beta}^{b-t\beta}&
\Big|u(t,x)-\big(S_tu(0)\big)(x)\Big|
\,dx\cr
&\leq
L\cdot\int_0^t
\left\{\liminf_{h\to 0+}
{1\over h}\int_{a+(\tau+h)\beta}
^{b-(\tau+h)\beta}\Big| u(\tau+h,x)-\big(S_h u(\tau)\big)(x)\Big|\,dx
\right\}d\tau\,.\cr}
\eqno(15.9)$$
To prove the identity (15.6) it 
thus suffices to show that the integrand on the right hand side of
(15.9) vanishes for all $\tau\in [0,T]$.

Fix any time $\tau\in [0,T]$
and let $\epsilon>0$ be given.
Since the total variation of $u(\tau,\cdot)$ is finite, we can
choose finitely many points 
$$a+\tau\beta=x_0<x_1<\cdots<x_N=b-\tau\beta$$
such that, for every $j=1,\ldots,N$,
$$\tv\big\{u(\tau,\cdot)\,;~]x_{j-1},\,x_j[\,\big\}<\epsilon\,.$$
By the necessity part of the theorem, which has been already
proved, the function $w(t,\cdot)\doteq S_{t-\tau}u(\tau)$ 
is itself a viscosity solution and hence it also
satisfies the estimates (15.4)-(15.5).
We now consider the midpoints $y_j\doteq(x_{j-1}+x_j)/2$.
Using the estimate (15.4) at each of the points $\xi=x_j$ and
the estimate (15.5) with $\xi\doteq y_j$ 
on each of the intervals $]x_{j-1},~x_j[\,$, 
taking $\beta>0$ sufficiently large we now compute
$$\eqalign{
\limsup_{h\to 0+}\,
&{1\over h}\int_{a+(\tau+h)\beta}
^{b-(\tau+h)\beta}\Big| u(\tau+h,x)-\big(S_h u(\tau)\big)(x)\Big|\,dx\cr
&\leq\sum_{j=1}^{N-1}
\limsup_{h\to 0+}
{1\over h}\int_{x_j-h\beta}
^{x_j+h\beta}
\bigg(
\Big| u(\tau+h,x)- U^\sharp
_{(u;\tau,x_j)}(\tau+h,~x)\Big|
\cr
&\qquad\qquad\qquad\qquad\qquad\qquad\qquad\qquad+\Big|U^\sharp
_{(u;\tau,x_j)}(\tau+h,~x)-\big(S_h u(\tau)\big)(x)
\Big|\bigg)
\,dx\cr
&\qquad
+\sum_{j=1}^{N}
\limsup_{h\to 0+}
{1\over h}\int_{x_{j-1}+h\beta}^{x_j-h\beta}
\bigg(
\Big| u(\tau+h,x)- U^\flat
_{(u;\tau,y_j)}(h,~x-x_j)\Big|
\cr
&\qquad\qquad\qquad\qquad\qquad\qquad\qquad\qquad +\Big|U^\flat
_{(u;\tau,y_j)}(h,~x-x_j)-\big(S_h u(\tau)\big)(x)
\Big|\bigg)
\,dx\cr
&\leq 0+\sum_{j=1}^N C\, \Big(\tv\big\{
u(\tau);~]x_{j-1},~x_j[\big\}\Big)^2\cr
&\leq C\, \epsilon\cdot\tv\big\{
u(\tau);~]\alpha+\tau\beta\,,~b-\tau\beta[\big\}\,.\cr}$$ 
Since $\epsilon>0$ was arbitrary, the integrand
on the right hand side of (15.9) must vanish at time $\tau$.  
This completes the proof of the lemma.
\endproof
\v
\n{\bf Remark 15.3.} From the proof of the sufficiency part, 
it is clear that the identity (15.6)
still holds if we require that the integral estimates (15.4)
hold only for $\tau$ outside a set of times $\N\subset [0,T]$
of measure zero.
By a well known result in the theory of BV functions [EG], 
any BV function of two variables $u=u(t,x)$
is either approximately continuous or has an approximate jump
discontinuity at every point $(\tau,\xi)$, with $\tau$ outside a
set $\N$ having zero measure. 
To decide whether a function $u$ is a viscosity solution, it thus
suffices to check (15.4) only at points of approximate jump,
where the Riemann problem is solved in terms of a single shock.
\v
Using Lemma 15.2, we now obtain at one stroke the uniqueness
of viscosity solutions and of vanishing viscosity limits:
\v
\n{\bf Completion of the proof of Theorem 1.} What remains to
be proved is that the whole family of viscous
approximations converges to a unique limit, i.e.
$$\lim_{\ve\to 0+} S^\ve_t\bar u=S_t\bar u\,,\eqno(15.10)$$ 
where the limit holds over all real values of $\ve$ and not only
along a particular sequence $\{\ve_m\}$.
If (15.10)
fails, we can find $\bar v$, $\tau$ and two different sequences
$\ve_m,\ve'_m\to 0$ such that
$$\lim_{m\to\infty}S^{\ve_m}_\tau\bar v\not= 
\lim_{m\to\infty}S^{\ve'_m}_\tau\bar v\,.
\eqno(15.11)$$
By extracting further subsequences, we can assume that
the limits
$$\lim_{m\to\infty} S^{\ve_m}_t\bar u=S_t\bar u\,,\qquad\qquad
\lim_{m\to\infty} S^{\ve'_m}_t\bar u=S'_t\bar u\,,\eqno(15.12)$$ 
exist in $\Ll$, for all $t\geq 0$ and $\bar u\in\U$.
By the analysis in Section 13, both $S$ and $S'$ are semigroups
of vanishing viscosity solutions. In particular, the necessity
part of
Lemma 14.2 implies that the map $t\mapsto v(t)\doteq S_t\bar v$
is a viscosity solution of (15.1), while the sufficiency part
implies $v(t)=S'_t v(0)$ for all $t\geq 0$. 
But this is in contradiction with (15.11), hence the unique limit
(15.10) is well defined.
\endproof
\v
\n{\bf Remark 15.4.} The above uniqueness result is obtained within the
family of vanishing viscosity limits of the form (1.13)$_{\ve}$,
with unit viscosity matrix.  
In the more general case (1.21)$_\ve$, if the system is not 
in conservation form,
we expect that the limit of solutions as $\ve\to 0$ will depend on the form
of the viscosity matrices $B(u)$. Indeed, by choosing
different matrices $B(u)$, one will likely alter the vanishing 
viscosity solutions of the Riemann problems (14.2).
In turn, this affects the definition of viscosity solution
at (15.4).
\vsk 
\n{\medbf 16 - Dependence on parameters and large time asymptotics}
\v
We wish to derive here a simple estimate on how the viscosity 
solution 
changes, depending on hyperbolic matrices $A(u)$.
\v
\n{\bf Corollary 16.1.} 
{\it
Assume that the two hyperbolic systems
$$\eqalign{u_t+A(u)u_x&=0,\cr
u_t+\Hat A(u) u_x&=0\,,\cr}
$$
both satisfy the hypotheses of Theorem 1.  Call $S,\Hat S$
the corresponding semigroups of viscosity solutions. 
Then, for every initial data $\bar u$ with small total variation, one has
the estimate}
$$\big\| S_t\bar u-\Hat S_t\bar u\big\|_{\L^1}=\O(1)\cdot
t\Big(\sup_u\big|\Hat A(u)-A(u)\big|\Big)
\cdot\tv\{\bar u\}\,.\eqno(16.1)$$
\v
\n{\bf Proof.}  Call $S^\ve,\Hat S^\ve$ the semigroups
of solutions to the corresponding viscous problems
$$\eqalign{u_t+A(u)u_x&=\ve\, u_{xx},\qquad\qquad
u_t+\Hat A(u) u_x=\ve\, u_{xx}\,.\cr}
$$
Let $L$ be the Lipschitz constant in (1.16) and call $w^\ve(t)\doteq
\Hat S^\ve_t\bar u$.
For every $t\geq 0$ we have the error estimate
$$\eqalign{\big\|\Hat S_t^\ve\bar u- S_t^\ve\bar u\big\|_{\L^1}
&=\big\|w^\ve(t)- S_t^\ve\bar u\big\|_{\L^1}\cr
&\leq
L\cdot\int_0^t \left\{\liminf_{h\to 0+}\,{
\big\|w^\ve(s+h)-S^\ve_h w^\ve(s)\big\|_{\L^1}\over h}\right\}\,ds\cr
&\leq L\cdot\int_0^t \int \Big|\Hat A\big( w^\ve(s,x)\big)-
A\big( w^\ve(s,x)\big)
\Big|\, \big|w^\ve_x(s,x)\big|\,dxds\cr
&\leq L \,\Big(\sup_u \big|\Hat A(u)- A(u)\big|\Big)\int_0^t 
\big\|w_x^\ve(s)\big\|_{\L^1}\,ds\cr
&=\O(1) \cdot t\,\Big(\sup_u \big|\Hat A(u)- A(u)\big|\Big)\
\tv\{\bar u\}
.\cr}$$
\endproof
\v
Next, we show that some semigroup trajectories are asymptotically
self-similar.  
\v
\n{\bf Corollary 16.2.} {\it 
Under the same assumption of Theorem 1,
consider an initial data $\bar u$ with small total variation, such that
$$\int_{-\infty}^0\big|\bar u(x)-u^-\big|\,dx+
\int_0^\infty\big|\bar u(x)-u^+\big|\,dx<\infty\,,\eqno(16.2)$$
for some states $u^-,u^+$.
Call $\omega(t,x)=\tilde\omega(x/t)$ the self-similar solution of the
corresponding Riemann problem (14.2).
Then the solution of the viscous Cauchy problem
$$u_t+A(u)u_x=u_{xx}\,,\qquad\qquad u(0,x)=\bar u(x)\eqno(16.3)$$
satisfies}
$$\lim_{\tau\to\infty} \int \big|u(\tau,\,\tau y)-\tilde\omega(y)
\big|\,dy=0\,.\eqno(16.4)$$
\v
\n{\bf Proof.} The assumption on $\omega$ implies that the limit
(14.25) holds.
For fixed $\tau$, call $\ve\doteq 1/\tau$ and consider the function
$v^\ve(t,x)\doteq u(\tau x,\,\tau t)$.
Clearly, $v^\ve$ satisfies the equation
$$v^\ve_t+A(v^\ve)v^\ve_x=\ve\,v^\ve_{xx}\,,\qquad\qquad
v^\ve(0,x)=\bar u(x/\ve).$$
Therefore,
$$\eqalign{\int\big|u(\tau,\,\tau y)-\tilde \omega(y)\big|\,dy&=
\int \big|v^\ve(1,x)-\omega(1,x)\big|\,dx\cr
&\leq \big\|S^\ve_1 v^\ve(0)-S^\ve_1\omega(0)\big\|_{\L^1}
+\big\|S^\ve_1\omega(0)-\omega(1)\big\|_{\L^1}\,.\cr}\eqno(16.5)
$$
Observing that 
$$\eqalign{\big\|S^\ve_1 v^\ve(0)-S^\ve_1\omega(0)\big\|_{\L^1}
&\leq L\cdot \big\|S^\ve_1 v^\ve(0)-S^\ve_1\omega(0)\big\|_{\L^1}
\cr
&=L\ve\,\left(
\int_{-\infty}^0\big|\bar u(x)-u^-\big|\,dx+
\int_0^\infty\big|\bar u(x)-u^+\big|\,dx\right),\cr}$$
and using (14.25), from 
(16.5) we obtain (16.4).
\endproof

\vsk
\c{\medbf Appendix A}
\v
We derive here the explicit form of
the evolution equations 
(6.1), for the variables
$v_i$ and $w_i$ defined by the decomposition
$$
u_x = \sum_i v_i \tilde r_i \bigl(u,v_i,\lambda_i^* -
\theta(w_i/v_i) \bigr),
\qquad\quad u_t = \sum_i \bigl( w_i  - \lambda_i^* v_i \bigr) \tilde r_i
\bigl( u, v_i, \lambda_i^* - \theta(w_i/v_i)\bigr).
\eqno(A.1)$$
By checking one by one all source terms, we then
provide an alternative proof of Lemma 6.1.
The computations are lengthy but straightforward:
one has to
rewrite the evolution equations
for $u_x$ and $u_t$ :
$$\left\{ \eqalign{ (u_x)_t+\big(A(u)u_x\big)_x -(u_x)_{xx}&=0\,,\cr
(u_t)_t+\big(A(u)u_t\big)_x -(u_t)_{xx}&=\big(u_x\bullet A(u)\big)u_t
-\big(u_t\bullet A(u)\big)u_x\,,
\cr}\right.\eqno(A.2)$$
in terms of $v_i,w_i$.
For convenience, we set
$
\theta_i \doteq \theta(w_i/v_i)
$.
The fundamental relation (4.23) can be written as
$$
v_i \tilde r_{i,u}\tr_i - A(u) \tilde r_i = - \tilde
\lambda_i \tilde r_i + \bigl( - \tilde \lambda_i + \lambda_i^* -
\theta_i \bigr) v_i \tilde r_{i,v}\,.\eqno(A.3)
$$
Differentiating (A.1) w.r.t.~$x$ and using (A.3) we obtain
$$\eqalign{
u_{xx} - A(u) u_x &=\sum_i v_{i,x} \tilde r_i + \sum_i v_i\tr_{i,x}-  
\sum_i A(u) v_i\tilde r_i\cr
&=
\sum_i v_{i,x} \tilde r_i + \sum_i v_i \Bigl[ v_i  
\tilde r_{i,u}\tr_i - A(u) \tilde r_i \Bigr] 
\cr 
&\qquad + \sum_i v_i 
\Bigl[ v_{i,x} \tilde r_{i,v} -
\theta_i' \bigl( (v_i w_{i,x} - w_i v_{i,x} ) / v_i^2\big) \tilde
r_{i,\sigma}  \Bigr]
+ \sum_{i \not= j} v_i v_j \tr_{i,u}\tr_j \cr
&=\sum_i (v_{i,x}-\tla_i v_i)\tr_i+\sum_i 
( - \tilde \lambda_i + \lambda_i^* -
\theta_i ) v_i^2 \tilde r_{i,v}\cr
&\qquad +
 \sum_i v_i 
\Bigl[ v_{i,x} \tilde r_{i,v} -
\theta_i' \bigl( (v_i w_{i,x} - w_i v_{i,x} ) / v_i^2\big) \tilde
r_{i,\sigma}  \Bigr] + \sum_{i \not= j} v_i v_j \tr_{i,u}\tr_j \cr
&= \sum_i \big( v_{i,x} - \tilde \lambda_i v_i \big) \Bigl[ \tilde
r_i + v_i \tilde r_{i,v} + \theta_i' \bigl( w_i / v_i \bigr) \tilde
r_{i,\sigma} \Bigr] + \sum_i \big( w_{i,x} - \tilde \lambda_i w_i
\big) \Big[-\theta'_i \tilde r_{i,\sigma}\Big] \cr
& \qquad + \sum_i v_i^2 \bigl( \lambda_i^* - \theta_i
\bigr) \tilde r_{i,v} + \sum_{i \not= j} v_i v_j \tilde r_{i,u}\tr_j\,,
\cr}\eqno(A.4)$$
$$\eqalign{
u_{tx} - A(u) u_t&=
\sum_i ( w_{i,x} - \lambda_i^* v_{i,x} ) \tilde r_i +
\sum_i ( w_i - \lambda_i^* v_i ) \tr_{i,x}-
\sum_i( w_i - \lambda_i^* v_i ) A(u)\tr_i\cr
&=\sum_i ( w_{i,x} - \lambda_i^* v_{i,x} ) \tilde r_i
\sum_i( w_i - \lambda_i^* v_i ) \big[v_i\tr_{i,u}\tr_i -A(u)\tr_i\big]\cr
&\qquad 
+\sum_i ( w_i - \lambda_i^* v_i ) \Bigl[ v_{i,x}
\tilde r_{i,v} - \theta_i' \bigl(( v_i w_{i,x} - w_i v_{i,x})/
v_i^2\big) \tilde r_{i,\sigma}  \Bigr] 
+ \sum_{i \not= j} \bigl( w_i - \lambda_i^* v_i \bigr) v_j 
\tilde r_{i,u}\tr_j \cr
&=\sum_i ( w_{i,x} - \tilde \lambda_i w_i ) \tilde r_i -
\sum_i \lambda_i^* ( v_{i,x} - \tilde \lambda_i v_i)
\tilde r_i +\sum_i 
(w_i-\lambda_i^* v_i)\big(-\tla_i+\lambda_i^*-\theta_i\big)
v_i\tr_{i,v}
\cr
&\qquad +\sum_i( w_i - \lambda_i^* v_i ) \Bigl[ v_{i,x}
\tilde r_{i,v} - \theta_i' \bigl(( v_i w_{i,x} - w_i v_{i,x})/
v_i^2\big) \tilde r_{i,\sigma}  \Bigr] 
+ \sum_{i \not= j} \bigl( w_i - \lambda_i^* v_i \bigr) v_j 
\tilde r_{i,u}\tr_j \cr
&= \sum_i \big( v_{i,x} - \tilde \lambda_i v_i \big) \Bigl[ w_i r_{i,v}
+ \theta'_i \bigl( w_i / v_i \bigr)^2 \tilde r_{i,\sigma} \Bigr]
+ \sum_i \big( w_{i,x} - \tilde \lambda_i w_i \big) \Bigl[ \tilde
r_i - \bigl( \theta'_i w_i/v_i \bigr) \tilde r_{i,\sigma} \Bigr] \cr
&\qquad + \sum_i v_i w_i \bigl(
\lambda_i^* - \theta_i \bigr) \tilde r_{i,v} + \sum_{i \not= j}
w_i v_j \tilde r_{i,u}\tr_j \cr
&\qquad - \sum_i \lambda_i^* \biggl\{ \Bigl( v_{i,x} - \tilde \lambda_i v_i
\Bigr) \Bigl[ \tilde r_i + v_i \tilde r_{i,v} + \bigl( \theta_i' w_i /
v_i \bigr) \tilde r_{i,\sigma} \Bigr] + \Bigl( w_{i,x} - \tilde
\lambda_i w_i \Bigr)\Big[-\theta'_i \tilde r_{i,\sigma}\Big] \cr
&\qquad\qquad \qquad \qquad + v_i^2 \bigl( \lambda_i^* - \theta_i
\bigr) \tilde r_{i,v} + \sum_{j \not= i} v_i v_j 
\tilde r_{i,u}\tr_j \biggr\}.
\cr}\eqno(A.5)$$
Differentiating (A.1) w.r.t.~$t$ one obtains
$$\eqalign{
u_{xt} &= \sum_i v_{i,t} \tilde r_i + \sum_i v_i r_{i,t}
\cr
&= \sum_i v_{i,t} \tilde r_i + \sum_i v_i \Big[  v_{i,t} \tilde r_{i,v} - 
\theta'_i \bigl( (w_{i,t} v_i - w_i v_{i,t} )/ v_i ^2\big) \tilde 
r_{i\sigma} \Big] + \sum_{i,j} v_i \bigl( w_j - \lambda_j^* v_j
\bigr) \tr_{i,u}\tr_j   \cr
&= \sum_i v_{i,t} \Bigl[ \tilde r_i + v_i \tilde
r_{i,v} + \bigl( \theta'_i w_i/v_i \bigr) \tilde r_{i,\sigma} \Bigr] 
+\sum_i w_{i,t} \Big[- \theta'_i \tilde r_{i,\sigma} \Big]+
\sum_{i,j} v_i \bigl( w_j - \lambda_j^* v_j \bigr)
\tr_{i,u}\tr_j\,,  
\cr}
\eqno(A.6)$$
$$\eqalign{ 
u_{tt} &= \sum_i \bigl( w_{i,t} - \lambda_i^* v_{i,t} \bigr) \tilde
r_i + \sum_i \bigl( w_i - \lambda_i^* v_i \bigr) \tr_{i,t}
\cr
&= \sum_i \bigl( w_{i,t} - \lambda_i^* v_{i,t} \bigr) \tilde r_i +
\sum_i \bigl( w_i - \lambda_i^* v_i \bigr) \Big[  v_{i,t} \tilde r_{i,v} - 
\theta'_i \bigl( (w_{i,t} v_i - w_i v_{i,t} )/ v_i^2\big) \tilde 
r_{i\sigma} \biggl]   \cr
&\qquad+ \sum_{i,j} \bigl( w_i - \lambda_i^* v_i
\bigr) \bigl( w_j - \lambda_j^* v_j \bigr)
\tr_{i,u}\tr_j   \cr
&= \sum_i v_{i,t} \Bigl[ w_i \tilde
r_{i,v} + \theta'_i \bigl( w_i/v_i \bigr)^2 \tilde r_{i,\sigma} \Bigr]
+ \sum_i w_{i,t} \Bigl[ \tilde r_i - \theta'_i \bigl( w_i/v_i \bigr)
\tilde r_{i,\sigma} \Bigr] + \sum_{i,j} w_i \bigl( w_j
- \lambda_j^* v_j \bigr) \tr_{i,u}\tr_j   \cr
&\qquad - \sum_i \lambda_i^* \biggl\{ v_{i,t} \Bigl[ \tilde r_i + v_i \tilde
r_{i,v} + \bigl( \theta_i' w_i/v_i \bigr) \tilde r_{i,\sigma} \Bigr]
- w_{i,t} \theta_i' \tilde r_{i,\sigma}  + \sum_j v_i
\bigl( w_j - \lambda_j^* v_j \bigr)\tr_{i,u}\tr_j
\biggr\}.  
\cr}\eqno(A.7)
$$
Differentiating again $u_{xx}-A(u)u_x$ 
and $u_{tx} - A(u) u_t$ w.r.t.~$x$, from (A.4) and (A.5) one finds
$$\eqalign{u_{tx}&=\bigl(u_x \bigr)_{xx} - \bigl( A(u) u_x \bigr)_x \cr
&= \sum_i \big( v_{i,xx} - (\tilde \lambda_i v_i
)_x \big) \Bigl[ \tilde r_i + v_i \tilde r_{i,v} +
\theta_i' ( w_i/v_i) \tr_{i,\sigma} \Bigr] + \sum_i
\big(
w_{i,xx}  - ( \tilde \lambda_i w_i )_x \big) \Bigl[ -
\theta_i' \tilde r_{i,\sigma} \Bigr] \cr
&\qquad + \sum_i ( v_{i,x} - \tilde \lambda_i v_i ) \biggl[ \sum_j
v_j \tr_{i,u}\tr_j + 2 v_{i,x} \tilde r_{i,v} + \Big(
-\theta_{i,x} + ( \theta_i' w_i/v_i )_x \Big) \tilde
r_{i,\sigma} + \sum_j v_j v_i \tr_{i,vu}\tr_j  
\cr
&\qquad \qquad + v_i v_{i,x} \tilde r_{i,vv} + \big(- v_i \theta_{i,x} +
\theta_i' v_{i,x}
w_i / v_i \big) \tilde r_{i,v\sigma}  + \sum_j  v_j \theta_i'
(w_i /v_i ) \tr_{i,\sigma u}\tr_j - \theta_{i,x}
\theta_i' (w_i /v_i ) \tilde r_{i,\sigma \sigma} \biggr]   \cr
&\qquad + \sum_i ( w_{i,x}  - \tilde \lambda_i w_i ) \biggl[ 
- \theta'_{i,x} \tilde r_{i,\sigma} - \sum_j
\theta_i' v_j \tr_{i,\sigma u}\tr_j - v_{i,x} \theta_i' \tilde r_{i,\sigma v}
+ \theta_i'\theta_{i,x} \tilde r_{i,\sigma \sigma} \biggr]   
\cr
&\qquad + \sum_i \big( v_i^2 ( \lambda_i^* - \theta_i
) \big)_x \tilde r_{i,v} + \sum_i
v_i^2 ( \lambda_i^* - \theta_i) \biggl[
\sum_j v_j  \tilde r_{i,v u}\tr_j
+ v_{i,x} \tilde r_{i,vv} - \theta_{i,x} \tilde r_{i,v\sigma} \biggr]   \cr
&\qquad + \sum_{i \not= j} \bigl( v_i v_j \bigr)_x \tr_{i,u}\tr_j +
\sum_{i \not= j} v_i v_j \biggl[ \sum_k v_k 
\big(\tr_{i,uu}(\tr_j\otimes\tr_k)
+\tr_{i,u}\tr_{j,u}\tr_k \big) \cr
&\qquad\qquad+ v_{j,x} \tilde
r_{i,u}\tr_{j,v} + v_{i,x} \tilde r_{i,uv} \tr_j - \theta_{j,x}
\tr_{i,u}\tilde r_{j,\sigma} - \theta_{i,x} \tilde
r_{i,u\sigma}\tr_j \biggr],
\cr}\eqno(A.8)
$$
$$\eqalign{
(&u_t )_{xx}- \bigl( A(u) u_t \bigr)_x  \cr
&= \sum_i \big( v_{i,xx} -
(\tilde \lambda_i v_i)_x \big) \Bigl[ w_i \tilde r_{i,v} +
\theta_i' \bigl( w_i/v_i \bigr)^2 \tilde r_{i,\sigma} \Bigr] + \sum_i
\big( w_{i,xx}  - ( \tilde \lambda_i w_i )_x \big) \Bigl[
\tilde r_i - \theta_i' (w_i/v_i) \tilde r_{i,\sigma} \Bigr]   \cr
&\qquad+ \sum_i ( v_{i,x} - \tilde \lambda_i v_i )
\biggl[ w_{i,x} \tilde r_{i,v} + \sum_j w_i v_j\, 
\tilde r_{i,vu}\tr_j + w_i v_{i,x} \,\tilde r_{i,vv} + 
\Big( -w_i \theta_{i,x} +
\theta_i' ( w_i/v_i )^2  v_{i,x} \Big) \tilde r_{i,v\sigma}
\cr
& \qquad \qquad\qquad+ \big( \theta_i' ( w_i/v_i )^2 \big)_x \tilde
r_{i,\sigma}  + \sum_j  v_j \theta_i' ( w_i/v_i)^2
 \tilde r_{i,\sigma u}\tr_j - \theta_i' ( w_i/v_i )^2
\theta_{i,x} \tilde r_{i,\sigma \sigma} \biggr]   \cr
&\qquad + \sum_i ( w_{i,x} - \tilde \lambda_i w_i ) \biggl[ \sum_j
v_j \tilde r_j \tr_{i,u} + v_{i,x}  \tilde r_{i,v} - \big(
\theta_{i,x} + ( \theta_i' w_i/v_i )_x \big) \tilde r_{i,
\sigma} - \sum_j v_j \theta_i' (w_i/v_i) \tilde r_{i,
\sigma u} \tr_j  
\cr
&\qquad \qquad \qquad-  v_{i,x} \theta_i' (w_i/v_i) \tilde r_{i,\sigma v}
+ \theta_{i,x} \theta_i' (w_i/v_i) \tilde
r_{i,\sigma\sigma} \biggr] + \sum_i \big( w_i v_i (
\lambda_i^* - \theta_i) \big)_x \tilde r_{i,v}   \cr
&\qquad + \sum_i
w_i v_i \bigl( \lambda_i^* -\theta_i \bigr) \biggl[ \sum_j v_j 
\tilde r_{i,v u}\tr_j + v_{i,x}
\tilde r_{i,vv} - \theta_{i,x} \tilde r_{i,v\sigma} \biggr] + \sum_{i
\not=j} ( w_i v_j )_x \tr_{i,u}\tr_j   \cr
&\qquad + \sum_{i \not= j} w_i v_j \biggl[ \sum_k v_k 
\big(\tr_{i,uu}(\tr_j\otimes\tr_k)
+\tr_{i,u}\tr_{j,u}\tr_k \big)
+ v_{j,x} 
\tr_{i,u}\tr_{j,v} + v_{i,x} \tilde r_{i,v u}\tr_j -
\theta_{j,x} \tr_{i,u} \tr_{j,\sigma}  - \theta_{i,x}
\tr_{i,u\sigma}\tr_j \biggr]   \cr
&\qquad - \sum_i \lambda_i^* \biggl\{ ( v_{i,x} - \tilde \lambda_i v_i
) \Bigl[ \tilde r_i + v_i \tilde r_{i,v} +  \theta_i'
(w_i/v_i) \tilde r_{i,\sigma} \Bigr] + \big( w_{i,x} - \tilde
\lambda_i w_i \big) \Bigl[ - \theta_i' \tilde r_{i,\sigma} \Bigr]
  \cr
&\qquad \qquad \qquad\qquad +v_i^2 \bigl(
\lambda_i^* - \theta_i \bigr) \tilde r_{i,v} + \sum_{j \not= i}
v_i v_j \tr_{i,u}\tr_j \biggr\}_x.  
\cr}\eqno(A.9)
$$
Substituting the expressions (A.6)--(A.9) inside (A.2) and observing that
$$\big(u_x\bullet A(u)\big)u_t
-\big(u_t\bullet A(u)\big)u_x=\sum_{j\not= i}
(w_i-\lambda_i^*v_i)v_j\Big[ \big(\tr_j\bullet A(u)\big)\tr_i-
\big(\tr_i\bullet A(u)\big)\tr_j\Big]\,,$$
we finally obtain an implicit system of $2n$ scalar equations,
describing the evolution of the components $v_i,w_i$:
$$\eqalign{
\sum_i &\Bigl( v_{i,t} + \bigl( \tilde \lambda_i v_i \bigr)_x -
v_{i,xx} \Bigr) \Bigl[ \tilde r_i + v_i \tilde r_{i,v} + \theta_i'
(w_i/v_i) \tilde r_{i,\sigma} \Bigr] + \sum_i  \Bigl( w_{i,t} + \bigl(
\tilde \lambda_i w_i \bigr)_x - w_{i,xx} \Bigr) \Bigl[ - \theta_i'
\tilde r_{i,\sigma} \Bigr] \cr
&= \sum_i \tr_{i,u}\tr_i \Bigl[ v_i \bigl( v_{i,x} -
\tilde \lambda_i v_i \bigr) - v_i \bigl( w_i - \lambda_i^* v_i
\bigr) \Bigr]   \cr
& \quad+ \sum_{i \not= j}
\tr_{i,u}\tr_j \Bigl[ \bigl( v_{i,x} - \tilde \lambda_i v_i
\bigr) v_j + \bigl( v_i v_j \bigr)_x - v_i \bigl( w_j - \lambda_j^*
v_j \bigr) \Bigr]   \cr
& \quad+ \sum_i \tilde r_{i,v} \Bigl[ 2 v_{i,x} \bigl( v_{i,x} - \tilde
\lambda_i v_i \bigr) + \big( v_i^2 (\lambda_i^*-\theta_i ) \big)_x \Bigr]  
\cr
& \quad
+ \sum_i
\tilde r_{i,\sigma} \Bigl[
( v_{i,x} - \tilde \lambda_i v_i) 
\big(- \theta_{i,x} + ( \theta_i' w_i / v_i )_x
\big) - \bigl(w_{i,x} - \tilde \lambda_i w_i \bigr) \theta'_{i,x}
\Bigr]   \cr
&\quad + \sum_i \tilde r_{i,vu}\tr_i \Bigl[  v_i^2
\bigl( v_{i,x} - \tilde \lambda_i v_i \bigr) + v_i^3
(\lambda_i^* -\theta_i) \Bigr]   \cr
& \quad + \sum_{i \not= j} \tilde r_{i,vu}\tr_j \Bigr[ v_i v_j
\bigl( v_{i,x} - \tilde \lambda_i v_i \bigr) + v_j v_i^2 
(\lambda_i^* -\theta_i)  \Bigr]   \cr
&\quad + \sum_i \tilde r_{i,vv} \Bigl[ v_i
v_{i,x} \bigl( v_{i,x} - \tilde \lambda_i v_i \bigr) + v_{i,x} 
v_i^2 (\lambda_i^* -\theta_i) \Bigr]   \cr
&\quad + \sum_i \tilde r_{i,v\sigma} \Bigl[ \bigl( v_{i,x} - \tilde \lambda_i
v_i \bigr) \Bigl( -v_i \theta_{i,x} + \theta_i' v_{i,x} w_i / v_i \Bigr) -
\bigl( w_{i,x} - \tilde \lambda_i w_i \bigr) v_{i,x} \theta_i' -
v_i^2 ( \lambda_i^* - \theta_i) \theta_{i,x}
\Bigr]   \cr
&\quad + \sum_i\tilde r_{i,\sigma u}\tr_i \Bigl[ \bigl( v_{i,x} -
\tilde \lambda_i v_i \bigr) \theta_i' w_i - \bigl( w_{i,x} - \tilde
\lambda_i w_i \bigr) v_i \theta_i' \Bigr]   \cr
& \quad + \sum_{i \not= j}  \tr_{i,\sigma u}\tr_j \Bigl[
\bigl( v_{i,x} - \tilde \lambda_i v_i \bigr) v_j \theta_i' w_i/v_i  -
\bigl( w_{i,x} - \tilde \lambda_i w_i \bigr) v_j \theta_i' \Bigr]   \cr
& \quad + \sum_i \tr_{i,\sigma \sigma} \Bigl[ -\bigl( v_{i,x} - \tilde
\lambda_i v_i \bigr) \theta_{i,x} \theta_i'w_i/v_i + \bigl( w_{i,x} -
\tilde \lambda_i w_i \bigr) \theta_i' \theta_{i,x} \Bigr]   \cr
& \quad + \sum_{i \not=j}
v_i v_j \biggl[ \sum_k v_k 
\big(\tr_{i,uu}(\tr_j\otimes\tr_k)
+\tr_{i,u}\tr_{j,u}\tr_k \big) 
+ v_{j,x} \tilde
r_{i,u}\tr_{j,v} + v_{i,x} \tilde r_{i,uv} \tr_j - \theta_{j,x}
\tr_{i,u}\tilde r_{j,\sigma} -\theta_{i,x} \tilde
r_{i,u\sigma}\tr_j \biggr]\cr
&\doteq\sum_i a_i(t,x),  
\cr}\eqno(A.10)
$$
\vfill\eject
$$\eqalign{
\sum_i & \Bigl( v_{i,t} + ( \lambda_i v_i)_x - v_{i,xx} \Bigr)
\Bigl[ w_i \tilde r_{i,v} + \theta_i' ( w_i / v_i)^2 \tilde
r_{i,\sigma} \Bigr] + \sum_i \Bigl( w_{i,t} + ( \lambda_i w_i
)_x - w_{i,xx} \Bigr) \Bigl[ \tilde r_i - \theta_i' ( w_i /
v_i ) \tilde r_{i,\sigma} \Bigr]   \cr
&- \sum_i  \lambda_i^* \biggl\{ \Bigl( v_{i,t} + \bigl( \tilde \lambda_i
v_i \bigr)_x -
v_{i,xx} \Bigr) \Bigl[ \tilde r_i + v_i \tilde r_{i,v} + \theta_i'
(w_i/v_i) \tilde r_{i,\sigma} \Bigr]
- \sum_i  \Bigl( w_{i,t} + \bigl(
\tilde \lambda_i w_i \bigr)_x - w_{i,xx} \Bigr) \Bigl[ \theta_i'
\tilde r_{i,\sigma} \Bigr]\biggr\} \cr
=& \sum_i \tilde r_i \tr_{i,u} \Bigl[ \bigl( w_{i,x} - \tilde
\lambda_i w_i \bigr) v_i - w_i \bigl( w_i - \lambda_i^* v_i \bigr)
\Bigr]   \cr
& + \sum_{i \not= j} \tilde
r_j \tr_{i,u} \Bigl[ \bigl( w_{i,x} - \tilde \lambda_i w_i
\bigr) v_j - w_i \bigl( w_j - \lambda_j^* v_j \bigr) + \bigl( w_i v_j
\bigr)_x \Bigr]   \cr
& + \sum_i \tilde r_{i,v} \Bigl[ \bigl( v_{i,x} - \lambda_i v_i
\bigr) w_{i,x} + \bigl( w_{i,x} - \tilde \lambda_i w_i \bigr) v_{i,x} +
\bigl( w_i v_i
(\lambda_i^* -\theta_i) \bigr)_x \Bigr]
\cr
& + \sum_i \tilde r_{i,\sigma} \Bigl[ \bigl( v_{i,x} - \tilde \lambda_i
v_i \bigr) \big( \theta_i' (w_i/v_i )^2 \big)_x - (
w_{i,x} - \tilde \lambda_i w_i ) \big( \theta_{i,x} + (
\theta_i' w_i / v_i )_x \big) \Bigr]   \cr
& + \sum_i \tilde r_{i,vu}\tr_i \Bigl[ \bigl( v_{i,x} -
\tilde \lambda_i v_i \bigr) w_i v_i +
w_i v_i^2 (\lambda_i^* -\theta_i)
\Bigr]   \cr
& + \sum_i \tilde r_{i,vv} \Bigl[ \bigl( v_{i,x} - \tilde \lambda_i
v_i \bigr) w_i v_{i,x} + w_i v_i v_{i,x} (\lambda_i^* -\theta_i) 
\Bigr]   \cr
& + \sum_{i \not= j} \tilde r_{i,v u}\tr_j \Bigl[
\bigl( v_{i,x} - \lambda_i v_i
\bigr) w_i v_j + w_i v_i v_j(\lambda_i^* -\theta_i) \Bigr]   
\cr
& + \sum_i \tilde r_{i,v\sigma} \Bigl[ ( v_{i,x} - \tilde \lambda_i
v_i ) \big( -w_i \theta_{i,x} + \theta_i' (w_i/v_i )^2
v_{i,x} \big) - ( w_{i,x} - \tilde \lambda_i w_i ) \theta_i'
v_{i,x} w_i/v_i - w_i v_i (\lambda_i^* -\theta_i) 
\theta_{i,x} \Bigr]   \cr
& + \sum_i \tilde r_{i,\sigma u} \tr_i\Bigl[ \bigl( v_{i,x}
- \tilde \lambda_i v_i \bigr) \theta_i' w_i^2 / v_i  -
\bigl( w_{i,x} - \tilde \lambda_i w_i \bigr) \theta_i' w_i \Bigr]  
\cr
& + \sum_{i \not= j} \tilde r_{i,\sigma u}
\tr_j \Bigl[
( v_{i,x} - \tilde \lambda_i v_i ) v_j \theta_i' ( w_i /
v_i )^2 - ( w_{i,x} - \tilde \lambda_i w_i ) v_j
\theta_i' w_i/v_i \Bigr]   \cr
& + \sum_i \tilde r_{i,\sigma \sigma} \Bigl[ -\bigl( v_{i,x} - \tilde
\lambda_i v_i \bigr) \theta_i' (w_i/v_i)^2 \theta_{i,x} +
\bigl( w_{i,x} - \tilde \lambda_i w_i \bigr) \theta' \theta_{i,x}
w_i/v_i \Bigr]   \cr
&
+ \sum_{i \not=j}
w_i v_j \biggl[ \sum_k v_k 
\big(\tr_{i,uu}(\tr_j\otimes\tr_k)
+\tr_{i,u}\tr_{j,u}\tr_k \big) 
+ v_{j,x} \tilde
r_{i,u}\tr_{j,v} + v_{i,x} \tilde r_{i,uv} \tr_j - \theta_{j,x}
\tr_{i,u}\tilde r_{j,\sigma} - \theta_{i,x} \tilde
r_{i,u\sigma}\tr_j \biggr]\cr
& + \sum_{i \not= j}(w_i-\lambda_i^*v_i) v_j \Big[\big( \tilde r_j
\bullet A(u)\big) \tilde r_i - \big(\tilde r_i
\bullet A(u)\big) \tilde r_j  \Big]  \cr
&
- \sum_i \lambda_i^*
a_i(t,x)   \cr
\doteq&~ \sum_i b_i(t,x) - \sum_i \lambda_i^* a_i(t,x).  
\cr}\eqno(A.11)
$$
Recalling the expression (5.10) for the differential $\partial\Lambda/
\partial(v,w)$, we recognize that
the equations (A.10)-(A.11) provide the explicit form of
the system (6.5).
The uniform invertibility of the differential of $\Lambda$
implies the estimates
$$\phi_j,~\psi_j=\O(1)\cdot\sum_i\big(|a_i|+|b_i|\big).$$
To prove Lemma 6.1, it thus suffices to show that
all the terms in the summations defining $a_i,b_i$ have the correct order of
magnitude.
\v
First of all, one checks that 
all those terms which involve a product of distinct
components $i\not= j$
can be bounded as
$$\O(1)\cdot \sum_{j\not= k}
\big(|v_jv_k|+|v_{j,x}v_k|+
|v_jw_k|+|v_{j,x}w_k|+|v_jw_{k,x}|+|w_jw_k|\big).
\eqno(A.12)$$
In most cases, this estimate is straightforward.
For the terms containing the factor $\theta_{i,x}$ or $\theta_{j,x}$
this is proved as follows.
Recalling the bounds (4.24)
we have, for example,
$$\eqalign{\theta_{j,x}\tr_{j,\sigma}
&=\O(1)\cdot v_j\,\theta'_j{w_{j,x}v_j-w_jv_{j,x}\over v_j^2}
=\O(1)\cdot\big(|w_{j,x}|+|v_{j,x}|\big)=\O(1)\cdot \delta_0^3\,,}$$
because of (5.24).
Hence
$$v_iv_j\theta_{j,x}\tr_{i,u}\tr_{j,\sigma}=\O(1)\cdot\delta_0^3|v_iv_j|\,.$$

Next, we look at each one of the remaining terms on the right hand side
of (A.10) and (A.11) and show that its size can
be bounded as claimed by Lemma 6.1. 
To appreciate the following computations, one should keep in mind that:
\v
\n{\bf 1.} By (6.16) there holds
$$v_{i,x}-(\tla_i-\lambda_i^*)v_i-w_i=\O(1)\cdot\delta_0\sum_{j\not= i}
|v_j|\,.$$
Therefore
$$\big(|v_i|+|w_i|+|v_{i,x}|+|w_{i,x}|\big)
\Big|v_{i,x}-(\tla_i-\lambda_i^*)v_i-w_i\Big|
=\O(1)\cdot\delta_0\sum_{j\not= i}
\big(|v_iv_j|+|w_iv_j|+|v_{i,x}v_j|+|w_{i,x}v_j|\big)\,.$$
\v
\n{\bf 2.} By (5.5) the cutoff functions satisfy
$\theta_i'=\theta_i''=0$ whenever $|w_i/v_i|\geq 3\delta_1$.
\v
\n{\bf 3.} By (4.24) we have
$\tr_{i,\sigma}/v_i\,,~\tr_{i,\sigma\sigma}/v_i\,,~
\tr_{i,\sigma u}/v_i=\O(1)$.
\v
\n{\bf 4.}
One can have $|w_i-\theta_iv_i|\not= 0$ only when $|w_i|>\delta_1|v_i|$.
In this case, (6.18) yields
$$v_i=\O(1)\cdot v_{i,x}+\O(1)\cdot\delta_0\sum_{j\not= i} |v_j|\,.$$
\v
What follows is a list of the various
terms, first those appearing in $a_i$, then the ones in $b_i$.
\v
\n Coefficients of $\tr_{i,u}\tr_i$ :
$$
v_i ( v_{i,x} - \tilde \lambda_i v_i) - v_i
( w_i - \lambda_i^* v_i ) =  v_i \bigr[
v_{i,x} - ( \tilde \lambda_i -
\lambda_i^* ) v_i - w_i \bigr],
$$
$$
v_i ( w_{i,x} - \tilde \lambda_i w_i) 
- w_i( w_i - \lambda_i^* v_i ) = \big[ v_i w_{i,x}
- v_{i,x} w_i \big] +  w_i
\big[ v_{i,x} - ( \tilde \lambda_i - \lambda_i^* ) v_i - w_i
\big].
$$
\v
\n Coefficients of $\tilde r_{i,v}$ :
$$\eqalign{
2 v_{i,x} ( v_{i,x} - \tilde \lambda_i v_i)& + \bigl( v_i^2 ( \lambda_i^*-
\theta_i ) \bigr)_x = 2
v_{i,x} \bigl[ v_{i,x} - ( \tilde \lambda_i - \lambda_i^*)
v_i -\theta_i v_i \bigr] + \theta_i' \big[ v_{i,x} w_i - v_i
w_{i,x} \big] \cr
&= 2 v_{i,x} \Bigl[ v_{i,x} - ( \tilde \lambda_i -
\lambda_i^* ) v_i - w_i \Bigr] + 2 v_{i,x} \big[ w_i -
\theta_i v_i \big] + \theta_i' \bigl[ v_{i,x} w_i - v_i w_{i,x}
\bigr],
\cr}
$$
$$\eqalign{
w_{i,x} & ( v_{i,x} - \tilde \lambda_i v_i ) + v_{i,x} ( w_{i,x}
- \tilde \lambda_i w_i ) + \bigl(
w_i v_i( \lambda_i^*-\theta_i) \bigr)_x \cr
&= 2 w_{i,x} \big[ v_{i,x} - ( \tilde
\lambda_i - \lambda_i^*)  v_i -w_i \big] + 2 w_{i,x}
\bigl[w_i - \theta_i v_i \bigr] + \bigl(
\lambda_i^*-\theta_i
-\tla_i+ \theta'_i w_i/v_i \bigr)
\bigl[ v_{i,x} w_i - v_i w_{i,x} \bigr].
\cr}
$$
\v
\n Coefficients of $\tilde r_{i,\sigma}/v_i$ :
$$
\eqalign{
v_i \bigl( v_{i,x} - \tilde \lambda_i v_i \bigr) \Bigl( -\theta_{i,x} + ( 
\theta'_i w_i / v_i )_x \Bigr) - v_i \bigl( w_{i,x} - 
\tilde \lambda_i w_i\bigr) \theta'_{i,x}  &= -\bigr( v_i w_{i,x} - w_i v_{i,x}
\bigr) \theta''_i \big( {w_i/ v_i} \big)_x \cr
&= -\theta_i'' \Big[v_i\,\big( {w_i/ v_i} \big)_x\Big]^2,
\cr}
$$
$$
v_i \bigl( v_{i,x} - \tilde \lambda_i v_i \bigr) \Bigl(
\theta'_i ( w_i / v_i )^2 \Bigr)_x - v_i \bigl( w_{i,x} - \tilde 
\lambda_i w_i \bigr) \Bigl( \theta_{i,x} + ( \theta'_i w_i / v_i
)_x \Bigl) = -\Big( \theta''_i (w_i/ v_i) + 2 \theta'_i
\Big) \Big[v_i\,\big( {w_i/ v_i} \big)_x\Big]^2.
$$
\v
\n Coefficients of $\tr_{i,vu}\tr_i$ :
$$
v_i^2 \bigl( v_{i,x} - \tilde \lambda_i v_i \bigr) +
v_i^3 (\lambda_i^* -\theta_i) = v_i^2 \Bigl[ v_{i,x} - (
\tilde \lambda_i - \lambda_i^* ) v_i - w_i \Bigr] + v_i^2 
\big[ w_i - \theta_i v_i \big],
$$
$$
v_i w_i \bigl( v_{i,x} - \tilde \lambda_i v_i \bigr) +
v_i^2 w_i (\lambda_i^*-\theta_i) = v_i w_i \Bigl[
v_{i,x} - ( \tilde
\lambda_i - \lambda_i^* ) v_i - w_i \Bigr] + 
v_i w_i \big[ w_i -\theta_i v_i \big].
$$
\v
\n Coefficients of $\tilde r_{i,vv}$ :
$$
v_i v_{i,x} ( v_{i,x} - \tilde \lambda_i v_i ) - v_{i,x} 
v_i^2 ( \lambda_i^*-\theta_i) = v_i v_{i,x} \bigl[
v_{i,x} -( \tilde \lambda_i - \lambda_i^* ) v_i
- w_i \bigr] + v_i v_{i,x} \bigl[ w_i -\theta_i v_i \big],
$$
$$
w_i v_{i,x} \bigl( v_{i,x} - \tilde \lambda_i v_i \bigr) - v_{i,x}  v_i 
w_i ( \lambda_i^*-\theta_i ) = w_i v_{i,x} \Bigl[ v_{i,x} -
\bigl( \tilde \lambda_i -
\lambda_i^* \bigr) v_i - w_i \Bigr] + w_i v_{i,x} \big[ w_i -
\theta_i v_i \big].
$$
\v
\n Coefficients of $\tilde r_{i,v\sigma}$ :
$$\eqalign{
(&v_{i,x} - \tilde \lambda_i v_i) \big(- v_i \theta_{i,x} +
\theta'_i v_{i,x} w_i/v_i  \big) - ( w_{i,x} - \tilde \lambda_i w_i
) \theta'_i v_{i,x} - v_i^2 ( \lambda_i^* -
\theta_i) \theta_{i,x} \cr
&= v_{i,x} \theta'_i \bigl(v_{i,x} w_i - v_i w_{i,x} \bigr) / v_i -
\theta_{i,x} v_i \Bigl( v_{i,x} - ( \tilde \lambda_i -
\lambda_i^*  + \theta_i) v_i \Bigr) \cr
&= 2 \theta_i' \left( v_{i,x} {w_i\over v_i}  - w_{i,x} \right)
\Big\{ \big[ v_{i,x} - ( \tilde \lambda_i - \lambda_i^* )
v_i - w_i \big] + [ w_i - 
\theta_i v_i] \Big\}
+\bigl(\tilde \lambda_i 
- \lambda_i^* +\theta_i\bigr) \theta'_i 
\bigl[ w_{i,x} v_i - w_i v_{i,x} \bigr],
\cr}
$$
$$
\eqalign{
\bigl(v_{i,x} - \tilde \lambda_i v_i \bigr)& \bigl(- w_i \theta_{i,x} +
\theta'_i
( w_i /v_i  )^2 v_{i,x} \bigr) - \bigl( w_{i,x} - \tilde
\lambda_i w_i
\bigr) \theta'_i v_{i,x} w_i / v_i - v_i w_i( \lambda_i^* -
\theta_i ) \theta_{i,x} \cr
&= 2  \theta'_i {w_i\over v_i} \left( v_{i,x} {w_i\over v_i} -
w_{i,x}
\right) \Big\{ \big[ v_{i,x} - ( \tilde \lambda_i -
\lambda_i^*) v_i - w_i \big] + [ w_i - \theta_i v_i ]
\Big\} \cr
&\qquad \qquad - \bigl( \tilde \lambda_i 
-\lambda_i^* +\theta_i\bigr) \theta'_i 
\bigl[ w_{i,x} v_i - w_i v_{i,x} \bigr]
w_i / v_i\,.\cr}
$$
\v
\n Coefficients of $\tilde r_{i,\sigma u}\tr_i/v_i$ :
$$
( v_{i,x} - \tilde \lambda_i v_i ) w_i v_i \theta'_i - (
w_{i,x} - \tilde \lambda_i w_i ) v_i^2 \theta'_i =
\theta'_i v_i \bigl[ v_{i,x} w_i - v_i w_{i,x} \bigr],
$$
$$
( v_{i,x} - \tilde \lambda_i v_i ) w_i^2 \theta'_i -
( w_{i,x} - \tilde
\lambda_i w_i ) w_i v_i \theta'_i = \theta'_i w_i \bigl[ v_{i,x}
w_i - v_i w_{i,x} \bigr].
$$
\v
\n Coefficients of $\tilde r_{i,\sigma\sigma}/v_i$ :
$$
-(v_{i,x} - \tilde \lambda_i v_i ) w_i \theta'_i \theta_{i,x} +
( w_{i,x} - \tilde \lambda_i w_i ) v_i \theta'_i \theta_{i,x}
= -(\theta'_i)^2  \Big[v_i\,\big( {w_i/ v_i} \big)_x\Big]^2,
$$
$$
-(v_{i,x} - \tilde \lambda_i v_i ) w_i \bigl( w_i / v_i \bigr)
\theta'_i \theta_{i,x} + \bigl( w_{i,x} - \tilde \lambda_i w_i \bigr)
w_i \theta'_i \theta_{i,x} = -( \theta'_i)^2 {w_i\over v_i} 
\Big[v_i\,\big( {w_i/ v_i} \big)_x\Big]^2.
$$
This completes our analysis, showing that all terms
in the summations that define $a_i,b_i$ have the correct order of
magnitude, as claimed by Lemma 6.1.
\endproof
\vsk
\c{\medbf Appendix B}
\v
We compute here the source terms $\hat\phi_i,\hat\psi_i$ in 
the equations (11.15) for the components of a first order perturbation,
and prove Lemma 11.4.  We recall that
$$
z = \sum_i \h_i \tilde r_i \bigl(u,v_i,\lambda_i^* -
\theta(\iot_i/\h_i) \bigr),
\qquad \io = \sum_i \bigl( \iot_i  - \lambda_i^* \h_i \bigr) \tilde r_i
\bigl( u, v_i, \lambda_i^* - \theta( \iot_i/\h_i)\bigr),\eqno(B.1)$$
$$
\hth_i \doteq  \theta \left(  {\iot_i\over \h_i} \right),
\qquad \hat r_i \doteq \tilde r_i (u,v_i,\,\lambda_i^*-\hth_i ),
\qquad\hla_i\doteq\la\hr_i,\,A(u)\hr_i\ra.
$$
As in (A.4)--(A.11), the computations are lengthy but straightforward:
one has to
rewrite the evolution equations
for $z$ and $\io$ :
$$\left\{ \eqalign{ 
z_t+\big(A(u)z\big)_x-z_{xx}&=
\big(u_x\bullet A(u)\big) z-\big(z\bullet A(u)\big) u_x\,,\cr
\io_t + \bigl( A(u) \io \bigr)_x - \io_{xx} &= \Big[ \big(
u_x \bullet A(u)
\big) z - \big( z \bullet A(u) \big) u_x \Big]_x 
- A(u) \Big[ \big( u_x \bullet A(u)
\big) z - \big( z \bullet A(u) \big) u_x \Big]
\cr &\qquad
+ \big( u_x \bullet A(u) \big) \io - \big( u_t \bullet A(u)
\big) z\,,\cr}\right.\eqno(B.2)$$
in terms of $\h_i,\iot_i$.
The fundamental relation (4.23) implies
$$A(u)\hr_i=\hla_i\hr_i+v_i\big(\hr_{i,u}\hr_i
+(\hla_i-\lambda_i^*+\hth_i)\hr_{i,v}\big)\,.\eqno(B.3)$$
Differentiating (B.1) w.r.t.~$x$ and
using (B.3) we obtain
$$\eqalign{
z_x - A(u) z &= \sum_i \h_{i,x} \hat r_i + \sum_i \h_i \Bigl[ v_i 
\hat r_{i,u}\tr_i - A(u) \hat r_i \Bigr] 
+ \sum_i \h_i \Bigl[ v_{i,x} \hr_{i,v} -
\big(\hth_i' ( \h_i \iot_{i,x} - \iot_i \h_{i,x} ) /
\h_i^2 \big)\hat r_{i,\sigma}  \Bigr] \cr
&~ + \sum_{i \not= j} \h_i v_j \hr_{i,u}\tr_j \cr
&= \sum_i \bigl( \h_{i,x} - \hat \lambda_i \h_i \bigr) \hat r_i + \sum_i
\h_i v_i r_{i,u}\bigl( \tilde r_i - \hat r_i \bigr) +
\sum_i \h_i \Bigl[ v_{i,x} - ( \hla_i-\lambda_i^*
+ \hth_i ) v_i \Bigr]
\hat r_{i,v} \cr
&~ - \sum_i \hth'_i \left[ \bigl( \iot_{i,x} - \hat \lambda_i
\iot_i \bigr) -  {\iot_i\over \h_i} \bigl( \h_{i,x} - \hat \lambda \h_i \bigr)
\right] \hat r_{i,\sigma} + \sum_{i \not= j} \h_i v_j\hat r_{i,u}\tr_j, \cr
&= \sum_i \bigl( \h_{i,x} - \hat \lambda_i \h_i \bigr) \Bigl[ \hat
r_i +\hth_i' (\iot_i / \h_i ) \hat
r_{i,\sigma} \Bigr] - \sum_i \bigl( \iot_{i,x} - \hat \lambda_i \iot_i
\bigr)\hth'_i \hat r_{i,\sigma} \cr
& +  \sum_i \h_i v_i  \tr_{i,u} \bigl( \tilde r_i - \hat r_i
\bigr)+
\sum_i \h_i \Bigl[  v_{i,x} - ( \hla_i-\lambda_i^*
+ \hth_i ) v_i 
\Bigr] \hat r_{i,v} 
+ \sum_{i \not= j} \h_i v_j r_{i,u}\tr_j\,,
\cr}
\eqno(B.4)
$$
$$
\eqalign{
\io_x - A(u) \io &= \sum_i \bigl( \iot_{i,x} - \lambda_i^* \h_{i,x}
\bigr) \hat r_i +
\sum_i \bigl( \iot_i - \lambda_i^* \h_i \bigr) \Bigl[ v_i  
\hr_{i,u}\tr_i - A(u) \hat r_i \Bigr] \cr
&~ + \sum_i \bigl( \iot_i - \lambda_i^* \h_i \bigr) \Bigl[ v_{i,x}
\hat r_{i,v} - \big(\hth_i' ( \h_i \iot_{i,x} - \iot_i \h_{i,x}
)/ \h_i^2\big) \hat r_{i,\sigma}  \Bigr]
+ \sum_{i \not= j} \bigl( \iot_i - \lambda_i^* \h_i \bigr) v_j
\hr_{i,u}\tr_j \cr
&= \sum_i \bigl( \iot_{i,x} - \hat \lambda_i \iot_i \bigr) \hat r_i -
\sum_i \lambda_i^* \bigl( \h_{i,x} - \hat \lambda_i \h_i \bigr)
\hr_i + \sum_i \bigl( \iot_i - \lambda_i^* \h_i \bigr) v_i
\hr_{i,u}( \tilde r_i - \hat r_i) \cr
&~ + \sum_i \bigl( \iot_i - \lambda_i^* \h_i \bigr) \Bigl[
v_{i,x} - (\hat \lambda_i - \lambda_i^* + \hth_i )
v_i \Bigr] \hat r_{i,v} \cr
&~ - \sum_i \hth'_i \left(  {\iot_i\over \h_i} - \lambda_i^* \right)
\left[ \bigl( \iot_{i,x} - \hat \lambda_i \iot_i \bigr) -  {\iot_i\over \h_i}
\bigl( \h_{i,x} - \hat \lambda_i \h_i \bigr) \right] \hat r_{i,\sigma}
+ \sum_{i \not= j} \bigl( \iot_i - \lambda_i^* \h_i \bigr) v_j 
\hr_{i,u}\tr_j \cr
&= \sum_i \big( \h_{i,x} - \hat \lambda_i \h_i \big) \Bigl[ 
\hth'_i \bigl( \iot_i / \h_i \bigr)^2 \hat r_{i,\sigma} \Bigr]
+ \sum_i \big( \iot_{i,x} - \hat \lambda_i \iot_i \big) \Bigl[ \hat
r_i - \hth'_i (\iot_i/\h_i ) \hat r_{i,\sigma} \Bigr] \cr
&~ + \sum_i \iot_i \Bigl[ v_{i,x} - ( \hat \lambda - 
\lambda_i^* + \hth_i) v_i \Bigr] \hat r_{i,v} + \sum_i
\iot_i v_i \hr_{i,u}\bigl( \tilde r_i - \hat r_i \bigr) +
\sum_{i \not= j}
\iot_i v_j\hr_{i,u}\tr_j \cr
&~ - \sum_i \lambda_i^* \biggl\{ \bigl( \h_{i,x} - \hat \lambda_i \h_i
\bigr) \Bigl[ \hat r_i + \hth_i'( \iot_i /
\h_i ) \hat r_{i,\sigma} \Bigr] -\bigl( \iot_{i,x} - \hat
\lambda_i \iot_i \bigr) \hth'_i \hat r_{i,\sigma}  \cr
&~ \quad \qquad \qquad + \h_i \Bigl[ v_{i,x} - ( \hat \lambda_i -
\lambda_i^* + \hth_i) v_i
\Bigr] \hat r_{i,v} + \h_i v_i  \hr_{i,u}\bigl( \tilde r_i - \hat r_i
\bigr) + \sum_{j \not= i} \h_i v_j
\hat r_{i,u}\tr_j \biggr\}.
\cr}
\eqno(B.5)
$$
Next, differentiating (B.1) w.r.t.~$t$ we obtain
$$\eqalign{
z_t &= \sum_i \h_{i,t} \hat r_i + \sum_i \h_i \bigl( \hr_{i,u} u_t+
v_{i,t} \hat r_{i,v} - \hth_{i,t} \hat r_{i,\sigma} \bigr) \cr
&= \sum_i \h_{i,t} \hat r_i + \sum_i \h_i \biggl[  v_{i,t} \hat r_{i,v} - 
\big(\hth'_i ( \iot_{i,t} \h_i - \iot_i \h_{i,t} )/\h_i
^2\big) \hat
r_{i,\sigma} \biggl] + \sum_{i,j} \h_i \bigl( w_j - \lambda_j^*v_j
\bigr) \hr_{i,u}\tr_j   \cr
&= \sum_i \h_{i,t} \Bigl[ \hat r_i + \hth'_i (\iot_i/\h_i
) \hat r_{i,\sigma} \Bigr] -
\sum_i \hth'_i\iot_{i,t}  \hat r_{i,\sigma} +
\sum_i \h_i v_{i,t} \hat r_{i,v} + 
\sum_{i,j} \h_i \bigl( w_j - \lambda_j^*v_j \bigr)
\hr_{i,u}\tr_j\,,  
\cr}\eqno(B.6)$$
$$\eqalign{
\io_t &= \sum_i \bigl( \iot_{i,t} - \lambda_i^* \h_{i,t} \bigr) \hat
r_i + \sum_i \bigl( \iot_i - \lambda_i^* \h_i \bigr) \Bigl(
\hat r_{i,u}u_t +
v_{i,t} \hat r_{i,v} - \hth_{i,t} \hat r_{i,\sigma} \Bigr) \cr
&= \sum_i \bigl( \iot_{i,t} - \lambda_i^* \h_{i,t} \bigr) \hat r_i +
\sum_i \bigl( \iot_i - \lambda_i^* \h_i \bigr) \biggl[  v_{i,t} \hat
r_{i,v} -\big(
\hth'_i ( \iot_{i,t} \h_i - \iot_i \h_{i,t} )/\h_i
^2 \big)\hat r_{i,\sigma} \biggl]   \cr
&\qquad + \sum_{i,j} \bigl( \iot_i - \lambda_i^* \h_i
\bigr) \bigl( w_j - \lambda_j^*v_j \bigr)  \hr_{i,u}\tr_j
\cr
&= \sum_i \h_{i,t} \Bigl[ \hth'_i ( \iot_i/\h_i )^2
\hat r_{i,\sigma} \Bigr]
+ \sum_i \iot_{i,t} \Bigl[ \hat r_i - \hth'_i ( \iot_i/\h_i )
\hat r_{i,\sigma} \Bigr] + \sum_i \iot_i v_{i,t} \hat r_{i,v} +
\sum_{i,j} \iot_i ( w_j
- \lambda_j^*v_j ) \hr_{i,u}\tr_j   \cr
&\qquad- \sum_i \lambda_i^* \biggl\{ \h_{i,t} \bigl[ \hat r_i + 
\hth_i' (\iot_i/\h_i) \hat r_{i,\sigma} \bigr]
- \hth_i'\iot_{i,t}  \hat r_{i,\sigma}  +
\h_i v_{i,t} \hat r_{i,v} + \sum_j \h_i
\bigl( w_j - \lambda_j^*v_j \bigr) \hr_{i,u}\tr_j
\biggr\}.  
\cr}\eqno(B.7)$$
Differentiating again $z_x-A(u)z$ and $\io_x-A(u)\io$ w.r.t.~$x$,
from (B.4) and (B.5) one finds
\vfill\eject
$$\eqalignno{
z_{xx} - \bigl( A(u) z \bigr)_x = &~ \sum_i \bigl( \h_{i,xx} -
( \hat \lambda_i \h_i)_x \bigr) \Bigl[ \hat r_i + 
\hth_i' (\iot_i/\h_i) \hat r_{i,\sigma} \Bigr] - \sum_i \bigl(
\iot_{i,xx}  - ( \hat \lambda_i \iot_i )_x \bigr) 
\hth_i' \hat r_{i,\sigma} \cr
&~ + \sum_i \bigl( \h_{i,x} - \hat \lambda_i \h_i \bigr) \biggl[ \sum_j
v_j \hr_{i,u}\tr_j + v_{i,x} \hat r_{i,v} + \bigl(
-\hth_{i,x} + (\hth_i' \iot_i/\h_i )_x \bigr) \hat
r_{i,\sigma}   \cr
&\quad \qquad + \hth_i' v_{i,x}
(\iot_i /\h_i) \hat r_{i,\sigma v}  + \sum_j  v_j \hth_i'
(\iot_i /\h_i ) \hat r_{i,u\sigma}\tr_j-
\hth_{i,x} \hth_i' (\iot_i /\h_i) \hat r_{i,\sigma
\sigma} \biggr]   \cr
&~ + \sum_i \bigl( \iot_{i,x}  - \hat \lambda_i \iot_i \bigr) \biggl[
-\hth_i'' (\iot_i/\h_i)_x \hat r_{i,\sigma} - \sum_j
\hth_i' v_j \hat r_{i,u\sigma}\tr_j - v_{i,x} \hth_i' \hat
r_{i,v\sigma} + \hth_i'
\hth_{i,x} \hat r_{i,\sigma \sigma} \biggr]   \cr
&~ + \sum_i \bigl( \h_i v_i \hr_{i,u}
( \tilde r_i - \hat r_i) 
\bigr)_x + \sum_i \Bigl( \h_i \big( v_{i,x} -( \hat \lambda_i -
\lambda_i^* + \hth_i) v_i
\bigr) \Bigr)_x \hat r_{i,v}   \cr
&~ + \sum_i \h_i \big( v_{i,x} - 
( \hat \lambda_i - \lambda_i^* + \hth_i) v_i \big) \biggl[
\sum_j v_j\hat r_{i,vu}\tr_j
+ v_{i,x} \hat r_{i,vv} - \hth_{i,x} \hat r_{i,v\sigma} \biggr]   \cr
&~ + \sum_{i \not= j} \bigl( \h_i v_j \bigr)_x \hr_{i,u}\tr_j +
\sum_{i \not= j} \h_i v_j \biggl[ \sum_k v_k  \bigl( 
\hr_{i,u}\tr_{j,u}\tr_k+\tr_{i,uu}(\tr_j\otimes\tr_k) \bigr) \cr
&\qquad \quad+ v_{j,x} \hat
r_{i,u}\tr_{j,v} + v_{i,x} \hr_{i,vu} \tr_j  
- \theta_{j,x}
\tr_{i,u}\tr_{j,\sigma} - \hth_{i,x}
\hr_{i,\sigma u}\tr_j \biggr], &(B.8)
\cr}
$$
$$
\eqalignno{
\io_{xx}& - \bigl( A(u) \io \bigr)_x ~=~ \sum_i \bigl( \h_{i,xx} -
( \hat \lambda_i \h_i )_x \bigr) \Bigl[ \hth_i'
( \iot_i/\h_i )^2 \hat r_{i,\sigma} \Bigr] + \sum_i
\bigl( \iot_{i,xx}  - ( \hat \lambda_i \iot_i )_x \bigr) \Bigl[
\hat r_i - \hth_i' (\iot_i/\h_i) \hat r_{i,\sigma} \Bigr]   \cr
&~ + \sum_i \bigl( \h_{i,x} - \hat \lambda_i \h_i \bigr)
\biggl[ \hth_i' ( w_i/v_i )^2  v_{i,x} \hat r_{i,v\sigma}
+ \bigl( \hth_i' ( \iot_i/\h_i )^2 \bigr)_x \hat
r_{i,\sigma}   \cr
&~ \qquad + \sum_j v_j \hth_i' \bigl( \iot_i/\h_i \bigr)^2
\hat r_{i,\sigma u}\hr_j - \hth_i' ( \iot_i/\h_i )^2
\hth_{i,x} \hat r_{i,\sigma \sigma} \biggr]   \cr
&~ + \sum_i \bigl( \iot_{i,x} - \tilde \lambda_i \iot_i \bigr) \biggl[ \sum_j
v_j \hr_{i,u}\tr_j + v_{i,x}  \hat r_{i,v} - \bigl(
\hth_{i,x} + ( \hth_i' \iot_i/\h_i )_x \bigr)
\hat r_{i,\sigma}   \cr
&~ \qquad - \sum_j v_j \hth_i' (\iot_i/\h_i) \hat r_{i,\sigma u}\tr_j
- v_{i,x} \hth_i' (\iot_i/\h_i ) \hat
r_{i,v\sigma} + \hth_{i,x} \hth_i' (\iot_i/\h_i)
\hat r_{i, \sigma \sigma} \biggr]   \cr
&~ + \sum_i \Bigl[ \iot_i \bigl( v_{i,x} - ( \hat \lambda_i -
\lambda_i^* +
\hth_i) v_i \bigr) \Bigr]_x \hat r_{i,v} + \sum_i \Bigl( \iot_i v_i
\hr_{i,u} ( \tilde r_i - \hat r_i ) \Bigr)_x   \cr
&~ + \sum_i
\iot_i \big( v_{i,x} - ( \hat \lambda_i - \lambda_i^* + \hth_i) v_i 
\big) \biggl[ \sum_j v_j  \hat r_{i,v}\tr_j + v_{i,x}
\hat r_{i,vv} - \hth_{i,x} \hat r_{i,v\sigma} \biggr] + \sum_{i
\not=j} ( \iot_i v_j )_x \hr_{i,u}\tr_j   \cr
&~ + \sum_{i \not= j} \iot_i v_j \biggl[ \sum_k v_k 
\bigl(\hr_{i,u}\tr_{j,u}\tr_k+\hr_{i,uu}(\tr_j\otimes\tr_k) \bigr) 
+ v_{j,x} \hr_{i,u} \tilde r_{j,v}
+ v_{i,x} \hat r_{i,vu}\tr_j -
\theta_{j,x} \hr_{i,u}\tilde r_{j,\sigma} - \hth_{i,x}
\hat r_{i,\sigma u}\tr_j \biggr]   \cr
&~ - \sum_i \lambda_i^* \biggl\{ \bigl( \h_{i,x} - \hat \lambda_i \h_i
\bigr) \Bigl[ \hat r_i + \hth_i'
(\iot_i/\h_i ) \hat r_{i,\sigma} \Bigr] -\bigl( \iot_{i,x} - \hat
\lambda_i \iot_i \bigr) \hth_i' \hat r_{i,\sigma} 
\cr
&\qquad \qquad + \h_i \bigl( v_{i,x} - ( \hat \lambda_i - 
\lambda_i^* + \hth_i) v_i \bigr) \hat r_{i,v} + \h_i v_i
\hr_{i,u}\bigl( \tilde r_i - \hat r_i \bigr) + \sum_{j \not= i}
v_i v_j \hr_{i,u}\tr_j \biggr\}_x.  
&(B.9)\cr}
$$
\vfill\eject
Substituting the expressions (B.6)--(B.9) inside (B.2)
we obtain an implicit system of $2n$ scalar equations
governing the evolution of the components $\h_i,\iot_i$:
$$\eqalign{
\sum_i & \Bigl( \h_{i,t} + ( \hat \lambda_i \h_i )_x -
\h_{i,xx} \Bigr) \Bigl[ \hat r_i + \hth_i'
(\iot_i/\h_i) \hat r_{i,\sigma} \Bigr] + \sum_i  \Bigl( \iot_{i,t} + \bigl(
\hat \lambda_i \iot_i \bigr)_x - \iot_{i,xx} \Bigr)\Big[-\hth_i'
\hat r_{i,\sigma} \Big]  
\cr
=~& \sum_i \hr_{i,u}\tr_i \Bigl[ v_i( \h_{i,x} -
\hat \lambda_i \h_i ) - \h_i ( w_i - \lambda_i^* )
\Bigr] 
\cr
&+ \sum_{i \not= j} 
\hr_{i,u}\tr_j \Bigl[ ( \h_{i,x} - \hat \lambda_i \h_i
) v_j + \big( \h_i v_j \big)_x - \h_i ( w_j - \lambda_j^*
v_j ) \Bigr]   
\cr
& + \sum_i \hat r_{i,v} \bigg[ ( \h_{i,x} - \hat \lambda_i \h_i
) v_{i,x} + \Bigl( \h_i \bigl( v_{i,x} - ( \hat
\lambda_i - \lambda_i^*  + \hth_i )v_i \bigr) \Bigr)_x -
\h_i v_{i,t} \bigg]   
\cr
& + \sum_i \hat r_{i,\sigma} \Bigl[
( \h_{i,x} - \hat \lambda_i \h_i
) \big( -\hth_{i,x} + ( \hth_i' \iot_i / \h_i )_x
\big) - ( \iot_{i,x} - \hat \lambda_i \iot_i )
\hth'_{i,x} \Bigr]   
\cr
& + \sum_i \hat r_{i,vu} \tr_i\Bigl[ \h_i v_i
\bigl( v_{i,x} - ( \hat \lambda_i - \lambda_i^*  +
\hth_i) v_i \bigr) \Bigr]   
\cr
& + \sum_{i \not= j}  \hat r_{i,vu}\tr_j \Bigr[ \h_i v_j
\bigl( v_{i,x} - ( \hat \lambda_i - \lambda_i^*  +
\hth_i) v_i \bigr)\Bigr]
\cr
&+\sum_i \hat r_{i,vv} \Bigl[ \h_i
v_{i,x} \bigl( v_{i,x} - ( \hat \lambda_i - \lambda_i^*
+ \hth_i) v_i \bigr) \Bigr]   
\cr
& + \sum_i \hat r_{i,v\sigma} \Bigl[ \bigl( \h_{i,x} - \hat \lambda_i
\h_i \bigr) \hth_i'  v_{i,x}\iot_i /\h_i -
\bigl( \iot_{i,x} - \hat \lambda_i \iot_i \bigr) v_{i,x} \hth_i' - \h_i
\bigl( v_{i,x} - ( \hat \lambda_i - \lambda_i^* +
\hth_i) v_i \big) \hth_{i,x}
\Bigr]   
\cr
& + \sum_i  \hat r_{i,\sigma u} \tr_i\Bigl[ \bigl( \h_{i,x} -
\hat \lambda_i \h_i \bigr) v_i \hth_i' \iot_i/\h_i - \bigl(
\iot_{i,x} - \hat \lambda_i \iot_i \bigr) v_i \hth_i' \Bigr]
\cr
& + \sum_{i \not= j} \hat r_{i,\sigma u}\tr_j \Bigl[
\bigl( \h_{i,x} - \hat \lambda_i \h_i \bigr) v_j \hth_i' \iot_i/\h_i  -
\bigl( \iot_{i,x} - \hat \lambda_i \iot_i \bigr) v_j \hth_i' \Bigr]   
\cr
& + \sum_i \hat r_{i,\sigma \sigma} \Bigl[- \bigl( \h_{i,x} - \hat
\lambda_i \h_i \bigr) \hth_i' \hth_{i,x} \iot_i/\h_i +
\bigl( \iot_{i,x} - \hat
\lambda_i \iot_i \bigr) \hth_i' \hth_{i,x} \Bigr] +
\sum_i \Bigl( \h_i v_i \hr_{i,u}( \tilde r_i - \hat r_i )\Bigr)_x   \cr
& + \sum_{i \not=j}\h_i v_j 
\bigg[\sum_k v_k \bigl( 
\hr_{i,u}\tilde r_{j,u}\tr_k +\hr_{i,uu}(\tr_j\otimes\tr_k)\bigr) 
+  v_{j,x}\hr_{i,u} \tilde
r_{j,v} +v_{i,x}\hr_{vu}\tr_j- \theta_{j,x} 
\hat r_{i,u}  \tr_{j,\sigma}-\hth_{i,x}\hr_{\sigma u}\tr_j\bigg]\cr
&
+ \sum_{i,j} \h_i v_j \Big[ \big(\tilde r_j \bullet A(u)\big)\hat r_i 
-\big(\hr_i\bullet A(u)\big)\tilde r_j \Big]   \cr
\doteq&~ \sum_i \hat a_i(t,x),  
\cr}
\eqno(B.10)
$$
\vfill\eject
$$
\eqalign{
& \sum_i \Bigl( \h_{i,t} + 
\bigl( \hat\lambda_i \h_i \bigr)_x - \h_{i,xx} \Bigr)
\Bigl[ \hth_i' ( \iot_i / \h_i)^2 \hat
r_{i,\sigma} \Bigr] + \sum_i \Bigl( \iot_{i,t} + \bigl( \hat \lambda_i \iot_i
\bigr)_x - \iot_{i,xx} \Bigr) \Bigl[ \hat r_i - \hth_i'(\iot_i /
\h_i ) \hat r_{i,\sigma} \Bigr]  \cr
& - \sum_i \lambda_i^* \biggl\{ \Bigl( \h_{i,t} + ( \hat \lambda_i
\h_i )_x -
\h_{i,xx} \Bigr) \Bigl[ \hr_i + \hth_i'
(\iot_i/\h_i) \hat r_{i,\sigma} \Bigr] + \sum_i  \Bigl( \iot_{i,t} + (
\hat \lambda_i \iot_i )_x - \iot_{i,xx} \Bigr) \Big[-\hth_i'
\hat r_{i,\sigma}\Big] \biggr\}  \cr
=~& \sum_i \hr_{i,u}\tr_i \Bigl[ \bigl( \iot_{i,x} - \hat
\lambda_i \iot_i \bigr) v_i - \iot_i \bigl( w_i - \lambda_i^* v_i
\bigr) \Bigr]
\cr
& + \sum_{i \not= j} 
\hr_{i,u}\tr_j \Bigl[ ( \iot_{i,x} - \hat \lambda_i \iot_i
) v_j - \iot_i ( w_j - \lambda_j^*v_j ) + ( \iot_i v_j
)_x \Bigr]   
\cr
& + \sum_i \hat r_{i,v} \bigg[ \Bigl( \iot_i \bigl( v_{i,x} - (
\hat \lambda_i
- \lambda_i^*  + \hth_i) v_i \bigr)
\Bigr)_x  + \bigl( \iot_{i,x} - \hat \lambda_i
\iot_i \bigr) v_{i,x} - \iot_i v_{i,t} \bigg]   
\cr
& + \sum_i \hat r_{i,\sigma} \Bigl[ \bigl( \h_{i,x} - \hat \lambda_i
\h_i \bigr) \bigl( \hth_i' ( \iot_i/\h_i )^2 \bigr)_x
- (\iot_{i,x} - \hat \lambda_i \iot_i ) \bigl( \hth_{i,x} + (
\hth_i' \iot_i / \h_i )_x \bigr) \Bigr]   
\cr
& + \sum_i \hat r_{i,vu}\tr_i \Bigl[ \bigl( v_{i,x} -
( \hat \lambda_i - \lambda_i^*  + \hth_i) v_i
\bigr) \iot_i v_i \Bigr]
\cr
& + \sum_i \hat r_{i,vv} \Bigl[ \bigl( v_{i,x} -
( \hat
\lambda_i - \lambda_i^*  + \hth_i) v_i 
\bigr) \iot_i v_{i,x} \Bigr]   
\cr
& + \sum_{i \not= j} \hat r_{i,vu}\tr_j \Bigl[
\bigl( v_{i,x} - ( \tilde \lambda_i - \lambda_i^*  +
\hth_i) v_i
\bigr) \iot_i v_j \Bigr]   
\cr
& + \sum_i \hat r_{i,v\sigma} \Bigl[ \bigl( \h_{i,x} - \hat \lambda_i
\h_i \bigr) \hth_i' (\iot_i/\h_i)^2
v_{i,x} - \bigl( \iot_{i,x} - \hat \lambda_i \iot_i \bigr) \hth_i'
v_{i,x} \iot_i/\h_i 
- \iot_i \bigl( v_{i,x} - (\hat \lambda_i -
\lambda_i^* + \hth_i) v_i \bigr)
\hth_{i,x} \Bigr]   
\cr
& + \sum_i \hat r_{i,\sigma u}\tr_i \Bigl[ \bigl( \h_{i,x}
- \hat \lambda_i \h_i \bigr) \hth_i' ( \iot_i/\h_i )^2
v_i  -
\bigl( \iot_{i,x} - \hat \lambda_i \iot_i \bigr) \hth_i' v_i
\iot_i/\h_i \Bigr]   
\cr
& + \sum_{i \not= j} \hat r_{i,\sigma u}\tr_j \Bigl[
\bigl( \h_{i,x} - \hat \lambda_i \h_i \bigr) v_j \hth_i' ( \iot_i /
\h_i)^2 - \bigl( \iot_{i,x} - \hat \lambda_i \iot_i \bigr) v_j
\hth_i' \iot_i/\h_i \Bigr]   
\cr
& + \sum_i \hat r_{i,\sigma \sigma} \Bigl[- \bigl( \h_{i,x} - \hat
\lambda_i \h_i \bigr) \hth_i' (\iot_i/\h_i)^2
\hth_{i,x} +
\bigl( \iot_{i,x} - \hat \lambda_i \iot_i \bigr) \hth'
\hth_{i,x} \iot_i/\h_i \Bigr] + \sum_i \Bigl( \iot_i v_i \hr_{i,u} (
\tilde r_i - \hat r_i ) \Bigr)_x   
\cr
& + \sum_{i \not= j}\iot_i v_j\bigg[\sum_k  v_k 
\bigl( \hr_{i,u}\tr_{j,u}\tr_k+ \hr_{i,uu}(\tr_j\otimes\tr_k)\bigr) 
+v_{j,x}\hr_{i,u}\tr_{j,v}  +v_{i,x}\hr_{i,vu}\tr_j- \theta_{j,x} 
\hr_{i,u} \tilde r_{j,\sigma}-\hth_{i,x}\hr_{i,\sigma u}\tr_j\bigg] 
\cr
&+ \sum_{i,j} \bigl( w_i \h_j - v_i \iot_j
\bigr) \big(\tr_i \bullet A(u)\big) 
\hat r_j + \sum_{i,j} 
\bigg[v_i \h_j\Big( \big( \tilde r_i \bullet A(u)\big) \hat r_j
- \big(\hat r_j \bullet A(u)\big) \tilde r_i \Big) \bigg]_x\cr
&
+ \sum_{i \not= j} ( \lambda_j^*- \lambda_i^* ) v_i \h_j
\,\big(\tilde r_i \bullet A(u)\big) \hat r_j  
+ \sum_{i,j} v_i \h_j A(u)\Big[\big(\tilde r_i \bullet A(u)\big) 
\hat r_j - \big(\hat
r_j \bullet A(u) \big)\tilde r_i \Big]\cr
& - \sum_i \lambda_i^* \hat a_i(t,x)
\cr
\doteq &~ \sum_i \hat b_i(t,x) - \sum_i \lambda_i^*  \hat a_i(t,x).  
\cr}\eqno(B.11)$$
\vfill\eject

Recalling the expression (11.11) for the differential $\partial\Hat\Lambda/
\partial(\h,\iot)$, we can write (B.10)-(B.11) in the more compact
form
$${\partial\Hat\Lambda\over\partial (\h,\iot)}\cdot
\pmatrix{\big[\h_{i,t} + ( \hat \lambda_i \h_i )_x -
\h_{i,xx}\big]\cr
\big[\iot_{i,t} + \bigl( \hat \lambda_i \iot_i
\bigr)_x - \iot_{i,xx}\big] \cr}=\sum_i\pmatrix{\hat a_i\cr
\hat b_i-\lambda_i^*\hat a_i\cr}\,.
$$
By the uniform invertibility of the differential of $\Hat\Lambda$,
to prove the estimates stated in Lemma 11.4, it suffices to 
show that, for every $i=1,\ldots,n$,
the four quantities
$$\hat a_i\,,\qquad \hat b_i\,,\qquad
\big((\tla_i-\hla_i)\h_i\big)_x\,,\qquad 
\big((\tla_i-\hla_i)\iot_i\big)_x\,,
$$ 
can all be bounded according to the right hand side 
of (11.16).
\v
We start by looking at all the terms in the expressions (B.10)-(B.11)
for  $\hat a_i$ and $ \hat b_i$.
First of all, one checks that 
all those terms
which involve a product of distinct
components $i\not= j$
can be bounded as 
$$\O(1)\cdot\sum_{j \not= k}\Big(  | \h_j \h_k |+| \h_j v_k| + 
| \h_{j,x} v_k | +| \h_j v_{k,x} |
+ | \h_j w_k |+ | \iot_j v_k | +  | \iot_{x,j} v_k | +|
\iot_j v_{k,x} | + | \iot_j w_k|\Big).
\eqno(B.12)$$
For convenience, quantities whose size is bounded 
as in (B.12) will be called ``transversal terms''.
More generally, quantities whose size is bounded 
according to the right hand side of (11.16) will be called ``admissible
terms''.  We denote by $\A$ the family of all admissible terms.
We now exhibit various additional
terms which are admissible.
\v
\n{\bf 1.} 
By (6.16) it follows
$$\eqalign{\big( |\h_i|+|\iot_i|+&|\h_{i,x}|+|\iot_{i,x}|\big)
\Big|v_{i,x}-(\tla_i-\lambda_i^*)v_i-w_i\Big|\cr
&=
\O(1)\cdot\delta_0
\sum_{j\not= i}
\big( |\h_iv_j|+|\iot_iv_j|+|\h_{i,x}v_j|+|\iot_{i,x}v_j|\big)~\in~
\A\,.\cr}\eqno(B.13)$$
\v
\n{\bf 2.} Two other other admissible terms are
$$\eqalign{\h_i[w_{i,x}v_i-w_iv_{i,x}]&=[\h_i w_{i,x}-w_i \h_{i,x}]v_i+
w_i[\h_{i,x}v_i-\h_iv_{i,x}]\in\A\,,\cr
\iot_i[w_{i,x}v_i-w_iv_{i,x}]&=[\iot_i w_{i,x}-w_i \iot_{i,x}]v_i+
w_i[\iot_{i,x}v_i-\iot_iv_{i,x}]\in\A\,.\cr}
\eqno(B.14)$$
\v
\n{\bf 3.} We now consider terms that involve the
difference between the speeds: $\hth_i-\theta_i$.
We claim that the following four quantities are admissible:
$$\h_iv_i(\hth_i-\theta_i),\quad\iot_iv_i(\hth_i-\theta_i),
\quad
\h_{i,x}v_i(\hth_i-\theta_i),\quad\iot_{i,x}v_i(\hth_i-\theta_i)\in\A\,.
\eqno(B.15)$$
Indeed, from the definitions and the bounds (4.24) it follows
$$|\hla_i-\tla_i|=\O(1)\cdot|\hr_i-\tr_i|=\O(1)\cdot v_i 
\,|\hth_i-\theta_i|=\O(1)\cdot \delta_0
\,|\hth_i-\theta_i|.\eqno(B.16)$$
Since $|\theta'|\leq 1$, one has
$$|\hth_i-\theta_i|\leq\big|(\iot_i/\h_i)-(w_i/v_i)\big|.$$
Using (6.16) and (11.12) we now obtain
$$\eqalign{|h_iv_i|\big|\hth_i-\theta_i\big|&\leq
|\iot_i v_i-w_i \h_i|\cr
&=\bigg| \Big( \h_{i,x}+(\hla_i-\lambda_i^*)\h_i+\O(1)\cdot
\delta_0\sum_{j\not= i}\big(|\h_j|+|v_j|\big)\Big)v_i
\cr
&\qquad - \Big( v_{i,x}+(\tla_i-\lambda_i^*)v_i+\O(1)\cdot\delta_0
\sum_{j\not= i}|v_j|\Big)\h_i\bigg|\cr
&=\left|(\h_{i,x}v_i-v_{i,x}h_i)+(\hla_i-\tla_i)v_ih_i+\O(1)\cdot
\delta_0\sum_{j\not= i}\big(|v_jv_i|+|h_j v_i|\big)\right|\cr
&\leq \big|\h_{i,x}v_i-v_{i,x}h_i\big|+\O(1)\cdot\delta_0 
\big|\hth_i-\theta_i\big||h_iv_i|+\O(1)\cdot
\delta_0\sum_{j\not= k}\big(|v_jv_k|+|h_j v_k|\big)
\,.\cr}$$
Hence
$$|\iot_i v_i-w_i \h_i|\leq 2 \big|\h_{i,x}v_i-v_{i,x}h_i\big|
+\O(1)\cdot
\delta_0\sum_{j\not= k}\big(|v_jv_k|+|h_j v_k|\big)\in\A\,,\eqno(B.17)$$
showing that the quantity
$h_iv_i(\hth_i-\theta_i)$ is admissible.

Observing that $\hth_i - \theta_i \not= 0$ only if either
$|\iot_i/\h_i| \leq 6 \delta_1$ and $|w_i/v_i| \leq 3 \delta_1$, or
else $|\iot_i/\h_i| \geq 6 \delta_1$ and $|w_i/v_i| \leq 3 \delta_1$, 
we can write
$$\eqalign{
\Bigl| \iot_i v_i ( \hth_i - \theta_i ) \Bigr| &\leq
|\iot_i/\h_i|\, \bigl| \h_i v_i ( \hth_i - \theta_i
) \bigr| \cdot\chi_{\big\{| \iot_i/\h_i| \leq 6 \delta_1\big\}}
+ 2 \delta_1 | \iot_i v_i | \cdot\chi_{ \big\{|
\iot_i / \h_i | \geq 6 \delta_1, ~| w_i/v_i | \leq 3
\delta_1\big\}} \cr
&\leq 6 \delta_1 |\iot_i/\h_i|\, \bigl| \h_i v_i ( \hth_i - \theta_i
) \bigr| + 4\delta_1 \bigl|
\iot_i v_i - w_i \h_i \bigr|.
\cr}$$
Hence $\iot_iv_i(\hth_i-\theta_i)\in\A$.
In turn, using (11.12) we obtain
$$
\h_{i,x} v_i ( \hth_i - \theta_i) = \iot_i v_i \bigl( \hth_i
- \theta_i \bigr) + ( \hla_i - \lambda_i^* \bigr) \h_i
v_i ( \hth_i - \theta_i) + \O(1)\cdot\delta_0\sum_{j\not= i}
\big(|v_iv_j|+|v_i h_j|\big),
$$
showing that the term $\h_{i,x} v_i ( \hth_i -
\theta_i )$ is also admissible. Finally, using (6.16) one can write
$$\eqalign{
\iot_{i,x} v_i ( \hth_i - \theta_i ) &= ( \hth_i -
\theta_i ) \bigl[ \iot_{i,x} v_i - v_{i,x} \iot_i \bigr] + \iot_i
v_{i,x} ( \hth_i - \theta_i ) \cr
&= ( \hth_i - \theta_i ) \bigl[ \iot_{i,x} v_i - v_{i,x}
\iot_i \bigr] + \iot_i w_i ( \hth_i - \theta_i ) \cr
&\qquad+
(\tla_i-\lambda_i^*) \iot_i v_i ( \hth_i - \theta_i )+
\O(1)\cdot\delta_0\sum_{j\not= i}|v_iv_j|\,.\cr}$$
To estimate the term $\iot_i w_i ( \hth_i - \theta_i )$, 
we observe that $\hth_i-\theta_i=0$ if $|w_i/v_i|$
and $|\iot_i/\h_i|$ are both $\geq 3\delta_1$.
Hence, using again (6.16), we can write
$$\eqalign{
\bigl| \iot_i w_i \bigl( \hth_i - \theta_i \bigr) \bigr| &= 3
\delta_1 \bigl| \iot_i v_i \bigl(
\hth_i - \theta_i \bigr) \bigr|\cdot \chi_{\big\{ | w_i/v_i|<
3 \delta_1 \big\}} + 3 \delta_1 \bigl| \h_i w_i ( \hth_i -
\theta_i \bigr) |\cdot \chi_{
\big\{ | \iot_i/\h_i| < 3 \delta_1 \big\}} \cr
&= 3 \delta_1 | \iot_i v_i |\big| \hth_i - \theta_i \big|+
\Big| \h_i \big( v_{i,x} - ( \tilde \lambda_i -
\lambda_i^* ) v_i \big) \Bigr| | \hth_i - \theta_i |
+\O(1)\cdot\sum_{j\not= i} |h_iv_j|\,.
\cr}$$
By the previous estimates, this shows that
$\iot_{i,x} v_i ( \hth_i - \theta_i )\in\A$,
completing the proof of (B.15).
By (B.16), the following terms are also admissible:
$$\h_i(\tla_i-\hla_i),\quad\iot_i(\tla_i-\hla_i),
\quad
\h_{i,x}(\tla_i-\hla_i),\quad\iot_{i,x}(\tla_i-\hla_i)\in\A.
\eqno(B.18)$$
\v
\n{\bf 4.} 
Next, 
we claim that
$$
h_i(\tr_i-\hr_i)_x\,,\quad \iot_i(\tr_i-\hr_i)_x\,,\quad
h_i(\tla_i-\hla_i)_x\,,\quad \iot_i(\tla_i-\hla_i)_x\in\A\,.
\eqno(B.19)$$
Indeed, one can write
$$\eqalign{
\h_i( \tr_i-\hr_i)_x &= \h_i v_i ( \hth_i - \theta_i) 
\biggl\{ \sum_j v_j
{( \tr_{i,u} - \hr_{i,u})\tr_j \over v_i ( \hth_i - \theta_i)}+
v_{i,x}  {\tr_{i,v}-\hr_{i,v}\over v_i ( \hth_i - \theta_i )} 
\biggr\} 
\cr
& \qquad+ \hth_i' ( \iot_i/\h_i ) 
\bigl[ v_{i,x} \h_i - \h_{i,x} v_i \bigr]
(\hr_{i,\sigma}/v_i) + \hth_i' 
\bigl[ v_i \iot_{i,x} - \iot_i v_{i,x} \bigr]
(\hr_{i,\sigma}/v_i) 
\cr
& \qquad+ (w_i/v_i) \theta_i' \bigl[ v_{i,x} \h_i - v_i \h_{i,x} \bigr]
( \tr_{i,\sigma} / v_i ) + \theta_i' \bigl[
\h_{i,x} w_i - \h_i
w_{i,x} \bigr] ( \tr_{i,\sigma}/v_i),
\cr}$$
$$\eqalign{
\iot_i ( \tr_i-\hr_i)_x
 &= \iot_i v_i ( \hth_i - \theta_i) \biggl\{ \sum_j v_j
{( \tr_{i,u} - \hr_{i,u})\tr_j \over v_i ( \hth_i - \theta_i)}+
v_{i,x}  {\tr_{i,v}-\hr_{i,v}\over v_i ( \hth_i - \theta_i )} 
\biggr\} 
\cr
& \qquad+ \hth_i' ( \iot_i/\h_i )^2 
\bigl[ v_{i,x} \h_i - \h_{i,x} v_i \bigr]
(\hr_{i,\sigma}/v_i) + \hth_i'(\iot_i/\h_i) 
\bigl[ v_i \iot_{i,x} - \iot_i v_{i,x} \bigr]
(\hr_{i,\sigma}/v_i) 
\cr
& \qquad+ (w_i/v_i) \theta_i' \bigl[ v_{i,x} \iot_i - v_i \iot_{i,x} \bigr]
( \tr_{i,\sigma} / v_i ) + \theta_i' \bigl[
\iot_{i,x} w_i - \iot_i
w_{i,x} \bigr] ( \tr_{i,\sigma}/v_i).
\cr}$$
By (4.24), the above expressions within braces are
uniformly bounded.  Hence the first two
quantities in (B.19) are admissible. To prove the admissibility
of the last two terms it suffices to repeat the above computation,
with $\tr_i$ and $\hr_i$ replaced by
by $\tla_i$ and $\hla_i$.
\v
In a similar way as in Appendix~A,
we are now ready to check one by one
all the (non-tranversal) terms
in the expressions of $\hat a_i,\hat b_i$
in (B.10)-(B.11), showing that all of them 
are admissible.
\v
\n Coefficients of $\hr_{i,u}\tr_i$ :
$$\eqalign{
v_i ( \h_{i,x} - &\hat \lambda_i \h_i
) - \h_i ( w_i - \lambda_i^* v_i )\cr 
&= \big[ v_i
\h_{i,x} - \h_i v_{i,x} \big] + \big[ v_i \h_i ( \tilde \lambda_i -
\hat \lambda_i ) \big] + \Bigl[ \h_i \bigr( v_{i,x} - ( \tilde \lambda_i -
\lambda_i^* ) v_i - w_i \bigr) \Bigr],
\cr}
$$
$$
\eqalign{
v_i ( \iot_{i,x} -& \hat \lambda_i \iot_i
) - \iot_i ( w_i - \lambda_i^* v_i )\cr 
&=
\bigl[ \iot_{i,x} v_i - v_{i,x} \iot_i
\bigr] + \bigl[ v_i \iot_i ( \tilde \lambda_i - \hat \lambda_i
) \bigr] 
+ \Bigl[ \iot_i \bigr( v_{i,x} - ( \tilde \lambda_i -
\lambda_i^* ) v_i - w_i \bigr) \Bigr].\cr}$$
\v
\n Coefficients of $\hr_{i,v}$ :
$$\eqalign{
v_{i,x} ( \h_{i,x} - \hat \lambda_i \h_i ) + &\Bigl( \h_i \bigl(
v_{i,x} - ( \hat \lambda_i - \lambda_i^* + \hth_i)
v_i \bigr) \Bigr)_x - \h_i v_{i,t}\cr
&= 2 \Bigl[
\h_{i,x} \bigl( v_{i,x} - ( \tilde \lambda_i - \lambda_i^* )
v_i - w_i \bigr) \Bigr] 
+ 2 \big[ \h_{i,x} ( w_i - \theta_i v_i ) \big]\cr
&\qquad  + \bigl[
\h_i v_{i,x} - \h_{i,x} v_i \bigr] \bigl(
\lambda_i^* - \hth_i - \hat \lambda_i + \hth'_i \iot_i/\h_i
\bigr) + \hth_i' \bigl[ v_i \iot_{i,x} - v_{i,x} \iot_i \bigr] 
\cr
&\qquad + 2 \bigl[ \h_{i,x} v_i ( \hat \lambda_i + \hth_i - \tilde \lambda_i
- \theta_i ) \bigr] + \Bigl[ \h_i \bigl( ( \tilde \lambda_i -
\hat \lambda_i )
v_i \bigr)_x \Bigr] - \h_i \phi_i\,,
\cr}
$$
$$\eqalign{
\Bigl( \iot_i \bigl( v_{i,x} - ( & \hat \lambda_i - \lambda_i^* 
+ \hth_i) v_i \bigr) \Bigr) _x + v_{i,x} ( \iot_{i,x}
- \hat \lambda_i \iot_i ) -
\iot_i v_{i,t} \cr
&= 2 \Bigl[ \iot_{i,x} \bigl( v_{i,x} - \bigl( \tilde
\lambda_i - \lambda_i^* \bigr) v_i - w_i \bigr) \Bigr] + 2 \Bigl[
\iot_{i,x} \bigl( w_i - \theta_i v_i \bigr) \Bigr]\cr
&\qquad  + \bigl(
\lambda_i^* -\hth_i- \hat \lambda_i + \hth'_i
\iot_i/\h_i \bigr)
\bigl[ v_{i,x} \iot_i - v_i \iot_{i,x} \bigr] 
+ \hth_i' ( \iot_i / \h_i )^2 \bigl[ v_i \h_{i,x} -
v_{i,x} \h_i \bigr] \cr
&\qquad + 2 \bigl[ \iot_{i,x} v_i ( \tla_i+\theta_i-\hla_i-\hth_i
) \bigr] + \Bigl[ \iot_i
\bigl(( \tilde \lambda_i - \hat \lambda_i )v_i \bigr)_x
\Bigr] - \iot_i \phi_i\,.
\cr}$$
\v
\n Coefficients of $\hat r_{i,\sigma}/v_i$ :
$$
v_i ( \h_{i,x} - \hat \lambda_i \h_i) \bigl(- \hth_{i,x} +( 
\hth'_i \iot_i / \h_i )_x \bigr) - v_i( \iot_{i,x} -
\hat \lambda_i \iot_i
) \hth'_{i,x}  = -v_i \left[ \hth_i'' \h_i \left(
 {\iot_i\over \h_i} \right)_x^2 \right],
$$
$$
v_i \bigl( \h_{i,x} - \hat \lambda_i \h_i \bigr) \Bigl(
\hth'_i \bigl( \iot_i / \h_i \bigr)^2 \Bigr)_x + v_i \bigl(
\iot_{i,x} - \hat
\lambda_i \iot_i \bigr) \Bigl( \hth'_{i,x} - \bigl(
\hth'_i \iot_i / \h_i \bigr)_x \Bigl) =- v_i \left[ \left(
\hth''_i  {\iot_i\over \h_i} + 2 \hth'_i
\right) \h_i \left(  {\iot_i\over \h_i} \right)_x^2 \right].
$$
\v
\n Coefficients of 
$\hat r_{i,vu}\tr_i$ :
$$\eqalign{
v_i \h_i \bigl( v_{i,x} - ( \hat \lambda_i -
\lambda_i^* + \hth_i) v_i \bigr)
&= \Bigl[ v_i \h_i \bigl( v_{i,x} - (
\tilde \lambda_i - \lambda_i^* ) v_i - w_i \bigr)
\Bigr] + \Bigl[ v_i \h_i \bigl( w_i - \theta_i v_i \bigr) \Bigr] \cr
&\qquad + \Bigl[ v_i^2 \h_i (\tla_i+\theta_i-\hla_i-\hth_i) \Bigr],
\cr}$$
$$\eqalign{
v_i \iot_i \bigl( v_{i,x} - ( \hat \lambda_i - \lambda_i^*  
+ \hth_i) v_i \bigr) &= \Bigl[ v_i \iot_i \bigl( v_{i,x} -
( \tilde
\lambda_i - \lambda_i^* ) v_i - w_i \bigr) \Bigr] +
\Bigl[ v_i \iot_i ( w_i - \theta_i v_i) \Bigr] \cr
&\qquad + \Bigl[  v_i^2 \iot_i (\tla_i+\theta_i-\hla_i-\hth_i)\Bigr].
\cr}$$
\v
\n Coefficients of 
$\hat r_{i,vv}$ :
$$\eqalign{
\h_i v_{i,x} \bigl( v_{i,x} - ( \hat \lambda_i - \lambda_i^* 
+ \hth_i) v_i \bigr) &= \Bigl[ \h_i v_{i,x}
\bigl( v_{i,x} - ( \tilde \lambda_i - \lambda_i^*
+ \theta_i)  v_i \bigr) \Bigr] + \Bigl[ \h_i v_i v_{i,x} (
\tla_i+\theta_i-\hla_i-\hth_i) \Bigr]
\cr
&= \Bigl[ v_i \h_{i,x} \bigl( v_{i,x} - ( \tilde \lambda_i -
\lambda_i^* ) v_i - w_i \bigr) \Bigr] + \Bigl[ \h_i v_i
v_{i,x} (\tla_i+\theta_i-\hla_i-\hth_i
) \Bigr] \cr
&\qquad + v_i\bigl[ \h_{i,x} ( w_i - \theta_i v_i ) \bigr] + \bigl(
v_{i,x} + ( \tilde \lambda_i - \lambda_i^*) v_i
- \theta_i v_i \bigr) \bigl[ \h_i v_{i,x} - v_i \h_{i,x} \bigr],
\cr}$$
$$\eqalign{
\iot_i v_{i,x} \bigl( v_{i,x} - ( \hat \lambda_i - \lambda_i^*
) v_i + \hth_i v_i \bigr) &= v_i \Bigl[ \iot_{i,x} \bigl(
v_{i,x} - ( \tilde \lambda_i - \lambda_i^* ) v_i -
w_i \bigr) \Bigr] + \big[ \iot_i v_i v_{i,x} 
( \tla_i+\theta_i-\hla_i-\hth_i) \big] \cr
&\qquad + v_i \big[ \iot_{i,x} ( w_i - \theta_i v_i ) \big] +
\bigl( v_{i,x} - ( \tilde \lambda_i - \lambda_i^* ) v_i
- \theta_i v_i \bigr) \big[ \iot_i v_{i,x} - v_i \iot_x \big]\,.
\cr}$$
\v
\n Coefficients of $\hat r_{i,v\sigma}$ :
$$\eqalign{
\bigl(\h_{i,x} &- \hat \lambda_i \h_i \bigr)  \hth'_i (\iot_i/h_i) v_{i,x}
- ( \iot_{i,x} - \hat \lambda_i \iot_i
) \hth'_i v_{i,x} - \h_i \bigl( v_{i,x} - ( \hat
\lambda_i - \lambda_i^*  +
\hth_i) v_i \bigr) \hth_{i,x} \cr
&= v_{i,x} \hth'_i ( \h_{i,x} \iot_i - \h_i \iot_{i,x} ) / \h_i 
-\hth_{i,x} \h_i \bigl( v_{i,x} - ( \hat \lambda_i -
\lambda_i^* + \hth_i) v_i \bigr) \cr
&= 2 \hth'_i \big( \h_{i,x}  (\iot_i/ \h_i)  -
\iot_{i,x} \big) \Bigl[
\bigl( v_{i,x} - ( \tilde \lambda_i - \lambda_i^* ) v_i -
w_i \bigr) + ( w_i - 
\theta_i v_i ) \Bigr]  \cr
&\qquad + \Bigl( 2 ( \tilde \lambda_i - \lambda_i^* + \theta_i
) - ( \hat \lambda_i
- \lambda_i^* + \hth_i ) \Bigr) \Big\{ \hth'_i \big[ v_{i,x}
\iot_i - v_i \iot_{i,x} \big] 
- \hth_i' ( \iot_i / \h_i )
\big[ v_{i,x} \h_i - v_i \h_{i,x} \big] \Big\},
\cr}$$
$$\eqalign{
\bigl(& \h_{i,x} - \hat \lambda_i \h_i \bigr)  \hth'_i
( \iot_i /\h_i )^2 v_{i,x} - \bigl( \iot_{i,x} - \hat \lambda_i
\iot_i
\bigr) \hth'_i v_{i,x}\iot_i / \h_i 
- \iot_i \bigl( v_{i,x} -
( \hat \lambda_i - \lambda_i^*  +
\hth_i) v_i \bigr) \hth_{i,x} \cr
&= 2 \theta'_i  {\iot_i\over \h_i} \left( \h_{i,x}
{\iot_i\over \h_i} - \iot_{i,x}
\right) \Bigl[ \bigl( v_{i,x} - ( \tilde \lambda_i - \lambda_i^*
) v_i - w_i \bigr) + ( w_i - \theta v_i ) \Bigr] \cr
&\qquad + \Bigl( 2 ( \tilde \lambda_i - \lambda_i^* + \theta_i) 
- ( \hat \lambda_i - \lambda_i^* + \hth_i ) \Bigr)
\Big\{ \hth_i' ( \iot_i / \h_i ) \Bigl[ v_{i,x} \iot_i - v_i
\iot_{i,x} \Bigr] 
- \hth_i' ( \iot_i/\h_i )^2 \Bigl[ v_{i,x}
\h_i - v_i \h_{i,x} \Bigr] \bigg\}.
\cr}$$
\v
\n Coefficients of 
$\hat r_{i,\sigma u}\tr_i/v_i$ :
$$
( \h_{i,x} - \hat \lambda \h_i ) \hth'_i  v_i^2 \iot_i/\h_i
- ( \iot_{i,x} - \hat
\lambda_i \iot_i ) v_i^2 \hth'_i = \hth'_i
v_i \bigl[
v_{i,x} \iot_i - v_i \iot_{i,x} \bigr] + \hth_i' v_i 
(\iot_i / \h_i ) \bigl[ v_i \h_{i,x} - \h_i v_{i,x} \bigr],
$$
$$\eqalign{
( \h_{i,x} - \hat \lambda \h_i ) &\hth'_i  v_i^2 
( \iot_i / \h_i) ^2
- ( \iot_{i,x} - \hat \lambda_i \iot_i ) \hth'_i 
v_i^2 ( \iot_i / \h_i)
\cr
& =\hth'_i v_i ( \iot_i / \h_i ) \Big\{ \bigl[ v_{i,x}
\iot_i - v_i
\iot_{i,x} \bigr] + ( \iot_i / \h_i ) \bigl[ v_i \h_{i,x} - \h_i
v_{i,x} \bigr] \Big\}.
\cr}$$
\v
\n Coefficients of 
$\hat r_{i,\sigma\sigma}/v_i$ :
$$
-\bigl( \h_{i,x} - \hat \lambda_i \h_i \bigr) v_i ( \iot_i / \h_i )
\hth'_i \hth_{i,x} +
\bigl( \iot_{i,x} - \hat \lambda_i \iot_i \bigr) v_i \hth'_i
\hth_{i,x}
=( \hth'_i)^2 v_i  \h_i \big[
 (\iot_i/ \h_i)_x\big]^2 ,
$$
$$
-\bigl( \h_{i,x} - \hat \lambda_i \h_i \bigr) v_i ( \iot_i / \h_i )^2
\hth'_i \hth_{i,x} + \bigl( \iot_{i,x} - \hat \lambda_i
\iot_i \bigr) \hth'_i \hth_{i,x} v_i ( \iot_i /
\h_i ) = ( \hth_i')^2 v_i  \iot_i\big[(\iot_i/ \h_i)_x\big]^2 .
$$
\v
There are a few remaining terms in (B.10)-(B.11)
which we now examine.
Recalling (B.14) we have
$$\eqalign{
\Bigl( \h_i v_i  \hr_{i,u}( \tilde r_i - \hat r_i )
\Bigr)_x =&\O(1)\cdot \h_{i,x}v_i^2(\theta_i-\hth_i)+
\O(1)\cdot \h_iv_iv_{i,x}(\theta_i-\hth_i)\cr
&\qquad +\O(1)\cdot \h_iv_i^2(\theta_i-\hth_i)
+\O(1)\cdot h_iv_i(\tr_i-\hr_i)_x\,,\cr}
$$
$$\eqalign{
\Bigl( \iot_i v_i  \hr_{i,u}( \tilde r_i - \hat r_i )
\Bigr)_x =&\O(1)\cdot \iot_{i,x}v_i^2(\theta_i-\hth_i)+
\O(1)\cdot \iot_iv_iv_{i,x}(\theta_i-\hth_i)\cr
&\qquad +\O(1)\cdot \iot_iv_i^2(\theta_i-\hth_i)
+\O(1)\cdot \iot_iv_i(\tr_i-\hr_i)_x\,,\cr}
$$
$$\eqalign{
\h_i v_i \Big[ \big(\tr_i \bullet A(u)\big)\hr_i 
-\big(\hr_i\bullet A(u)\big)\tr_i \Big]
&=\O(1)\cdot \h_i v_i(\tr_i-\hr_i)\,,\cr
\h_i v_i A(u)\Big[ \big(\tr_i \bullet A(u)\big)\hr_i 
-\big(\hr_i\bullet A(u)\big)\tr_i \Big]
&=\O(1)\cdot \h_i v_i(\tr_i-\hr_i)\,.
\cr
( w_i \h_i - v_i \iot_i)
\big(\tr_i \bullet A(u)\big) 
\hat r_i &=\O(1)\cdot |w_i \h_i-\iot_i v_i|\,,\cr}$$
$$\eqalign{\bigg[v_i \h_i\Big( \big( \tilde r_i \bullet A(u)\big) \hat r_i
- \big(\hat r_i \bullet A(u)\big) \tilde r_i \Big) \bigg]_x
&=\O(1)\cdot \big(|v_{i,x}\h_i|+|v_i\h_{i,x}|\big)
|\tr_i-\hr_i|+\O(1)\cdot v_i\h_i (\tr_i-\hr_i)_x\,,\cr}$$
These terms are all admissible because of (B.15)--(B.19).

We have thus completed the analysis of all terms in (B.10)-(B.11),
showing that the quantities $\hat a_i$, $\hat b_i$ are admissible.
The admissibility of the terms $\bigl( ( \tla_i - \hla_i ) \h_i \bigr)_x$
and $\bigl( ( \tla_i - \hla_i ) \iot_i \bigr)_x$ follows immediately
from (B.18) and (B.19).  This completes the proof of
Lemma 11.4.
\endproof
\vsk

\c{\medbf Appendix C}
\v
Aim of this section is to derive
energy estimates for the components $\h_i$, $\iot_i$
and prove the bounds (11.33)-(11.34).
We write the evolution equations (11.15)
for the components $\h_i,\iot_i$ in the form
$$
\left\{\eqalign{\h_{i,t} + ( \tla_i\h_i )_x -
\h_{i,xx} &= \hat\phi_i\,,\cr
\iot_{i,t} + ( \tla_i\iot_i )_x -
\iot_{i,xx} &= \hat\psi_i\,.\cr}\right.
\eqno(C.1)
$$
For convenience, we define
$\hvt_i\doteq\vth(\iot_i/\h_i)$.
Multiplying the first equation
in (C.1) 
by $\h_i \hvt_i$ and
integrating by parts, we obtain
$$\eqalign{\int\hvt_i\h_i\hat\phi_i\,dx&=
\int\Big\{\hvt_i\h_i \h_{i,t}+\hvt_i
\h_i(\tla_i\h_i)_x-\hvt_i\h_i\h_{i,xx}
\Big\}\,dx\cr
&=\int\Big\{\hvt_i
(\h_i^2/ 2)_t
-\hvt_i\tla_i \h_i\h_{i,x}-\hvt_{i,x}\tla_i\h_i^2+\hvt_i \h_{i,x}^2+
\hvt_{i,x}\h_{i,x}\h_i\Big\}\,dx\cr
&=\int\Big\{\big(\hvt_i\h_i^2/2 \big)_t+(\tla_i\hvt_i)_x(\h_i^2/2)
-\big(\hvt_{i,t}
+2\tla_i\hvt_{i,x}-\hvt_{i,xx}\big)(\h_i^2/2) +\hvt_i \h_{i,x}^2
\Big\}\,dx\,.\cr}$$
Therefore
$$\eqalign{ \int \hvt_i \h_{i,x}^2\,dx&=
-{d\over dt}\left[\int\hvt_i \h_i^2/2\,dx\right]
+\int \big(\hvt_{i,t}
+\tla_i\hvt_{i,x}-\hvt_{i,xx}\big)(\h_i^2/2)\,dx\cr
&\qquad -\int \tla_{i,x}\hvt_i(\h_i^2/2)\,dx+
\int\hvt_i\h_i\hat\phi_i\,dx\,.\cr}
\eqno(C.2)$$
As in (9.14),
a direct computation yields
$$\hvt_{i,t}
+\tla_i\hvt_{i,x}-\hvt_{i,xx}
=\hvt'_i \left( {\hat\psi_i\over \h_i} - {\iot_i\over \h_i}
{\hat\phi_i\over \h_i} \right) + 2 \hvt'_i {\h_{i,x}\over \h_i} \left(
{\iot_i\over \h_i} \right)_x
- \hvt''_i \left( {\iot_i\over \h_i} \right)_x^2\,.
\eqno(C.3)$$
Since $\tla_{i,x}=(\tla_i-\lambda_i^*)_x$,
integrating by parts and using the second estimate in (11.13) one obtains
$$\eqalign{\left|\int \tla_{i,x}\hvt_i(\h_i^2/2)\,dx\right|&=
\left|\int (\tla_i-\lambda_i^*) \big(\hvt_{i,x}\h_i^2/2
+\hvt_i\h_i\h_{i,x}\big)\,dx
\right|\cr
&\leq \|\tla_i-\lambda_i^*\|_{\L^\infty}\cdot\bigg\{
{1\over 2}\int
\big|\hvt_i'\big|\,|\iot_{i,x}\h_i-\iot_i \h_{i,x}|\,dx+{5\over 2\delta_1} 
\int \hvt_i \h_{i,x}^2\,dx\cr
&\qquad+\O(1)\cdot
\delta_0\sum_{j\not= i}\int\big(|\h_iv_j|+|\h_i\h_j|\big)\,dx\bigg\}\cr
&\leq \int |\iot_{i,x}\h_i-\h_{i,x}\iot_i|\,dx
+{1\over 2}\int \hvt_i \h_{i,x}^2\,dx
+\delta_0\sum_{j\not= i}\int\big(|\h_iv_j|+|\h_i\h_j|\big)\,dx
\,,\cr}
\eqno(C.4)$$
because
$$|\tla_i-\lambda_i^*|=\O(1)\cdot\delta_0<\!<\delta_1\leq 1\,.$$
Using (C.3) and (C.4) in (C.2), we now obtain
$$\eqalign{
{1\over 2}\int \hvt_i\,\h_{i,x}^2 \, dx &\leq -
{d\over dt}\left[ \int {\hvt_i 
\h_i^2\over 2} \,dx\right] 
+{1\over 2}\int \big|\hvt_i'\big|\big(|\h_i\hat\psi_i|
+|\iot_i\hat\phi_i|\big)\,dx
+\int \left| \hvt_i' \h_i\h_{i,x}\left(\iot_i\over \h_i\right)_x
\right|\,dx
\cr
&\qquad+{1\over 2}\int\left|\hvt_i''\h_i^2
\left(\iot_i\over \h_i\right)_x^2\right|\,dx
+\int
\big| \iot_{i,x}
\h_i -\iot_i \h_{i,x} \big|\, dx\cr
&\qquad
+\delta_0
\sum_{j\not=i}\int
\big(|\h_iv_j|+|\h_i\h_j|\big)
\,dx
+\int |\h_i\hat\phi_i|\,dx\,.\cr}\eqno(C.5)$$
Recalling the definition of $\hvt_i$,
on regions where $\hvt_i'\not=0$ one has $|\iot_i/\h_i|\leq 4\delta_1/5$,
hence the bounds (11.14) hold. In turn, they imply
$$\eqalign{
\left|\hvt_i' \h_i\h_{i,x}\left(\iot_i\over \h_i\right)_x
\right|&\leq
{5\over 2\delta_1}\left| \hvt_i' \h_i^2\left(\iot_i\over \h_i\right)_x
\right|+\O(1)\cdot
\delta_0\sum_{j\not= i}\left| \hvt_i' \h_i\left(\iot_i\over \h_i\right)_x
\right|\big(|v_j|+|h_j|\big)
\cr
&=\O(1) \cdot |\iot_{i,x}\h_i-\iot_i\h_{i,x}|
+\O(1)\cdot
\delta_0\sum_{j\not= i}\big(|v_j\iot_{i,x}|+|v_j\h_{i,x}|+
|\h_j\iot_{i,x}|+|\h_j\h_{i,x}|\big)\,.\cr}
\eqno(C.6)$$
Using (C.6) and then the bounds (11.18), (11.26), (11.28), 
(11.30) and (11.31), from (C.5) we
conclude
$$\eqalign{
\int_{\hat t}^T \!\int  \hvt_i\,\h_{i,x}^2 \, dxdt
&\leq \int\hvt_i \h_i^2(\hat t,x)\,dx
+\O(1)\cdot\int_{\hat t}^T\!\int \big(|\h_i\hat\psi_i|
+|\iot_i\hat\phi_i|\big)\,dxdt
\cr
&\qquad
+\O(1)\cdot\int_{\hat t}^T\!\int|\iot_{i,x}\h_i-\iot_i\h_{i,x}|\,dxdt
+\O(1)\cdot\int_{\hat t}^T\!\int 
\big|\h_i
(\iot_i/ \h_i)_x\big|^2\,dxdt\cr
&\qquad
+\O(1)\cdot\delta_0\sum_{j\not= i}\int_{\hat t}^T\!\int 
\big(|v_j\iot_{i,x}|+|v_j\h_{i,x}|+
|\h_j\iot_{i,x}|+|\h_j\h_{i,x}|\big)\,dxdt
\cr
&\qquad
+\delta_0\sum_{j\not=i}\int_{\hat t}^T\!\int
\big(|\h_iv_j|+|\h_i\h_j|\big)\,dxdt
+2\int_{\hat t}^T\!\int |\h_i\hat\phi_i|\,dxdt\cr
&=\O(1)\cdot \delta_0^2\,,\cr}\eqno(C.7)$$
proving the estimate (11.33)
\v
We now perform a similar computation for $\iot_{i,x}^2$.
Multiplying the second equation in (C.1) by $\hvt_i\iot_i$ and
integrating by parts, one obtains
$$\int\hvt_i\iot_i\hat\psi_i\,dx
=\int\Big\{\big(\hvt_i\iot_i^2/2 \big)_t+(\tla_i\hvt_i)_x(\iot_i^2/2)
-\big(\hvt_{i,t}
+2\tla_i\hvt_{i,x}-\hvt_{i,xx}\big)(\iot_i^2/2) +\hvt_i \iot_{i,x}^2
\Big\}\,dx\,.$$
Therefore, the identity (C.2) still holds, with
$\h_i,\hat\phi_i$ replaced by $\iot_i,\hat\psi_i$, respectively:
$$\eqalign{ \int \hvt_i \iot_{i,x}^2\,dx&=
-{d\over dt}\left[\int\hvt_i \iot_i^2/2\,dx\right]
+\int \big(\hvt_{i,t}
+\tla_i\hvt_{i,x}-\hvt_{i,xx}\big)(\iot_i^2/2)\,dx\cr
&\qquad -\int \tla_{i,x}\hvt_i(\iot_i^2/2)\,dx+
\int\hvt_i \iot_i\hat\psi_i\,dx\,.\cr}
\eqno(C.8)$$
The equality (C.3) can again be used.  To obtain a suitable
replacement for (C.4) we observe that, if $\hvt_i\not= 0$,
then (11.13) implies
$$|\iot_i\iot_{i,x}|\leq 2|\h_{i,x}\iot_{i,x}|+\O(1)\cdot
\delta_0 \sum_{j\not= i}
\big(|v_j\iot_{i,x}|+|\h_j\iot_{i,x}|\big)$$
and hence
$$|\iot_i\iot_{i,x}|\leq 
\h^2_{i,x} + \iot_{i,x}^2
+\O(1)\cdot\delta_0 \sum_{j\not= i}
\big(|v_j\iot_i|+|\h_j\iot_i|\big)
\,.$$
Integrating by parts we thus obtain
$$\eqalign{&\left|\int \tla_{i,x}\hvt_i(\iot_i^2/2)\,dx\right|=
\left|\int (\tla_i-\lambda_i^*) \big(\hvt_{i,x}\iot_i^2/2+
\hvt_i\iot_i\iot_{i,x}
\big)\,dx
\right|\cr
&\qquad\qquad\leq \|\tla_i-\lambda_i^*\|_{\L^\infty}\cdot\bigg\{\int
\big|\hvt_i'\big|\,|\iot_{i,x}\h_i-\iot_i\h_{i,x}|\,
\left|\iot_i^2\over \h_i^2
\right|\,dx+
\int \hvt_i \h_{i,x}^2\,dx\cr
&\qquad\qquad\qquad+ 
\int \hvt_i \iot_{i,x}^2\,dx
+\O(1)\cdot\delta_0\sum_{j\not= i}\int\big(|v_j \iot_i|+|\h_j
\iot_i|\big)\,dx\bigg\}\cr
&\qquad\qquad\leq
\int |\iot_{i,x}\h_i-\h_{i,x}\iot_i|\,dx
+{1\over 2}\int \hvt_i \h_{i,x}^2\,dx
+ {1\over 2}
\int \hvt_i \iot_{i,x}^2\,dx
+\delta_0\sum_{j\not= i}\int\big(|v_j \iot_i|+|\h_j
\iot_i|\big)\,dx
\,.\cr}
\eqno(C.9)$$
Using (C.3) and (C.9) in (C.8) and observing that
$|\iot_i^2/ \h_i^2|\leq \delta_1^2$ on the region where 
$\hvt_i'\not=0$, 
we now obtain an estimate similar to (C.5):
$$\eqalign{
{1\over 2}\int \hvt_i\,w_{i,x}^2 \, dx &\leq -
{d\over dt}\left[ \int {\hvt_i 
\iot_i^2\over 2} \,dx\right] 
+{\delta_1^2\over 2}\int 
|\hvt_i'|\big(|\h_i\hat\psi_i|+|\iot_i\hat\phi_i|\big)\,dx
+\delta_1^2\int \left| \hvt_i' \h_i\h_{i,x}\left(\iot_i\over \h_i\right)_x
\right|\,dx
\cr
&\qquad +{\delta_1^2\over 2}\int\left|\hvt_i''\h_i^2
\left(\iot_i\over \h_i\right)_x^2\right|\,dx
+\int| \iot_{i,x}
\h_i -\iot_i \h_{i,x} |\, dx+ {1\over 2}\int \hvt_i \h_{i,x}^2\,dx\cr
&\qquad 
+\delta_0\sum_{j\not= i}\int\big(|v_j\iot_i|+|h_j\iot_i|\big)\,dx
+\int |\iot_i\hat\psi_i|\,dx\,.\cr}\eqno(C.10)$$
Using (C.6) and then the bounds
(C.7), (11.18), (11.26), (11.28), 
(11.30) and (11.31), from (C.10) we
conclude
$$\eqalign{
\int_{\hat t}^T \!&\int \hvt_i\,\iot_{i,x}^2 \, dxdt
\leq \int\hvt_i \iot_i^2(\hat t,x)\,dx
+\O(1)\cdot\int_{\hat t}^T\!\int 
\big(|\h_i\hat\psi_i|+|\iot_i\hat\phi_i|\big)\,dxdt
\cr
&\qquad
+\O(1)\cdot\int_{\hat t}^T\!\int|\iot_{i,x}\h_i-\iot_i\h_{i,x}|\,dxdt
+\O(1)\cdot\int_{\hat t}^T\!\int 
\big|\h_i
(\iot_i/ \h_i)_x\big|^2\,dxdt\cr
&\qquad
+\O(1)\cdot\delta_0\sum_{j\not= i}\int_{\hat t}^T\!\int 
\big(|v_j\iot_{i,x}|+|v_j\h_{i,x}|+
|\h_j\iot_{i,x}|+|\h_j\h_{i,x}|\big)\,dxdt
\cr
&\qquad+\int_{\hat t}^T\!\int \hvt_i\,\h_{i,x}^2 \, dxdt
+\delta_0\sum_{j\not= i}\int\big(|v_j \iot_i|+|\h_j
\iot_i|\big)\,dx
+2\int_{\hat t}^T\!\int |\h_i\hat\phi_i|\,dxdt\cr
&=\O(1)\cdot \delta_0^2\,,\cr}\eqno(C.11)$$
proving the estimate (11.34).
\vfill\eject
\c{\medbf Appendix D}
\v
We derive here the two estimates (12.9)-(12.10), used in the proof
of Lemma 12.1.
$$
\eqalign{
\|A &\|_{L^\infty} \int_0^t\!\!\int \bigl| G_x(t-s,x-y) \bigr|
E(s,y) dyds \cr
&= \| A \|_{L^\infty} \int_0^t\!\!\int  {|x -
y|\over 4 (t-s) \sqrt{\pi(t-s)}} B(s)
\cr
&\qquad\qquad\qquad
\cdot\exp \biggl\{ -  {(x-y)^2\over 4(t-s)} + 4 \| DA \|_{L^\infty}
\int_0^s \bigl\| u_x(\sigma) \bigr\|_{L^\infty} d\sigma + s - y
\biggr\} dyds \cr
&\leq \| A \|_{L^\infty} \exp \left\{ 4 \| DA \|_{L^\infty}
\int_0^t \bigl\| u_x(\sigma) \bigr\|_{L^\infty} d\sigma
+ t - x \right\} \cr
&\qquad \qquad\qquad\qquad\cdot
\int_0^t  {B(s)\over 4 (t-s)\sqrt{\pi(t-s)}} \left(\int |x-y|
\exp \left\{ -  {(y + 2(t-s) -x)^2\over 4(t-s)} \right\} dy\right)ds \cr
&= \exp \left\{ 4 \| DA \|_{L^\infty}
\int_0^t \bigl\| u_x(\sigma) \bigr\|_{L^\infty} d\sigma
+ t - x \right\} \int_0^t  {\| A \|_{L^\infty}
B(s)\over \sqrt{\pi(t-s)}} \left(\int \bigl| \zeta - \sqrt{t-s} \bigr|
e^{ - \zeta^2} d\zeta\right) ds \cr
&\leq \exp \left\{ 4 \| DA \|_{L^\infty}
\int_0^t \bigl\| u_x(\sigma) \bigr\|_{L^\infty} d\sigma
+ t - x \right\}
\int_0^t \|A\|_{\L^\infty}  \biggl( {1\over \sqrt{t-s}} +
\sqrt{\pi} \biggr) B(s)\,ds \cr
&\leq\exp \left\{ 4 \| DA \|_{L^\infty}
\int_0^t \bigl\| u_x(\sigma) \bigr\|_{L^\infty} d\sigma
+ t - x \right\} \left(  {B(t)\over 2} -  {1\over 2} \right) \cr
&=  {1\over 2} E(t,x) - 
 {1\over 2} \exp \left\{ 4 \| DA \|_{L^\infty}
\int_0^t \bigl\| u_x(\sigma) \bigr\|_{L^\infty} d\sigma
+ t - x \right\} \cr
& \leq  {1\over 2} E(t,x)-  {1\over 2} e^{t-x},
\cr}
$$
\v
$$
\eqalign{ 2 \|DA\|_{L^\infty}&
\int_0^t \bigl\| u_x(s) \bigr\|_{L^\infty} \left(\int G(t-s,x-y)
E(s,y) dy\right)ds \cr
&= 2 \|DA\|_{L^\infty} \int_0^t B(s) \cdot \exp \left\{ 4
\|DA\|_{L^\infty}
\int_0^s \bigl\| u_x(\sigma) \bigr\|_{L^\infty} d\sigma +
s \right\}\cr
&\qquad\qquad\qquad\qquad\cdot  {\bigl\| u_x(s) \bigr\|_{L^\infty}\over 2
\sqrt{\pi(t-s)}}
\left(\int \exp \left\{ -  {(x-y)^2\over 4(t-s)} - y \right\}
dy\right)ds \cr
&\leq B(t) e^{t-x} \int_0^t 2 \|DA\|_{L^\infty}
\bigl\| u_x(s) \bigr\|_{L^\infty} \exp \left\{ 4
\|DA\|_{L^\infty}
\int_0^s \bigl\| u_x(\sigma) \bigr\|_{L^\infty} d\sigma \right\}
ds \cr
&=   B(t) e^{t-x} \bigg[{1\over 2}\exp \left\{ 4
\|DA\|_{L^\infty}
\int_0^t \bigl\| u_x(\sigma) \bigr\|_{L^\infty} d\sigma \right\}
-{1\over 2}\bigg]\cr
& \leq
{1\over 2} E(t,x)-  {1\over 2} e^{t-x}.
\cr}$$
\vfill\eject
\n{\bf Acknowledgements.}  This research was supported in part
by the Italian M.U.R.S.T.
The work has been completed while the first author was visiting
Academia Sinica in Taipei and Max Planck Institute in Leipzig.
\vsk
\c{\medbf References}
\v
\i{[AM]} F.~Ancona and A.~Marson, Well posedness for general
$2\times 2$ systems of conservation laws, {\it Amer. Math. Soc. Memoir},
to appear.
\v
\i{[BaJ]} P.~Baiti and H.~K.~Jenssen, On the front tracking 
algorithm,
{\it J. Math. Anal. Appl.} {\bf 217} (1998), 395-404.
\v
\i{[BLFP]} P.~Baiti, P.~LeFloch and B.~Piccoli,
Uniqueness of classical and nonclassical solutions for nonlinear
hyperbolic systems,
{\it J. Differential Equations}  {\bf 172} (2001), 59-82.         
\v
\i{[BiB1]}  S.~Bianchini and A.~Bressan, BV estimates for a class of
viscous hyperbolic systems,  {\it Indiana Univ. Math. J.} {\bf 49},
(2000), 1673-1713.
\v
\i{[BiB2]} S.~Bianchini and A.~Bressan, On a Lyapunov 
functional relating shortening curves and conservation laws, 
{\it Nonlinear Anal., T.M.A.}, to appear.
\v
\i{[BiB3]} S.~Bianchini and A.~Bressan, A case study in vanishing viscosity,
{\it Discr. Cont. Dyn. Sys.} {\bf 7} (2001), 449-476.
\v
\i{[BiB4]} S.~Bianchini and A.~Bressan, 
A center manifold technique for tracing viscous waves, 
{\it Comm. Pure Applied Analysis}, to appear.
\v
\i{[B1]} A.~Bressan, Contractive metrics for nonlinear 
hyperbolic systems, {\it Indiana Univ. J. Math.} {\bf 37} (1988), 409-421.
\v
\i{[B2]} A.~Bressan, Global solutions to systems of conservation laws by 
wave-front tracking, {\it J. Math. Anal. Appl.} {\bf 170} (1992),
414-432.
\v
\i{[B3]} A.~Bressan, The unique limit of the Glimm scheme,
{\it Arch. Rational Mech. Anal.} {\bf 130} (1995), 205-230.
\v
\i{[B4]} A.~Bressan, A locally contractive metric for systems of 
conservation laws, {\it Ann. Scuola Norm. Sup. Pisa} {\bf IV - 22}
(1995), 109-135.
\v
\i{[B5]} A.~Bressan,
{\it Hyperbolic 
Systems of Conservation Laws. The One Dimensional Cauchy Problem}.
Oxford University Press, 2000.
\v
\i{[BC]} A.~Bressan and R.~M.~Colombo, The semigroup generated
by $2\times 2$ conservation laws, {\it Arch. Rational Mech. Anal.}
{\bf 133} (1995), 1-75.
\v
\i{[BCP]} A.~Bressan, G.~Crasta and B.~Piccoli, Well posedness of the
Cauchy problem for $n\times n$ conservation laws, {\it
Amer. Math. Soc. Memoir} {\bf 694} (2000).
\v
\i{[BG]} A.~Bressan and P.~Goatin, Oleinik type estimates and 
uniqueness for $n\times n$ conservation laws, {\it J. Diff.
Equat.} {\bf 156} (1999), 26-49.
\v
\i{[BLF]} A.~Bressan and P.~LeFloch, Uniqueness of weak solutions to
systems of conservation laws, {\it Arch. Rat. Mech. Anal.} {\bf 140} 
(1997), 301-317.
\v
\i{[BLe]} A.~Bressan and M.~Lewicka, A uniqueness condition for
hyperbolic systems of conservation laws, {\it Discr.
Cont. Dynam. Syst.} {\bf 6} (2000), 673-682.
\v
\i{[BLY]} A.~Bressan, T.~P.~Liu and T.~Yang,  $ L^1$ stability
estimates for $n\times n$ conservation laws,
{\it Arch. Rational Mech. Anal.} {\bf 149} (1999), 1-22.
\v
\i{[Cl]} F.~H.~Clarke, {\it Optimization and Nonsmooth Analysis},
Wiley - Interscience, 1983.
\v
\i{[Cr]} M.~Crandall, The semigroup approach to first-order quasilinear 
equations in several space variables, {\it Israel J. Math.} {\bf 12}
(1972), 108-132.
\v
\i{[D]} C.~Dafermos, {\it Hyperbolic Conservation Laws in Continuum
Physics}, Springer-Verlag, Berlin 1999.
\v
\i{[DP1]} R.~DiPerna, Global existence of solutions to nonlinear hyperbolic 
systems of conservation laws, { \it J. Diff. Equat.}  {\bf 20}
(1976), 187-212.
\v
\i{[DP2]} R.~DiPerna, Convergence of approximate solutions to conservation
laws, {\it Arch. Rational Mech. Anal.} {\bf 82} (1983), 27-70.
\v
\i{[EG]} L.~C.~Evans and R.~F.~Gariepy, {\it 
Measure theory and fine properties of functions}, CRC Press, Boca Raton (1992).
\v
\i{[F]} L.~Foy, Steady-state solutions of hyperbolic systems of 
conservation laws with viscosity terms. {\it Comm. Pure Appl. Math.} {\bf
17} (1964), 177-188.
\v
\i{[G]} J.~Glimm, Solutions in the large for nonlinear hyperbolic systems 
of equations, {\it Comm. Pure Appl. Math.} {\bf 18} (1965), 697-715.
\v
\i{[GX]} J.~Goodman and Z.~Xin, Viscous limits for piecewise smooth
solutions to systems of conservation laws, {\it Arch. Rational Mech. Anal.}
{\bf 121} (1992), 235-265.
\v
\i{[Kr]} S.~Kruzhkov, First order quasilinear equations with several space 
variables, {\it Math. USSR Sbornik} {\bf 10} (1970), 217-243.
\v
\i{[Lx]} P.~Lax, Hyperbolic systems of conservation laws II, {\it
Comm. Pure Appl. Math.} {\bf 10} (1957), 537-566.
\v
\i{[L1]} T.~P.~Liu, The entropy condition and the admissibility of shocks,
{\it J. Math. Anal. Appl.} {\bf 53} (1976), 78-88.
\v
\i{[L2]} T.~P.~Liu, Admissible solutions of hyperbolic
conservation laws, {\it Amer. Math. Soc. Memoir} {\bf 240} (1981).
\v
\i{[L3]} T.~P.~Liu, Nonlinear stability of shock waves,
{\it Amer. Math. Soc. Memoir} {\bf 328} (1986).
\v
\i{[LY]} T.~P.~Liu and T.~Yang, Weak solutions of general systems of
hyperbolic conservation laws, preprint.
\v
\i{[MP]} A.~Majda and R.~Pego, Stable viscosity matrices for 
systems of conservation laws, {\it J. Differential Equations} {\bf 56}
(1985), 229-262.
\v
\i{[Ri]} N.~H.~Risebro, A front-tracking alternative to the random choice
method, {\it Proc. Amer. Math. Soc.} {\bf 117} (1993), 1125-1139.
\v 
\i{[Se]} D.~Serre, {\it Systems of Conservation Laws I, II},
Cambridge University Press, 2000.
\v
\i{[SX]} A.~Szepessy and Z.~Xin, Nonlinear stability abd viscous shocks,
{\it Arch. Rational Mech. Anal.} {\bf 122} (1993), 53-103.
\v
\i{[SZ]} A.~Szepessy and K.~Zumbrun, Stability of rarefaction waves in
viscous media, {\it Arch. Rational Mech. Anal.} {\bf 133} (1996), 249-298.
\v
\i{[V]} A.~Vanderbauwhede, Centre manifolds, normal forms and elementary
bifurcations, {\it Dynamics Reported, Vol.~2} (1989), 89-169.
\v
\i{[Yu]} S.~H.~Yu, Zero-dissipation limit of solutions with shocks for
systems of hyperbolic conservation laws, {\it Arch. Rational
Mech. Anal.}
{\bf 146} (1999), 275-370.
\bye